\documentclass[10pt, final]{amsart}
\usepackage{amsmath}
\usepackage{amsfonts}
\usepackage{amssymb}
\usepackage{amscd}
\usepackage{amsthm}
\usepackage[mathscr]{eucal}
\theoremstyle{plain}% default
\theoremstyle{definition}
\newtheorem{theorem}{Theorem}[section]
\newtheorem{lemma}[theorem]{Lemma}
\newtheorem{proposition}[theorem]{Proposition}
\newtheorem{corollary}[theorem]{Corollary}
\newtheorem{definition}[theorem]{Definition}
\newtheorem{remark}[theorem]{Remark}
\newtheorem{remarks}[theorem]{Remarks}
\newtheorem{examples}[theorem]{Examples}

\newtheorem*{claim}{Claim}
\newcommand{\mc}[1]{{\mathcal #1}}
\newcommand{\df}[2]{\displaystyle\frac{#1}{#2}}
\newcommand{\seq}[2]{\ensuremath{({#1}_{#2})_{#2\in\mathbb{N}}}}
 \DeclareMathOperator{\supp}{supp}

\DeclareMathOperator{\conv}{conv}
\DeclareMathOperator{\convq}{conv_{\mathbb{Q}}}
\DeclareMathOperator{\spann}{span}
 \DeclareMathOperator{\ran}{ran}

\DeclareMathOperator{\dist}{dist} \DeclareMathOperator{\ind}{ind}
\DeclareMathOperator{\sgn}{sgn} \DeclareMathOperator{\Ext}{Ext}
\newcommand{\R}{\mathbb{R}}
\newcommand{\N}{\mathbb{N}}
\newcommand{\Q}{\mathbb{Q}}
\newcommand{\e}{\varepsilon}
\newcommand{\co}{\ensuremath{c_{00}(\N)\;}}
\newcommand{\eqs}{{\mathfrak X}}
\newcommand{\norm}[1][\cdot]{\Vert #1\Vert}
\theoremstyle{plain}% default

\begin{document}

\title[\emph{The attractors method and Hereditarily James Tree
 Spaces}]{\bf Saturated extensions, the attractors method   and Hereditarily James Tree
 Spaces}
 \dedicatory{Dedicated
to the memory of R.C. James}

\author{Spiros A. Argyros}
\address[S.A. Argyros]{Department of Mathematics,
National Technical University of Athens}
\email{sargyros@math.ntua.gr}

\author{Alexander D. Arvanitakis}
\address[A. D. Arvanitakis]{Department of Mathematics, National
Technical University of Athens} \email{aarva@math.ntua.gr}

\author{Andreas G. Tolias}\thanks{Research partially
supported by EPEAEK program ``PYTHAGORAS"}
\address[A. Tolias]{Department of Applied Mathematics,
University of Crete} \email{atolias@tem.uoc.gr}

\maketitle

\tableofcontents

 \setcounter{section}{-1}
 \section{Introduction}
The purpose of the present work is to provide examples of HI
Banach spaces with no reflexive subspace and study their
structure. As is well known W.T. Gowers \cite{G1} has constructed
a Banach space $\mathfrak{X}_{gt}$  with a boundedly complete
basis $(e_n)_n$, not containing $\ell_1,$ and such that all of its
infinite dimensional subspaces have non separable dual. We shall
refer to this space as the Gowers Tree space. The predual
$(\mathfrak{X}_{gt})_*,$ namely the space generated by the
biorthogonal of the basis,  also has the property that it does not
contain $c_0$ or a reflexive subspace. It remains unknown whether
$\mathfrak{X}_{gt}$ is HI and moreover the structure of
$\mathcal{L}(\mathfrak{X}_{gt})$ is unclear. Notice that Gowers
dichotomy \cite{G2} yields that $\mathfrak{X}_{gt}$ and
$(\mathfrak{X}_{gt})_*$ contain HI subspaces. The structure of
$\mathfrak{X}_{gt}^*$ also remains unclear. The main obstacle for
understanding the structure of $\mathfrak{X}_{gt}$ or
$\mathcal{L}(\mathfrak{X}_{gt})$ is the use of a probabilistic
argument for establishing the existence of vectors with
 certain properties.

Our approach in constructing HI spaces with no reflexive subspace,
is different from Gowers' one. In particular we avoid the use of
any probabilistic argument and thus we are able to control the
structure of the spaces as well as the structure of the spaces of
bounded linear operators acting on them. Moreover we are able to
provide examples of spaces $X$ exhibiting a vast difference
between the  structures of $X$ and $X^*.$

The following are the highlight of our results:
\begin{itemize}
\item There exists a HI Banach space $X$ with a shrinking basis
and with no reflexive subspace. Moreover
every $T:X\to X$ is of the form $\lambda I + W$ with $W$ weakly compact
(and hence strictly singular).
\end{itemize}
 The absence of reflexive subspaces in $X$ in
 conjunction with the property that every
  strictly singular operator is weakly compact is  evidence supporting the
existence of Banach spaces such that every non Fredholm operator
is compact.
\begin{itemize}
\item The dual $X^*$ of the previous $X$ is HI and reflexively
saturated and the dual of every subspace $Y$ of $X$ is also HI.
\end{itemize}
This shows a strong divergence between the structure of $X$ and
$X^*.$ We recall that in \cite{AT2} a reflexive HI space $X$ is
constructed whose dual $X^*$ is  unconditionally saturated. The
analogue of this in the present setting
is the following one:%\\[0.1cm]
\begin{itemize}
\item There exists a HI Banach space $Y$ with a shrinking basis
and with no reflexive subspace, such that the
dual space $Y^*$ is reflexive and unconditionally saturated.%\\[0.2mm]
\end{itemize}
The definition of the  space $Y$ requires an adaptation of  the
methods of \cite{AT2} within the present framework of building
spaces with no reflexive subspace.
\begin{itemize}
\item There exists a partition of the basis $(e_n)_n$ of the
previous $X$ into two sets $(e_n)_{n \in L_1},$ $(e_n)_{n\in L_2}$
such that setting $X_{L_1}=\overline{\spann}\{e_n:\;n\in L_1\}$,
$X_{L_2}=\overline{\spann}\{e_n:\;n\in L_2\}$, both
 $X_{L_1}^*$, $X_{L_2}^*$  are HI with no reflexive
subspace.
\end{itemize}
The pairs $X_{L_i}, X_{L_i}^*$ for $i=1,2$ share similar
properties with the pair $(\mathfrak{X}_{gt})_*$ and
$\mathfrak{X}_{gt}.$
\begin{itemize}
\item The space $X^{**}$ is non separable and every  $w^*$-closed
subspace of $X^{**},$  is either   isomorphic to $\ell_2$ or is
non-separable and  contains $\ell_2.$ Therefore every quotient of
$X^*$ has a further quotient isomorphic to $\ell_2.$ Moreover
$X^{**}/X$ is isomorphic to $\ell_2(\Gamma).$
\end{itemize}
It seems also possible that $\mathfrak{X}_{gt}^*$ satisfies a
similar to the above property although this is not easily shown.
Further $X^*$ is the first example of a HI space with the
following property: $X^*/Y$ is HI whenever $Y$ is $w^*$-closed
(this is equivalent to say that for every subspace $Z$ of $X,$
$Z^*$ is HI) and also every quotient of $X^*$ has a further
quotient which is isomorphic to $\ell_2.$
\begin{itemize}
 \item There exists a non separable HI Banach space $Z$
not containing a reflexive subspace and such that every bounded
linear operator $T:Z\to Z$
 is of the form $T= \lambda I+ W$ with $W$ a
weakly compact (hence strictly singular) operator with separable
range.
\end{itemize}
This is an extreme construction resulting from a variation of the
methods used in the construction of the space $X$ involved in the
previous results. The  fact of the matter is that these methods
are not stable. Thus some minor changes in the initial data could
produce spaces with entirely different structure. Notice that the
space $Z$ is of the form $Y^{**}$ with $Y$ and $Y^*$ sharing
similar properties with the pair $X, X^*$ appearing in the
previous statements.

We shall proceed to a more detailed presentation of the results of
the paper and also of the methods used for constructing the
spaces, which are interesting on their own. We have divided the
rest of the introduction in three subsections. The first concerns
the structure of Banach spaces not containing $\ell_1, c_0$ or
reflexive subspace. The second is devoted to  saturated extensions
and in the third we explain the method of attractors which permits
the construction of dual pairs $X, X^*$ with strongly divergent
structure.

\subsection{Hereditarily James Tree spaces}

Separable spaces like Gowers Tree space undoubtedly have peculiar
structure. Roughly speaking, in every subspace one can find a
structure similar to the James tree basis. Next we shall attempt
to be more precise. Thus we shall define the Hereditarily James
Tree spaces, making more transparent their structure. We begin by
recalling some of the fundamental characteristics of James'
paradigm.

In the sequel we shall denote by $(\mc{D}, \prec)$ the dyadic tree
and by $[\mc{D}],$ the set of all branches (or the body) of
$\mc{D}.$ As usual we would consider that the nodes of $\mc{D}$
consist of finite sequences of $0$'s and $1$'s and $a\prec b$ iff
$a$ is an initial part of $b.$ The lexicographic order of
$\mc{D},$ denoted by $\prec_{lex}$ defines a well ordering which
is consistent with the tree order (i.e. $a\prec b$ yields that $a
\prec_{lex} b$).

{\bf The space $JT.$}

The James Tree space $JT$ (\cite{J}) is defined as the completion
of $(c_{00}(\mc{D}), \| \cdot \|_{JT})$ where for $x\in
c_{00}(\mc{D}),$ $\|x\|_{JT}$ is defined as follows:
\[ \| x\|_{JT} = \sup \big\{ \big(
\sum_{i=1}^n \big( \sum_{n\in s_i} x(n)\big)^2\big)^{1/2}
:\; (s_i)_{i=1}^n \text{ pairwise disjoint segments }\big\}.\]

The main properties of the space $JT,$ is that  does not contain
$\ell_1$ and has nonseparable dual.

Next, we list some properties of $JT$ related to our
consideration.

\begin{itemize}
\item The Hamel basis $(e_a)_{a\in\mc{D}}$ of $c_{00}(\mc{D})$
ordered with the lexicographic order defines a (conditional)
boundedly complete basis of $JT.$

\item For every branch $b$ in $[\mc{D}],$ $b=(a_1 \prec a_2 \prec
\cdots \prec a_n \cdots )$ the sequence $(e_{a_n})_n$ is non
trivial weak Cauchy and moreover $b^* = w^*-\sum_{n=1}^\infty
e_{a_n}^*$ defines a norm one functional in $JT^*.$

\item The biorthogonal functionals of the basis
$(e_a^*)_{a\in\mc{D}}$ generate the predual $JT_*$ of $JT$ and
they satisfy the following property.

For every segment $s$ of $\mc{D}$ setting $s^* = \sum_{a\in s}
e_a^*$ we have that $\|s^*\| =1.$
\end{itemize}

It is worth pointing out an alternative definition of the norm of
$JT.$ Thus we consider the following subset of $c_{00}(\mc{D}),$
\[ G_{JT} = \big\{ \sum_{i=1}^n \lambda_i s_i^* :
(s_i)_{i=1}^n \text{ are disjoint finite segments and }
\sum_{i=1}^n \lambda_i^2 \le 1 \big\} \]

Here $s_i^*$ are defined as before. It is an easy exercise to see
that the norm induced by the set $G_{JT}$ on $c_{00}(\mc{D})$
coincides with the norm of $JT.$

{\bf The James Tree properties.}

Let $X$ be a space with a Schauder basis $(e_n)_n.$ A block
subspace $Y$ of $X$ has the boundedly complete (shrinking) James
tree property if there exists a seminormalized block (in the
lexicographical order $\prec_{lex}$ of $\mc{D}$) sequence
$(y_a)_{a\in \mc{D}}$ in $Y$ and a $c>0$ such that the following
holds.
\begin{enumerate}
\item {\bf(boundedly complete)} There exists a bounded family
$(b^*)_{b\in[\mc{D}]}$  in $X^*$,
 such that for each $b\in [\mc{D}]$,
  $b= (a_1, a_2, \dots, a_n, \dots)$ the sequence $(y_{a_n})_n$
is  non trivial weakly Cauchy with $\lim b^*(y_{a_n})
>c$ and $\lim b_1^*(y_{a_n}) =0$ for all $b_1 \neq b.$

\item {\bf(shrinking)} For all finite segments $s$ of $\mc{D}$,
 $\| \sum\limits_{a\in s} y_a\| \le c$.
\end{enumerate}

Let's observe that $(e_a)_{a\in \mc{D}}$ in $JT$ satisfies the
boundedly complete James Tree property while $(e_a^*)_{a\in
\mc{D}}$ in $JT_*$ satisfies the shrinking one. Also, if the
initial space $X$ has a boundedly complete basis only the
boundedly complete James Tree property could occur. A similar
result holds if $X$ has a shrinking basis. Finally if $Y$ has the
boundedly complete James Tree property, then $Y^*$ is non
separable and if $X$ has a shrinking basis and $Y$ has the
(shrinking) James Tree property, then $Y^{**}$ is non separable.

For simplicity, in the sequel we shall consider that the initial
space $X$ has either a boundedly complete or a shrinking basis.
Thus if a block subspace has the James Tree property, then it will
be determined as either boundedly complete or shrinking according
to the corresponding property of the initial basis.

\begin{definition}
Let $X$ be a Banach space with a Schauder basis.

\begin{enumerate}
\item[(a)] A family $\mc{L}$ of block subspaces of $X$ has the
James Tree property, provided every $Y$ in $\mc{L}$ has that
property.

\item[(b)] The space $X$ is said to be Hereditarily James Tree
(HJT) if it does not contain $c_0,$ $\ell_1$ and every block
subspace $Y$ of $X,$ has the James Tree property.
\end{enumerate}

\end{definition}

It follows from Gowers' construction that the Gowers Tree space
$\mathfrak{X}_{gt},$ and its predual $(\mathfrak{X}_{gt})_*$ are
HJT spaces.

One of the results of the present work is that HJT property is not
preserved under duality. Namely, there exists a HJT space $X$ with
a shrinking basis, such that $X^*$ is reflexively (even
unconditionally) saturated. However, in the same example there
exists a subspace $Y$ of $X$ with $Y^*$ also an HJT space.

One of the basic ingredients in our approach to building HJT
spaces is the following space:
\begin{proposition}
There exists a space $JT_{\mc{F}_2}$ with a boundedly complete
basis $(e_n)_n$ such that the following hold:
\begin{enumerate}
\item[(i)] The space $JT_{\mc{F}_2}$ is $\ell_2$ saturated.

\item[(ii)] The basis $(e_n)_n$ is normalized weakly null and for
every $M\in [\N]$ the subspace $X_M =\overline{\spann}\{e_n:\;n\in
M\}$ has the James tree property.
\end{enumerate}
\end{proposition}

It is clear that none subsequence $(e_n)_{n\in M}$ is
unconditional. Thus the basis of $JT_{\mc{F}_2}$ shares similar
properties with the classical Maurey Rosenthal example \cite{MR}.
We shall return to this space in the sequel explaining more about
its structure and its difference from Gowers' space.

{\bf Codings and tree structures.} As  is well known, every
attempt to impose tight (or conditional) structure in some Banach
space, requires the definition of the conditional elements which
in turn results from the existence of special sequences defined
with the use of a coding. What  is less well known is that the
codings induce a tree structure in the special sequences. As we
shall explain shortly, the James tree structure in the subspaces
of HJT spaces, like $\mathfrak{X}_{gt},$ $(\mathfrak{X}_{gt})_*$
or the spaces presented in this paper, are directly related to
codings.

Let's start with a general definition of a coding, and the
obtained special sequences. Consider a collection $(F_j)_j$ with
each $F_j$ a countable family of elements of $c_{00}(\N).$ To make
more transparent the meaning of our definitions, let's assume that
each $F_j=\{\frac{1}{m_j^2} \sum_{k\in F}e_k^*:\; F\subset \N,\;
\# F\le n_j\}$ where $(m_j), (n_j)$ are appropriate fast
increasing sequences of natural numbers. Notice that the elements
of the family $\mc{T}=\cup_j F_j$ and the combinations of them
will play the role of functionals belonging to a norming set. This
explains the use of $e_k^*$ instead of $e_k.$ For simplicity, we
also assume that the families $(F_j)_j$ are pairwise disjoint.
This happens in the aforementioned example although it is not
always true. Under this additional assumption to each $\phi \in
\cup_j F_j$ corresponds a unique index by the rule $\text{ind}
(\phi) = j$ iff $\phi \in F_j.$ Further for a finite block
sequence $s=(\phi_1, \dots, \phi_d)$ with each $\phi_i \in \cup_j
F_j,$ we define $\text{ind} (s) = \{ \text{ind}(\phi_1), \dots,
\text{ind} (\phi_d)\}.$

{\bf The $\sigma$-coding:} Let $\Omega_1, \Omega_2$ be a partition
of $\N$ into two infinite disjoint subsets. We denote by $\mc{S}$
the family of all block sequences $s=(\phi_1 < \phi_2 < \cdots <
\phi_d)$ such that $\phi_i\in \cup_j F_j,$ $\text{ind}(\phi_1) \in
\Omega_1, $ $\{\text{ind}(\phi_2) < \cdots < \text{ind}(\phi_d) \}
\subset \Omega_2.$ Clearly $\mc{S}$ is countable, hence there
exists an injection
\[ \sigma: \mc{S} \to \Omega_2 \]
satisfying $\sigma(s) > \text{ind}(s)$ for every $s\in\mc{S}.$

{\bf The $\sigma$-special sequences:} A sequence $s=(\phi_1 <
\phi_2 < \cdots < \phi_n)$ in $\mc{S}$ is said to be a
$\sigma$-special sequence iff for every $1\le i < n$ setting $s_i
= (\phi_1 < \cdots < \phi_i)$ we have that \[ \phi_{i+1} \in
F_{\sigma(s_i)}.\]

The following tree-like interference holds for $\sigma$-special
sequences.

Let $s, t$ be two $\sigma$-special sequences with $s=(\phi_1,
\dots, \phi_n),$ $t=(\psi_1, \dots, \psi_m).$ We set
\[ i_{s, t} = \max \{i: \text{ind}(\phi_i) = \text{ind}(\psi_i) \} \]
 if the later set is non empty. Otherwise
we set $i_{s,t} = 0.$  Then the following are easily checked.
\begin{enumerate}
\item[(a)] For every $i< i_{s,t}$ we have that $\phi_i = \psi_i.$
\item[(b)] $\{\text{ind}(\phi_i) : i > i_{s, t}\} \cap
\{\text{ind}(\psi_j) : j> i_{s,t}\} = \emptyset.$
\end{enumerate}

These two properties immediately yield that the set $\mc{T} \cup_jF_j$ endowed with the partial order $\phi\prec_{\sigma}
\psi$ iff there exists a $\sigma$-special sequence $(\phi_1,
\dots, \phi_n)$ and $1\le i < j\le n$ with $\phi=\phi_i$ and
$\psi= \phi_j$ is a tree.

Now for the given tree structure $(\mc{T}, \prec_\sigma)$ we will
define norms similar to the classical James tree norm mentioned
above.

{\bf The space $JT_{\mc{F}_2}$:} For the first application the
family $(F_j)_j$ is the one defined above.

For  a $\sigma$-special sequence $s= (\phi_1, \dots, \phi_n)$ and
an interval $E$ of $\N$ we set $s^*= \sum_{k=1}^n \phi_k$ and let
 $Es^*$ be
the restriction of $s^*$ on $E$ (or the pointwise product $s^*
\chi_E$). A $\sigma$-special functional $x^*$ is any element
$Es^*$ as before.

Also, for a $\sigma$-special functional $x^*= Es^*,$ $s=(\phi_1,
\dots,\phi_n)$, we let $\text{ind}(x^*) =\{\text{ind}(\phi_k):
\supp \phi_k \cap E\neq \emptyset\}.$ We consider the following
set
\[ \begin{split}
\mc{F}_2 = &\{\pm e_n^*:\;n\in\N\} \cup \{
 \sum_{i=1}^d a_i x_i^*:\;a_i\in\Q,\;\sum\limits_{i=1}^d
  a_i^2 \le 1,\; (x_i^*)_{i=1}^d\text{ are} \\
&\text{ $\sigma$-special functionals with $(\ind(x_i^*))_{i=1}^d$
pairwise disjoint,}\;d\in\N \} \end{split} \]

The space $JT_{\mc{F}_2}$ is the completion of $(c_{00},
\|.\|_{\mc{F}_2})$ where for $x\in c_{00},$
\[ \|x\|_{\mc{F}_2} = \sup \{ \phi(x) : \phi \in \mc{F}_2 \}.\]

Comparing the norming set $\mc{F}_2$ with the norming set $G_{JT}$
of $JT$ one observes that $\sigma$-special functionals in
$\mc{F}_2$ play the role of the functionals $s^*$ defined by the
segments of the dyadic tree $\mc{D}.$ As we have mentioned in
Proposition 0.2, the space $JT_{\mc{F}_2}$ , like $JT$, is
$\ell_2$ saturated, but for every $M\in [\N],$ the subspace $X_M \overline{\spann}\{e_n:\;n\in M\}$ has non separable dual. The
later is a consequence of the fact that the tree structure
$(\mc{T}, \prec_\sigma)$ is  richer than that of the dyadic tree
basis in $JT$. Indeed, it is easy to check that for every $M\in
[\N]$ we can construct a block sequence $(\phi_a)_{a\in \mc{D}}$
such that
\begin{enumerate}
\item[(i)] $\phi_a = \frac{1}{m_{j_a}^2} \sum_{k\in F_a}e_k^*$
where $\# F_a = n_{j_a}$ and $F_a\subset M$, while $F_a < F_\beta$
if $a\prec_{lex} \beta$.

\item[(ii)] For a branch $b=(a_1 \prec a_2 \prec \cdots \prec a_n
\prec \cdots)$  of $\mc{D}$ and for every $n\in \N$ we have that
$(\phi_{a_1}, \dots, \phi_{a_n})$ is a $\sigma$-special sequence.
\end{enumerate}

Defining now $x_a=\frac{m_{j_a}^2}{n_{j_a}}\sum_{k\in F_a}e_k,$
the family $(x_a)_{a\in \mc{D}}$ provides the James tree structure
of $X_M$.

{\bf The Gowers Tree space.} The definition of $\mathfrak{X}_{gt}$
uses similar ingredients with the corresponding of $JT_{\mc{F}_2}$
although structurally the two spaces are entirely different. The
norming set $G_{gt}$ of Gowers space is saturated under the
operations $(\mc{A}_{n_j},\frac{1}{m_j})_j.$ We recall that a
subset $G$ of $c_{00}$ is closed (or saturated) for the operation
$(\mc{A}_n,\frac{1}{m})$ if for every $\phi_1 < \phi_2 < \cdots <
\phi_k,$ $k\le n$ with $\phi_i \in G,$ $i=1, \dots, k,$ the
functional $\phi= \frac{1}{m}\sum_{i=1}^k \phi_i$ belongs to $G$.

The norming set $G_{gt}$ is the minimal subset of $c_{00}$
satisfying the following conditions:

\begin{enumerate}
 \item[(i)] $\{\pm e_k^*:k\in\N\}\subset G_{gt},$ $G_{gt}$ is
 symmetric and closed under the operation of restricting  elements to the
 intervals.
 \item[(ii)] $G_{gt}$ is closed in the
 $(\mc{A}_{n_j},\frac{1}{m_j})_j$ operations. We also set
 \[K_j= \{\phi\in G_{gt}: \phi \text{ is the result of a }
 (\mc{A}_{n_j},\frac{1}{m_j})
 \text{ operation}\}\]
 \item[(iii)] $G_{gt}$ contains the  set
 \begin{eqnarray*}
  && \big\{\sum_{i=1}^d a_i x_i^*:\;a_i\in\Q,\;\sum\limits_{i=1}^d
  a_i^2 \le 1,\; (x_i^*)_{i=1}^d,\;\sigma\text{-special functionals} \\
 &&\text{with $(\ind(x_i^*))_{i=1}^d$
pairwise disjoint,}\;d\in\N \big\}
\end{eqnarray*}

 \item[(iv)] $G_{gt}$ is rationally convex.

\end{enumerate}

 We explain briefly condition (iii).
 For a coding $\sigma$, the $\sigma$-special sequences
 $(\phi_1,\ldots,\phi_n)$ are defined as in the case of
 $\mc{F}_2$.
 Here the set $K_j$ plays the role of the corresponding $F_j$ in
 $\mc{F}_2$.
  The $\sigma$-special functionals $x^*$, are defined
  as in the case of
 $\mc{F}_2$.

Let's observe that $G_{gt}$ is almost identical with $\mc{F}_2$,
although the spaces  defined by them are entirely different. The
essential difference between $\mc{F}_2$ and $G_{gt}$ is that in
the case of $\mc{F}_2$ each $F_j,$ $j\in \N$ does not norm any
subspace of $JT_{\mc{F}_2},$ while in $\mathfrak{X}_{gt}$ each
$K_j$ defines an equivalent norm on $\mathfrak{X}_{gt}.$ The later
means that in every subspace $Y$ of $\mathfrak{X}_{gt},$ the
families $K_j$, $j\in\N$ as well as  $\{x^*: x^* \text{ is a
$\sigma$-special functional}\}$ and $\{\sum\limits_{i=1}^d
\lambda_i x_i^*:\; \sum\limits_{i=1}^n \lambda_i^2 \le 1,\;
(\ind(x_i^*))_{i=1}^d \text{ are pairwise disjoint}\}$ define
equivalent norms making it difficult to distinguish the action of
them on the elements of $Y.$ Thus, while the spaces of the form
$JT_{\mc{F}_2}$ can be studied in terms of the classical theory,
the space $\mathfrak{X}_{gt}$ requires advanced tools, like Gowers
probabilistic argument, which do not permit a complete
understanding of its structure.

\subsection{Saturated extensions}

The method of HI extensions appeared in the  Memoirs monograph
\cite{AT1} and was used to derive the following two results:

\begin{itemize}
\item Every separable Banach space $Z$ not containing $\ell_1$ is
a quotient of a separable HI space $X,$ with the additional
property that $Q^* Z^*$ is a complemented subspace of $X^*.$ (Here
$Q$ denotes the quotient map from $X$ to $Z.$)

\item There exists a nonseparable HI Banach space.
\end{itemize}

Roughly speaking, the method of HI extensions provides a tool to
connect a given norm, usually defined through a norming set $G$
with a HI norm. The resulting new norm will preserve some of the
ingredients of the initial norm and will also be HI. To some
extent, HI extensions, have similar goals with HI interpolations
(\cite{AF}) and some of the results could be obtained with both
methods. However it seems that the method of extensions is very
efficient when we want to construct dual pairs $X, X^*$ with
divergent structure. This actually requires the combination of
extensions with the method of attractors, which appeared in
\cite{AT2} where a reflexive HI space $X$ is constructed with
$X^*$ unconditionally saturated.

In the sequel we shall provide a general definition of  saturated
extensions which include several forms of extensions which
appeared elsewhere (c.f. \cite{AT1, AT2, ArTo})

Let $\mc{M}$ be a compact family of finite subsets of $\N.$ For
the purposes of the present paper, $\mc{M}$ will be either some
$\mc{A}_n = \{ F\subset \N:\; \#F\le n\},$ or some $\mc{S}_n,$ the
$n^{th}$ Schreier family. For $0< \theta < 1,$ the $(\mc{M},
\theta)$-operation on $c_{00}$ is a map which assigns to each
$\mc{M}$-admissible block sequence $(\phi_1 < \phi_2 < \cdots <
\phi_n),$ the functional $\theta \sum_{i=1}^n \phi_i.$ (We recall
that $\phi_1, \phi_2, \dots, \phi_n$ is $\mc{M}$-admissible if $\{
\min\supp \phi_i:\;i=1,\ldots,n\}$ belongs to $\mc{M}.$) A subset
$G$ of $c_{00}$ is said to be closed in the $(\mc{M},
\theta)$-operation, if for every $\mc{M}$-admissible block
sequence $\phi_1, \dots, \phi_n,$ with each $\phi_i \in G,$ the
functional $\theta \sum_{i=1}^n \phi_i$ belongs to $G$. When we
refer to saturated norms we shall mean that there exists a norming
set $G$ which is closed under certain $(\mc{M}_j,\theta_j)_j$
operations.

Let $G$ be a subset of $c_{00}.$ The set $G$ is said to be a
ground set if it is symmetric, $\{e_n^*:\;n\in\N\}$ is contained
in $G,$ $\|\phi\|_\infty \le 1,$ $\phi(n) \in \mathbb{Q}$ for all
$\phi\in G$  and $G$ is closed under the restriction of its
elements to  intervals of $\N.$ A ground norm, $\| \cdot \|_G$ is
any norm induced on $c_{00}$ by a ground set $G.$ It turns out
that for every space $(X, \|\cdot\|_X)$ with a normalized Schauder
basis $(x_n)_n$ there exists a ground set $G_X$ such that the
natural map $e_n \mapsto x_n$ defines an isomorphism between $(X,
\|\cdot \|_X)$ and $\overline{(c_{00}, \| \cdot \|_{G_X})}.$

{\bf Saturated extensions of a ground set $G$.} Let $G$ be a
ground set, $(m_j)_j$ an appropriate sequence of natural numbers
and $(\mc{M}_j)_j$ a sequence of compact families such that
$(\mc{M}_j)_j$ is either $(\mc{A}_{n_j})_j$ or $(S_{n_j})_j.$

Denote by $E_G$ the minimal subset of $c_{00}$ such that
\begin{enumerate}
\item[(i)] The ground set $G$ is a subset of $E_G.$

\item[(ii)] The set $E_G$ is closed in the  $(\mc{M}_j,
\frac{1}{m_j})$ operation.

\item[(iii)] The set $E_G$ is rationally convex.
\end{enumerate}

\begin{definition}
A subset $D_G$ of $E_G$ is said to be a saturated extension of the
ground set $G$ if the following hold:

\begin{enumerate}
\item[(i)] The set $D_G$ is a subset of $E_G,$ the ground set $G$
is contained in $D_G$ and $D_G$ is closed under restrictions of
its elements to  intervals.

\item[(ii)] The set $D_G$ is closed under even operations
$(\mc{M}_{2j}, \frac{1}{m_{2j}})_j.$

\item[(iii)] The set $D_G$ is rationally convex.

\item[(iv)] Every $\phi \in D_G$ admits a tree analysis
$(f_t)_{t\in T}$ with each $f_t \in D_G.$
\end{enumerate}
\end{definition}

Denoting by $\|\cdot \|_{D_G}$ the norm on $c_{00}$ induced by
$D_G$ and letting $X_{D_G}$ be the space $\overline{(c_{00},
\|\cdot \|_{D_G})},$ we call $X_{D_G}$ a {\em saturated extension}
of the space $X_G = \overline{(c_{00}, \| \cdot\|_G)}.$

Let's point out that the basis $(e_n)_n$ of $c_{00}$ is a
bimonotone boundedly complete Schauder basis of $X_{D_G}$ and that
the identity $I: X_{D_G} \to X_G$ is a norm one operator. Observe
also that we make no assumption concerning the odd operations. As
we will see later making several assumptions for the odd
operations, we will derive saturated extensions with different
properties.

A last comment on the definition of $D_G,$ is related to the
condition (iv). The tree analysis $(f_t)_{t\in T}$ of a functional
$f$ in $E_G$ describes an inductive procedure for obtaining $f$
starting from elements of the ground set $G$ and either applying
operations $(\mc{M}_j, \frac{1}{m_j})$ or taking rational convex
combinations. This tree structure is completely irrelevant to the
tree structures discussed in the previous subsection. Its role is
to help estimate upper bounds of the norm of vectors in $X_{D_G}.$

{\em Properties and variants of Saturated extensions.}

As we have mentioned, for $x\in c_{00},$ $\|x\|_G \le \|
x\|_{D_G}.$ This is an immediate consequence of the fact that
$G\subset D_G.$ On the other hand, there are cases of ground sets
$G$ such that $D_G$ does not add more information beyond  $G$
itself. Such a case is when $G$ defines a norm $\|\cdot\|_G$
equivalent to the $\ell_1$ norm. The measure of the fact that
$\|\cdot \|_{D_G}$ is strictly greater than $\|\cdot \|_G$ on a
subspace $Y$ of $X_{D_G}$ is that the identity operator $I:
X_{D_G} \to X_G$ restricted to $Y$ is a strictly singular one. If
$I: X_{D_G} \to X_G$ is strictly singular, then we refer to
strictly singular extensions. The first result we want to mention
is that strictly singular extensions are reflexive ones. More
precisely the following holds:

\begin{proposition}
Let $Y$ be a closed subspace of $X_{D_G}$ such that $I|_Y : Y \to
X_G$ is strictly singular. Then Y is reflexively saturated. In
particular $X_{D_G}$ is reflexively saturated whenever it is a
strictly singular extension.
\end{proposition}

Next we proceed to specify the odd operations and to derive
additional information on the structure of $X_{D_G}$ whenever
$X_{D_G}$ is a strictly singular extension.

{\em (a) Unconditionally saturated extensions.}

This is the case where $D_G = E_G = D^u_G.$ In this case the
following holds:

\begin{proposition}
Let $Y$ be a closed subspace of $X_{D_G^u}$ such that $I|_Y: Y \to
X_G$ is strictly singular. Then $Y$ is unconditionally (and
reflexively) saturated.
\end{proposition}

{\em (b) Hereditarily Indecomposable  extensions.}

HI extensions, are the most important ones. In this case the
norming set $D_G^{hi}$ is defined as follows. $D_G^{hi}$ is the
minimal subset of $c_{00}$ satisfying the following conditions
\begin{enumerate}
 \item[(i)] $\{e_n^*:\;n\in\N\}\subset D_G^{hi}$, $D_G^{hi}$ is
 symmetric and closed under restriction of its elements to
 intervals.
 \item[(ii)] $D_G^{hi}$ is closed under
 $(\mc{M}_{2j},\frac{1}{m_{2j}})_j$ operations.
 \item[(iii)] For each $j$, $D_G^{hi}$ is closed under
 $(\mc{M}_{2j-1},\frac{1}{m_{2j-1}})$ operation on $2j-1$ special sequences.
 \item[(iv)] $D_G^{hi}$ is rationally convex.
 \end{enumerate}

 The $2j-1$ special sequences are defined through a coding $\sigma$
 and satisfy the following conditions.
 \begin{enumerate}
 \item[(a)] $(f_1,\ldots,f_d)$ is $\mc{M}_{2j-1}$ admissible
 \item[(b)] For $i\le i\le d$ there exists some $j\in \N$
 such that $f_i\in K_{2j}=\big\{\frac{1}{m_{2j}}\sum\limits_{l=1}^k\phi_l:\;
 \phi_1<\cdots<\phi_k\mbox{ is } \mc{M}_{2j} \mbox{ admissible, }
 \phi_l\in D_{G}^{hi}\}$ and if $i>1$ then
 $2j=\sigma(f_1,\ldots,f_{i-1})$.
 \end{enumerate}Notice that in the definition of $D_G^{hi}$ we do
 not refer to the tree analysis. The reason is that the existence
 of a tree analysis follows from the minimality of $D_G^{hi}$ and
 the conditions involved in its definition.

The analogue of the previous results also holds for HI extensions.

\begin{proposition}
Let $Y$ be a closed subspace of $X_{D_G^{hi}}$ such that $I|_Y : Y
\to X_G$ is strictly singular. Then $Y$ is a HI space. In
particular strictly singular and HI  extensions yield HI spaces.
\end{proposition}

The above three propositions indicate that if we wish to have
additional structure on $X_{D_G},$ $X_{D_G^u},$ $X_{D_G^{hi}}$ we
need to consider strictly singular extensions. As is shown in
\cite{AT1}, this is always possible. Indeed, for every ground set
$G$ such that the corresponding space $X_G$ does not contain
$\ell_1$
 there exists a
family $(\mc{M}_j, \frac{1}{m_j})_j$  such that the saturated
extension of $G$ by this family
 is a strictly singular one. Thus the following is
proven (\cite{AT1}).

\begin{theorem}
Let $X$ be a  Banach space with a normalized Schauder basis
$(x_n)_n$ such that $X$ contains no isomorphic copy of $\ell_1$.
Then there exists a HI space $Z$ with a normalized basis $(z_n)_n$
such that the map $z_n \mapsto x_n$ has a linear extension  to a
bounded operator $T: Z\to X.$
\end{theorem}

This theorem in conjunction with the following one  yields that
every separable Banach space $X$ not containing $\ell_1$ is the
quotient of a  HI space.

\begin{theorem}[\cite{AT1}]
Let $X$ be a separable Banach space not containing $\ell_1.$ Then
there exists a space $Y$ not containing $\ell_1$, with a
normalized Schauder basis $(y_n)_n$ and a bounded linear operator
$T: Y\to X$ such that $(Ty_n)_n$ is a dense subset of the unit
sphere of $X$.
\end{theorem}

 {\bf The predual $(X_{D_G^{hi}})_*$.}
 As we have mentioned before the basis \seq{e}{n} of
 $X_{D_G^{hi}}$ is boundedly complete,
 hence the space $(X_{D_G^{hi}})_*$, which is the subspace of $X_{D_G^{hi}}^*$
 norm generated by the biorthogonal functionals
 $(e_n^*)_{n\in\N}$, is a predual of $X_{D_G^{hi}}$.
  In many cases it is shown
 that $(X_{D_G^{hi}})_*$ is also a HI space. This requires some
 additional information concerning the weakly
 null block sequences in $X_G$, which is stronger than the
 assumption that
  the identity map $I:X_{D^{hi}_G}\to  X_G$ is strictly singular.
  For example in \cite{AT1}, for extensions using the operations
  $(\mc{S}_{n_j},\frac{1}{m_j})_j$,  had been assumed   that the
  ground set
  $G$ is $\mc{S}_2$ bounded. In the present paper for the
  operations $(\mc{A}_{n_j},\frac{1}{m_j})_j$ we introduce the
  concept of strongly strictly singular extension which yields
  that $(X_{D^{hi}_G})_*$ is HI. It is also worth pointing out that
 $(X_{D^{hi}_G})_*$ is not necessarily reflexively saturated as
 happens for the strictly singular extensions $X_{D_G}$
 $X_{D^{hi}_G}$. This actually will be a key point in our approach
 for constructing HI spaces with no reflexive subspace.

\subsection{The attractors method}

Let's return to our initial goal, namely constructing HI spaces
with no reflexive subspace. It is clear from our preceding
discussion that HI extensions of ground sets $G$ such that $X_G$
does not contain $\ell_1$ yield reflexively saturated HI spaces.
Therefore there is no hope to obtain HI spaces with no reflexive
subspace as a result of a HI extension of a ground set $G$. As
mentioned in \cite{ArTo}  saturation and HI methods share common
metamathematical ideas with the forcing method in set theory. In
particular the fact that  HI extensions are reflexively saturated
is similar to the well known collapse phenomena in the extensions
of models of set theory via the forcing method. An illustrating
example of such phenomena in HI extensions  is that
$\mathfrak{X}_{gt}$ is a quotient of a HI and reflexively
saturated space. In spite of all these discouraging observations
we claim that HI extensions could help to yield HI spaces with no
reflexive subspace, and this is closely related to the structure
of $(X_{D_G^{hi}})_*$. Evidently from the initial stages of HI
theory, (\cite{GM1},\cite{GM2},\cite{AD}) and for many years, a
norming set $D$ was defined, using saturation methods and codings,
in such a way as to impose certain properties in the space $(X,
\|\cdot \|_D).$ In \cite{AT2} the norming set $D$ was designed to
impose divergent properties in $(X, \|\cdot \|_D)$ and $(X,
\|\cdot \|_D)^*$. The method used for this is the attractors
method, not so named in \cite{AT2}, which will also be used in the
present work.

The general principle of the attractors method is the following:

We are interested in designing a ground set $G$ and a HI extension
$D_G^{hi}$ such that the following two divergent properties hold:

\begin{enumerate}
\item[(a)] $X_{D_G^{hi}}$ is a strictly singular extension of
$X_G.$ In other words every subspace $Y$ of $X_{D_G^{hi}}$
contains a further subspace $Z$ on which the $G$-norm becomes
negligible.

\item[(b)] The set $G$ is asymptotic in $(X_{D_G^{hi}})_*.$ This
means that there exists $c>0$ such that for every  subspace $Y$ of
$(X_{D_G^{hi}})_*$ and every $\e >0$ there exists $\phi \in G$
with $\|\phi \|_{(X_{D_G^{hi}})_*} \ge c$ and $\dist(\phi, Y) <
\e.$
\end{enumerate}

In other words, we want $G$ to be small, as a norming set for the
space $X_{D_G^{hi}}$ and large as a subset of $(X_{D_G^{hi}})_*.$
Notice that such a relation between $G$ and $D_G^{hi}$ requires
the two sets to be built with similar materials, and moreover to
impose certain special functionals in $D_G^{hi}$ (we call these
attractor functionals) which will allow us to attract in every
subspace of $(X_{D_G^{hi}})_*$ part of the structure of the set
$G.$

Let us be more precise explaining how we define the corresponding
sets $G$ and $D_G^{hi}$ to obtain a HI extension $X_{D_G^{hi}},$
such that $(X_{D_G^{hi}})_*$ is also HI and does not contain
reflexive subspaces.

{\bf The ground set $\mc{F}_2$ and the norming set
$D_{\mc{F}_2}$.}
 We start by defining  the following family $(F_j)_j$.
We shall use the sequence of positive integers $(m_j)_j,$
$(n_j)_j$ recursively defined as follows:

 \begin{itemize}
 \item  $m_1=2$  and $m_{j+1}=m_j^5$.
 \item $n_1=4$,  and $n_{j+1}=(5n_j)^{s_j}$  where  $s_j=\log_2 m_{j+1}^3$.
 \end{itemize}
 We set $F_0=\{\pm e_n^*:\;n\in\N\}$ and for
  $j=1,2,\ldots$ we set
 \[ F_j=\big\{\frac{1}{m_{4j-3}^2}\sum\limits_{i\in I}\pm e^*_i:\;
  \#(I)\le \frac{n_{4j-3}}{2} \big\}\cup\big\{0\big\}. \]

Using the family $(F_j)_j$ and a coding $\sigma_\mc{F},$ we define
the ground set $\mc{F}_2$ in the same manner as in the first
subsection.

Next we define the set $D_{\mc{F}_2}$ which is a HI extension of
$\mc{F}_2$ with attractors as follows:

The set $D_{\mc{F}_2}$ is a minimal subset of $c_{00}$ satisfying
the following properties:

\begin{enumerate}
 \item[(i)] ${\mc{F}_2}\subset D_{\mc{F}_2}$,
 $D_{\mc{F}_2}$ is  symmetric (i.e. if $f\in D_{\mc{F}_2}$ then  $-f\in D_{\mc{F}_2}$) and
  $D_{\mc{F}_2}$ is closed under the restriction of its elements
 to intervals of $\N$ (i.e. if $f\in D_{\mc{F}_2}$ and $E$ is an interval of
 $\N$ then $Ef\in D_{\mc{F}_2}$).
 \item[(ii)] $D_{\mc{F}_2}$ is closed under  $(\mc{A}_{n_{2j}},\frac{1}{m_{2j}})$
 operations, i.e. if $f_1<f_2<\cdots<f_{n_{2j}}$ belong to $D_{\mc{F}_2}$
 then the functional
 $f=\frac{1}{m_{2j}}(f_1+f_2+\cdots+f_{n_{2j}})$  belongs also to
 $D_{\mc{F}_2}$.
 \item[(iii)] $D_{\mc{F}_2}$ is closed under  $(\mc{A}_{n_{4j-1}},\frac{1}{m_{4j-1}})$
 operations on special sequences i.e. for every $n_{4j-1}$ special sequence
 $(f_1,f_2,\ldots,f_{n_{4j-1}})$   the functional
 $f=\frac{1}{m_{4j-1}}(f_1+f_2+\cdots+f_{n_{4j-1}})$
 belongs to $D_{\mc{F}_2}$.
 In this case we say that $f$ is
 a {\bf special functional}.
  \item[(iv)]
 $D_{\mc{F}_2}$ is closed under  $(\mc{A}_{n_{4j-3}},\frac{1}{m_{4j-3}})$
 operations on attractor sequences i.e. for every $4j-3$ attractor sequence
 $(f_1,f_2,\ldots,f_{n_{4j-3}})$   the functional
 $f=\frac{1}{m_{4j-3}}(f_1+f_2+\cdots+f_{n_{4j-3}})$
 belongs to $D_{\mc{F}_2}$.
 In this case we say that $f$  is
 an {\bf attractor}.
  \item[(v)] The set $D_{\mc{F}_2}$ is rationally convex.
 \end{enumerate}

 In the above definition, the special functionals
 and the attractors, defined in (iii) and (iv) respectively,
 require some more explanation. First, the $n_{4j-1}$ special sequences $(f_1, \dots, f_{n_{4j-1}})$ are
 defined through a coding $\sigma$ as in the case of the aforementioned HI extensions. Thus each $f_i,$ $1\le i
 \le n_{4j-1}$ belongs to some
 \[ K_{2j} = \{ \frac{1}{m_{2j}} \sum_{l=1}^{n_{2j}} \phi_l:\;
  \phi_1 < \cdots < \phi_{n_{2j}},\; \phi_l \in D_{\mc{F}_2} \} \]
 where $2j$ is equal to  $\sigma(f_1, \dots, f_{i-1})$ whenever $1<i.$

 The special functionals will determine the HI property of the extension $D_{\mc{F}_2}.$

 Each $4j-3$ attractor sequence $f_1< \cdots < f_{n_{4j-3}}$
  is of the following form. All the odd members of the
 sequence are elements of $\cup_j K_{2j}$
 while the even members are $f_{2i} = e_{\ell_{2i}}^*$ and furthermore the
 sequence $f_1, \dots, f_{n_{4j-3}}$
  is determined by the coding $\sigma$ in a similar manner to the $n_{4j-1}$
 special sequence. Let us observe that for every $j\in \N$ there exist
  many $P \subset \N$ with the following properties. First $\# P  \frac{n_{4j-3}}{2}$, hence
  $\frac{1}{m_{4j-3}^2} \sum_{\ell \in P} e_\ell^* \in F_j$
 and also there exists an attractor sequence
 $(f_1, \dots, f_{4j-3})$  with $\{e_\ell^*: \ell \in P\}$
 coinciding with the even terms of
 the sequence. The purpose of the attractors
  is to make the family $\cup_j F_j$ asymptotic in the space
 $(X_{D_{\mc{F}_2}})_*.$ In particular using attractors, the following is proved.

 For every subspace $Y$ of $(X_{D_{\mc{F}_2}})_*$
  and every $j\in \N$ there exist $\phi_j \in Y$ and
  $\psi_j=\frac{1}{m_{4j-3}^2} \sum_{\ell \in P} e_\ell^* \in
  F_j$,
     such that
 \begin{enumerate}
 \item[(a)] $\|\phi_j + \psi_j\| > \frac{1}{144}.$
 \item[(b)] $\|\phi_j - \psi_j\| \le \frac{1}{m_{4j-3}}.$
 \end{enumerate}

 This shows that indeed $\cup_j F_j$ is asymptotic
  and furthermore we can copy a complete dyadic block subtree of
 $(\mc{T} = \cup_j F_j, \prec_{\sigma_\mc{F}})$
  into the subspace $Y$ which yields that $Y$ is not reflexive.

  The following summarizes the main steps in our approach to
  constructing HI spaces with no reflexive subspace.

  \begin{enumerate}
  \item[(1)] For two appropriately chosen sequences $(m_j)_j,$ $(n_j)_j$
  we set $F_j = \{ \frac{1}{m_{4j-3}^2} \sum_{k\in F}e_k^* : \#
  F\le \frac{n_{4j-3}}{2} \}$ and for the family $(F_j)_j$ we
  construct the norming set $\mc{F}_2$ and the James Tree space
  $JT_{\mc{F}_2}.$

  \item[(2)] The space $JT_{\mc{F}_2}$ does not contain $\ell_1$
  and for
  every  weakly null sequence $(x_n)_n$ in $JT_{\mc{F}_2}$ with $\|
  x_n\| \le C,$ $\lim \|x_n\|_\infty =0$ and every $m\in \N$ there
  exists $L\in [\N]$ such that for every $y^* \in \mc{F}_2$
  \begin{equation}  \# \{n\in L: |y^*(x_n)| \ge \frac{1}{m} \} \le
  66 m^2 C^2. \label{star} \end{equation}

  \item[(3)] We consider the HI extension with attractors
  $D_{\mc{F}_2}^{hi}$ of $\mc{F}_2$ defined by the operations
  $(\mc{A}_{n_j}, \frac{1}{m_j}),$ and we denote by
  $\mathfrak{X}_{\mc{F}_2}$ the space $\overline{(c_{00}, \|\cdot
  \|_{D_{\mc{F}_2}^{hi}})}.$

  \item[(4)] Inequality \eqref{star} yields that
  $\mathfrak{X}_{\mc{F}_2}$ is a strongly strictly singular
  extension of $JT_{\mc{F}_2}.$ Therefore:

  \begin{enumerate}
  \item[(i)] The space $\mathfrak{X}_{\mc{F}_2}$ is HI and
  reflexively saturated.

  \item[(ii)] The predual $(\mathfrak{X}_{\mc{F}_2})_*$ is HI.
  \end{enumerate}

  \item[(5)] Using the attractor functionals, we copy into every
  subspace of $(\mathfrak{X}_{\mc{F}_2})_*$
   a complete dyadic subtree of $(\mc{T}, \prec_\mc{F})$
  which shows that $(\mathfrak{X}_{\mc{F}_2})_*$ is a Hereditarily
  James Tree space (HJT) and hence it does not contain a reflexive
  subspace.
  \end{enumerate}

 Notice that $(\mathfrak{X}_{\mc{F}_2})_*$ shares with the
 space $X$, in the statements presented at the beginning of the
 introduction,  most of the properties stated there. However
 for some of the properties  a variation is required.
  In fact there exists a complete subtree
 $(\mc{T}',\prec_{\sigma_F})$ of $(\mc{T},\prec_{\sigma_F})$ such
 that for the corresponding space
 $\mathfrak{X}_{\mc{F}_2'}$ we have that
 $(\mathfrak{X}_{\mc{F}_2'})^*/(\mathfrak{X}_{\mc{F}_2'})_*=\ell_2(\Gamma)$
 with $\#\Gamma=2^{\omega}$. The space
 $(\mathfrak{X}_{\mc{F}_2'})_*$ coincides with $X$ in the
 aforementioned statements.

 The construction of a nonseparable HI space $Z$
  not containing reflexive subspaces requires changing the framework
   with the operations $(\mc{S}_{n_j}, \frac{1}{m_j})_j$ instead of
 $(\mc{A}_{n_j}, \frac{1}{m_j})_j$.
  In this framework  the set $\mc{F}_s$  and the space
  $JT_{\mc{F}_s}$ are defined. More precisely $\mc{F}_s$ is defined in a
  similar manner as $\mc{F}_2$ based on the families
  \[ F_j = \{ \frac{1}{m_{4j-3}^2} \sum_{i\in I}\pm  e_i^*:\; I \in
  \mc{S}_{n_{4j-3}} \}. \]
  Using the coding $\sigma_{\mc{F}}$, we define the special
  functionals and their indices as in the $\mc{F}_2$ case. Finally
  we set
  \begin{eqnarray*}
   \mc{F}_s& =&  \{ \pm e_n^*:\; n\in\N\}
    \cup \{ \sum_{i=1}^d x_i^* :\; \min
  \supp x_i^*\ge d,\; i=1,\ldots,n,\\
  & &  (\ind (x_i^*))_{i=1}^d \text{ are
    pairwise disjoint} \}
  \end{eqnarray*}

 The HI extension with attractors is defined similarly to the
 $\eqs_{\mc{F}_2}$ case. Then $\mathfrak{X}_{\mc{F}_s}$ is an asymptotic
 $\ell_1$ and reflexively saturated HI space and also
 $(\mathfrak{X}_{\mc{F}_s})_*$ is HI while not containing any reflexive
 subspace. Passing to a complete subtree $(T',
 \prec_{\sigma_{\mc{F}}})$ of $(T, \prec_{\sigma_\mc{F}})$ and to
 the corresponding $\mc{F}_s',$ $ \mathfrak{X}_{\mc{F}_s'},$ we
 obtain the additional property that
 $\mathfrak{X}_{\mc{F}_s'}^*/(\mathfrak{X}_{\mc{F}_s'})_* \cong
 c_0(\Gamma)$ with $\# \Gamma = 2^\omega$. As  is shown in
 \cite{AT1}, this yields that $\mathfrak{X}_{\mc{F}_s'}^*$ is HI
 and since it contains a subspace $((\mathfrak{X}_{\mc{F}_s'})_*)$
 with no reflexive subspace, the space
 $\mathfrak{X}_{\mc{F}_s'}^*$ has the same property.

 Let's mention also that   an HI
 asymptotic $\ell_1$ Banach space $X$ not containing a reflexive
 subspace, with nonseparable dual $X^*$ which is also HI not
 containing any reflexive subspace, has been constructed in \cite{AGT}.
  This space is the analogue of
 $\eqs_{gt}$ in the frame of the operations
 $(\mc{S}_{n_j},\frac{1}{m_j})_j$.

 The last variant we present, concerns the HI space $Y$ with a
 shrinking basis,  not containing a reflexive subspace, such that
 the dual $Y^*$ is    unconditionally and reflexively saturated.

 For this, starting with the set $\mc{F}_2$ we pass to an
 extension only with attractors and additionally we subtract  a
 large portion of the conditional structure of the attractors. This
 permits us to show that the extension space
 $\mathfrak{X}_{\mc{F}_2}^{us}$ is
 unconditionally saturated.
  The remaining part of the conditional structure of the
 attractors, forces the predual $(\mathfrak{X}_{\mc{F}_2}^{us})_*$
 to be HI and not to contain any reflexive subspace.

 The paper is organized as follows. The first two sections are
 devoted to the presentation of the strictly singular and strongly
 strictly singular (HI) extensions with attractors of a ground space $Y_G$. We
 shall denote these as $X_G$. For the results presented in these two sections
  the attractors play no role. Thus all statements remain valid
  whether we consider extensions with attractors or not.
   The strictly singular extension, as
 they have defined before, provide information about $X_G$ and
 $\mc{L}(X_G)$. In particular the following is proven:

  \begin{theorem} If $X_G$ is a strictly singular extension (with or
 without attractors) then
 the natural basis of $X_G$ is boundedly complete,
 the space $X_G$ is HI, reflexively saturated
  and every $T$ in
 $\mc{L}(X_G)$ is of the form $T=\lambda I+ S$ with $S$ a strictly
 singular operator.
 \end{theorem}

  Strongly strictly singular extensions concern the HI
 property of $(X_G)_*$ and the structures of $\mc{L}(X_G)$,
 $\mc{L}((X_G)_*)$. The following theorem includes the main
 results of Section 2.
  \begin{theorem} If $X_G$ is a strongly
  strictly singular extension (with or
 without attractors) then in addition to the above we have the
 following
 \begin{enumerate}
 \item[(i)] The predual $(X_G)_*$ is HI.
 \item[(ii)] Every strictly singular $S$ in
 $\mc{L}(X_G)$ is  weakly compact.
 \item[(iii)] Every $T$ in $\mc{L}((X_G)_*)$ is of the form
 $T=\lambda I+S$ with $S$ being a strictly singular and weakly compact
 operator.
 \end{enumerate}
 \end{theorem}
 The following result concerning the quotients of $X_G$ is also
 proved in Section 2.
  \begin{theorem} If $X_G$ is a strongly
  strictly singular extension (with or
 without attractors) and $Z$ is a $w^*$ closed subspace of $X_G$
 then the quotient $X_G/Z$ is HI.
 \end{theorem}
 The dual form of the above theorem
  is the following. For
 every subspace $W$ of $(X_G)_*$ the dual space $W^*$ is HI.
 Notice also that  the additional assumption that $Z$ is $w^*$ closed can not
 be dropped,  as the results of Section 5 indicate.

  In Section 3  we study  the
 spaces $JT_{\mc{F}_2}$.
  We are mainly concerned with proving the
 aforementioned \eqref{star} yielding that the extension $\eqs_{\mc{F}_2}$
 is a strongly strictly singular one.
 In Section 4 using  attractors we prove that
 $(\eqs_{\mc{F}_2})_*$ is a Hereditarily James Tree (HJT) space
 and hence it does not contain any reflexive subspace.
 Section 5  is devoted to the study of
 $(\eqs_{\mc{F}_2})^*$. It is shown that for every subspace $Y$ of
 $(\eqs_{\mc{F}_2})_*$, the space $\ell_2$ is isomorphic to a subspace of
 the nonseparable space $Y^{**}$. We also describe the definition
 of $\eqs_{\mc{F}_2'}$ which has the additional property that
 $(\eqs_{\mc{F}_2'})^*/(\eqs_{\mc{F}_2'})_*$ is isomorphic to
 $\ell_2(\Gamma)$. Section 6 and Section 7 contain the variants
 $\eqs_{\mc{F}_s}$, $\eqs_{\mc{F}_2}^{us}$ mentioned before. We
 have also included two appendices. In Appendix A we present a
 proof of a form of the basic inequality used for  estimating
 upper bounds for the action of functionals on certain vectors.
 In Appendix B we proceed to a systematic study of the James  Tree
 spaces
 $JT_{\mc{F}_2}$, $JT_{\mc{F}_s}$ and $JT_{\mc{F}_{2,s}}$.
 We actually show that $JT_{\mc{F}_2}$ is $\ell_2$ saturated while
 $JT_{\mc{F}_s}$ and $JT_{\mc{F}_{2,s}}$ are $c_0$ saturated.
  The study of James Tree spaces in Section 3 and Appendix B
  is not related
 to HI techniques and
 uses classical Banach space theory with Ramsey' s
 theorem also playing a key role.

 \section{Strictly singular extensions  with attractors}

In this section we introduce the ground sets $G$ and then we
define the extensions $X_G = T[G, (\mc{A}_{n_j}, \frac{1}{m_j})_j,
\sigma]$ with low complexity saturation methods. Attractors are
also defined. We provide conditions yielding the HI property of
the extension $X_G$ and we study the space of the operators
$\mc{L}(X_G).$ The results and the techniques are analogue to the
corresponding of \cite{AT1} where extensions using higher
complexity saturation methods are presented. We refer the reader
to \cite{ArTo} for an exposition of  low complexity extensions. We
also point out that the attractors in the present and next section
will be completely neutralized. Their role will be revealed in
Section 4 where we study the structure of $(\eqs_{\mc{F}_2})_*$.

  \begin{definition}\label{def5}
 {\bf (ground sets)}
 A set $G\subset\co$ is said to be ground if the following
 conditions are satisfied
 \begin{enumerate}
 \item[(i)] $e_n^*\in G$ for $n=1,2,\ldots$, $G$ is symmetric (i.e.
 if $g\in G$ then $-g\in G$), and closed under restriction of its
 elements to intervals of $\N$ (i.e. if $g\in G$ and $E$ is an
 interval of $\N$ then $Eg\in G$).
 \item[(ii)] $\|g\|_{\infty}\le 1$ for every $g\in G$ and
 $g(n)\in \Q$ for every $g\in G$ and $n\in\N$.
 \item[(iii)] Denoting by   $\|\;\|_G$ the ground norm on \co  defined
 by the rule $\|x\|_G=\sup\{g(x):\; g\in G\}$,  the ground
 space $Y_G$, which  is the completion of $(\co,\|\;\|_G)$
 contains no isomorphic copy of $\ell_1$.
 \end{enumerate}
 \end{definition}

 It follows readily that the standard basis $(e_n)_n$ of $Y_G$ is
 a bimonotone Schauder basis. The converse is also true. Namely if
 $Y$ has a bimonotone basis $(y_n)_n$ and $\ell_1$ does not embed
 into $Y$ then there exists a ground set $G$ such that $Y_G$ is
 isometric to $Y$. Also, as is well known, every space $(Y,\|\;\|)$
 with a basis \seq{e}{n}
 admits an equivalent norm $|||\cdot|||$ such that $\seq{e}{n}$ is
 a bimonotome basis for $(Y,|||\;|||)$.

 \begin{definition}\label{def13}
 {\bf (HI extensions with attractors)}
  We fix two strictly increasing sequences of even positive integers
 $\seq{m}{j}$ and $\seq{n}{j}$ defined as follows:
 \begin{itemize}
 \item  $m_1=2$  and $m_{j+1}=m_j^5$.
 \item $n_1=4$,  and $n_{j+1}=(5n_j)^{s_j}$  where  $s_j=\log_2 m_{j+1}^3$.
 \end{itemize}
 We let  $D_G$   be the minimal
 subset of \co satisfying the following conditions:
 \begin{enumerate}
 \item[(i)] $G\subset D_G$,
 $D_G$ is  symmetric (i.e. if $f\in D_G$ then  $-f\in D_G$) and
  $D_G$ is closed under the restriction of its elements
 to intervals of $\N$ (i.e. if $f\in D_G$ and $E$ is an interval of
 $\N$ then $Ef\in D_G$).
 \item[(ii)] $D_G$ is closed under  $(\mc{A}_{n_{2j}},\frac{1}{m_{2j}})$
 operations, i.e. if $f_1<f_2<\cdots<f_{n_{2j}}$ belong to $D_G$
 then the functional
 $f=\frac{1}{m_{2j}}(f_1+f_2+\cdots+f_{n_{2j}})$   also belongs to
 $D_G$. In this case we say that the functional $f$ is the result of
 an $(\mc{A}_{n_{2j}},\frac{1}{m_{2j}})$ operation.
 \item[(iii)] $D_G$ is closed under  $(\mc{A}_{n_{4j-1}},\frac{1}{m_{4j-1}})$
 operations on special sequences, i.e. for every $n_{4j-1}$ special sequence
 $(f_1,f_2,\ldots,f_{n_{4j-1}})$,   the functional
 $f=\frac{1}{m_{4j-1}}(f_1+f_2+\cdots+f_{n_{4j-1}})$
 belongs to $D_G$.
 In this case we say that $f$ is a result of an
 $(\mc{A}_{n_{4j-1}},\frac{1}{m_{4j-1}})$ operation and that $f$ is
 a {\bf special functional}.
  \item[(iv)]
 $D_G$ is closed under  $(\mc{A}_{n_{4j-3}},\frac{1}{m_{4j-3}})$
 operations on attractor sequences, i.e. for every $4j-3$ attractor sequence
 $(f_1,f_2,\ldots,f_{n_{4j-3}})$,   the functional
 $f=\frac{1}{m_{4j-3}}(f_1+f_2+\cdots+f_{n_{4j-3}})$
 belongs to $D_G$.
 In this case we say that $f$ is a result of an
 $(\mc{A}_{n_{4j-3}},\frac{1}{m_{4j-3}})$ operation and that $f$ is
 an {\bf attractor}.
  \item[(v)] The set $D_G$ is rationally convex.
 \end{enumerate}
   The   space $X_G=T[G,(\mc{A}_{n_j},\frac{1}{m_j})_{j},\sigma]$,
 which is the completion of the space $(\co,\|\;\|_{D_G})$, is called a {\bf strictly
 singular extension with attractors}
 of the space $Y_G$,
 provided that the identity operator $I:X_G\to Y_G$
 is strictly singular.

 The norm satisfies the following implicit formula.
 \begin{eqnarray*}
  \|x\|& =& \max\bigg\{\|x\|_G,\;\sup\limits_j\{\sup \frac{1}{m_{2j}}
 \sum\limits_{i=1}^{n_{2j}}\|E_ix\|\},\\
 &&\sup\{\phi(x):\;\phi\text{
 special functional}\},\;
      \sup\{\phi(x):\;\phi\text{ attractor}\}\bigg\}
 \end{eqnarray*}
 where the inside supremum in the second term is taken over all choices
 $(E_i)_{i=1}^{n_{2j}}$ of successive intervals of $\N$.
  \end{definition}

 We next complete the definition of the norming set $D_G$ by giving
 the precise definition of  special functionals and  attractors.

 From the minimality of $D_G$ it follows that each $f\in D_G$ has one
 of the following forms.
 \begin{enumerate}
 \item[(i)] $f\in G$. We then say that $f$ is of  type 0.
  \item[(ii)] $f=\pm Eh$ where $h$ is a result of an
  $(\mc{A}_{n_{j}},\frac{1}{m_{j}})$ operation and $E$ is an interval.
   In this case we say that $f$ is of  type I. Moreover we
   say that the integer $m_j$ is a  weight of $f$ and we write
   $w(f)=m_j$. We notice that an $f\in D_G$ may have many weights.
 \item[(iii)] $f$ is a rational convex combination of type 0 and
 type I functionals. In this case we say that $f$ is of  type
 II.
 \end{enumerate}

 \begin{definition}\label{def10}
 {\bf ($\sigma$ coding,  special sequences and attractor sequences)}
  Let $\Q_s$ denote the set of all finite
 sequences $(\phi_1,\phi_2,\ldots,\phi_d)$ such that
 $\phi_i\in \co$, $\phi_i\neq 0$ with $\phi_i(n)\in \Q$ for all $i,n$ and
 $\phi_1<\phi_2<\cdots<\phi_d$.
 We fix a pair $\Omega_1,\Omega_2$ of disjoint infinite subsets of
 $\N$.
 From the fact that $\Q_s$ is
 countable we are able to define a Gowers-Maurey type injective
 coding function
 $\sigma:\Q_s\to \{2j:\;j\in\Omega_2\}$ such that
 $m_{\sigma(\phi_1,\phi_2,\ldots,\phi_d)}>\max\{\frac{1}{|\phi_i(e_l)|}:\;l\in\supp
 \phi_i,\;i=1,\ldots,d\}\cdot\max\supp \phi_d$.
 Also, let $(\Lambda_i)_{i\in\N}$ be a sequence of pairwise
 disjoint infinite subsets of $\N$ with $\min \Lambda_i>m_i$.
 \begin{enumerate}
 \item[(A)]
 A finite  sequence $(f_i)_{i=1}^{n_{4j-1}}$ is said to be a
 {\bf $n_{4j-1}$ special sequence}  provided that
 \begin{enumerate}
 \item[(i)]  $(f_1,f_2,\ldots,f_{n_{4j-1}})\in\Q_s$ and $f_i\in D_G$
 for $i=1,2,\ldots,n_{4j-1}$.
 \item[(ii)] $w(f_1)=m_{2k}$ with $k\in\Omega_1$, $m_{2k}^{1/2}>n_{4j-1}$
  and for each $1\le i<n_{4j-1}$,
  $w(f_{i+1})=m_{\sigma(f_1,\ldots,f_{i})}$.
 \end{enumerate}
 \item[(B)]
 A finite  sequence $(f_i)_{i=1}^{n_{4j-3}}$ is said to be a
 {\bf $n_{4j-3}$ attractor sequence}  provided that
 \begin{enumerate}
 \item[(i)]  $(f_1,f_2,\ldots,f_{n_{4j-3}})\in\Q_s$ and $f_i\in D_G$
 for $i=1,2,\ldots,n_{4j-3}$.
 \item[(ii)] $w(f_1)=m_{2k}$ with $k\in\Omega_1$, $m_{2k}^{1/2}>n_{4j-3}$
  and $w(f_{2i+1})=m_{\sigma(f_1,\ldots,f_{2i})}$
 for each $1\le i<n_{4j-3}/2$.
 \item[(iii)]  $f_{2i}=e^*_{l_{2i}}$ for some
  $l_{2i}\in\Lambda_{\sigma(f_1,\ldots,f_{2i-1})}$, for
  $i=1,\ldots,n_{4j-3}/2$.
 \end{enumerate}
 \end{enumerate}
  \end{definition}
 The definition of the special functionals and the attractors completes the
 definition of the norming set $D_G$ of the space $X_G$.

 \begin{remarks}\label{rem6}
 \begin{enumerate}
 \item[(i)] Since the sequence $(\frac{n_{2j}}{m_{2j}})_j$ increases
 to infinity and the norming set $D_G$ is closed in the
 $(\mc{A}_{n_{2j}},\frac{1}{m_{2j}})$ operations we get that the
 Schauder basis \seq{e}{n} of $X_G$ is boundedly complete.
 \item[(ii)] The special
  sequences, as in previous constructions
  (see for example \cite{GM1},\cite{G2},\cite{AD},\cite{AT1}), are responsible for
  the absence of unconditionality in every subspace of $X_G$.
 \item[(iii)] The attractors do not effect the results of the present
 section. Their role is to attract the structure of $G$, for
 certain $G$, in every  subspace of the predual
 $(X_G)_*=\overline{\spann}\{e_n^*:\;n\in\N\}$ of the
 space $X_G$. So, if we discard condition (iv) in the definition
 of the norming set $D_G$ (Definition \ref{def13}) then the
 corresponding space $X_G$, which we call a {\bf strictly singular
 extension} provided that the identity operator $I:X_G\to Y_G$ is strictly
 singular,
 shares all the properties we shall prove in this
 section with those $X_G$'s which are strictly singular extensions with attractors.
 \end{enumerate}
 \end{remarks}

 \begin{definition}\label{def7}
 {\bf (Rapidly increasing sequences)}
 A block sequence $(x_{k})$ in $X_{G}$ is said to be a
 $(C,\varepsilon)$ rapidly
 increasing sequence  (R.I.S.), if $\|x_k\|\le C$,
  and
  there exists a strictly
 increasing sequence $(j_{k})$ of positive integers such that

 (a) $(\max\supp x_{k})\frac{1}{m_{j_{k+1}}}<\varepsilon$.

 (b) For every $k=1,2,\ldots$ and every $f\in D_G$ with
 $w(f)=m_i$, $i<j_{k}$  we have that $|f(x_{k})| \le
 \frac{C}{m_i}$.
 \end{definition}

 \begin{remark}\label{rem2}
 A subsequence of a $(C,\e)$ R.I.S. remains a
 $(C,\e)$ R.I.S. while a sequence which is a $(C,\e)$ R.I.S.
 is also  a $(C',{\e}')$
 R.I.S. if  $C'\ge C$ and ${\e}'\ge \e$.
 \end{remark}

 \begin{proposition}\label{prop13}
 Let $(x_k)_{k=1}^{n_{j_0}}$ be a $(C,\e)$ R.I.S.  with
 $\e\le\frac{2}{m_{j_0}^2}$ such that
 for every $g\in G$, $\#\{k:\;|g(x_k)|>\e\}\le n_{j_0-1}$. Then

1) For every $f\in D_G$ with $w(f)=m_i$,
\[
|f(\frac{1}{n_{j_0}}\sum\limits_{k=1}^{n_{j_0}}x_{k})| \leq
\begin{cases}
\frac{3C}{m_{j_0}m_i}\,,\quad &\text{if}\,\,\, i<j_0\\
\frac{C}{n_{j_0}}+\frac{C}{m_i}+C\e\,,&\text{if}\,\,\, i\ge j_0
\end{cases}\]
 In particular
$\|\frac{1}{n_{j_0}}\sum\limits_{k=1}^{n_{j_0}}x_{k}\|
\leq\frac{2C}{m_{j_0}}$.

 2) If $(b_k)_{k=1}^{n_{j_0}}$ are scalars with $|b_k|\le 1$  such
 that
 \begin{equation}\label{eq16}
 |h(\sum\limits_{k\in E} b_{k}x_{k})| \leq C(\max_{k\in E}
 |b_{k}|+\varepsilon\sum\limits_{k\in E}|b_{k}|)
 \end{equation}
 for every interval $E$ of positive integers
 and every $h\in D_G$ with $w(h)=m_{j_0}$, then
 \[
 \|\frac{1}{n_{j_0}}\sum\limits_{k=1}^{n_{j_0}}b_{k}x_{k}\|
 \leq\frac{4C}{m_{j_0}^{2}}. \]
 \end{proposition}
 The proof of the above proposition
 is based on what we call the basic inequality
 (see also \cite{AT1}, \cite{ALT}).
 Its proof is presented in  Appendix A.

  \begin{remark}
 The validity of  Proposition \ref{prop13} is independent of the
 assumption that the operator $I:X_G\to Y_G$ is strictly singular.
  \end{remark}

 In the present section we shall prove several properties of the
 space $X_G$ provided that the space $X_G$ is a strictly singular
 extension of $Y_G$.

 \begin{definition}\label{def19}
{\bf (exact pairs)}
  A pair $(x,\phi)$ with $x\in X_{G}$
 and $\phi\in D_G$ is said to be a $(C,j,\theta)$ exact pair
 (where $C\ge 1$,  $j\in\N$, $0\le\theta\le 1$)
 if the  following conditions are satisfied:
 \begin{enumerate}
 \item[(i)] $1\le \|x\|\le C$, for every
 $\psi\in D_G$ with $w(\psi)=m_i$, $i\neq j$ we have
 that $|\psi(x)|\le\frac{2C}{m_i}$ if $i<j$, while $|\psi(x)|\le\frac{C}{m_j^2}$
 if $i>j$ and $\|x\|_{\infty}\le \frac{C}{m_j^2}$.
 \item[(ii)] $\phi$ is of type I with $w(\phi)=m_j$.
 \item[(iii)] $\phi(x)=\theta$
  and $\ran x=\ran \phi$.
 \end{enumerate}
 \end{definition}

 \begin{definition}
  {\bf ($\ell_1^k$ averages)}
  Let $k\in\mathbb{N}$. A  finitely supported vector $x\in
  X_{G}$    is said to be a
  $C-\ell_{1}^{k}$ average if $\|x\|> 1$ and there
  exist $x_{1}<\ldots<x_{k}$ with $\norm[x_{i}]\leq C$ such that
  $x=\frac{1}{k}\sum\limits_{i=1}^{k}x_{i}$.
 \end{definition}

 \begin{lemma}\label{lem10}
 Let $j\in\N$ and $\e>0$.
 Then every block subspace of
 $X_{G}$ contains a vector $x$ which is a $2-
 \ell_1^{n_{2j}}$ average.
 If $X_G$ is a strictly singular extension (with or without attractors)
  then we may
 select $x$ to additionally  satisfy $\|x\|_G<\e$.
 \end{lemma}
 \begin{proof}[\bf Proof.]
 If the identity operator $I:X_G\to Y_G$ is  strictly
 singular we may pass to a further block subspace on
 which the restriction of  $I$ has norm
 less than $\frac{\e}{2}$.
For the remainder of the proof in this case, and the proof in the
general case, we refer to
 \cite{GM1} (Lemma 3) or to \cite{AM}.
 \end{proof}

 \begin{lemma}\label{lem9}
 Let   $x$ be a
  $C-\ell_{1}^{k}$ average.  Then for
 every $n\le k$ and
  every sequence of intervals
  $E_{1}<\ldots<E_{n}$,
 we have that $ \sum\limits_{i=1}^{n}\norm[E_{i}x]\leq
 C(1+\frac{2n}{k})$. In particular
 if $x$ is a $C-\ell_{1}^{n_j}$ average then for every
 $f\in D_G$ with $w(f)=m_i$, $i<j$ then $|f(x)|\le
 \frac{1}{m_{i}}C(1+\frac{2n_{j-1}}{n_j})\le
 \frac{3C}{2}\frac{1}{m_i}$.
  \end{lemma}
 We refer to \cite{S} or to
  \cite{GM1}   (Lemma 4)
   for a proof.

  \begin{remark}\label{rem1}
 Let $(x_k)_k$ be  a block sequence in $X_{G}$ such that each
 $x_k$ is a $\frac{2C}{3}-\ell_1^{n_{j_k}}$ average and
 let $\varepsilon>0$ be such that
 $\#(\ran(x_{k}))\frac{1}{m_{j_{k+1}}}<\varepsilon$.
  Then Lemma \ref{lem9}
 yields that condition (b) in the  definition of R.I.S.
 (Definition \ref{def7})
 is also
 satisfied hence $(x_k)_k$ is a $(C,\varepsilon)$ R.I.S. In this case we
 shall call $(x_k)_k$ a $(C,\varepsilon)$  R.I.S. of
 $\ell_1$  averages. From this observation and Lemma \ref{lem10}
 it follows  that if $X_G$ is a strictly singular extension
 of $Y_G$ then for every $\e>0$, every block subspace of $X_{G}$
 contains a $(3,\e)$ R.I.S. of
 $\ell_1$ averages \seq{x}{k} with $\|x_k\|_G<\e$.
  \end{remark}

  \begin{proposition}\label{prop20}
   Suppose that $X_G$ is a strictly singular
  extension of $Y_G$ (with or without attractors).
 Let $Z$ be a block subspace of $X_{G}$, let $j\in\N$ and let $\e>0$.
 Then there exists a $(6,2j,1)$ exact pair $(x,\phi)$ with $x\in Z$
 and $\|x\|_G<\e$.
 \end{proposition}
 \begin{proof}[\bf Proof.] From the fact that the identity
 operator $I:X_G\to Y_G$ is strictly singular
  we may
 assume, passing to a block subspace of $Z$, that
 $\|z\|_{G}< \frac{\e}{6}\|z\|$ for every $z\in Z$.
 We choose a $(3,\frac{1}{n_{2j}})$ R.I.S. of $\ell_1$ averages in
 $Z$,   $(x_k)_{k=1}^{n_{2j}}$. For $k=1,2,\ldots,n_{2j}$ we
 choose $\phi_k\in D_G$ with $\ran\phi_k=\ran x_k$ such that
 $\phi_k(x_k)>1$.
 We set $\phi=\frac{1}{m_{2j}}\sum\limits_{k=1}^{n_{2j}}\phi_k$.
 We have that
 $\eta=\phi(\frac{m_{2j}}{n_{2j}}\sum\limits_{k=1}^{n_{2j}}x_k)>1$.
 On the other hand Proposition \ref{prop13}
 yields that
 $\|\frac{m_{2j}}{n_{2j}}\sum\limits_{k=1}^{n_{2j}}x_k\|\le 6$.
 We set
 $x=\frac{1}{\eta}\frac{m_{2j}}{n_{2j}}\sum\limits_{k=1}^{n_{2j}}x_k$.

 We have that $1=\phi(x)\le\|x\|\le 6$,  hence also
 $\|x\|_{G}\le \e$, while $\ran\phi=\ran x$.
 From Proposition \ref{prop13} it follows that for every $\psi\in
 D_G$ with $w(\psi)=m_{i}$, $i\neq 2j$ we have that
 $|\psi(x)|\le \frac{9}{m_i}$ if $i<2j$ while
 $|\psi(x)|\le m_{2j}(\frac{3}{n_{2j}}
 +\frac{3}{m_i}+\frac{3}{n_{2j}})\le
 \frac{1}{m_{2j}^2}$ if $i>m_{2j}$.
 Finally $\|x\|_{\infty}\le
 \frac{m_{2j}}{n_{2j}}\max\limits_k\|x_k\|_{\infty}\le
 \frac{3m_{2j}}{n_{2j}}< \frac{1}{m_{2j}^2}$.

 Therefore $(x,\phi)$ is a $(6,2j,1)$ exact pair with $x\in Z$ and $\|x\|_G<\e$.
 \end{proof}

 \begin{definition}\label{depseq}
 {\bf (dependent sequences and attracting sequences)}
 \begin{enumerate}
  \item[(A)]
 A double sequence $(x_k,x_k^*)_{k=1}^{n_{4j-1}}$ is said to be a
 $(C,4j-1,\theta)$ dependent sequence (for $C>1$,
 $j\in\mathbb{N}$,     and
 $0\le\theta\le 1$) if there exists a
 sequence $(2j_k)_{k=1}^{n_{4j-1}}$ of even integers such that the
 following conditions are fulfilled:
 \begin{enumerate}
 \item[(i)] $(x^*_k)_{k=1}^{n_{4j-1}}$ is a $4j-1$ special
 sequence with
 $w(x^*_{k})=m_{2j_{k}}$ for each $k$.
 \item[(ii)] Each $(x_{k},x_{k}^*)$ is a
 $(C,2j_{k},\theta)$ exact pair.
 \end{enumerate}
 \item[(B)]
 A double sequence $(x_k,x_k^*)_{k=1}^{n_{4j-3}}$ is said to be a
 $(C,4j-3,\theta)$ attracting sequence (for $C>1$,
 $j\in\mathbb{N}$,     and
 $0\le\theta\le 1$) if there exists a
 sequence $(2j_k)_{k=1}^{n_{4j-3}}$ of even integers such that the
 following conditions are fulfilled:
 \begin{enumerate}
 \item[(i)] $(x^*_k)_{k=1}^{n_{4j-3}}$ is a $4j-3$ attractor sequence
 with  $w(x^*_{2k-1})=m_{2j_{2k-1}}$ and
 $x^*_{2k}=e^*_{l_{2k}}$ where $l_{2k}\in \Lambda_{2j_{2k}}$
  for all $k\le n_{4j-3}/{2}$.
 \item[(ii)] $x_{2k}=e_{l_{2k}}$.
 \item[(iii)] Each $(x_{2k-1},x_{2k-1}^*)$ is a
 $(C,2j_{2k-1},\theta)$ exact pair.
 \end{enumerate}
 \end{enumerate}
 \end{definition}

 \begin{remark}\label{rem3} If $(x_k,x_k^*)_{k=1}^{n_{4j-1}}$ is
 a  $(C,4j-1,\theta)$ dependent sequence
 (resp. $(x_k,x_k^*)_{k=1}^{n_{4j-3}}$ is a $(C,4j-1,\theta)$ attracting sequence)
  then the sequence
 $(x_k)_k$ is  a $(2C,\frac{1}{n_{4j-1}^2})$
 R.I.S. (resp. a $(2C,\frac{1}{n_{4j-3}^2})$ R.I.S.).
 Let examine this for a $(C,4j-3,\theta)$ attracting sequence (the
 proof for a dependent sequence is similar).
 First $\|x_k\|=\|e_{l_k}\|=1$ if $k$ is even
 while $\|x_k\|\le C$ if $k$ is odd,
 as follows from the fact that $(x_k,x_k^*)$ is a
 $(C,2j_{k},\theta)$ exact pair.

   Second,   the growth condition
 of the coding function $\sigma$ in Definition \ref{def10}
 and condition (ii) in the same definition yield
 that     for each $k$ we have
 that $(\max\supp x_k)\frac{1}{m_{2j_{k+1}}}=\max\supp
 x_k^*\cdot\frac{1}{m_{\sigma(x_1^*,\ldots,x_k^*)}}$\\
 $<\min\{|x_i^*(e_l)|:\;l\in\supp x_i^*,\;i=1,\ldots,k\}
 \le \frac{1}{m_{2j_1}}<\frac{1}{n_{4j-3}^2}$.

 Finally, if $f\in D_G$ with $w(f)=m_i$, $i<2j_k$ then
 $|f(x_k)|=|f(e_{l_k})|\le \|f\|_{\infty}\le\frac{1}{m_i}$ if
 $k$ is even, while
 $|f(x_k)|\le\frac{2C}{m_i}$ if $k$ is odd,
 since in this case $(x_k,x_k^*)$ is a
 $(C,2j_{k},\theta)$ exact pair.
 \end{remark}

 \begin{proposition} \label{prop11}
 \begin{enumerate}
 \item[(i)]
 Let $(x_k,x_k^*)_{k=1}^{n_{4j-1}}$ be a $(C,4j-1,\theta)$ dependent
 sequence such that  $\|x_{k}\|_G\le \frac{2}{m_{4j-1}^2}$
 for $1\le k\le n_{4j-1}$. Then
  we have that
 \[\|\frac{1}{n_{4j-1}}\sum\limits_{k=1}^{n_{4j-1}}(-1)^{k+1}x_{k}\|
 \le\frac{8C}{m_{4j-1}^2}.\]
 \item[(ii)] If $(x_k,x_k^*)_{k=1}^{n_{4j-3}}$ is a
 $(C,4j-3,\theta)$ attracting  sequence with
 $\|x_{2k-1}\|_G\le \frac{2}{m_{4j-3}^2}$
 for $1\le k\le n_{4j-3}/2$ and for  every $g\in G$ we
 have that $\#\{k:\;|g(x_{2k})|>\frac{2}{m_{4j-3}^2}\}\le
 n_{4j-4}$ then
 \[\|\frac{1}{n_{4j-3}}\sum\limits_{k=1}^{n_{4j-3}}(-1)^{k+1}x_{k}\|
 \le\frac{8C}{m_{4j-3}^2}.\]
 \item[(iii)]
 If $(x_k,x_k^*)_{k=1}^{n_{4j-1}}$ is a $(C,4j-1,0)$ dependent
 sequence with $\|x_{k}\|_G\le \frac{2}{m_{4j-1}^2}$
 for $1\le k\le n_{4j-1}$ then
  we have that
 \[\|\frac{1}{n_{4j-1}}\sum\limits_{k=1}^{n_{4j-1}}x_{k}\|
 \le\frac{8C}{m_{4j-1}^2}.\]

 \end{enumerate}
 \end{proposition}
 \begin{proof}[\bf Proof.]
 The conclusion  will follow from Proposition \ref{prop13} 2) after
 showing that the required conditions are fulfilled.
 We shall only show (i); the proof of (ii) and (iii) is similar.

 From the previous remark the sequence $(x_k)_{k=1}^{n_{4j-1}}$ is
 a  $(2C,\frac{1}{n_{4j-1}^2})$  R.I.S.
 hence
 it  is a $(2C,\frac{2}{m_{4j-1}^2})$ R.I.S.
 (see Remark \ref{rem2}).
 It remains to show that for $f\in D_G$ with $w(f)=m_{4j-1}$ and for every
 interval $E$ of positive integers we have that
 \[  |f\big(\sum\limits_{k\in E}(-1)^{k+1}x_{k}\big)|
 \le 2C(1+\frac{2}{m_{4j-1}^2}\#(E)).\]

 Such an $f$ is of the form
 $f=\frac{1}{m_{4j-1}}
 (Fx_{t-1}^*+x_t^*+\cdots+x_{r}^*+f_{r+1}+\cdots+f_d)$
 for some $4j-1$ special sequence
 $(x_1^*,x_2^*,\ldots,x_{r}^*,f_{r+1},\ldots,f_{n_{4j-1}})$
  where $x^*_{r+1}\neq f_{r+1}$ with
 $w(x^*_{r+1})=w(f_{r+1})$,  $d\le n_{4j-1}$ and $F$ is an interval of
 the  form $[m,\max\supp x_{t-1}^*]$.

 We estimate the value $f(x_{k})$ for each $k$.
 \begin{itemize}
 \item If $k<t-1$ we have that $f(x_k)=0$.
 \item If $k=t-1$ we get $|f(x_{t-1})|=\frac{1}{m_{4j-1}}|Fx^*_{t-1}(x_{t-1})|
        \le \frac{1}{m_{4j-1}}\|x_{t-1}\|  \le \frac{C}{m_{4j-1}}.$
 \item If $k\in\{t,\ldots,r\}$ we have that
   $f(x_{k})=\frac{1}{m_{4j-1}}x^*_{k}(x_k)=\frac{\theta}{m_{4j-1}}$.
 \item If $k> r+1$, then the injectivity of the coding function
 $\sigma$ and the definition of special functionals yield that
 $w(f_{i})\neq m_{2j_{k}}$ for all
 $i\ge r+1$. Using the fact that
 $(x_{k},x_{k}^*)$ is a  $(C,2j_{k},\theta)$ exact pair
 and taking into account that $n_{4j-1}^2<m_{2j_1}\le \sqrt{m_{2j_{k}}}$
   we get that
 \begin{eqnarray*}
 |f(x_{k})|&= & \frac{1}{m_{4j-1}}|(f_{r+1}+\ldots+f_d)(x_{k})| \\
 &\le&   \frac{1}{m_{4j-1}}
 \Big(\sum\limits_{w(f_{i})<m_{2j_{k}}}|f_{i}(x_{k})|
      + \sum\limits_{w(f_{i})>m_{2j_{k}}}|f_{i}(x_{k})|   \Big)\\
 &\le &  \frac{1}{m_{4j-1}}
 \Big(\sum\limits_{4j-1<l<2j_{k}}\frac{2C}{m_l}+
  n_{4j-1}\frac{C}{m_{2j_{k}}^2}  \Big)\\
 &\le & \frac{C}{m_{4j-1}^2}
 \end{eqnarray*}
  \item For $k=r+1$, using a similar argument to the previous case
  we get   that  $|f(x_{r+1})|\le \frac{C}{m_{4j-1}}
                 +\frac{C}{m_{4j-1}^2}<\frac{2C}{m_{4j-1}}$.
 \end{itemize}

 Let  $E$ be an interval. From the previous estimates we get
 that
 \begin{eqnarray*}
 |f\big(\sum\limits_{k\in E}(-1)^{k+1}x_{k}\big)|
   &\le&  |f(x_{t-1})|+
 \big|\sum\limits_{k\in E \cap[t,r]}
 \frac{\theta}{m_{4j-1}}(-1)^{k+1}\big|\\
 & & +|f(x_{r+1})|+
 \sum\limits_{k\in E \cap(r+1,n_{4j-1}]}|f(x_{k})|\\
 &\le&
 \frac{C}{m_{4j-1}}+\frac{1}{m_{4j-1}}+\frac{C+1}{m_{4j-1}}+\frac{C}{m_{4j-1}^2}\#(E)\\
 &&< 2C(1+\frac{2}{m_{4j-1}^2}\#(E)).
  \end{eqnarray*}
 The proof of the proposition is complete.
 \end{proof}

 \begin{theorem}\label{theo1} If the space $X_G$ is a strictly
 singular extension, (with or without attractors)
 then it is Hereditarily
 Indecomposable.
 \end{theorem}
 \begin{proof}[\bf Proof.]
  Let $Y$ and $Z$ be a pair of block subspaces of $X_{G}$ and let
 $\delta>0$. We choose $j\in\N$ with $m_{4j-1}>
 \frac{48}{\delta}$. Using Proposition
 \ref{prop20}   we inductively construct a
 $(6,4j-1,1)$ dependent sequence $(x_k,x_k^*)_{k=1}^{n_{4j-1}}$ with
 $x_{2k-1}\in Y$, $x_{2k}\in Z$ and
 $\|x_{k}\|_G<\frac{1}{m_{4j-1}^2}$ for all  $k$.
 From Proposition \ref{prop11} (i) we get that
 $\|\frac{1}{n_{4j-1}}\sum\limits_{k=1}^{n_{4j-1}}(-1)^{k+1}x_{k}\|
 \le\frac{48}{m_{4j-1}^2}$. On the other hand the functional
 $x^*=\frac{1}{m_{4j-1}}\sum\limits_{k=1}^{n_{4j-1}}x_k^*$ belongs
 to $D_G$ and the estimate of $x^*$ on the vector
 $\frac{1}{n_{4j-1}}
 \sum\limits_{k=1}^{n_{4j-1}}x_{k}$ gives that
 $\|\frac{1}{n_{4j-1}}
 \sum\limits_{k=1}^{n_{4j-1}}x_{k}\|\ge \frac{1}{m_{4j-1}}$.

 Setting
  $y=\sum\limits_{k=1}^{n_{4j-1}/2}x_{2k-1}$
 and $z=\sum\limits_{k=1}^{n_{4j-1}/4}x_{2k-1}$ we
 get that $y\in Y$, $z\in Z$ and $\|y-z\|<\delta\|y+z\|$.
 Therefore the space $X_{G}$ is Hereditarily Indecomposable.
 \end{proof}

 \begin{proposition}   \label{prop14}
 If $X_G$ is a strictly singular extension (with or without
 attractors) then
 the dual $X_G^*$ of the space
 $X_G=T[G,(\mc{A}_{n_j},\frac{1}{m_j})_{j},\sigma]$ is the norm closed
 linear span of the $w^*$ closure of $G$.
 \[   X_G^*=\overline{\spann}(\overline{G}^{w^*}).\]
 \end{proposition}
 \begin{proof}[\bf Proof.] Assume the contrary. Then setting
 $Z=\overline{\spann}(\overline{G}^{w^*})$ there exist $x^*\in
 X_G^*\setminus Z$ with $\|x^*\|=1$
  and $x^{**}\in X_G^{**}$ such that
 $Z\subset\ker x^{**}$,
  $\|x^{**}\|=2$ and $x^{**}(x^*)=2$.
 The space $X_G$ contains no isomorphic copy of $\ell_1$, since
 $X_G$ is a HI space, thus from the Odell-Rosenthal theorem there
 exist a sequence \seq{x}{k} in $X_G$ with $\|x_k\|\le 2$ such that
 $x_k\stackrel{w^*}{\longrightarrow}x^{**}$.
 Since each $e_n^*$ belongs to $Z$
 we get that $\lim\limits_ke_n^*(x_k)=0$ for all $n$, thus,
 using a  sliding hump argument,
 we may assume that \seq{x}{k} is a block sequence.
 Since also $x^*(x_k)\to x^{**}(x^*)=2$
 we may also assume that  $1<x^*(x_k)$ for all $k$.
 Let's observe that every convex combination of \seq{x}{k} has
 norm greater than 1.

 Considering each $x_k$ as a continuous function
 $x_k:\overline{G}^{w^*}\to\R$ we have that the sequence
 \seq{x}{k} is uniformly bounded and tends pointwise to 0,
 hence it is a weakly null sequence in $C(\overline{G}^{w^*})$.
 Since $Y_G$ is isometric to a subspace of $C(\overline{G}^{w^*})$
 we get that $x_k \stackrel{w}{\longrightarrow}0$ in $Y_G$ thus
 there exists a convex block sequence \seq{y}{k} of \seq{x}{k}
 with $\|y_k\|_G\to 0$. We may thus assume that
 $\|y_k\|_G<\frac{\e}{2}$ for all $k$, where $\e=\frac{1}{n_4}$.
 We may construct a block sequence \seq{z}{k} of \seq{y}{k} such
 that \seq{z}{k} is a $(3,\e)$ R.I.S. of $\ell_1$ averages while
 each $z_k$ is an average of \seq{y}{k} with $\|z_k\|_G<\e$
 (see  Remark \ref{rem1}).
 Proposition \ref{prop13}
 yields that the vector
 $z=\frac{1}{n_4}\sum\limits_{k=1}^{n_4}z_k$
 satisfies $\|z\| \le  \frac{2\cdot 3}{m_4}<1$.
 On the other hand,   the vector $z$, being a convex combination
 of \seq{x}{k}, satisfies $\|z\|>1$.
  This
  contradiction completes the proof of the proposition.
  \end{proof}
 \begin{remark} The content of the above proposition is that the
 strictly singular   extension
(with or without attractors)
 $X_G=T[G,(\mc{A}_{n_j},\frac{1}{m_j})_{j},\sigma]$ of the space $Y_G$
 is actually a
 {\bf reflexive extension}. Namely if $\overline{G}^{w^*}$ is a
 subset of \co then a consequence of Proposition \ref{prop14}
 is that the space $X_G$ is reflexive. Furthermore, if $X_G$ is
 nonreflexive then the quotient space $X_G^*/(X_G)_*$ is norm
 generated by the classes of the elements of the set
 $\overline{G}^{w^*}$. Related to this is also the next.
 \end{remark}

 \begin{proposition}\label{prop15}
 The   strictly singular extension (with or without
 attractors)
 $X_G$
   is reflexively saturated (or somewhat reflexive).
 \end{proposition}
 \begin{proof}[\bf Proof.]
 Let $Z$ be a block subspace of $X_G$. From the fact that the
 identity operator $I:X_G\to Y_G$ is strictly singular we may
 choose a normalized block sequence \seq{z}{n} in $Z$,
 with $\sum\limits_{n=1}^{\infty}\|z_n\|_G<\frac{1}{2}$. We claim
 that the space $Z'=\overline{\spann}\{z_n:\;n\in\N\}$ is a
 reflexive subspace of $Z$.

 It is enough to show that the Schauder basis \seq{z}{n} of $Z'$
 is boundedly complete and shrinking. The first follows from the
 fact that \seq{z}{n} is a block sequence of the boundedly
 complete basis \seq{e}{n} of $X_G$. To see that \seq{z}{n} is
 shrinking it is enough to show that
 $\|f|_{\overline{\spann}\{z_i:\;i\ge
 n\}}\|\stackrel{n\to\infty}{\longrightarrow}0$ for every $f\in
 X_G^*$. From Proposition \ref{prop14} it is enough to prove it
 only for $f\in \overline{G}^{w^*}$. Since
 $\sum\limits_{n=1}^{\infty}\|z_n\|_G<\frac{1}{2}$ the conclusion
 follows.
 \end{proof}

 \begin{proposition}\label{prop16} Let $Y$ be an infinite
 dimensional closed subspace of $X_G$. Every bounded linear
 operator $T:Y\to X_G$ takes the form $T=\lambda I_Y+S$
  with $\lambda\in\R$ and $S$ a strictly singular
 operator ($I_Y$ denotes the inclusion map from $Y$ to $X_G$).
 \end{proposition}

 The proof of Proposition \ref{prop16}
 is similar to the corresponding result for the
 space of Gowers and Maurey (Lemmas 22 and 23 of \cite{GM1})
 and is based on the following lemma.

 \begin{lemma}\label{lem11}
 Let $Y$ be a  subspace of $X_G$ and let $T:Y\to X_G$ be a
 bounded linear operator.
 Let \seq{y}{l} be a block
 sequence of $2-\ell_1^{n_l}$ averages with increasing
 lengths in $Y$
 such that $(Ty_l)_{l\in\N}$ is also a block sequence and
 $\lim\limits_l\|y_l\|_G=0$.
 Then $\lim\limits_l \dist(Ty_l,\R y_l)=0$.
 \end{lemma}
 \begin{proof}[\bf Proof of Proposition \ref{prop16}]
 Assume that $T$ is not strictly singular. We shall determine a
 $\lambda\neq 0$ such that $T-\lambda I_Y$ is strictly singular.

 Let $Y'$ be an infinite dimensional closed subspace of $Y$ such
 that $T:Y'\to T(Y')$ is an isomorphism. By standard perturbation
 arguments and using the fact that $X_G$ is a strictly singular
 extension of $Y_G$, we may assume, passing to a subspace, that
 $Y'$ is a block subspace of $X_G$ spanned by a normalized block sequence
 $(y'_n)_{n\in\N}$ such that
 $(Ty'_n)_{n\in\N}$ is also a block sequence and
 $\sum\limits_{n=1}^{\infty}\|y'_n\|_G<1$.
 From  Lemma  \ref{lem10} we may choose a block sequence \seq{y}{n}
 of $2-\ell_1^{n_i}$  averages of increasing lengths in
 $\spann\{y'_n:\;n\in\N\}$ with
 $\|y_n\|_G\to 0$. Lemma \ref{lem11} yields that
 $\lim\limits_n \dist(Ty_n,\R y_n)=0$. Thus there exists a
 $\lambda\neq 0$ such that $\lim\limits_n\|Ty_n-\lambda y_n\|=0$.

 Since the restriction of $T-\lambda I_Y$ to any finite
 codimensional subspace of $\overline{\spann}\{y_n:\;n\in\N\}$ is
 clearly not an isomorphism and since also $Y$ is a HI space,
 it follows from Proposition 1.2 of \cite{AT1} that the operator
 $T-\lambda I_Y$ is strictly singular.
  \end{proof}

  \section{Strongly strictly singular extensions}
  It is not known whether the
   predual $(X_G)_*=\overline{\spann}\{e_n^*:\;n\in\N\}$ of the
 strictly singular extension $X_G$
  is in general a Hereditarily
 Indecomposable space. In this section we
 introduce the concept of strongly strictly singular extensions
 which permit us to ensure  the HI property for the space $(X_G)_*$ and
 to obtain additional information for this space  as well as for
 the spaces $\mc{L}(X_G),$ $\mc{L}((X_G)_*).$
 We also study the quotients of $X_G$ with $w^*$ closed subspaces
 $Z$  and we show that these quotients are HI.

\begin{definition} \label{aa11} Let $G$ be a ground set and $X_G$ be an
extension of the space $Y_G.$ The space $X_G$ is said to be a
strongly strictly singular extension provided the following
property holds:

For every $C>0$ there exists $j(C)\in \N$ such that for every
$j\ge j(C)$ and every $C$-bounded block sequence $\seq{x}{n}$ in
$X_G$ with $\|x_n\|_\infty\to 0$ and $\seq{x}{n}$ a weakly null
sequence in $Y_G,$ there exists $L\in[\N]$ such that for every
$g\in G$
\[\#\{n\in L:\; |g(x_n)|> \frac{2}{m_{2j}^2}\}\le n_{2j-1}.\]
\end{definition}

\begin{remark}\label{rem12} Let's observe, for later use, that if $\seq{x}{n}$ is
a R.I.S. of $\ell_1$ averages (Remark \ref{rem1}), then
$\|x_n\|_\infty \to 0.$ Therefore if $X_G$ is a strongly singular
extension of $Y_G,$ there exists a subsequence
$(x_{l_n})_{n\in\N}$ such that the sequence $y_n= x_{l_{2n-1}}-
x_{l_{2n}}$ is weakly null and satisfies the above stated
property.
\end{remark}

\begin{proposition}\label{propsss}
 If $X_G$ is a strongly strictly singular extension (with or
without attractors) of $Y_G,$ then the identity map $I: X_G\to
Y_G$ is a strictly singular operator.
\end{proposition}
\begin{proof}[\bf Proof.] In any block subspace of $X_G$ we may
consider a block sequence $(x_n)_{n\in\N}$ with $1\le
\|x_n\|_{X_G}\le 2,$ $\|x_n\|_\infty \to 0$ and $(x_n)_{n\in\N}$
being weakly null. Passing to a subsequence $(x_{l_n})_{n\in\N}$
for $j\ge j(2)$ we obtain that
\[ \|\frac{1}{n_{2j}} \sum_{i=1}^{n_{2j}} x_{l_i}\|_{X_G} \ge
\frac{1}{m_{2j}} \] and on the other hand \[ \|\frac{1}{n_{2j}}
\sum_{i=1}^{n_{2j}} x_{l_i}\|_{Y_G} \le \frac{2}{m_{2j}^2} +
\frac{2n_{2j-1}}{n_{2j}}<\frac{3}{m_{2j}^2}\] which yields that
$I$ is not an isomorphism in any block subspace of $X_G.$
\end{proof}

 Definition \ref{aa11} is the analogue of the definition of
 $\mc{S}_2$ bounded or $\mc{S}_{\xi}$ bounded sets (see
 \cite{AT1} where the norming sets are defined with the use
 of saturation methods
 of the form  $(\mc{S}_{\xi_j},\frac{1}{m_j})_j$)
  in the context of saturation methods of low
 complexity, i.e. of the form $(\mc{A}_{n_j},\frac{1}{m_j})_j$.
 As we have noticed earlier the assumption of a strongly strictly
 singular extension (with or without attractors) is required in
 order to prove that the predual space $(X_G)_*$ is Hereditarily
 Indecomposable.

 The HI property of the dual space $X_G^*$ essentially depends
 on the internal structure of the set $G$. Thus we shall see
 examples of strongly strictly singular extensions (with or
 without attractors) such that $X_G^*$ is either HI or contains
 $\ell_2(\N)$.

    \begin{definition}\label{def11}
 {\bf ($c_0^k$ vectors)}
 Let $k\in\N$. A finitely supported vector $x^*\in (X_G)_*$ is said to be a
 $C- c_{0}^k$  vector if there exist
 $x_1^*<\cdots<x_k^*$ such that $\|x_i^*\|>C^{-1}$,
 $x^*=x_1^*+\cdots+x_k^*$ and $\|x^*\|\le 1$.
 \end{definition}

 \begin{remark}\label{rem4}  The fact that the norming set $D_G$
 is rationally convex yields that $D_G$ is pointwise dense in
 the unit ball  of the space
 $B_{X_G^*}$.
 Since also the norming set $D_G$
  is closed in  $(\mc{A}_{n_{2j}},\frac{1}{m_{2j}})$
 operations we get that for every $j$, if
 $f_1,f_2,\ldots,f_{n_{2j}}$ is a block sequence in
 $X_G^*$ with $\|f_i\|\le 1$ then
 $\|\frac{1}{m_{2j}}\sum\limits_{i=1}^{n_{2j}}f_i\|\le 1$.
 \end{remark}

 \begin{lemma}\label{lem16}
 Let $(x_\ell^*)_{\ell\in\N}$ be a block sequence in $(X_G)_*$.
 Then for every $k$ there exists a
 $x^*\in\spann\{x_\ell^*:\;\ell\in\N\}$ which is a $2-c_0^k$ vector.
 \end{lemma}
 The proof is based on Remark \ref{rem4} and
 can be found in \cite{AT2} (Lemma 5.4).

 \begin{lemma}\label{lem17}
 For every  $2-c_0^k$ vector $x^*$
 and every $\e>0$ there exists a $2-c_0^k$ vector $f$ with $f\in
 D_G$, $\ran f=\ran x^*$ and  $\|x^*-f\|<\e$.
 \end{lemma}
 \begin{proof}[\bf Proof.]
 This follows from the fact that the norming set $D_G$ is pointwise
 dense in $B_{X_G^*}$.
 \end{proof}

 \begin{lemma}\label{lem18} If $x^*$ is a $C-c_0^k$ vector
 then there exists a
 $C-\ell_1^k$ average $x$ with $\ran(x)=\ran(x^*)$ and $x^*(x)>1$.
 \end{lemma}
 \begin{proof}[\bf Proof.]
 Let $x^*=\sum\limits_{i=1}^kx_i^*$ where
 $x_1^*<\cdots<x_k^*$, $\|x_i^*\|> C^{-1}$ and $\|x^*\|\le 1$.
 For $i=1,\ldots,k$ we choose $x_i\in X_G$ with $\|x_i\|\le
 1$,  $x_i^*(x_i)>C^{-1}$ and $\ran(x_i)= \ran(x_i^*)$.
 We set $x=\frac{1}{k}\sum\limits_{i=1}^k(Cx_i)$.
 Then $\|Cx_i\|\le C$ for $i=1,\ldots,k$, while
 $\|x\|\ge x^*(x)>1$. Also, since
  $\ran(x)=\ran(x^*)$,  $x$ is the desired $C-\ell_1^k$
 average.
 \end{proof}

 \begin{proposition}\label{prop29}
 Let $Z$ be a block subspace of $(X_G)_*$ and let $k\in\N$,
 $\e>0$.
 Then there exists  a $2-\ell_1^k$ vector $y$ and $y^*\in D_G$ such
 that $y^*(y)>1$, $\ran(y^*)=\ran(y)$ and $\dist(y^*,Z)<\e$.
 \end{proposition}
 \begin{proof}[\bf Proof.]
 From Lemma \ref{lem16} we can choose
  a $2-c_0^k$ vector $x^*=\sum\limits_{i=1}^kx_i^*$ in $Z$.
 Lemma \ref{lem18} yields the existence of a
 $2-\ell_1^k$ average $y$ with $\ran(y)=\ran(x^*)$ and $x^*(y)>1$.
 Applying Lemma \ref{lem17} we can find $y^*\in D_G$ with
 $\ran(y^*)=\ran(x^*)$ and
 $\|y^*-x^*\|<\min\{\e,\frac{x^*(y)-1}{2}\}$.
 It is clear that $y$ and $y^*$ satisfy the desired conditions.
 \end{proof}

 \begin{lemma}\label{lem25}
 Suppose that $X_G$ is a strongly strictly singular extension
 (with or without attractors)
  and let  $Z$ be a block subspace of $(X_G)_*$.
 Then for every $j>1$ and $\e>0$ there exists a $(18,2j,1)$ exact
 pair $(z,z^*)$ with $\dist(z^*,Z)<\e$ and
 $\|z\|_G\le\frac{3}{m_{2j}}$.
 \end{lemma}
 \begin{proof}[\bf Proof.]
 Using Proposition \ref{prop29} we may select  a block
 sequence \seq{y}{l} in $X_G$ and a sequence $(y^*_l)_{l\in\N}$
 such that
 \begin{enumerate}
 \item[(i)] Each $y_l$ is a $2-\ell_1^{n_{i_l}}$ average where
 \seq{i}{l} is an increasing sequence of integers.
 \item[(ii)] $y_l^*\in D_G$ for all $l$ and  $\sum\limits_{l=1}^{\infty}\dist(y_l^*,Z)<\e$.
 \item[(iii)] $y_l^*(y_l)>1$ and $\ran y_l^*=\ran y_l$.
 \end{enumerate}
 From Remark \ref{rem1} we may assume (passing, if necessary, to a subsequence)
 that \seq{y}{l} is a $(3,\e)$ R.I.S. Since $Y_G$  contains no
 isomorphic copy of $\ell_1$ we may assume, passing again to a
 subsequence, that \seq{y}{l} is a weakly Cauchy sequence in $Y_G$.
 Setting $x_l=\frac{1}{2}(y_{2l-1}-y_{2l})$ it is clear that
 $\|x_l\|_{\infty}\to 0$   while
 \seq{x}{l} is a weakly null sequence in $Y_G$. From the fact that
 $X_G$ is a strongly strictly singular extension (with or without
 attractors) it follows that there exists $M\in[\N]$
 such that for every $g\in G$ the set $\{l\in
 M:\;|g(x_l)|>\frac{2}{m_{2j}^2}\}$ has at most $n_{2j-1}$
 elements (notice that $66m_{2j}^4<n_{2j-1}$).
 We may assume that $M=\N$. Also \seq{x}{l} is a
 $(3,\e)$ R.I.S.
  We set
 \begin{equation*}
 z^*=\frac{1}{m_{2j}}\big(\sum\limits_{l=1}^{n_{2j}-1}y_{2l-1}^*-y^*_{2l}\big).
 \end{equation*}
  From Proposition \ref{prop13} we get that
 $\|\frac{1}{n_{2j}}\sum\limits_{l=1}^{n_{2j}}x_l\|\le
 \frac{6}{m_{2j}}$ while
 $z^*(\frac{m_{2j}}{n_{2j}}\sum\limits_{l=1}^{n_{2j}}x_l)>1$ hence
 there exists $\eta$ with $\frac{1}{6}\le \eta <1$ such that
 $z^*(\eta\frac{m_{2j}}{n_{2j}}\sum\limits_{l=1}^{n_{2j}}x_l)=1$.
 We set
 \begin{equation*}
 z=\eta\frac{m_{2j}}{n_{2j}}\sum\limits_{l=1}^{n_{2j}}x_l.
 \end{equation*}

 It follows easily from Proposition \ref{prop13} that $(z,z^*)$ is
 a $(18,2j,1)$ exact pair. From condition (ii) we get that
 $\dist(z^*,Z)<\e$. Finally we have that
 $\|z\|_G\le\frac{3}{m_{2j}}$.
 Indeed, let $g\in G$.
 Since $\#\{l:\;|g(x_l)|>\frac{2}{m_{2j}^2}\}\le n_{2j-1}$
 and $|g(x_l)|\le \|x_l\|\le 2$ for all $l$ we have that
\[
 |g(z)|\le
 \frac{m_{2j}}{n_{2j}}\sum\limits_{l=1}^{n_{2j}}|g(x_l)|\le
 \frac{m_{2j}}{n_{2j}}(\frac{2}{m_{2j}^2}n_{2j}+2n_{2j-1})<\frac{3}{m_{2j}}.
 \]
 Therefore $\|z\|_G\le\frac{3}{m_{2j}}$.
 \end{proof}

 \begin{lemma}\label{lem26}
 Let $X_G$ be a strongly strictly singular extension (with or without
 attractors)
  and let $Y,Z$ be
 a pair of block subspaces of $(X_G)_*$. Then for every
 $\e>0$ and $j>1$ there exists a $(18,4j-1,1)$ dependent sequence
 $(x_k,x_k^*)_{k=1}^{n_{4j-1}}$ with $\sum
 \dist(x_{2k-1}^*,Y)<\e$,
  $\sum \dist(x_{2k}^*,Z)<\e$ and $\|x_k\|_G\le\frac{2}{n_{4j-1}^2}$
 for all $k$.
 \end{lemma}
 \begin{proof}[\bf Proof.]
 This is an immediate consequence of Lemma \ref{lem25}.
 \end{proof}

 \begin{theorem}\label{th6.1}
 If $X_G$ is a strongly strictly singular extension (with or without attractors),
  then the
 predual space $(X_G)_*$ is Hereditarily Indecomposable.
 \end{theorem}
 \begin{proof}[\bf Proof.]
 Let $Y,Z$ be a pair of block subspaces of $(X_G)_*$. For every
 $j>1$   using Lemma \ref{lem26}
  we may select
 a $(18,4j-1,1)$ dependent sequence
 $(x_k,x_k^*)_{k=1}^{n_{4j-1}}$ with $\sum \dist(x_{2k-1}^*,Y)<1$
 and $\sum \dist(x_{2k}^*,Z)<1$ and $\|x_k\|_G\le\frac{2}{m_{4j-1}^2}$
 for all $k$.

   The functional
   $x^*=\frac{1}{m_{4j-1}}\sum\limits_{k=1}^{n_{4j-1}}x_k^*$
   belongs to the norming set $D_G$ hence $\|x^*\|\le 1$.
   From Proposition \ref{prop11} we get that
   $\|\frac{1}{n_{4j-1}}\sum\limits_{k=1}^{n_{4j-1}}(-1)^{k+1}x_k\|\le
   \frac{144}{m_{4j-1}^2}$.

   We set
   \[h_Y=\frac{1}{m_{4j-1}}\sum\limits_{k=1}^{n_{4j-1}/2}x_{2k-1}^*\mbox{ and }
     h_Z=\frac{1}{m_{4j-1}}\sum\limits_{k=1}^{n_{4j-1}/2}x_{2k}^*.        \]
  Estimating
  $h_Y-h_Z$ on the vector
 $\frac{1}{n_{4j-1}}\sum\limits_{k=1}^{n_{4j-1}}(-1)^{k+1}x_k$
 yields that $\|h_Y-h_Z\|\ge \frac{m_{4j-1}}{144}$
 while we obviously have that $\|h_Y+h_Z\|=\|x^*\|\le 1$.

 From the fact that $\dist(h_Y,Y)<1$ and  $\dist(h_Z,Z)<1$ we may select
 $f_Y\in Y$ and $f_Z\in Z$ with $\|h_Y-f_Y\|<1$ and
 $\|h_Z-f_Z\|<1$. From the above estimates we conclude that
 $\|f_Y-f_Z\|\ge (\frac{m_{4j-1}}{432}-\frac{2}{3})\|f_Y+f_Z\|$.
 Since we can find such $f_Y$ and $f_Z$ for arbitrary large $j$
 it  follows that $(X_G)_*$ is Hereditarily Indecomposable.
 \end{proof}

 The next two theorems concern the structure of $\mc{L}(\mathfrak
 X_G),$ $\mc{L}((\mathfrak X_G)_*)$.
 We start with the following lemmas. The first is the analogue of
 Lemma \ref{lem11} for strongly strictly singular extensions.

 \begin{lemma} \label{newA}
 Assume that $X_G$ is a strongly strictly singular extension.
 Let $Y$ be a subspace of $X_G$ and let $T:
 Y\to X_G$ be a bounded linear operator. Let $(y_\ell)_{\ell\in\N}$ be a
 block sequence in $Y$ of $C$-$\ell_1^{j_k}$ averages with $\lim
 j_k=\infty.$ Furthermore assume that $(Ty_\ell)_{\ell\in\N}$ is also a block
 sequence. Then
 \[\lim \dist (Ty_\ell, \R y_\ell)=0.\]
 \end{lemma}

 \begin{proof}[\bf Proof.] Assume that the conclusion fails. We
 may  assume, passing to a subsequence
  that there exists $\delta >0$ such that for every
 $\ell\in\N,$ $\dist (Ty_\ell, \R y_\ell)>\delta$ and moreover
 that $(y_{\ell})_\ell$ is a R.I.S.
  Next for each
 $\ell\in\N,$ we choose $\phi_\ell$ such that $\supp \phi_\ell
 \subset \ran (y_\ell \cup Ty_\ell),$ $\phi_\ell \in D_G,$
 $\phi_\ell(Ty_\ell)>\frac{\delta}{2}$ and $\phi_\ell(y_\ell)=0.$
 From Remark \ref{rem12}
 and since $X_G$ is a strictly singular extension of $Y_G,$
 for every $j\in\N,$ $j>j(C)$ we can find a subsequence
 $(y_{\ell_k})_{k\in\N}$ such that the sequence
 $w_k=(y_{\ell_{2k-1}}-y_{\ell_{2k}})/2$ is weakly null and for
 every $g\in G,$
 \[ \#\Big\{ k\in \N : |g(w_k)| > \frac{2}{m_{2j}^2}\Big\} \le
 n_{2j-1} .\]
 This yields that for every $j>j(C)$ there exists $w_{k_1} <
 w_{k_2} < \cdots < w_{k_{n_{2j}}}$ and $\phi_{k_1} < \phi_{k_2} <
 \cdots < \phi_{k_{n_{2j}}}$ such that setting $w \frac{m_{2j}}{n_{2j}} \sum\limits_{i=1}^{n_{2j}} w_{k_i}$ and $\phi  \frac{1}{m_{2j}}\sum\limits_{i=1}^{n_{2j}} \phi_i,$ we have that
 \[ \|w\|\le 6C,\quad \phi\in D_G, \quad \phi(Tw)> \frac{\delta}{2}\quad
\phi(w)=0, \quad \text{and} \quad \|w\|_G<
 \frac{3}{m_{2j}^2}.\] In particular $(w, \phi)$ is $(6C, 2j,
 0)$ exact pair with $\|w\|_G< \frac{3}{m_{2j}^2}.$ The remaining
 part of the proof follows the arguments of Lemmas 22 and 23 of
 \cite{GM1} using Proposition \ref{prop11}, (iii).
 \end{proof}

 The next lemma is easy and its proof is included in the proof of
 Theorem 9.4 of \cite{AT1}.

 \begin{lemma} \label{newB} Let $X$ be a Banach space with a
 boundedly complete basis \seq{e}{n} not containing $\ell_1.$
 Assume that $T: X\to X$ is a bounded linear non weakly compact
 operator. Then there exist two block sequences \seq{x}{n}
 \seq{y}{n} and $y$ in $X$ such that the following hold:
 \begin{enumerate} \item[(i)] $(x_n)_{n\in\N}$ is normalized,
 $x_n\stackrel{w^*}{\to} x^{**} \in X^{**}\setminus X.$
 \item[(ii)] $(y_n)_{n\in\N}$ is bounded, $y_n\stackrel{w^*}{\to}
 y^{**}\in X^{**}\setminus X.$
 \item[(iii)] $\|Tx_n-(y+y_n)\|\to 0.$
 \end{enumerate}
 \end{lemma}

 \begin{theorem}\label{th7}
 If $X_G$ is a strongly strictly singular extension (with or without
 attractors),
 then every bounded
 linear operator $T:X_G\to X_G$ takes the form $T=\lambda I+S$
 with $S$ a strictly singular and weakly compact operator.
 \end{theorem}
 \begin{proof}[\bf Proof.]
 We already know from Proposition \ref{prop16} that
 every bounded
 linear operator $T:X_G\to X_G$ is of the form $T=\lambda I+S$
 with $S$ a strictly singular operator so it remains to show that
 every strictly singular operator $S:X_G\to X_G$ is weakly
 compact.

 Assume now that there exists a strictly singular
 $T\in\mc{L}(X_G)$ which is not weakly compact. Then from Lemma
 \ref{newB}, there exist $(x_n)_{n\in\N},$ $(y_n)_{n\in\N},$ $y$ in $X_G$
 satisfying the conclusions of Lemma \ref{newB}. It follows that there
 exists a subsequence $(x_n)_{n\in L}$ such that setting
 $Z=\overline{\spann}\{x_n:\;n\in L\}$
 there exists a compact perturbation of $T|Z$
 denoted by $\tilde{T}$ such that for $n\in L,$ we have that
 $\tilde{T}(x_n)= y+z_n.$ For simplicity of notation assume that $L=\N$.
 Since $x_n \stackrel{w^*}{\to}
 x^{**} \in X_G^{**}\setminus X_G$ and $y_n \stackrel{w^*}{\to} y^{**}
 \in X_G^{**} \setminus X_G$ we may assume that every convex
 combination of $(x_n)_{n\in\N}$
 has norm greater
 than $\delta >0.$

 Choose $(z_k)_{k\in\N},$
  $z_k= \frac{1}{n_{j_k}}\sum\limits_{i\in F_k}\frac{1}{\delta}x_i$ with
 $\# F_k =n_{j_k}$ and $F_k < F_{k+1}$.
  Then setting $w_k=\frac{z_{2k-1}-z_{2k}}{\|z_{2k-1}-z_{2k}\|}$,
  Lemma \ref{newA} (actually its proof)
 yields that
 \[ \lim\dist (\tilde{T}w_k, \R w_k)=0.\]
 From this we conclude that for some subsequence $(w_k)_{k\in L}$,
 $\tilde{T}|\overline{\spann}\{w_k:\;k\in L\}$ is an isomorphism contradicting our
 assumption that $T$ is strictly singular.
 \end{proof}

\begin{theorem} \label{aa12} Let $X_G$ be a strongly singular extension
(with or without attractors) of $Y_G.$ Then every bounded
linear
operator $T: (X_G)_* \to (X_G)_*$ is of the form $T= \lambda I+ S$
with $S$ strictly singular.
\end{theorem}

 The proof of this result follows the lines of Proposition 7.1 in
 [AT2]. We first state two auxiliary lemmas.

 \begin{lemma} Let $X$ be a HI space with a Schauder basis
 $(e_n)_{n\in\N}$. Assume that $T: X\to X$ is a bounded linear operator
 which is not of the form $T=\lambda I+ S$ with $S$ strictly singular. Then
 there exists $n_0$ and $\delta >0$ such that for every $z \in
 X_{n_0}= \overline{\spann}\{e_n:\; n\ge n_0\},$ $\dist(Tz, \R z) \ge \delta \|z\|.$
 \end{lemma}
 \begin{proof}[\bf Proof.] If not, then there exists a normalized block
 sequence $(z_n)_{n\in\N}$ such that $\dist(Tz_n, \R z_n)\le \frac{1}{n}.$
 Choose $\lambda \in \R$ such that $\|Tz_n-\lambda z_n\|_{n\in
 L}\to 0$ for a subsequence $(z_n)_{n\in L}$. Then for a further
 subsequence $(z_n)_{n\in M}$ we have that
  $T-\lambda I|_{\overline{\spann}\{z_n:\;n\in M\}}$
 is a compact operator. The HI property of $X$ easily yields
 that $T-\lambda I$ is a strictly singular operator, contradicting
 our assumption.
 \end{proof}

 \begin{lemma} \label{aa10} Let $T:(X_G)_* \to (X_G)_*$ be a bounded linear
 operator with $\|T\| = 1.$ Assume that for some $\delta >0,$ and
 $n_0\in \N,$ $\dist(Tf, \R f) \ge \delta \|f\|$ for all
 $f\in (X_G)_*$ with $n_0< \supp f.$ Then for every
 $k\in \N$ and every block subspace $Z$ of $(X_G)_*$ there exist a
 $z^* \in Z$ with $\|z^*\| \le 1$ and a
 $\frac{2}{\delta}$-$\ell_1^k$ average $z$  such that $z^*(z) =0$,
  $Tz^*(z) >1$ and $\ran z \subset \ran z^* \cup \ran Tz^*.$
 \end{lemma}

 \begin{proof}[\bf Proof.] By Lemma \ref{lem16} there exists
 a 2-$c_0^k$ vector  $z^*=\sum\limits_{i=1}^k z_i^*$ in $Z$ with $n_0 <
 \min \supp z^*.$ Since $\dist (Tz_i^*, \R z_i^*) \ge \delta
 \|z_i^*\| > \frac{\delta}{2}$ we may choose for each $i=1, \dots,
 k$ a vector $z_i\in X_G$ with $\|z_i\|< \frac{2}{\delta}$ and $\supp
 z_i \subset \ran z_i^* \cup \ran Tz_i^*$ satisfying $Tz_i^* (z_i)
 >1$ and $z_i^*(z_i) = 0.$ We set $z = \frac{1}{k}
 \sum\limits_{i=1}^k
 z_i.$ It is easy to check that $z$ is the desired vector.
 \end{proof}

 \begin{proof}[\bf Proof of Theorem \ref{aa12}] On the contrary assume that
 there exists $T\in \mc{L}((X_G)_*)$ which is not of the desired form.
 Assume further that $\|T\|=1$ and $Te_n^*$ is finitely supported
 with $\lim \min \supp T e_n^*=\infty.$ (We may assume the later
 conditions  from
 the fact that the basis $(e_n^*)_{n\in\N}$ of $(X_G)_*$ is weakly null.)
 In particular for every  block sequence $(z_n^*)_{n\in\N}$ in $(X_G)_*$
 there exists a subsequence $(z_n^*)_{n\in L}$ such that $(\ran
 z_n^* \cup \ran Tz_n^*)_n$ is a sequence of successive subsets of
 $\N.$

 Let $\delta >0$ and $n_0 \in \N$ be as in Lemma \ref{aa10} and
 let
 $j(\frac{2}{\delta})$ be the corresponding index such that for all
 $j\ge j(\frac{2}{\delta})$ the conclusion of Definition \ref{aa11}
 holds for $\frac{2}{\delta}$-bounded block sequences of $X_G.$
 Using arguments similar to those of Lemma \ref{lem25} for $j\ge
 j(\frac{2}{\delta})$ we can find an $(\frac{18}{\delta},  2j, 0)$
 exact pair $(z, z^*),$ with $Tz^*(z) > 1$ and $\|z\|_G \le
 \frac{3}{m_{2j}}.$
 Then for every $j\in\N$ there exists  a $(\frac{18}{\delta}, 4j-1, 0)$
 dependent sequence $(z_k,
 z_k^*)_{k=1}^{n_{4j-1}}$, such that
 $z_k^*(z_k)=0$,
 $Tz_k^*(z_k) >1$,
 $\|z_k\|_G \le \frac{1}{m_{4j-1}^2}$,
 $(\ran z_k^* \cup \ran Tz^*_k)_{k=1}^{n_{4j-1}}$ are successive subsets of $\N$
 and $\ran z_k \subset I_k$ where $I_k$ is the
 minimal interval of $\N$ containing $\ran z_k^* \cup \ran
 Tz_k^*$.

 Proposition \ref{prop11} yields that
 \[ \|\frac{1}{n_{4j-1}} \sum_{k=1}^{n_{4j-1}} z_k\| \le
 \frac{144}{m_{4j-1}^2 \delta}. \]

 Finally $\|\frac{1}{m_{4j-1}} \sum\limits_{k=1}^{n_{4j-1}}Tz_k^*\| \le
 1$ (since $\|T\|\le 1$ ) and also
 \[ 1 \ge \|\frac{1}{m_{4j-1}} \sum_{k=1}^{n_{4j-1}} Tz_k^*\|
 \ge\frac{m_{4j-1}^2 \delta}{144 m_{4j-1}} \frac{1}{n_{4j-1}}
 \sum_{k=1}^{n_{4j-1}} Tz_k^*(z_k) \ge
 \frac{m_{4j-1}\delta}{144}.\]
 This yields a contradiction for sufficiently large $j\in\N.$
 \end{proof}

 The following lemma is similar to a corresponding result used by
 V. Ferenczi \cite{Fe} in order to show that every quotient
 of the space constructed by W.T Gowers and B. Maurey remains
 Hereditarily Indecomposable.

 \begin{lemma}\label{lem30}
 Suppose that $X_G$ is a strictly singular extension of $Y_G$
 (with or without attractors). Let $Z$ be $w^*$ closed subspace
 of $X_G$ and let $Y$ be a closed subspace of $X_G$ with $Z\subset Y$
 such that the quotient space $Y/Z$ is infinite dimensional. Then
 for every $m,N\in\N$ and $\e>0$ there exists
 $x\in\spann\{e_i:\;i\ge m\}$ which is a $2-\ell_1^N$ average with $\dist(x,Y)<\e$
 and there exists
 $f\in B_{(X_G)_*}$ with $\dist(f,Z_{\bot})<\e$ such that $f(x)>1$.
 \end{lemma}
 \begin{proof}[\bf Proof.]
 We recall that from the fact that $Z$ is  $w^*$ closed the quotient space $X_G/Z$
 may be identified with the dual of the annihilator $Z_{\bot}=\{f\in
 (X_G)_*:\; f(z)=0\;\;\;\forall z\in Z\}$, i.e.  $X_G/Z=(Z_{\bot})^*$.
 Pick a normalized sequence
 $(\hat{y}_n')_{n\in\N}$ in $Y/Z$ with $\hat{y}_n'\stackrel{w^*}{\to}
 0$. From W.B. Johnson's and H.P. Rosenthal's  work on $w^*-$ basic sequences (\cite{JR})
 and their $w^*$ analogue of the
 classical  Bessaga - Pelczynski theorem, we may assume, passing
 to a subsequence, that $(\hat{y}_n')_{n\in\N}$ is a $w^*$ basic
 sequence.  Hence there exists  a bounded
 sequence $(z_n^*)_{n\in\N}$ in $Z_{\bot}$ such that $(z_n^*,\hat{y}_n')_{n\in\N}$
 are biorthogonal ($z_n^*(\hat{y}_m')=\delta_{nm}$) and
 $\sum\limits_{i=1}^n \hat{y}(z_n^*)\hat{y}_i'\to \hat{y}$ for
 every $\hat{y}$ in  the weak${}^*$ closure of the linear span of
 the sequence $(\hat{y}_n')_{n\in\N}$.

 Since $(X_G)_*$ contains no isomorphic copy of $\ell_1$ (as it a
 space with separable dual) we may assume, passing to a
 subsequence, that $(z^*_n)_{n\in\N}$ is weakly Cauchy, hence
 $(z^*_{2n-1}-z^*_{2n})_{n\in\N}$ is weakly null.
 Using a sliding hump argument and passing to a subsequence we may
 assume that with an error up to $\e$ this sequence is a block
 sequence with respect to the standard basis  of $(X_G)_*$.

 We set $y_n^*=z^*_{2n-1}-z^*_{2n}$ and $\hat{y}_n=\hat{y}_{2n-1}'$
 for $n=1,2,\ldots$. Then $(y_n^*,\hat{y}_n)_{n\in\N}$
 are biorthogonal, $(y_n^*)_{n\in\N}$ is a weakly null block
 sequence in $(X_G)_*$ with $y_n^*\in Z_{\bot}$, while
 $(\hat{y_n})$ is a normalized $w^*-$ basic sequence in $Y/Z$.

 We choose $k,j\in\N$  such that $2^k>m_{2j}$ and $(2N)^k\le
 n_{2j}$. We set
 \[\mc{A}_1=\big\{L\in [\N],\;\; L=\{l_i,\;i\in\N\}:\;
 \|\frac{1}{2N}\sum\limits_{i=1}^{2N}(-1)^{i+1}\hat{y}_{l_i}\|>\frac{1}{2}\big\}\]
 \[  \mbox{ and }\qquad\mc{B}_1=[\N]\setminus \mc{A}_1\]
 From Ramsey's theorem we may find a homogenous set $L$ either in
 $\mc{A}_1$ or in $\mc{B}_1$. We may assume that $L=\N$.

 Suppose first that the homogenous set is in $\mc{A}_1$, i.e.
 $\|\frac{1}{2N}\sum\limits_{i=1}^{2N}(-1)^{i+1}\hat{y}_{l_i}\|>\frac{1}{2}$
 for every $l_1<l_2<\ldots <l_{2N}$ in $\N$. For each $n$ we may
 choose $y_n\in Y\subset X_G$ with $\|y_n\|=1$ and
 $Q(y_n)=\hat{y}_n$.
 Passing to a subsequence we may assume (again with an error up to $\e$)
 that the sequence
 $x_n=y_{2n-1}-y_{2n}$ is a weakly null block sequence in $X_G$
 with $\min\supp x_i\ge m$. We set
 $x=\frac{1}{N}\sum\limits_{i=1}^Nx_i$.
 It is clear that $x$ is a 2-$\ell_1^N$ average while
 since $Qx\in Y/Z\subset X_G/Z=(Z_{\bot})^*$ and $\|Qx\|>1$ there
 exists $f\in Z_{\bot}$ with $\|f\|\le 1$ such that $f(x)>1$.

 On the other hand if the homogenous set is in $\mc{B}_1$ then we
 may assume, passing again to a subsequence that there exists
 $a_1\ge 2$ such that setting
 $\hat{y}_{2,n}=a_1\cdot\frac{1}{2N}\sum\limits_{i=(n-1)(2N)+1}^{n(2N)}(-1)^{i+1}\hat{y}_{i}$
 for $i=1,2,\ldots$, $(\hat{y}_{2,n})_{n\in\N}$ is a normalized
 sequence in $Y/Z$.  We may again apply Ramsey's theorem defining
 $\mc{A}_2$, $\mc{B}_2$ as before, using the sequence
 $(\hat{y}_{2,n})_{n\in\N}$ instead of $(\hat{y}_{n})_{n\in\N}$.
 If the homogenous set is in $\mc{A}_2$ the proof finishes as
 before while if it is in $\mc{B}_2$ we continue defining
 $(\hat{y}_{3,n})_{n\in\N}$ $\mc{A}_3$, $\mc{B}_3$ and so on.

 If in none of the first $k$ steps we  arrived at a homogenous set
 in some
 $\mc{A}_i$  then there exist $a_1,a_2,\ldots,a_k\ge
 2$ and $l_1<l_2<\cdots<l_{(2N)^k}$ in $\N$ such that the vector
 \[\hat{y}=a_1a_2\cdot\ldots\cdot
 a_k\frac{1}{(2N)^{k}}\sum\limits_{i=1}^{(2N)^{k}}(-1)^{i+1}\hat{y}_{l_i}\]
 satisfies $\|\hat{y}\|=1$.

 But then the functional
 $y^*=\frac{1}{m_{2j}}\sum\limits_{i=1}^{(2N)^{k}}(-1)^{i+1}y_{l_i}^*$
 belongs to $Z_{\bot}$ and satisfies $\|y^*\|\le 1$ (as
 $(y_{l_i}^*)_i$ is a block sequence with
 $\|y_{l_i}^*\|\le 1$  and $(2N)^{k}\le n_{2j}$).
 Therefore, taking into account the biorthogonality, we ge that
 \[ 1=\|\hat{y}\|\le
 y^*(\hat{y})=\frac{1}{m_{2j}}\frac{2^k}{(2N)^k}(2N)^k=\frac{2^k}{m_{2j}}\]
 which contradicts  our choice of $k$ and $j$.
 \end{proof}

 \begin{lemma}\label{lem31}
 Suppose that $X_G$ is a strongly strictly singular extension of $Y_G$
 (with or without attractors) and let $Y$ and $Z$ be as in Lemma \ref{lem30}.
 Then for every $j>1$ and every $\e>0$ there exists a $(18,2j,1)$
 exact pair $(y,f)$ with
 $\dist(y,Y)<\e$, $\|y\|_G<\frac{3}{m_{2j}}$ and
 $\dist(f,Z_{\bot})<\e$.
 \end{lemma}
 \begin{proof}[\bf Proof.]
 Let \seq{\e}{i} be a sequence of positive reals with
 $\sum\limits_{i=1}^{\infty}\e_i<\e$.
 Using Lemma \ref{lem30} we may inductively construct a block sequence \seq{x}{i}
 in $X_G$, a sequence \seq{\phi}{i} in $B_{(X_G)_*}$ and a sequence of
 integers $t_1<t_2<\cdots$ such that the following are satisfied:
 \begin{enumerate}
 \item[(i)] $\dist(x_i,Y)<\e_i$, $\dist(\phi_i,Z_{\bot})<\e_i$ and
 $\phi_i(x_i)>1$.
 \item[(ii)] The sequence \seq{x}{i} is a  sequence of $2-\ell_1$
 averages of increasing length and $\min\supp x_i\ge t_i$.
 \item[(iii)] The restriction of the functional $\phi_i$ to the space
 $\overline{\spann}\{e_n:\;n\ge t_{i+1}\}$ has norm at most $\e_i$.
 \end{enumerate}

 Passing  to a subsequence we may assume that the sequences
 \seq{\phi}{i} and \seq{x}{i} are weakly Cauchy. Thus the sequence
 $(-\phi_{2n-1}+\phi_{2n})_{n\in\N}$ is weakly null; so we
 may assume, passing again to a subsequence, that it is a
 block sequence and that,
 since  $(x_{2n})_{n\in\N}$ is a weakly Cauchy sequence, the
 sequence $y_n=x_{4n-2}-x_{4n}$ is a weakly null block sequence in
 $X_G$  and thus also in
 $Y_G$, therefore, from the fact that $X_G$ is a
 strongly strictly singular extension of $Y_G$ we may assume,
 passing to a subsequence, that for every $g\in G$ the set
 $\{n\in\N:\; |g(y_n)|>\frac{2}{m_{2j}^2}\}$ contains at most
 $n_{2j-1}$ elements and also that \seq{y}{n} is a $(3,\e)$ R.I.S.
 of $\ell_1$ averages.

 We set $f_n=\frac{1}{2}(-\phi_{4n-3}+\phi_{4n-2})$ for
 $n=1,2,\ldots$, and we may assume that $\max(\supp f_n\cup\supp
 y_n)<\min (\supp f_{n+1}\cup\supp y_{n+1})$ for all n.
 Finally we set
 \[ y'=\frac{m_{2j}}{n_{2j}}\sum\limits_{i=1}^{n_{2j}}y_i\mbox{ and }
 f=\frac{1}{m_{2j}}\sum\limits_{i=1}^{n_{2j}}f_i. \]

 As in the proof of Lemma \ref{lem25} we obtain that $(y,f)$,
 where $y$ is a suitable scalar multiple of $y'$, is the desired exact pair.
 \end{proof}

 Using Lemma \ref{lem31} we  prove the following:

 \begin{theorem}\label{prop30}
 If $X_G$ is a strongly strictly singular extension of $Y_G$
 (with or without attractors) and $Z$ is a $w^*$ closed subspace
 of $X_G$ of infinite codimension then the quotient space $X_G/Z$
 is  Hereditarily Indecomposable.
 \end{theorem}
 \begin{proof}[\bf Proof.] Let $Y_1$ and $Y_2$ be subspaces of $X_G$
 with $Z\hookrightarrow Y_1\cap Y_2$ such that $Z$ is of infinite
 codimension in each $Y_i,$ $i=1,2.$ Then for every $\e >0$
 and $j\in \N$ we may select an $(18, 4j-1, 1)$ dependent sequence
 $\chi =(x_k, x_k^*)_{k=1}^{n_{4j-1}}$ such that
 \begin{enumerate}\item[(i)] $\|x_k\|_G < \frac{1}{m_{4j-1}^2},$
 $k=1, \dots, n_{4j-1}.$
 \item[(ii)] $\dist (x_{2k-1}, Y_1) <\e,$ $\dist (x_{2k},
 Y_2) < \e.$
 \item[(iii)] $\dist (x_k^*, Z_\perp) < \e.$
 \end{enumerate}

 Let $Q: X_G\to X_G/Z$ be the quotient map.
 If $\e >0$ is sufficiently small
 Proposition \ref{prop11}
 easily yields that
 \[ \dist (S_{Q(Y_1)}, S_{Q(Y_2)}) < \frac{C}{m_{4j-1}} \]
 where $C$ is a constant independent of $j$. The proof is complete.
 \end{proof}

  \section{The James tree space    $JT_{\mc{F}_{2}}.$}

In this section we define a class of James Tree-like spaces. These
spaces share some of the main properties of the classical $JT$
space. Namely they do not contain an isomorphic copy of the space
$\ell_1.$ Furthermore they have a bimonotone basis. In particular
their norming set is a ground set and a specific example of this
form will be the ground set for our final constructions. The
principal goal is to prove  the inequality in Proposition
\ref{aa4} yielding that that the ground set $\mc{F}_2$ defined in
the next section admits a strongly strictly singular extension. In
Appendix B we present a systematic study of $JT_{\mc{F}_2}$ spaces
and of some variants of them.

 \begin{definition}\label{def1}
  {\bf (JTG families)}
 A family  $\mc{F}=(F_j)_{j=0}^{\infty}$  of subsets of
 $c_{00}(\N)$ is said to be a {\bf James Tree Generating family}
  (JTG
 family)  provided it  satisfies the following conditions:
 \begin{enumerate}
 \item[(A)] $F_0=\{\pm e_n^*:\;n\in\N\}$
 and each $F_j$ is nonempty, countable, symmetric, compact in the
 topology of pointwise convergence and closed under restrictions to
 intervals of $\N$.
 \item[(B)] Setting ${\tau}_j=\sup\limits\{\|f\|_{\infty}:\; f\in
 F_j\}$, the sequence \seq{\tau}{j}
 is strictly
 decreasing and $\sum\limits_{j=1}^{\infty}{\tau}_j\le 1$.
 \end{enumerate}
  \end{definition}

 \begin{definition}\label{def4}
 {\bf (The $\sigma_{\mc{F}}$ coding)}
 Let $(F_j)_{j=0}^{\infty}$ be a JTG family.
 We fix a pair $\Xi_1,\Xi_2$ of
 disjoint infinite subsets of $\N$.
  Let $W=\{(f_1,\ldots,f_d):\; f_i\in \cup_{j=1}^\infty F_j,\;f_1<\cdots<f_d,\;d\in\N\}$.
  The set $W$ is countable so we
 may select an 1--1 coding function $\sigma_{\mc{F}}:W\to \Xi_2$ such that for
 every $(f_1,\ldots,f_d)\in W$,
 \[ \sigma_{\mc{F}}(f_1,\ldots,f_d)>\max\big\{k:\; \exists
 i\in\{1,\ldots,d\}\;\mbox{ with } f_i\in F_k\big\}.  \]

 A finite or infinite block
 sequence $(f_i)_i$ in
 $\bigcup\limits_{j=1}^{\infty}F_j\setminus\{0\}$
  is
 said to be a $\sigma_{\mc{F}}$ special sequence provided
 $f_1\in\bigcup\limits_{l\in \Xi_1}F_l$
 and
 $f_{i+1}\in F_{\sigma_{\mc{F}}(f_1,\ldots,f_i)}$ for all $i$.
 A  $\sigma_{\mc{F}}$ special functional $x^*$
  is any functional of the form
 $x^*=E\sum\limits_i f_i$ with $(f_i)_i$ a $\sigma_{\mc{F}}$
 special sequence (when the sum $\sum\limits_i f_i$ is infinite it
 is considered in the pointwise topology) and $E$ an  interval of $\N$.
 If the interval $E$ is finite then $x^*$ is said to be a finite
 $\sigma_{\mc{F}}$ special functional.
 We denote by $\mathscr{S}$ the set of all finite $\sigma_{\mc{F}}$
 special functionals.
 Let's observe that $\overline{\mathscr{S}}^{w^*}$ is the set of
 all $\sigma_{\mc{F}}$ special functionals.
 \end{definition}

 \begin{definition}\label{def222}
 \begin{enumerate}
 \item[(A)]
 Let $s=(f_i)_i$ be a $\sigma_{\mc{F}}$ special sequence.
 Then for for each $i$ we define the $\ind_s(f_i)$ as follows.
 $\ind_s(f_1)=\min\{j:\;f_1\in F_j\}$ while for $i=2,3,\ldots$
 $\ind_s(f_i)=\sigma_{\mc{F}}(f_1,\ldots,f_{i-1})$.
 \item[(B)]
 Let  $s=(f_i)_i$ be a $\sigma_{\mc{F}}$ special sequence and let $E$ be an interval.
 The set of indices of the $\sigma_{\mc{F}}$ special functional
 $x^*=E\sum\limits_if_i$ is the set
 $\ind_s(x^*)=\{\ind_s(f_i):\;Ef_i\neq 0\}$.
 \item[(C)] A (finite or infinite) family of $\sigma_{\mc{F}}$ special functionals $(x_k^*)_k$
 is said to be  disjoint if  for each $k$
 there exists a $\sigma_{\mc{F}}$
  special sequence $s_k=(f^k_i)_i$  and  interval $E_k$ such that
 $x_k^*=E_k\sum\limits_if^k_i$
 and $(\ind_{s_k}(x_k^*))_k$ are pairwise disjoint.
 \end{enumerate}
 \end{definition}

 \begin{remark}
 \begin{enumerate}
 \item[(a)] Our definition of $\ind_s(f_i)$ and $\ind_s(x^*)$, which
 is
 rather technical, is required by the fact that we did not assume
 $(F_i\setminus\{0\})_i$ to be pairwise disjoint, hence the same $f$ could
 occur in several different $\sigma_{\mc{F}}$ special sequences.
 \item[(b)] Let's observe  that for every family
 $(x_i^*)_{i=1}^d$ of disjoint $\sigma_{\mc{F}}$ special
 functionals,
 $\|\sum\limits_{i=1}^dx_i^*\|_{\infty}\le 1$
 (recall that
 $\sum\limits_{j=1}^{\infty}\tau_j\le 1$).
 \item[(c)] Let $s_1=(f_i)_i$, $s_2=(h_i)_i$ be two
  distinct $\sigma_{\mc{F}}$ special
 sequences.
 Then $\ind_{s_1}(f_i)\neq\ind_{s_2}(h_j)$ for $i\neq j$
 while there exists $i_0$ such that $f_i=h_i$ for all $i<i_0$ and
  $\ind_{s_1}(f_i)\neq \ind_{s_2}(h_i)$
 for  $i>i_0$.
 \item[(d)] For every family $(s_i)_{i=1}^d$ of infinite $\sigma_{\mc{F}}$ special
 sequences there exists $n_0$ such that $(Es_i^*)_{i=1}^d$ are disjoint, where
 $E=[n_0,\infty)$ and $s_i^*$ denotes the $\sigma_{\mc{F}}$ special functional
 defined by the $\sigma_{\mc{F}}$ special sequence $s_i$.
 \end{enumerate}
 \end{remark}

 \begin{definition}\label{def17}
  {\bf (The norming set   $\mc{F}_2$).}
  Let $(F_j)_{j=0}^{\infty}$ be a JTG family.
 We set
 \begin{eqnarray*}
 \mc{F}_2 & = & F_0\cup \big\{\sum\limits_{k=1}^da_kx^*_k:\; a_k\in\Q,\;
 \sum\limits_{k=1}^da_k^2\le 1,\;
 \mbox{ and } \\
   & & (x_k^*)_{k=1}^d\mbox{ is a family of disjoint finite $\sigma_{\mc{F}}$
   special functionals}
    \big\}
 \end{eqnarray*}

 The space $JT_{\mc{F}_2}$ is defined as the
 completion of the space $(\co,\|\;\|_{\mc{F}_{2}})$
 where $\|x\|_{\mc{F}_2}
 =\sup\{f(x):\; f\in \mc{F}_2\}$ for $x\in\co$.
 \end{definition}

 \begin{remark}\label{rem11}
 Let's observe that the standard basis \seq{e}{n} of \co is a
normalized bimonotone
 Schauder basis of the space $JT_{\mc{F}_2}$.
  \end{remark}

 \begin{theorem} \label{th8}
 \begin{enumerate}
 \item[(i)]
 The space $JT_{\mc{F}_2}$ does not contain $\ell_1$.
  \item[(ii)]
  $JT_{\mc{F}_2}^*=\overline{\spann}(\{e_n^*:\;n\in\N\}\cup\{b^*:\;b\in
  \mc{B}  \})$ where $\mc{B}$ is the set of all infinite $\sigma_{\mc{F}}$
  special sequences.
 \end{enumerate}
 \end{theorem}
 The proof of the above theorem is almost identical with the
 proofs of Propositions 10.4 and  10.11 of \cite{AT1}. We proceed
 to a short description of the basic steps.

 Let's start by observing the following.
 \begin{eqnarray*}
 \overline{\mc{F}_2}^{w^*} & = &
  F_0\cup \big\{\sum\limits_{k=1}^{\infty}a_kx^*_k:\; a_k\in\Q,\;
 \sum\limits_{k=1}^{\infty}a_k^2\le 1,\;
 \mbox{ and } \\
   & & (x_k^*)_{k=1}^{\infty}
   \mbox{ is a family of disjoint $\sigma_{\mc{F}}$ special functionals}
    \big\}
 \end{eqnarray*}
 Also for a disjoint family $(x_i^*)_{i=1}^{\infty}$ of special
 functionals and $(a_i)_{i=1}^{\infty}$ in $\R$, we have that\\
 $\|\sum\limits_{i=1}^{\infty}a_ix_i^*\|_{JT_{\mc{F}_2}^*}\le
 \big(\sum\limits_{i=1}^{\infty}a_i^2\big)^{1/2}$.
 The above observations yield the following:

 \begin{lemma}\label{lem29}
  $\overline{\mc{F}}_2^{w^*}
\subset \overline{\spann}(\{e_n^*:\;n\in\N\}\cup\{b^*:\;b\in\mc{B}
  \})$ where $\mc{B}$ is the set of all infinite $\sigma_{\mc{F}}$
  special sequences.
 \end{lemma}
 Observe also that $\overline{\mc{F}_2}^{w^*}$ is $w^*$ compact and
 1-norming hence contains the set $\Ext(B_{JT_{\mc{F}_2}^*})$.
 Rainwater's theorem and the above results yield that a bounded
 sequence \seq{x}{k} is weakly Cauchy if and only if
 $\lim\limits_k e_n^*(x_k)$ and  $\lim\limits_k b^*(x_k)$ exist
 for all $n$ and  infinite special sequences $b$. This is
 established by the following.

 \begin{lemma}\label{lem27}
 Let \seq{x}{k} be a bounded sequence in $JT_{\mc{F}_2}$ and let
 $\e>0$. Then there exists a finite family $x_1^*,\ldots,x_d^*$ of
 disjoint special functionals and an $L\in[\N]$ such that
 \[\limsup\limits_{k\in L} |x^*(x_k)|\le \e\]
 for every special functional $x^*$ such that the family
 $x^*,x_1^*,\ldots,x_d^*$ is disjoint.
 \end{lemma}
 For a proof we refer the reader to the proof of Lemma 10.5 \cite{AT1}.

 \begin{lemma}\label{lem28} Let \seq{x}{k} be a bounded sequence in
 $JT_{\mc{F}_2}$. There exists an $M\in[\N]$
 such that for every special functional $x^*$ the sequence
 $(x^*(x_k))_{k\in M}$ converges.
 \end{lemma}
 Also for the proof of this we refer the reader to the proof of
 Lemma 10.6 of \cite{AT1}

 \begin{proof}[\bf Proof of Theorem \ref{th8}]
 (i)\; Let \seq{x}{k} be a bounded sequence in
 $JT_{\mc{F}_2}$. By an easy diagonal argument we may assume that
 for every $n\in\N$, $\lim\limits_n e_n^*(x_k)$ exists. Lemma
 \ref{lem28} also yields that there exists a subsequence
 $(x_{l_k})_{k\in\N}$ such that for every special sequence $b$,
 $\lim\limits_k b^*(x_{l_k})$ also exists. As we have mentioned
 above
 Lemma \ref{lem29} yields that
 $(x_{l_k})_{k\in\N}$ is weakly Cauchy. \\
 (ii)\; Since $\Ext(B_{JT_{\mc{F}_2}^*})\subset\overline{\mc{F}_2}^{w^*}$
 and $\ell_1$ does not embed into $JT_{\mc{F}_2}$
 Haydon's theorem \cite{Ha} yields that
 $\overline{\mc{F}_2}^{w^*}$  norm generates
 $JT_{\mc{F}_2}^*$. Lemma \ref{lem29} yields the desired result.
 \end{proof}

The remaining part of this section concerns the proof of
Proposition \ref{aa4}, stated below. This will be used in the next
section to show that a specific ground set $\mc{F}_2$ admits a
strongly strictly singular extension.

 \begin{definition} \label{aa1} Let $(x_n)_n$ be a bounded block
sequence in $JT_{\mc{F}_2}$ and $\e >0.$ We say that $(x_n)_n$ is
{\em $\e$-separated} if for every $\phi\in \cup_{j\in
\mathbb{N}}F_j$
\[ \# \{n: |\phi(x_n)|\ge\e \}\le 1. \]

In addition, we say that $(x_n)_n$ is {\em separated} if for every
$L\in [\N]$ and $\e >0$ there exists an $M\in [L]$ such that
$(x_n)_{n\in M}$ is $\e$-separated.
\end{definition}

\begin{lemma} \label{aa2} Let $(x_n)_n$ be a
bounded separated   sequence in $JT_{\mc{F}_2}$ such that for
every infinite $\sigma_{\mc{F}}$ special functional $b^*$ we have
that $\lim\limits_n b^*(x_n)=0$. Then for every $\e
>0,$ there exists an $L\in [\N]$ such that for all $y^* \in
\overline{\mathscr{S}}^{w^*},$
\[\# \{ n\in L:\; |y^*(x_n)| \ge \e \}\le 2. \]
\end{lemma}

\begin{proof}[\bf Proof] Assume the contrary and fix an $\e >0$ such
that the statement of the lemma is false. Define
\[ A= \big\{ (n_1<n_2< n_3) \in [\N]^3:\; \exists y^* \in
\overline{\mathscr{S}}^{w^*},\; |y^*(x_{n_1})|, |y^*(x_{n_2})|,
|y^*(x_{n_2})|\ge \e \big\}\] and $B= [\N]^3 \setminus A.$ Then
Ramsey's Theorem yields that there exists an $L\in [\N]$ such that
either $[L]^3\subset A$ or $[L]^3 \subset B.$ Our assumption
rejects the second case, so we conclude that for all $n_1 < n_2 <
n_3 \in L,$ there is a $y_{n_1, n_2, n_3}^* \in
\overline{\mathscr{S}}^{w^*}$ such that $|y_{n_1, n_2,
n_3}^*(x_{n_i})| \ge \e,$ $i=1,2,3.$

Since $(x_n)_{n\in L}$ is separated, we may assume by passing to a
subsequence that for $\e' = \frac{\e}{8},$ $(x_n)_{n\in L}$ is
$\e'$-separated. For reasons of simplicity in the notation we may
moreover and do assume that $(x_n)_{n\in \N}$ has both properties.

For all triples $(1<n< k),$ let $y_{n,k}^*$ denote an element in
$\overline{\mathscr{S}}^{w^*}$ such that $|y_{n,k}^*(x_i)|\ge \e,$
$i=1, n, k.$ Moreover, let $y_{n,k}^* = E_{n,k} \sum_{i=1}^\infty
\phi_{n,k}^i$ where $(\phi_{n,k}^i)_{i\in\N}$ is a
$\sigma_{\mc{F}}$ special sequence and $E_{n,k}\subset \N$ is an
interval. For $1<n <k$ we define the number $[n,k]$ as follows:
\[ [n,k] = \min \{ i\in\N:\; \max \supp \phi_{n,k}^i \ge \min \supp
x_k\}. \] Also, let $A= \{(n<k)\in [\N\setminus\{1\}]^2:\;
|\phi_{n,k}^{[n,k]}(x_n)|\le \e'\}$ and
$B=[\N\setminus\{1\}]^2\setminus A.$

Again,  using Ramsey's theorem and passing to a subsequence, we
may and do assume that $[\N\setminus\{1\}]^2\subset A$ or
$[\N\setminus\{1\}]^2 \subset B.$ Notice in the second case, that
since $(x_n)_n$ is $\e'$-separated, we have that for all $1<n<k$,
$|\phi_{n,k}^{[n,k]}(x_k)| \le \e'. $ We set
\[ s_{n,k}= \begin{cases} (\phi_{n,k}^1, \dots,
\phi_{n,k}^{[n,k]-1}),  & \text{ if } [\N\setminus\{1\}]^2\subset A\\
(\phi_{n,k}^1, \dots, \phi_{n,k}^{[n,k]}), &\text{ if }
[\N\setminus\{1\}]^2\subset B \end{cases}. \]

\begin{claim} There is an $M>0$ such that for all $k\in \N,$
\[\# \{s_{n,k} :\; 2\le n \le k-1\} \le M.\]\end{claim}

Let $(x_n)_n$ be bounded by some $c>0.$ Next fix any $k\in \N$ and
consider the following two cases:

The first case is $[\N\setminus\{1\}]^2\subset B.$ In this case
$\phi_{n,k}^{[n,k]} \in s_{n,k}$ and if $s_{n_1, k}\neq s_{n_2,
k}$ then $\phi_{n_1, k}^{[n_1, k]}$ is incomparable to $\phi_{n_2,
k}^{[n_2, k]}$ in the sense,that every two special functionals,
extending $s_{n_1, k}$ and  $s_{n_2, k}$  respectively, have
disjoint sets of indices.

So let $s_{n_j, k},$ $1\le j\le N$  all be different from each
other and consider the $\sigma_{\mc{F}}$ special functionals
$z_{n_j}^* = E_{n_j,k}y_{n_j,k}^*, $ $1\le j\le N$ where
$E_{n_j,k} = (\max \supp \phi_{n_j, k}^{[n_j,k]}, \infty).$
According to the previous observation these functionals have
pairwise disjoint indices. Moreover \begin{equation} \label{aa2.1}
|z_{n_j}^*(x_k)| = |y_{n_j,k}^*(x_k) - \phi_{n_j, k}^{[n_j,
k]}(x_k)| \ge \e-\e'\end{equation} since $[\N\setminus\{1\}]^2
\subset B$ and $(x_n)_n$ is $\e'$-separated.

Inequality \eqref{aa2.1} yields that \[
\Big(\sum_{j=1}^N(z_{n_j}^* (x_k))^2\Big)^{1/2} \ge (\e-\e')
N^{1/2}.\] Therefore there are $(a_j)_{j=1}^N$ with $\sum_{j=1}^n
a_j^2 \le 1$ such that \[ \sum_{j=1}^N a_j z_{n_j}^*(x_k) \ge (\e
- \e') N^{1/2}.\] On the other hand, by the definition of the norm
on $JT_{\mc{F}_2},$ $\sum_{j=1}^N a_j z_{n_j}^*(x_k) \le
\|x_k\|\le c.$ It follows that $N\le (\frac{c}{\e-\e'})^2$ and
this is the required upper estimate for $N.$

The second case is $[\N\setminus\{1\}]^2\subset A.$ As in the
first case, if $1< n_1< n_2 <k$ and $s_{n_1, k} \neq s_{n_2, k}$
then $\phi_{n_1, k}^{[n_1, k]}$ and $\phi_{n_2, k}^{[n_2, k]}$ are
incomparable and since $s_{n_1, k} \neq s_{n_2, k}$ they also have
different indices. As in the first case let $s_{n_j, k},$ $1\le j
\le N$ all be different from each other and set $z_{n_j}^* E_{n_j,k}y_{n_j,k}^*, $ $1\le j\le N$ where in this case $E_{n,k}
= [\min \supp \phi_{n_j, k}^{[n_j,k]}, \infty).$ By our previous
observation it follows that these $\sigma_{\mc{F}}$ special
functionals have pairwise disjoint indices. Notice also that
$|z_{n_j}^*(x_k)| = |y_{n_j, k}^*(x_k)|\ge \e.$ Therefore exactly
as in the first case we obtain an upper estimate for $N$
independent of $k$ and this finishes the proof of the claim.

In the case where $[\N\setminus\{1\}]^2 \subset B,$
$|s_{n,k}^*(x_n)| = |y_{n,k}^* (x_n)| \ge \e.$

In the case where $[\N\setminus\{1\}]^2\subset A,$
$|s_{n,k}^*(x_n)| = |y_{n,k}^*(x_n) - \phi_{n,k}^{[n, k]}(x_n)|
\ge \e -\e'.$

In any case we have that $|s_{n,k}^*(x_n)| \ge \e-\e'
>0.$

Combining this with the previous claim, we get that for any $k\ge
3$ there are $z_{1, k}^*, \dots, z_{M,k}^* \in
\overline{\mathscr{S}}^{w^*}$ such that for any $1<n<k$ there is
$i\in [1, M]$ so that $|z_{i,k}^* (x_n) | \ge \e-\e'.$

Since $\overline{\mathscr{S}}^{w^*}$ is weak-* compact we can pass
to an $L\in [\N]$ such that $(z_{i,k}^*)_{k\in L}$ is weak-*
convergent to some $z_i^* \in\overline{\mathscr{S}}^{w^*}.$ It is
easy to see that in this case, for any $n\in\N$ there is an $i\in
[1, M]$ so that $|z_i^*(x_n)|\ge \e-\e'$. Therefore there exists
an infinite subset $P$ of $\N$ and $1\le i_0\le M$ such that
$|z_{i_0}^*(x_n)|\ge \e-\e'$ for every $n\in P$.
 It  also follows that $z_{i_0}^*$ is an infinite
 $\sigma_{\mc{F}}$ special functional.
  These contradict our assumption that $\lim\limits_n
  b^*(x_n)=0$ for
  every infinite $\sigma_{\mc{F}}$ special functional $b^*$.
\end{proof}

We now prove the following lemma about $JT_{\mc{F}_2}:$

\begin{lemma} \label{aa3} Let $x\in JT_{\mc{F}_2}$ with finite
support and $\e >0.$ There exists $n\in \N$ such that if $y^*\sum_{k=1}^da_ky_k^*\in \mc{F}_2$ with $\max\{|a_k|: 1\le k\le d\}
<\frac{1}{n},$ then $|y^*(x)|<\e.$
\end{lemma}

\begin{proof}[\bf Proof] Let $\delta= \frac{\e}{2\sum_{n\in \N}|x(n)|}$
and $\tau_j=\sup\{\|f\|_\infty : f\in F_j\}.$ Since
$\sum_{j=1}^\infty \tau_j\le 1,$ by the definition of a $JTG$
family, there is a $j_0 \in\N$ such that $\sum_{j=j_0+1}^\infty
\tau_j< \delta.$ Let $n$ be such that $\frac{1}{n}< \frac{\e}{2j_0
\|x\|}.$

Assume that $y^* = \sum_{k=1}^d a_ky_k^* \in \mc{F}_2$ with $\max
\{|a_k|:\; 1\le k\le d\}< \frac{1}{n}.$ For every $k\in [1, d]$
let $y_k^* = y_{k,1}^* + y_{k,2}^*$ with $\ind(y_{k,1}^*) \subset
\{1, \dots, j_0\}$ and $\ind(y_{k,2}^*) \subset \{j_0+1, j_0+2,
\dots \}.$ So we may write $y^*= \sum_{k=1}^d a_k y_{k,1}^* +
\sum_{k=1}^d a_k y_{k,2}^*.$ Notice now that for any $n\in\N,$
$|\sum_{k=1}^d a_k y_{k,2}^*(n)| \le \sum_{k=1}^d
\|y_{k,2}^*\|_\infty$ and since $(\ind (y_{k,2}^*))_{k=1}^d$ are
pairwise disjoint and all greater than $j_0$ we get that
$\sum_{k=1}^d \|y_{k,2}^*\|_\infty < \delta.$ Therefore
$\|\sum_{k=1}^d a_k y_{k,2}^* \|_\infty < \delta $ and it follows
that
\begin{equation} \label{aa3.1} \Big| \sum_{k=1}^d a_k y_{k,2}^*(x)
\Big| \le \sum_{n\in \N} \delta |x(n)| = \frac{\e}{2}.
\end{equation}

On the other hand since  $(y_{k,1}^*)_{k=1}^d$ have pairwise
disjoint indices,
 at most $j_0$
 of them are non-zero  and
$|y_{k,1}^*(x)| \le \|x\|.$ Therefore $|\sum_{k=1}^d a_k y_{k,1}^*
(x)| \le j_0 \frac{1}{n} \|x\| < \frac{\e}{2}.$ Combining this
with \eqref{aa3.1} we get that $|y^*(x)|<\e$ as required.
\end{proof}

We combine now Lemmas \ref{aa2} and \ref{aa3} to prove the
following:

\begin{proposition} \label{aa4} Let $(x_n)_n$ be a weakly null
separated sequence in $JT_{\mc{F}_2}$ with $\|x_n\|_{\mc{F}_2}\le
C$ for all $n$. Then for all $m\in \N,$ there is $L\in [\N]$ such
that for every $y^*\in \mc{F}_2,$
\[\#\{ n\in L: |y^*(x_n)| \ge \frac{1}{m} \} \le 66 m^2C^2.\]
\end{proposition}

\begin{proof}[\bf Proof.] We may and do assume that
$\{x_n : n\in \N\}$ is normalized. We set $\delta_1= \frac{1}{4m}$
and we find $L_1\in [\N]$ such that
 Lemma \ref{aa2} is valid for $\e=\delta_1$. Then we set $n_1=\min
 L_1$ and using Lemma \ref{aa3} we find $n=r_1\in \N$ such
that the conclusion of Lemma \ref{aa3} is valid for $\e =\delta_1$
and $x= x_{n_1}$.  Then, after setting $\delta_2
=\min\{\frac{1}{8mr_1^2},\delta_1\}$ we find $L_2 \in
[\N\setminus\{n_1\}]$ such that Lemma \ref{aa2} is valid for
$\e=\delta_2$. We set $n_2 =\min L_2$ and we find $n=r_2\in\N$
with $r_2>r_1$ such that the conclusion of Lemma \ref{aa3} is
valid for $\e =\delta_2$ and $x\in\{x_{n_1},x_{n_2}\}$.

 Recursively, having defined $\delta_1\ge \delta_2\ge \cdots \ge
 \delta_{p-1}$,
$L_1\supset L_2 \supset L_3 \supset \cdots \supset L_{p-1}$,
$n_1<n_2<\cdots<n_{p-1}$ and
 $r_1 <r_2 < \cdots < r_{p-1}$, we set
 $\delta_{p}=\min\{\frac{1}{4m 2^{p-1}r_{p-1}^2},\delta_{p-1}\}$
 and
 we find $L_p \in [L_{p-1}\setminus\{n_{p-1}\}]$ such that Lemma \ref{aa2}
is valid for $\e=\delta_1$. We set  $n_p=\min L_p$
 and we find $n=r_p>r_{p-1}$
 that the conclusion of Lemma \ref{aa3} is satisfied for $\e =\delta_{p}$
and  $x\in\{x_{n_1}, x_{n_2}, \dots, x_{n_p}\}$.
 At the end we
consider the set $L= \{n_1 < n_2 < \cdots <n_p<\cdots\}.$

The crucial properties of this construction are the following:
\begin{enumerate}
 \item If $\sum_{k=1}^\ell a_k y_k^* \in \mc{F}_2$ and $|a_k| <
 \frac{1}{r_{p}},$ for $k=1, \dots, \ell$ then we have that $|\sum_{k=1}^\ell
 a_k y_k^* (x_{n_i})|< \delta_{p}$ for all $i=1, \dots, p.$
 \item For every  $x^* \in \mathscr{S}$ we have that $\# \{i\ge p :\; |x^*(x_{n_i})|
 \ge \delta_p \} \le 2.$
\end{enumerate}

We will make use of these two properties to prove the proposition.

Let $d= 66 m^2.$ It suffices to prove that if $n_{\ell_1}<
n_{\ell_2}< \cdots < n_{\ell_d}$ and $y^* = \sum_{k=1}^\ell a_k
y_k^* \in \mc{F}_2,$ then there is an $1\le i\le d$ such that
$|y^*(x_{n_{\ell_i}})| <\frac{1}{m}.$

We set \[ \begin{split} A_1 & = \{ k\in [1,\ell]:\; |a_k| \ge
\frac{1}{r_{\ell_1}}\},\\
A_p & = \{k\in [1, \ell]:\; \frac{1}{r_{\ell_p}}\le |a_k|<
\frac{1}{r_{\ell_{p-1} } }\}, \quad \text{for $1< p < d$}\\A_d &=
\{k\in [1, \ell] :\; \frac{1}{r_{\ell_d}} > |a_k| \}. \end{split}
\]

Observe that for $p<d$, we have that $\#A_p\le r_{\ell_p}^2$. By
property (1) we have that for any $p\in [1, d),$
\begin{equation} \label{aa4.1} \Big|\sum_{k\in \cup_{j>p}A_j} a_k
y_k^*(x_{n_{\ell_p}})\Big| < \delta_{\ell_p} \le \delta_1 \frac{1}{4m}. \end{equation}

Next, for $1\le j<p\le d$ we set \[B_{j,p}=\{k\in A_j:\;
|y_k^*(x_{n_{\ell_p}})| \ge\delta_{{\ell_j}+1}\}\] and then for
$p\in (1,d]$ we define
\[  B_p=\bigcup\limits_{j<p}B_{j,p}.\]

Since for $p>j$ we have that $\ell_p\ge \ell_j+1$, property (2)
yields that for every $k\in A_j$ there exist at most two
$B_{j,p}$'s
 containing $k$. Hence every $k\in\{1,\ldots,\ell\}$ belongs to
at most two $B_p$'s.

 Next we shall estimate the term $\sum_{k\in \cup_{j<p}A_j
\setminus B_p}| a_k y_k^* (x_{n_{\ell_p}})|.$ Let $p\in (1,d]$. We
have  that
\begin{equation} \label{mistake} \begin{split}
\sum_{k\in \cup_{j<p}A_j\setminus B_p} |a_ky_k^*(x_{n_{\ell_p}})|
 & = \sum\limits_{j=1}^{p-1}\sum\limits_{k\in A_j\setminus B_{j,p}}
 |a_ky_k^*(x_{n_{\ell_p}})| \\
 & \le \sum\limits_{j=1}^{p-1}\#(A_j) \delta_{\ell_j+1}\le
 \sum\limits_{j=1}^{p-1} r_{\ell_j}^2\frac{1}{4m
 2^{\ell_j}r_{\ell_j}^2}<\frac{1}{4m}.
\end{split} \end{equation}

We now argue  that for at least $\frac{d}{2} +1$ many of
$\{x_{n_{\ell_1}}, \dots, x_{n_{\ell_d}}\},$ we have that
$|\sum_{k\in A_p} a_k y_k^*( x_{n_{\ell_p}})| < \frac{1}{4m}.$ If
this is not the case, then for at least $\frac{d}{2}$ many,
$|\sum_{k\in A_p} a_k y_k^* (x_{n_{\ell_p}})|\ge \frac{1}{4m}$ and
therefore $\sum_{k\in A_p} a_k^2 \ge \frac{1}{16 m^2}.$ Thus \[
\sum_{k=1}^\ell a_k^2 \sum_{p=1}^d \sum_{k\in A_p} a_k^2 \ge
\frac{d}{2} \cdot \frac{1}{16m^2} = \frac{d}{32 m^2}
\] which is a contradiction since $d= 66 m^2$ and $\sum_{k=1}^\ell
a_k^2\le 1.$

Now we shall prove that for at least $\frac{d}{2} +1$ many of
$\{x_{n_{\ell_1}}, \dots, x_{n_{\ell_d}}\},$ $|\sum_{k\in B_p}a_k
y_k^* (x_{n_{\ell_p}})| < \frac{1}{4m}.$ Again if this is not the
case, then for at least $\frac{d}{2}$ many $|\sum_{k\in
B_p}a_ky_k^*(x_{n_{\ell_p}})|\ge\frac{1}{4m}$ and therefore
$\sum_{k\in B_p} a_k^2 \ge \frac{1}{16 m^2}.$ Since every $k$
appears in at most two $B_p$'s, we have that
\[ 2\ge 2\sum_{k=1}^\ell a_k^2 \ge \sum_{p=1}^d \sum_{k\in
B_p}a_k^2 \ge \frac{d}{2} \cdot \frac{1}{16m^2} = \frac{d}{32 m^2}
\] which is a contradiction.

These last two observations show that there exists at least one
$p\in [1, d]$ such that both $|\sum_{k\in
A_p}a_ky_k^*(x_{n_{\ell_p}})|<\frac{1}{4m}$ and $|\sum_{k\in B_p}
a_k y_k^*(x_{n_{\ell_p}})|< \frac{1}{4m}.$

Combining this with \eqref{aa4.1} and \eqref{mistake}, we get that
for this particular $p,$
\[ \begin{split} \Big|\sum_{k=1}^\ell a_k
y_k^*(x_{n_{\ell_p}})\Big| \le & \Big|\sum_{k\in \cup_{j>p}A_j}
a_k y_k^*(x_{n_{\ell_p}})\Big| + \Big|\sum_{k\in A_p} a_k
y_k^*(x_{n_{\ell_p}})\Big| \\ & + \Big|\sum_{k\in B_p} a_k
y_k^*(x_{n_{\ell_p}})\Big| + \Big| \sum_{k\in \cup_{j<p}A_j
\setminus B_p} a_k y_k^*(x_{n_{\ell_p}})\Big|\\ < & \frac{1}{m}
\end{split} \] as required.
\end{proof}

\section{The space $(\eqs_{\mc{F}_{2}})_*$ and the space
of the operators $\mc{L}((\mathfrak X_{\mc{F}_2})_*)$}\label{secr}

In this section we proceed to construct a HI space not containing
a reflexive subspace. This space is  $(\mathfrak X_{\mc{F}_2})_*$
where $\mathfrak X_{\mc{F}_2}$ is the strongly strictly singular
HI extension (Sections 1 and 2) of the set $\mc{F}_2.$ The set
$\mc{F}_2$ is defined from a family $F=(F_j)_j$ as in Section 3.
The proof that $(\mathfrak X_{\mc{F}_2})_*$ does not contain a
reflexive subspace, uses the method of attractors and the key
ingredient is the attractor functional and the attracting
sequences introduced in Section 1. The structure of the quotients
of $\mathfrak X_{\mc{F}_2}$ is also investigated.

{\bf The family $F = (F_j)_{j\in \N}$}

We shall use the sequence of positive integers $(m_j)_j,$
$(n_j)_j$ introduced in Definition \ref{def13} of strictly
singular extensions which for convenience we recall:

 \begin{itemize}
 \item  $m_1=2$  and $m_{j+1}=m_j^5$.
 \item $n_1=4$,  and $n_{j+1}=(5n_j)^{s_j}$  where  $s_j=\log_2 m_{j+1}^3$.
 \end{itemize}
 We set $F_0=\{\pm e_n^*:\;n\in\N\}$ and for
  $j=1,2,\ldots$ we set
 \[ F_j=\big\{\frac{1}{m_{4j-3}^2}\sum\limits_{i\in I}\pm e^*_i:\;
  \#(I)\le \frac{n_{4j-3}}{2} \big\}\cup\big\{0\big\}. \]

In the sequel we shall denote by $\mathfrak X_{\mc{F}_2}$ the HI
extension of $JT_{\mc{F}_2}$ with ground set  $\mc{F}_2$ defined
by  the aforementioned  family \seq{F}{j} as in Definition
\ref{def17}.

\begin{proposition} \label{propa1} The space $\mathfrak
X_{\mc{F}_2}$ is a strongly strictly singular extension of
$JT_{\mc{F}_2}(=Y_{\mc{F}_2}).$
\end{proposition}
\begin{proof}[\bf Proof.]
Let $C>0$. We select $j(C)$ such that
$\frac{33}{2}m_{2j}^4C^2<n_{2j-1}$ for every $j\ge j(C)$ and we
shall show that the integer $j(C)$ satisfies the conclusion of
Definition \ref{aa11}.

Let \seq{x}{n} be a block sequence in $\eqs_{\mc{F}_2}$ such that
$\|x_n\|\le C$ for all $n$, $\|x_n\|_{\infty}\to 0$ and \seq{x}{n}
is a weakly null sequence in $JT_{\mc{F}_2}$.
 It suffices to show that the sequence \seq{x}{n} is separated
 (Definition \ref{aa1}). Indeed, then Proposition \ref{aa4} and
 our choice of $j(C)$ yield that for every $j\ge j(C)$ there exists $L\in [\N]$
 such that for every $y^*\in\mc{F}_2$ we have that
 $\#\{n\in L:\;|y^*(x_n)|>\frac{2}{m_{2j}^2}\}\le 66
 (\frac{m_{2j}^2}{2})^2C^2<n_{2j-1}$.

 In order to show that the sequence \seq{x}{n} is separated we start with
 the following easy observations:
 \begin{enumerate}
 \item[(i)] If $m_{4j_0-3}^2>\frac{C}{\e}\#\supp(x)$ and $\|x\|\le
 C$ then for every $\phi\in\bigcup\limits_{j\ge j_0}F_j$ we have that
 $|\phi(x)|\le \e$.
  \item[(ii)] If $\|x\|_{\infty}<\frac{2\e}{n_{4j_0-3}}$  and
  $\phi\in\bigcup\limits_{j\le j_0}F_j$ then $|\phi(x)|\le \e$.
 \end{enumerate}

  Let $L\in [\N]$ and $\e>0$. Using (i) and (ii) we may
  inductively select $1=j_0<j_1<j_2<\cdots$ in  $\N$ and
  $k_1<k_2<\cdots$ in $L$ such that for each $i$ and
  $\phi\in \bigcup\limits_{j\not\in [j_{i-1},j_i)} F_j$ we have
  that $|\phi(x_{k_i})|<\e$. Setting $M=\{k_1,k_2,\ldots\}$ we
  have that the sequence $(x_n)_{n\in M}$ is $\e$-separated.
  Therefore the sequence \seq{x}{n} is separated.
\end{proof}

A consequence of the above proposition and the results of Sections
1 and 2 is the following:

\begin{theorem}
\begin{enumerate}
\item[(a)] The space $\mathfrak X_{\mc{F}_2}$ is HI and
reflexively saturated.

\item[(b)] The predual $(\mathfrak X_{\mc{F}_2})_*$ is HI.

\item[(c)] Every bounded linear operator $T: \mathfrak
X_{\mc{F}_2} \to \mathfrak X_{\mc{F}_2}$ is of the form $T=\lambda
I +S,$ with $S$ strictly singular and weakly compact.

\item[(d)] Every bounded linear operator $T: (\mathfrak
X_{\mc{F}_2})_* \to (\mathfrak X_{\mc{F}_2})_*$ is of the form
$T=\lambda I +S,$ with $S$ strictly singular.

\end{enumerate}
\end{theorem}

\begin{proof}[\bf Proof.]
 All the above properties are consequences of the fact that
 $\eqs_{\mc{F}_2}$ is a strongly strictly singular extension of
 $JT_{\mc{F}_2}$. In particular (a) follows from Proposition \ref{prop15}
 and Theorem \ref{theo1}, (b) follows from Theorem \ref{th6.1},
 (c) follows from Theorem \ref{th7} while (d) follows from Theorem
 \ref{aa12}.
\end{proof}

 \begin{proposition}\label{prop18}
 Let $(x_k,x_k^*)_{k=1}^{n_{4j-3}}$ be a $(18,4j-3,1)$ attracting
 sequence in $\eqs_{\mc{F}_{2}}$ such that
 $\|x_{2k-1}\|_{\mc{F}_{2}}\le\frac{2}{m_{4j-3}^2}$ for
 $k=1,\ldots,n_{4j-3}/2$. Then
 \[
 \|\frac{1}{n_{4j-3}}\sum\limits_{k=1}^{n_{4j-3}}(-1)^{k+1}x_k\|
 \le \frac{144}{m_{4j-3}^2}.  \]
 \end{proposition}
 \begin{proof}[\bf Proof.] The conclusion follows by an
 application of Proposition
 \ref{prop11} (ii) after checking  that for every $g\in \mc{F}_{2}$ we have that
 $|g(x_{2k})|>\frac{2}{m_{4j-3}^2}$ for at most $n_{4j-4}$
 $k$'s. From the fact that
 each $x_{2k}$ is of the form $e_l$ it suffices to show that
 for every $g\in \mc{F}_{2}$, the cardinality of the set
 $\{l:\;|g(e_l)|>\frac{2}{m_{4j-3}^2}\}$ is at most $n_{4j-4}$.

 Let $g\in \mc{F}_{2}$, $g=\sum\limits_{i=1}^da_ig_i$ where
 $\sum\limits_{i=1}^da_i^2\le 1$ and $(g_i)_{i=1}^d$ are $\sigma_{\mc{F}}$
 special  functionals with disjoint indices.
 For each $i$ we divide the functional $g_i$ into two parts,
  $g_i=y_i^*+z_i^*$, with
 $\ind(y_i^*)\subset \{1,\ldots,j-1\}$
 and $\ind(z_i^*)\subset \{j,j+1,\ldots\}$.
 For $l\not\in\bigcup\limits_{i=1}^d\supp(y_i^*)$
 we have that $|g(e_l)|\le \sum\limits_{i=1}^d|z_i^*(e_l)|
 <\sum\limits_{r=j}^{\infty}\frac{1}{m_{4r-3}^2}<\frac{2}{m_{4j-3}^2}$.
 Since
 $\#\big(\bigcup\limits_{i=1}^d\supp(y_i^*)\big)\le
  \frac{n_1}{2}+\frac{n_5}{2}+\cdots\frac{n_{4j-7}}{2}<n_{4j-5}$
  the conclusion  follows.

  The proof of the proposition is complete.
 \end{proof}

 \begin{definition}\label{def12}
 Let $\chi=(x_k,x_k^*)_{k=1}^{n_{4j-3}}$ be a $(18,4j-3,1)$
 attracting
 sequence, with $\|x_{2k-1}\|_{\mc{F}_{2}}\le \frac{2}{m_{4j-3}^2}$ for
 $1\le k\le n_{4j-3}/2$. We set
 \begin{align*}
 g_{\chi} &
 =\frac{1}{m_{4j-3}^2}\sum\limits_{k=1}^{n_{4j-3}/2}x_{2k}^*\\
  F_{\chi} &
 =-\frac{1}{m_{4j-3}^2}\sum\limits_{k=1}^{n_{4j-3}/2}x_{2k-1}^*\\
  d_{\chi} &
  =\frac{m_{4j-3}^2}{n_{4j-3}}\sum\limits_{k=1}^{n_{4j-3}}(-1)^kx_{k}
 \end{align*}
 \end{definition}

 \begin{lemma}\label{lem21} If $\chi$ is a $(18,4j-3,1)$
 attracting
 sequence,  $\chi=(x_k,x_k^*)_{k=1}^{n_{4j-3}}$,
 with $\|x_{2k-1}\|_{\mc{F}_{2}}\le \frac{2}{m_{4j-3}^2}$ for
 $1\le k\le n_{4j-3}/2$ then
 \begin{enumerate}
 \item[(1)] $\|g_{\chi}-F_{\chi}\|\le \frac{1}{m_{4j-3}}$.
 \item[(2)] $\frac{1}{2}=g_{\chi}(d_{\chi})\le \|d_{\chi}\|\le
 144$, and hence $\|g_{\chi}\|\ge \frac{1}{288}$.
%\item[(3)] $F_{\chi}(d_{\chi})=\frac{1}{2}$.
 \end{enumerate}
 \end{lemma}
 \begin{proof}[\bf Proof.]
 (1) We have $g_{\chi}-F_{\chi}
 =\frac{1}{m_{4j-3}}
 \big(\frac{1}{m_{4j-3}}\sum\limits_{k=1}^{n_{4j-3}}x_{k}^*\big)$.
 Since $(x_k^*)_{k=1}^{n_{4j-3}}$ is a special sequence of length
 $n_{4j-3}$, the functional
 $\frac{1}{m_{4j-3}} \sum\limits_{k=1}^{n_{4j-3}}x_{k}^*$
 belongs to $D_G$ and hence to $B_{\eqs_{\mc{F}_{2}}^*}$
 The conclusion follows.

 (2) It is straightforward from Definitions \ref{depseq} and
 \ref{def12} that
 $g_{\chi}(d_{\chi})=\frac{1}{2}$. Since
 $\|x_{2k-1}\|_{\mc{F}_{2}}\le \frac{2}{m_{4j-3}^2}$ for
 $k=1,\ldots, n_{4j-3}/2$, Proposition \ref{prop18} yields that
 $\|d_{\chi}\|\le 144$. Thus
 $\|g_{\chi}\|\ge \frac{g_{\chi}(d_{\chi})}{\|d_{\chi}\|}
 \ge\frac{\frac{1}{2}}{144}=\frac{1}{288}$.
 %(3) It follows readily from the fact that $x_{2k-1}^*(x_{2k-1})=1$
 %for $1\le k\le n_{4j-3}/2$, since $(x_k,x_k^*)_{k=1}^{n_{4j-3}}$
 %is a $(18,4j-3,1)$ attracting sequence.
 \end{proof}

\begin{lemma}\label{lem20}
 Let $Z$ be a block subspace of $(\eqs_{\mc{F}_{2}})_*$. Also, let
 $\e>0$ and  $j>1$.
  There exists a $(18,4j-3,1)$ attracting sequence
 $\chi=(x_k,x_k^*)_{k=1}^{n_{4j-3}}$
 with $\sum\limits_{k=1}^{n_{4j-3}/2}\|x_{2k-1}\|_{\mc{F}_{2}}<\frac{1}{n_{4j-3}^2}$
  and $\dist(F_{\chi},Z)<\e$.
 \end{lemma}
 \begin{proof}[\bf Proof.]
 We select an integer $j_1$ such that
 $m_{2j_1}^{\frac{1}{2}}>n_{4j-3}$.
 From Lemma \ref{lem25} we may select a $(18,2j_1,1)$ exact pair
 $(x_1,x_1^*)$  with $\dist(x_1^*,Z)<\frac{\e}{n_{4j-3}}$ and
 $\|x_1\|_{\mc{F}_{2}}\le\frac{2}{m_{2j_1}}$.
 Let $2j_2=\sigma(x_1^*)$. We select  $l_2\in\Lambda_{2j_2}$ and
 we set  $x_2=e_{l_2}$ and
 $x_2^*=e_{l_2}^*$.

 We then set $2j_3=\sigma(x_1^*,x_2^*)$ and we select, using Lemma
 \ref{lem25}, a
 $(18,2j_3,1)$ exact pair
 $(x_3,x_3^*)$  with $x_2<x_3$, $\dist(x_3^*,Z)<\frac{\e}{n_{4j-3}}$ and
 $\|x_3\|_{\mc{F}_{2}}\le\frac{2}{m_{2j_3}}$.

 It is clear that we may inductively construct a
 $(18,4j-3,1)$ attracting sequence
 $\chi=(x_k,x_k^*)_{k=1}^{n_{4j-3}}$ such that
 $\sum\limits_{k=1}^{n_{4j-3}/2}\|x_{2k-1}\|_{\mc{F}_{2}}
 \le\sum\limits_{k=1}^{n_{4j-3}/2}\frac{2}{m_{2j_{2k-1}}}<\frac{1}{n_{4j-3}^2}$
 and $\dist(x_{2k-1}^*,Z)<\frac{\e}{n_{4j-3}}$
 for $1\le k\le n_{4j-3}/2$.
 It follows that
 $\dist(F_{\chi},Z)\le \frac{1}{m_{4j-3}^2}
 \sum\limits_{k=1}^{n_{4j-3}/{2}}
 \dist(x_{2k-1}^*,Z)<\e$.
 \end{proof}

 \begin{theorem}\label{th3}
 The space $(\eqs_{\mc{F}_{2}})_*$ is a Hereditarily James Tree
 (HJT)  space. In particular it does not contain
 any reflexive subspace and
 every infinite dimensional  subspace $Z$ of
 $(\eqs_{\mc{F}_{2}})_*$ has nonseparable second dual $Z^{**}$.
 \end{theorem}
 \begin{proof}[\bf Proof.]
  Since each subspace of  $(\eqs_{\mc{F}_{2}})_*$ has a further
  subspace
 isomorphic to a block subspace it is enough to consider a block
 subspace $Z$ of $(\eqs_{\mc{F}_{2}})_*$ and to show that
 $Z$ has the James Tree property.

 We select a  $j_{\emptyset}\in \Xi_1$ with $j_{\emptyset}\ge 2$.
 We shall inductively
 construct a family $(\chi_a)_{a\in\mc{D}}$ of attracting
 sequences and a family
 $(j_a)_{a\in\mc{D}}$ of integers such that
 \begin{enumerate}
 \item[(i)] If $a<_{lex}\beta$ then $d_{\chi_a}<d_{\chi_\beta}$.
 \item[(ii)] For every $\beta\in\mc{D}$,
 $\chi_{\beta}
 =\big(x^{\beta}_k,(x^{\beta}_k)^*\big)_{k=1}^{n_{4j_{\beta}-3}}$
 is a  $(18,4j_{\beta}-3,1)$ attracting sequence with
 $\dist(F_{\chi_{\beta}},Z)<\frac{1}{m_{4j_{\beta}-3}}$
 and $\|x_{2k-1}^{\beta}\|_{\mc{F}_{2}}\le\frac{2}{m_{4j_{\beta}-3}^2}$
  for
 $k=1,\ldots,n_{4j_{\beta}-3}/2$.
 \item[(iii)] If  $\beta\in\mc{D}$ with $\beta\neq\emptyset$
  then
 $j_{\beta}=\sigma_{\mc{F}}((g_{\chi_{a}})_{a<\beta})$.
 \end{enumerate}

 The induction runs on the lexicographical ordering of $\mc{D}$.
 In the first step, i.e. for $\beta=\emptyset$, we select
 a  $(18,2j_{\emptyset}-1,1)$ attracting sequence with
 $\dist(F_{\chi_{\emptyset}},Z)<\frac{1}{m_{4j_{\emptyset}-3}}$
 and $\|x_{2k-1}^{\emptyset}\|_{\mc{F}_{2}}
 \le\frac{2}{m_{4j_{\emptyset}-3}^2}$ for
 $k=1,\ldots,n_{4j_{\emptyset}-3}/2$.
 In the general inductive step, we assume that
 $(j_a)_{a<_{lex}\beta}$ and $(\chi_a)_{a<_{lex}\beta}$ have been
 constructed for some $\beta\in\mc{D}$.
 Since
 $\{a\in\mc{D}:\;a<\beta\}\subset\{a\in\mc{D}:\;a<_{lex}\beta\}$,
 the attracting sequences $(\chi_a)_{a<\beta}$ have already been
 constructed so we may set
 $j_{\beta}=\sigma_{\mc{F}}((g_{\chi_{a}})_{a<\beta})$.
 Denoting by $\beta^-$ the immediate predecessor of $\beta$ in the
 lexicographical ordering, we select, using Lemma \ref{lem20},
 a $(18,2j_{\beta}-1,1)$ attracting sequence
 $\chi_{\beta}
 =\big(x^{\beta}_k,(x^{\beta}_k)^*\big)_{k=1}^{n_{4j_{\beta}-3}}$
 with $d_{\chi_{\beta^-}}<d_{\chi_{\beta}}$ such that
 $\dist(F_{\chi_{\beta}},Z)<\frac{1}{m_{4j_{\beta}-3}}$
 and $\|x_{2k-1}^{\beta}\|_{\mc{F}_{2}}\le\frac{2}{m_{4j_{\beta}-3}^2}$ for
 $1\le k\le n_{4j_{\beta}-3}/2$. The inductive construction is
 complete.

 For each branch $b$ of the dyadic tree
 the sequence  $(g_{\chi_{a}})_{a\in b}$ is a
 $\sigma_{\mc{F}}$ special sequence. Thus the series
 $\sum\limits_{a\in\mc{D}}g_{\chi_{a}}$
 converges in the $w^*$ topology
 to  a $\sigma_{\mc{F}}$ special functional $g_b\in
 \overline{G}^{w^*}\subset \overline{D_G}^{w^*}= B_{\eqs_{\mc{F}_{2}}^*}$.

 For each $\beta\in \mc{D}$ we select a
 $z^*_\beta\in Z$ such that
 $\|z_\beta^*-F_{\chi_{\beta}}\|<\frac{1}{m_{4j_{\beta}-3}}$.
 Then Lemma \ref{lem21} (1) yields that
 $\|z_\beta^*-g_{\chi_{\beta}}\|\le \|z_\beta^*-F_{\chi_{\beta}}\|+
 \|F_{\chi_{\beta}}-g_{\chi_{\beta}}\|<\frac{2}{m_{4j_{\beta}-3}}$.
 Now let $b$ be a branch of the dyadic tree.
 Since
 $\sum\limits_{a\in b}\|z^*_{a}-g_{\chi_{a}}\|<
  \sum\limits_{a\in b}\frac{2}{m_{4j_{a}-3}}<
  \frac{3}{m_{4j_{\emptyset}-3}}<\frac{1}{1152}$
 it follows that the series
 $\sum\limits_{a\in\mc{D}}z^*_{a}$
 is also $w^*$ convergent and its $w^*$ limit $z^*_b\in Z^{**}$ satisfies
 $\|z^*_b-g_b\|<\frac{1}{1152}$. This actually yields that the
 block sequence $(z^*_a)_{a\in\mc{D}}$ defines a James Tree
 structure in the subspace $Z$.

 The family $\{z_b^*:\; b\mbox{ a branch of }\mc{D}\}$ is a family
 in $Z^{**}$ with the cardinality of the continuum. We complete the proof
 of the theorem by showing that for $b\neq b'$ we have that
 $\|z_b^*-z_{b'}^*\|\ge \frac{1}{576}$.
 Let $b\neq b'$ be two branches of the dyadic tree.
 We select $a\in\mc{D}$ with $a\in b\setminus b'$ (i.e. $a$ is an
 initial part of $b$ but not of $b'$).
 Then our construction and Lemma \ref{lem21} (2) yield that
 \begin{eqnarray*}
  \|z_b^*-z_{b'}^*\|&\ge&\|g_b-g_{b'}\|-\|z_b^*-g_b\|-\|z_{b'}^*-g_{b'}\|\\
 &>&
 \frac{(g_b-g_{b'})(d_{\chi_a})}{\|d_{\chi_a}\|}-\frac{1}{1152}-\frac{1}{1152}\\
 &  \ge & \frac{g_{\chi_a}(d_{\chi_a})}{144}-\frac{1}{576}
 =\frac{\frac{1}{2}}{144}-\frac{1}{576}=\frac{1}{576}.
 \end{eqnarray*}
 \end{proof}

\begin{proposition}\label{prop28}
 For every
 block subspace
 $Y=\overline{\spann}\{y_n:\;n\in\N\}$
 of  $\eqs_{\mc{F}_{2}}$ there exist a further
 block subspace
 $Y'=\overline{\spann}\{y'_n:\;n\in\N\}$
 and a block subspace $Z=\overline{\spann}\{z_k:\;k\in\N\}$
 of  $\eqs_{\mc{F}_{2}}$ such that the following are satisfied.
 The space $Z$ is reflexive, the spaces $Y'$ and $Z$ are
 disjointly supported (i.e. $\supp z_k \cap \supp y_n'=\emptyset$
 for all $n,k$) and the space
 $X=\overline{\spann}(\{z_k:\;k\in\N\}\cup\{y_n':\;n\in\N\})$ has
 nonseparable dual.
 \end{proposition}
  \begin{proof}[\bf Proof.]
  The proof is similar to that of Theorem \ref{th3}.
  Let $Y$  be a block subspace
  of $\eqs_{\mc{F}_{2}}$. Using Proposition \ref{prop20}
  we may inductively construct (the induction runs on the
  lexicographic order of the dyadic tree $\mc{D}$) a family
   $(\chi_a)_{a\in\mc{D}}$ of attracting
 sequences and a family
 $(j_a)_{a\in\mc{D}}$ of integers such that the following
 conditions are satisfied:
 \begin{enumerate}
 \item[(i)] If $a<_{lex}\beta$ then $d_{\chi_a}<d_{\chi_\beta}$.
 \item[(ii)] For every $\beta\in\mc{D}$,
 $\chi_{\beta}
 =\big(x^{\beta}_k,(x^{\beta}_k)^*\big)_{k=1}^{n_{4j_{\beta}-3}}$
 is a  $(18,2j_{\beta}-1,1)$ attracting sequence with
 $x_{2k-1}^{\beta}\in Y$
 and $\|x_{2k-1}^{\beta}\|_{\mc{F}_{2}}
 \le\frac{2}{m_{4j_{\beta}-3}^2}$ for
 $k=1,\ldots,n_{4j_{\beta}-3}/2$.
 \item[(iii)] $j_{\emptyset}\in \Xi_1$ with  $j_{\emptyset}\ge 2$,
   while if
   $\beta\in\mc{D}$ with $\beta\neq\emptyset$,
  then
 $j_{\beta}=\sigma_{\mc{F}}((g_{\chi_{\emptyset}})_{a<\beta})$.
 \end{enumerate}

 For each $a\in \mc{D}$ we set
 $z_a=\frac{2m_{4j_a-3}}{n_{4j_a-3}}
 \sum\limits_{k=1}^{n_{4j_a-3}/2}x^a_{2k}$ and we consider the space
 $Z=\overline{\spann}\{z_a:\;a\in\mc{D}\}$.

 We first observe that for each $a\in \mc{D}$
 the functional
 $f_a=\frac{1}{m_{4j_a-3}}\sum\limits_{k=1}^{n_{4j_a-3}}(x^a_k)^*$
 belongs to $D_G\subset B_{\eqs_{\mc{F}_{2}}}$, hence
  $\|z_a\|\ge f_a(z_a)=1$. On the other hand we have that
 $\|z_a\|_{\mc{F}_{2}}\le \frac{2}{m_{4j_a-3}}$. Indeed, let
 $g=\sum\limits_{i=1}^da_ig_i\in \mc{F}_{2}$
 (i.e. $\sum\limits_{i=1}^da_i^2\le 1$ while $(g_i)_{i=1}^d$ and
 $\sigma_{\mc{F}}$ special functionals with pairwise disjoint indices).
 For each $i=1,\ldots,d$ let $g_i=y_i^*+z_i^*$ with
 $\ind(y_i^*)\subset\{1,\ldots,j_a-1\}$ and
  $\ind(z_i^*)\subset\{j_a,j_a+1,\ldots\}$. Then
  \begin{eqnarray*}
  |g(z_a)|&\le& \sum\limits_{i=1}^d |y_i^*(z_a)|+
 \sum\limits_{i=1}^d |z_i^*(z_a)|\\
 &\le & \frac{2m_{4j_a-3}}{n_{4j_a-3}}(\frac{n_1}{2}+\cdots+
  \frac{n_{4j_a-7}}{2})+
 \frac{2m_{4j_a-3}}{n_{4j_a-3}}\sum\limits_{r=j_a}^{\infty}
 \frac{n_{4j_a-3}}{2}\frac{1}{m_{4r-3}^2}\\
 &\le &\frac{2}{m_{4j_a-3}}.
 \end{eqnarray*}
  It follows that $\sum\limits_{a\in\mc{D}}
 \frac{\|z_a\|_{\mc{F}_{2}}}{\|z_a\|} \le
 \sum\limits_{a\in\mc{D}}\frac{2}{m_{4j_a-3}}<\frac{1}{2}$ which
 yields that the space $Z=\overline{\spann}\{z_a:\;a\in\mc{D}\}$
 is reflexive (see Proposition \ref{prop15}).

 For every branch $b$ of the dyadic tree, the
 functional $g_b$ which is defined to be the $w^*$ sum of the
 series
 $\sum\limits_{a\in\mc{D}}g_{\chi_{a}}$ belongs to
 $\overline{\mc{F}_{2}}^{w^*}\subset B_{\eqs_{\mc{F}_{2}}}$.
 The family $\{g_b|_X:\;b \mbox{ a branch of }\mc{D}\}$ is
 a family of $X^*$ with the cardinality of the continuum.
 For $b\neq b'$, selecting $a\in b\setminus b'$ the vector
 $d_{\chi_a}=\frac{m_{4j_a-3}^2}{n_{4j_a-3}}
 \sum\limits_{k=1}^{n_{4j_a-3}}(-1)^kx^a_k$ belongs to $Y+Z$ while
 from Lemma \ref{lem21} we have that $\|d_{\chi_a}\|\le 144$.
 Thus $\|g_b|_X-g_{b'}|_X\|_{X^*}\ge
  \frac{g_b(d_{\chi_a})-g_{b'}(d_{\chi_a})}{\|d_{\chi_a}\|}
 \ge \frac{\frac{1}{2}-0}{144}=\frac{1}{288}$.

 Therefore $X^*$ is nonseparable.
 \end{proof}

 \begin{lemma}\label{lem23}
 If $S:(\eqs_{\mc{F}_{2}})_*\to
 (\eqs_{\mc{F}_{2}})_*$ is a strictly singular operator then its
 conjugate operator  $S^*:\eqs_{\mc{F}_{2}}\to \eqs_{\mc{F}_{2}}$
 is also strictly singular.
 \end{lemma}
 \begin{proof}[\bf Proof.]
 From Theorem \ref{th7} the operator $S^*$ takes the form
 $S^*=\lambda I_{\eqs_{\mc{F}_{2}}}+W$ with $\lambda\in\R$ and
 $W:\eqs_{\mc{F}_{2}}\to \eqs_{\mc{F}_{2}}$ a strictly singular and weakly
 compact operator. We have to show that $\lambda=0$.

 The operator $W^*:\eqs^*_{\mc{F}_{2}}\to \eqs^*_{\mc{F}_{2}}$ is also weakly
 compact, while $W^*=S^{**}-\lambda I_{\eqs_{\mc{F}_{2}}^*}$ which
 yields that $W^*((\eqs_{\mc{F}_{2}})_*)\subset (\eqs_{\mc{F}_{2}})_*$.
 These facts, in conjunction to the fact that $(\eqs_{\mc{F}_{2}})_*$
 contains no reflexive subspace (Theorem \ref{th3}), imply that
 the restriction $W^*|(\eqs_{\mc{F}_{2}})_*$ is strictly singular.
 Thus, since
 $\lambda I_{(\eqs_{\mc{F}_{2}})_*}=S-W^*|(\eqs_{\mc{F}_{2}})_*$
 with both $S,W^*|(\eqs_{\mc{F}_{2}})_*$ being strictly singular, we get
 that $\lambda=0$.
 \end{proof}

 \begin{corollary}\label{cor6}
 Every bounded linear operator
 $T:(\eqs_{\mc{F}_{2}})_*\to(\eqs_{\mc{F}_{2}})_*$ takes the form $T=\lambda
 I+W$ with $\lambda\in\R$ and $W$ a weakly compact operator.
 \end{corollary}
 \begin{proof}[\bf Proof.]
 We know from Theorem \ref{aa12} that $T=\lambda I+W$ with
 $W$ a strictly singular operator. Lemma \ref{lem23} yields that
 $W^*$ is also strictly singular.  From Theorem
 \ref{th7} we get that $W^*$ is weakly compact, hence $W$ is
 weakly compact.
 \end{proof}

 \begin{theorem}\label{th9}
 Let $Z$ be a $w^*$ closed subspace of $\eqs_{\mc{F}_2}$ of infinite
 codimension such that
 for every $i=1,2,\ldots$ we have that
  \begin{equation}\label{eq18}
  \liminf\limits_{k\in \Lambda_i}\dist(e_k,Z)=0
 \end{equation}
  (\seq{\Lambda}{i}
 are the sets appearing in Definition \ref{def10}). Then
   every
 infinite dimensional  subspace of $\eqs_{\mc{F}_2}/Z$ has
 nonseparable  dual.
 \end{theorem}
 \begin{proof}[\bf Proof.]
  We denote by $Q$ the quotient operator $Q:\eqs_{\mc{F}_2}\to
  \eqs_{\mc{F}_2}/Z$ and we recall that since $Z$ is $w^*$ closed, $Z_{\bot}$ 1-norms
  $\eqs_{\mc{F}_2}/Z$.
   Let $Y$ be a
 closed subspace of $\eqs_{\mc{F}_2}$ with $Z\subset Y$ such that $Y/Z$ is
 infinite dimensional; we shall show that $(Y/Z)^{*}$ is nonseparable.

  For a given $j\in\N$  using
 Lemma \ref{lem31} and our assumption \eqref{eq18} are able to construct a $(18,4j-3,1)$
 attracting sequence $\chi=(x_k,x_k^*)_{k=1}^{n_{4j-3}}$
 such that each one of the sums $\sum\dist(x_{2k-1},Y)$,
 $\sum\dist(x_{2k-1}^*,Z_{\bot})$,  $\sum \|x_{2k-1}\|_{\mc{F}_2}$, $\sum\dist(x_{2k},Z)$
 is as small as we wish.
   Setting  $d_{\chi}^1
  =\frac{m_{4j-3}^2}{n_{4j-3}}\sum\limits_{k=1}^{n_{4j-3}/2}x_{2k-1}$
 we get that $Qd_{\chi}$ is almost equal to $Qd_{\chi}^1$ which
 almost belongs to $Y/Z$. Also $F_{\chi}$ almost belongs to
 $Z_{\bot}$, while
 $F_{\chi}(d_{\chi})=\frac{1}{2}$ and
  $\|F_{\chi}-g_{\chi}\|\le
 \frac{1}{m_{4j-3}}$.

 Using these estimates we are able to construct a
 dyadic tree
 $(\chi_a)_{a\in\mc{D}}$ of attracting
 sequences and a family
 $(j_a)_{a\in\mc{D}}$ of integers satisfying
 \begin{enumerate}
 \item[(i)] If $a<_{lex}\beta$ then $d_{\chi_a}<d_{\chi_\beta}$.
 \item[(ii)] For every $\beta\in\mc{D}$,
 $\chi_{\beta}
 =\big(x^{\beta}_k,(x^{\beta}_k)^*\big)_{k=1}^{n_{4j_{\beta}-3}}$
 is a  $(18,4j_{\beta}-3,1)$ attracting sequence with
 $\dist(F_{\chi_{\beta}},Z_{\bot})<\frac{1}{m_{4j_{\beta}-3}}$,
 $\dist(Qd_{\chi_{\beta}},Y/Z)<\frac{1}{m_{4j_{\beta}-3}}$
 and $\|x_{2k-1}^{\beta}\|_{\mc{F}_{2}}\le\frac{2}{m_{4j_{\beta}-3}^2}$
  for
 $k=1,\ldots,n_{4j_{\beta}-3}/2$.
 \item[(iii)] If  $\beta\in\mc{D}$ with $\beta\neq\emptyset$
  then
 $j_{\beta}=\sigma_{\mc{F}}((g_{\chi_{a}})_{a<\beta})$, while $j_{\emptyset}\in \Xi_1$ with $j_{\emptyset}\ge
 3$.
 \end{enumerate}
 For every $\beta\in\mc{D}$ we select $H_{\beta}\in Z_{\bot}$ with
 $\|H_{\beta}-F_{\chi_{\beta}}\|<\frac{1}{m_{4j_{\beta}-3}}$ and
 then for every branch $b$ of $\mc{D}$ we denote by $h_b$  the
 $w^*$ limit of the series $\sum\limits_{\beta\in b}H_{\beta}$.

 Using the above estimates, and arguing similarly to the proof of Theorem
 \ref{th3},
 we obtain that $\{h_b|Y:\; b\;\mbox{ is a
 branch of } \mc{D}\}$ is a discrete family in $(Y/Z)^*$ and therefore
 $(Y/Z)^{*}$ is nonseparable.
 \end{proof}

 \begin{remark}
 Actually it can be shown that
 the space $\eqs_{\mc{F}_2}/Z$ is HJT.
 \end{remark}

 \begin{corollary}\label{cor7}
 There exists a partition of the basis $(e^*_n)_{n\in\N}$ of
 $(\eqs_{\mc{F}_2})_*$ into
 two sets $(e^*_n)_{n \in L_1},$ $(e^*_n)_{n\in L_2}$ such that setting
$X_{L_1}=\overline{\spann}\{e^*_n:\;n\in L_1\}$,
$X_{L_2}=\overline{\spann}\{e^*_n:\;n\in L_2\}$ both
 $X_{L_1}^*$, $X_{L_2}^*$  are HI with no reflexive
 subspace.
 \end{corollary}
 \begin{proof}[\bf Proof.]
 We choose $L_1\in [\N]$ such that the sets $\Lambda_i\cap L_1$
 and
 $\Lambda_i\setminus L_1$ are infinite for each $i$ and we set
 $L_2=\N\setminus L_1$.
 The spaces $X_{L_i}=\overline{\spann}\{e_n^*:\;n\in L_i\}$, $i=1,2$
 satisfy the desired properties.  Indeed, since  $X_{L_1}^*$ is
 isometric to $\eqs_{\mc{F}_2}/\overline\spann\{e_n:\;n\in L_2\}$,
 Theorem \ref{th9} yields that $X_{L_1}^*$  has no reflexive subspace
 while from Theorem \ref{prop30} we get that it is an HI space.
 For  $X_{L_2}^*$ the  proof is completely analogous.
 \end{proof}

 \section{The structure of $\mathfrak X^*_{\mc{F}_2}$
 and a variant of $\mathfrak X_{\mc{F}_2}$}

 In the present section the structure of $\mathfrak
 X_{\mc{F}_2}^*$ is studied. This space is not HI since for
 every  subspace $Y$ of $(\mathfrak X_{\mc{F}_2})_*$ the space
 $\ell_2$  embeds into $Y^{**}.$ We also present a variant of
 $\mathfrak X_{\mc{F}_2}$, denoted $\mathfrak X_{\mc{F}_2'}$,
 such that $\mathfrak
 X_{\mc{F}_2'}^*/(\mathfrak X_{\mc{F}_2'})_*$ is isomorphic to
 $\ell_2(\Gamma)$ which yields some peculiar results on the
 structure of $\mathfrak X_{\mc{F}_2'}$ and $\mathfrak
 X_{\mc{F}_2'}^*$.
   Another variant of $\mathfrak X_{\mc{F}_2}$
 yielding a
 HI dual not containing reflexive subspace is also discussed.
 It is well known that, in $JT$ (James tree space) the
 quotient space $JT^*/JT_*$ is isometric to $\ell_2(\Gamma)$.
 It seems
 unlikely to have the same property for $\eqs_{\mc{F}_2}$.
 The main difficulty
 concerns the absence of biorthogonality between disjoint $\sigma_{\mc{F}}$ special
 functionals.
 However the next Proposition indicates that in some cases
 phenomena analogous to those in $JT$  also occur.

 \begin{proposition} \label{aa30} Let $(b_n^*)_n$ be a disjoint family of
 $\sigma_{\mc{F}}$ special functionals each one defined by an infinite special
 sequence $b_n = (f_1^n, \dots, f_k^n, \dots ).$ Assume
 furthermore
 that for each $(n,k)$ there exists a $(18, 4j_{(n,k)}-3, 1)$
 attracting sequence
 $\chi_{(n,k)}= (x_\ell^{(n,k)},
 (x_\ell^{(n,k)})^*)_{\ell=1}^{n_{4j_{(n,k)}-3}},$
  (Definition \ref{depseq}) with $\|x_{2\ell-1}^{(n,k)}\|_{\mc{F}_2}<
 \frac{1}{n_{4j_{(n,k)}-3}^2}$ and $f_n^k=g_{\chi_{(n,k)}},$
 (Definition \ref{def12}).

 Then $(b_n^*)_n$ is equivalent to the standard $\ell_2$-basis.
 \end{proposition}

 Let's provide a short description of the proof. We start with the
 following lemma:

 \begin{lemma}\label{lem32}
 Let $\chi=(x_k, x_k^*)_{k=1}^{n_{4j-3}}$ be a
 $(18, 4j-3, 1)$ attracting sequence such that
 $\|x_{2k-1}\|_{\mc{F}_2}< \frac{1}{n_{4j-3}^2}$ for all $k$. Then for every
 $\phi \in \mc{F}_2$ of the form
  $\phi=\sum\limits_{i=1}^d a_i \phi_i$ with
 $j\not\in \cup_{i=1}^d \ind (\phi_i)$ we have that $|\phi(d_\chi)|<
 \frac{1}{m_{4j-3}}$. (Recall that $d_\chi
 =\frac{m_{4j-3}^2}{n_{4j-3}}\sum\limits_{k=1}^{n_{4j-3}}(-1)^kx_{k}$,
 see Definition  \ref{def12}).
 \end{lemma}
 \begin{proof}[\bf Proof.]
 We set $d_{\chi}^1
  =\frac{m_{4j-3}^2}{n_{4j-3}}\sum\limits_{k=1}^{n_{4j-3}/2}x_{2k-1}$
and $d_{\chi}^2
  =\frac{m_{4j-3}^2}{n_{4j-3}}\sum\limits_{k=1}^{n_{4j-3}/2}x_{2k}$.
 From  our assumption that $\|x_{2k-1}\|_{\mc{F}_2}<
 \frac{1}{n_{4j-3}^2}$ for every $k$, we get that
 $|\phi(d_{\chi}^1)|\le
 \frac{m_{4j-3}^2}{n_{4j-3}}\frac{n_{4j-3}}{2}\frac{1}{n_{4j-3}^2}$.

 If  $f\in F_i$ for some $i<j$  we have that $|f(d_{\chi}^2)|\le
 \frac{1}{m_{4i-3}^2}\frac{m_{4j-3}^2}{n_{4j-3}}\frac{n_{4i-3}}{2}$,
 while for $f\in F_i$ with $i>j$ we have that
 $|f(d_{\chi}^2)|\le
 \frac{1}{m_{4i-3}^2}\frac{m_{4j-3}^2}{n_{4j-3}}\frac{n_{4j-3}}{2}$.

 Therefore
 \begin{eqnarray*}
 |\phi(d_{\chi})| & \le & |\phi(d_{\chi}^1)|+|\sum_{i=1}^d a_i\phi_i(d_{\chi}^2)|\le
 |\phi(d_{\chi}^1)|+\sum\limits_{i=1}^d|\phi_i(d_{\chi}^2)| \\
  & \le & |\phi(d_{\chi}^1)|
  +\sum\limits_{i<j}\sup\{|f(d_{\chi}^2)|:\; f\in F_i\}
  +\sum\limits_{i>j}\sup\{|f(d_{\chi}^2)|:\; f\in F_i\}\\
 & \le &
 \frac{m_{4j-3}^2}{2n_{4j-3}^2}+\frac{m_{4j-3}^2}{n_{4j-3}}\sum\limits_{i<j}\frac{n_{4i-3}}{2m_{4i-3}^2}+
 \frac{m_{4j-3}^2}{2}\sum\limits_{i>j}\frac{1}{m_{4i-3}^2}\\
 & &<\frac{1}{m_{4j-3}}.
 \end{eqnarray*}
 \end{proof}

 The content of the above lemma is that each
 $b^*$, defined by an infinite $\sigma_{\mc{F}}$ special sequence $b$ as in the
 previous proposition, is almost biorthogonal to any other $(b')^*$
 which is disjoint from $b.$

 Next we describe the main steps in the proof of  Proposition
 \ref{aa30}.

 \begin{proof}[\bf Proof of Proposition \ref{aa30}:] The proof
 follows the main lines of the proof of Lemma 11.3 of \cite{AT1}.
 Given $(a_i)_{i=1}^d$ with $a_i\in \Q$ such that $\sum_{i=1}^d a_i^2 =1,$ we
 have that $\sum\limits_{i=1}^d a_i b_i^* \in \mc{F}_2$ hence $\|\sum\limits_{i=1}^d a_i
 b_i^*\| \le 1.$

 In order to complete the proof we shall  show that
  \begin{equation} \label{aa30.1}
 \frac{1}{1000} \le \|\sum\limits_{i=1}^d a_i b_i^*\| \end{equation}
 which yields the desired result.

 To establish \eqref{aa30.1},  we choose $k\in\N$
 with $\frac{(5n_{2k-1})^{\log_2(m_{2k})}}{n_{2k}}<\frac{\e}{4d}$
 and then we choose
 $\{l^i_t:\; 1\le i\le d,\; 1\le t\le n_{2k}\}$ such that
 setting $x_{(t,i)} = d_{\chi_{(i, l^i_t)}}$
 the following conditions are satisfied.
 First, the  sequence
 $(x_{(t,i)})_{1\le i\le d,\; 1\le t\le n_{2k}}$ ordered
 lexicographically (i.e. $(t,i)<_{lex}(t',i')$ iff $t<t'$ or $t=t'$ and
 $i<i'$) is a $(144,\e)$ R.I.S. with associated  sequence
 $4j_{(t,i)}'-3:=4j_{(i,l^i_t)}-3$ while $m_{4j_{(1,1)}'-3}>\frac{2d n_{2k}}{\e}$.

% the
% following conditions are satisfied:
% \begin{enumerate}
% \item[(i)] $m_{4j_{(1,l^1_1)}-3}^{1/2}>n_{2k}$,
% $\supp d_{\chi_{(1,l^1_t)}}< \supp d_{\chi_{(2,l^2_t)}}<\cdots <\supp
% d_{\chi_{(d,l^d_t)}}<d_{\chi_{(1,l^1_{t+1})}}$ and
% $j_{(1,l^1_t)}<j_{(2,l^2_t)}<\cdots
% <j_{(d,l^d_t)}<j_{(1,l^1_{t+1})}$ for each $t$.
% \item[(ii)] $\max\supp
% d_{\chi_{(i,l^i_t)}}\frac{1}{m_{4j_{(i+1,l^{i+1}_t)}< 5
%\item[(iii)]
% \end{enumerate}
%
% sufficiently large and
% inductively we define $\{x_{(i, t)}: 1\le i\le d, 1\le t\le k\}$
% such that $\supp x_{(i,t)} < \supp x_{(q,s)}$ when $(i,t) <_{lex}
% (q,s)$ and each $x_{(i,t)} = d_{\chi_{(i, n_t)}}$ (Definition
% \ref{def12}.)

 We set $z_i = \frac{1}{n_{2k}}\sum_{t=1}^{n_{2k}}
 x_{(t,i)}$ for $i=1,\ldots,d$.
 In order to prove \eqref{aa30.1} it is enough to show that
 \begin{enumerate}
 \item[(i)] $(\sum\limits_{r=1}^d a_r b_r^*)(\sum\limits_{i=1}^d a_i
 z_i)>\frac{1}{2}-\e$.
 \item[(ii)] $\|\sum\limits_{i=1}^d a_iz_i\|\le 288$.
 \end{enumerate}

  (i) is an easy consequence of Lemma \ref{lem32}. Indeed
 \begin{eqnarray*}
 (\sum\limits_{r=1}^d a_r b_r^*)(\sum\limits_{i=1}^d a_i
 z_i) & = & \sum\limits_{r=1}^d a_r^2
 b_r^*(z_r)+\sum\limits_{i=1}^d\sum\limits_{r\neq
 i}a_rb_r^*(z_i) \\
 & \ge &
 \frac{1}{2}-\frac{1}{n_{2k}}
 \sum\limits_{i=1}^d|a_i|\cdot|(\sum\limits_{r\neq i}a_rb_r^*)(\sum\limits_{t=1}^{n_{2k}}x_{(t,i)})|\\
    & \ge &
    \frac{1}{2}-\frac{1}{n_{2k}}\sum\limits_{i=1}^d|a_i|(\sum\limits_{t=1}^{n_{2k}}\frac{1}{4m_{4j_{(t,i)}'-3}})\\
 &    \ge & \frac{1}{2}-\frac{1}{n_{2k}}\frac{2d}{m_{4j_{(1,1)}'-3}}>\frac{1}{2}-\e.
 \end{eqnarray*}

 For each $(t,i)$ we set $k_{(t,i)}=\min\supp x_{(t,i)}$ we set and
 $s_i=\{k_{(t,i)}:\; t=1,2,\ldots, n_{2k}\}$.
 We consider the set
 \begin{eqnarray*}
  \mc{H}_2  = \{e_n^*:\;n\in\N\}& \cup &\big\{
 \sum\limits_{i=1}^d\sum\limits_j\lambda_{i,j}s_{i,j}^*:\;\lambda_{i,j}\in\Q,\;
 \sum\limits_{i=1}^d\sum\limits_j\lambda_{i,j}^2\le 1,\text{ where}   \\
   & & (s_{i,j})_j \mbox{ are disjoint subintervals of }s_i\big\}
 \end{eqnarray*}
 and the norming set $D'$ of space
 $T[\mc{H}_2,(\mc{A}_{5n_j},\frac{1}{m_j})_{j\in\N}]$.

 We also set $\tilde{z}_i = \frac{1}{n_{2k}}\sum\limits_{t=1}^{n_{2k}}
 e^*_{k_(t,i)}$ for $i=1,\ldots,d$.

\begin{claim} For every $f\in D_{\mc{F}_2}$ (where $D_{\mc{F}_2}$ is the
 norming set of the space $\eqs_{\mc{F}_2}$) there exist an $h\in
 D'$ with nonnegative coordinates such that
 $|f(\sum\limits_{i=1}^d a_iz_i)|\le 288 h(\sum\limits_{i=1}^d
 |a_i|\tilde{z}_i)+ \e$.
\end{claim}
 The proof of the above claim is obtained using similar methods to the proof of the basic inequality
 (Proposition \ref{prop17}).

 Arguing  in a similar manner to the corresponding part of
 Lemma 11.3 of \cite{AT1}
 we shall show that $h(\sum\limits_{i=1}^d
 |a_i|\tilde{z}_i)\le 1+\e$.
 We may assume that the functional $h$ admits a tree
 $T_h=(h_a)_{a\in\mc{A}}$ (see Definition \ref{tree}) such that each $h_a$ is either of type 0
 (then $h_a\in \mc{H}_2$) or of type $I$, and moreover that the coordinates of each $h_a$ are
 nonnegative.
 Let $(g_{a_s})_{s=1}^{s_0}$ be
 the functionals corresponding to the maximal elements of
  the tree $\mc{A}$. We denote by $\preceq$
 the ordering  of the tree $\mc{A}$.
 Let
 \[ A=\left\{s\in\{1,2,\ldots,s_0\}:\;\; \prod\limits_{\gamma\prec a_s}
 \df{1}{w(h_\gamma)}\le\df{1}{m_{2k}}  \right\}  \]
 \[ B=\{1,2,\ldots, s_0\}\setminus A  \]
 and set $h_A=h|_{\bigcup\limits_{s\in A}\supp g_{a_s}}$,
 $h_B=h|_{\bigcup\limits_{s\in B}\supp g_{a_s}}$.

 We have that $h_A(\tilde{z}_i)\le \df{1}{m_j}$ for each $i$ thus
 \begin{equation}\label{qu1}
 h_A(\sum\limits_{i=1}^d|a_i|\tilde{z}_i)\le
   \df{1}{m_{2k}}\sum\limits_{i=1}^d|a_i|\le \df{d}{m_j}
   <\df{\varepsilon}{2}.
 \end{equation}

 It remains to estimate the value $h_B(\sum\limits_{i=1}^d|a_i|\tilde{z}_i)$.
 We observe that
 \[   \sum\limits_{i=1}^d|a_i|\tilde{z}_i
  =\sum\limits_{t=1}^{n_{2k}}\frac{1}{n_{2k}}(\sum\limits_{i=1}^d|a_i|e_{k_{(t,i)}}).\]
 We set
 \begin{eqnarray*}
  E_1 &= & \Big\{t\in\{1,2,\ldots,n_{2k}\}:\;\mbox{ the set }
 \{k_{(t,1)},k_{(t,2)},\ldots,k_{(t,d)}\} \mbox{ is contained}
  \\
 &  &\mbox{in }\ran g_{a_s} \mbox{ for some $s\in B$
   or does not intersect any } \ran g_{a_s},\;  s\in B \Big\}\\
  E_2 & = & \{1,2,\ldots,n_{2k}\}\setminus E_1 .
  \end{eqnarray*}
 For each $s=1,2,\ldots ,s_0$ set
 $\theta_s
 =\frac{1}{n_{2k}}\#\Big\{t:\;\{l^1_t,l^2_t,\ldots,l^d_t\}\subset\ran g_{a_s}\Big\}$
 and observe that $\sum\limits_{s\in B}\theta_s\le 1$.

 We first estimate the quantity $g_{a_s}
 (\sum\limits_{t\in E_1}\frac{1}{n_{2k}}(\sum\limits_{i=1}^d|a_i|e_{k_{(t,i)}}))$
 for $s\in B$.
 Each $g_{a_s}$ being in $\mc{H}_2$ takes the form
 $g_{a_s}=\sum\limits_i\sum\limits_j \lambda_{i,j}s_{i,j}^*$.
 For  $1\le i'\le d$ we get that
 $(\sum\limits_j \lambda_{i',j}s_{i',j}^*)
 (\sum\limits_{t\in E_1}\frac{1}{n_{2k}}(\sum\limits_{i=1}^d|a_i|e_{k_{(t,i)}}))
  |\tilde{a}_{i'}| (\sum\limits_j \lambda_{i',j}s_{i',j}^*)
                 (\sum\limits_{t\in E_1}\frac{1}{n_{2k}}e_{k_{(t,i)}})$
                 $\le
      |\tilde{a}_{i'}| (\max\limits_{j} \lambda_{i',j}) \theta_s.$
      Thus
 \begin{eqnarray*}
  g_{a_s}(\sum\limits_{t\in
 E_1}\frac{1}{n_{2k}}(\sum\limits_{i=1}^d|a_i|e_{k_(t,i)}))
 & \le & \theta_s \sum\limits_{i=1}^d |\tilde{a}_i|\max\limits_j
 \lambda_{i,j}\\
 & \le & \theta_s (\sum\limits_{i=1}^d \max\limits_j \lambda_{i,j}^2)^{\frac{1}{2}}
   (\sum\limits_{i=1}^d |\tilde{a}_i|^2)^{\frac{1}{2}}\le\theta_s.
   \end{eqnarray*}
 Therefore
 \begin{equation} \label{qu2}
  h_B(\sum\limits_{t\in
  E_1}\frac{1}{n_{2k}}(\sum\limits_{i=1}^d|a_i|e_{k_{(t,i)}}))
 \le \sum\limits_{s\in B} g_{a_s}
 (\sum\limits_{t\in
 E_1}\frac{1}{n_{2k}}(\sum\limits_{i=1}^d|a_i|e_{k_{(t,i)}}))
 \le \sum\limits_{s\in B}\theta_s \le 1.
 \end{equation}

  From the definition
 of the set $E_2$, the set
 $\{k_{(t,1)},k_{(t,2)},\ldots,k_{(t,d)}\}$,  for each $t\in E_2$,
  intersects at least one
 but is not contained  in any $\ran g_{a_s}$, $s\in B$.
 Also as in the proof of Lemma \ref{lem14}
 we get that $\#(B)\le (5n_{2k-1})^{\log_2(m_{2k})}$.
 These yield that $\#(E_2)\le 2 (5n_{2k-1})^{\log_2(m_{2k})}$.
 Therefore from our choice of $k$ we derive that
 \begin{equation}\label{qu3}
 h_B(\sum\limits_{t\in
 E_2}\frac{1}{n_{2k}}(\sum\limits_{i=1}^d|a_i|e_{k_{(t,i)}}))
 \le (\sum\limits_{t\in E_2}\frac{1}{n_{2k}})(\sum\limits_{i=1}^d|a_i|)
 <  \frac{2 (5n_{2k-1})^{\log_2(m_{2k})}}{n_{2k}}<\df{\varepsilon}{2}.
 \end{equation}

 From (\ref{qu1}),(\ref{qu2}) and (\ref{qu3}), we conclude that
 \begin{eqnarray*}
 h(\sum\limits_{i=1}^d|a_i|e_{k_{(t,i)}}) & \le &
  h_A(\sum\limits_{i=1}^d|a_i|\tilde{z}_i)+
  h_B(\sum\limits_{t\in E_1}\frac{1}{n_{2k}}(\sum\limits_{i=1}^d|a_i|e_{k_{(t,i)}}))\\
  && +
  h_B(\sum\limits_{t\in E_2}\frac{1}{n_{2k}}(\sum\limits_{i=1}^d|a_i|e_{k_{(t,i)}}))
   \le  \df{\varepsilon}{2}+1+\df{\varepsilon}{2}=1+\varepsilon.
 \end{eqnarray*}
 \end{proof}

 As a consequence we obtain the following:
 \begin{theorem} For every infinite dimensional subspace $Y$ of
 $(\mathfrak X_{\mc{F}_2})_*,$ the space $\ell_2$ is isomorphic to a
 subspace of $Y^{**}.$
 \end{theorem}

 {\bf A variant of $\mathfrak X_{\mc{F}_2}$}

 Next we shall indicate how we can obtain a space $\mathfrak X_{\mc{F}_2'}$
 similar to $\mathfrak X_{\mc{F}_2}$ satisfying the additional
 property that $\mathfrak X_{\mc{F}_2'}^* /(\mathfrak X_{\mc{F}_2'})_*$ is isomorphic to
 $\ell^2(\Gamma).$ Notice that such a space has the following
 peculiar property:

 \begin{proposition}
 Granting that $\mathfrak X_{\mc{F}_2'}^* /(\mathfrak X_{\mc{F}_2'})_*$ is isomorphic to
 $\ell^2(\Gamma)$, every  infinite dimensional $w^*$-closed
 subspace $Z$ of $\mathfrak X_{\mc{F}_2'}^*$ is either  nonseparable or
 isomorphic to $\ell_2.$
 \end{proposition}

 \begin{proof}[\bf Proof] Let $Q: \mathfrak X_{\mc{F}_2'}^* \to
 \mathfrak X_{\mc{F}_2'}^*/(\mathfrak X_{\mc{F}_2'})_*$ be the quotient map.
 There are two cases.
 If there exists a subspace
 $Z' \hookrightarrow Z$ of finite codimension such
 that $Q|_{Z'}$ is an isomorphism, then $Z$ is isomorphic
 to $\ell_2$.
 If not  then
 there exists a normalized block sequence $(v_n)_{n\in\N}$ in
 $(\mathfrak X_{\mc{F}_2'})_*$ such that
 $\sum\limits_{n=1}^{\infty}\dist(v_n,Z)<\frac{1}{3456}$.
 Setting $V=\overline{\spann}\{v_n:\;n\in\N\}$ we
 observe that $\dist(S_{V},Z)\le
 \frac{1}{3456}$ hence, since
 $Z$ is $w^*$-closed,
 \begin{equation}\label{equ1}
 \dist(S_{\overline{V}^{w^*}},Z)\le
 \frac{1}{3456}.
 \end{equation}

 As in the proof of Theorem \ref{th3} we consider a James Tree structure
 $(w_a)_{a\in\mc{D}}$ in $V$
 such that the corresponding family $\{w_b:\; b\in [\mc{D}]\}$
 satisfies the following properties:
 \begin{enumerate}
 \item[(i)] $\|w_b\|\le 2$ for every $b\in
 [\mc{D}]$.
 \item[(ii)] For $b\neq b'$ in $[\mc{D}]$ we have that
 $\|w_b-w_{b'}\|\ge \frac{1}{576}$.
 \end{enumerate}
 The above (i) and \eqref{equ1} yield that for every
 $b\in [\mc{D}]$ there exists $z_b\in Z$ such that
  \begin{equation}\label{equ2}
 \|z_b-w_b\|\le \frac {1}{1728}
 \end{equation}

 From \eqref{equ2} and the above (ii) we conclude that
 for $b\neq b'$ in $[\mc{D}]$ we have that
 $\|z_b-z_{b'}\|\ge \frac{1}{1728}$ which yields that $Z$ is
 nonseparable.
 \end{proof}

 The following summarizes some of the properties
 of the space $\eqs_{\mc{F}_2'}$.

 \begin{corollary} There exists a separable Banach space
  $\eqs_{\mc{F}_2'}$ such that
 \begin{enumerate}
 \item[(i)] The space $\eqs_{\mc{F}_2'}$ is HI and reflexively
 saturated.
 \item[(ii)] Every quotient of $\eqs_{\mc{F}_2'}$ has a further quotient
 isomorphic to $\ell_2.$
 \item[(iii)] Every quotient of $\eqs_{\mc{F}_2'}$ either has nonseparable dual
 or it is isomorphic to $\ell_2.$
 \item[(iv)] There exists a quotient of $\eqs_{\mc{F}_2'}$ not containing
 reflexive subspaces.
 \end{enumerate}
 \end{corollary}

 Before presenting the definition of the norming set
 $D_{\mc{F}_2'}$ let's explain our motivation. First we observe
 that Proposition \ref{aa30} yields that for a sequence
 $(b_n^*)_n$ satisfying the assumptions, the sequence
 $([b_n^*])_n$ in the quotient space $W=\mathfrak X_{\mc{F}_2}^*/
 (\mathfrak X_{\mc{F}_2})_*$ is equivalent to the $\ell_2$ basis.
This in particular yields that $W$ contains copies of
$\ell_2(\Gamma)$ with $\# \Gamma$ equal to the continuum. Our
intention is to define $\mc{F}_2' \subset \mc{F}_2$ and
$D_{\mc{F}_2'} \subset D_{\mc{F}_2}$ such that every infinite
$\sigma_{\mc{F}}$-special sequence $b=(f_1, f_2, \dots, f_n,
\dots)$ satisfies the requirements of Proposition \ref{aa30} with
respect to the norm induced by the set $D_{\mc{F}_2'}.$ Clearly if
this is accomplished, then granting Proposition \ref{aa30}, the
quotient $\mathfrak X_{\mc{F}_2'}^*/(\mathfrak X_{\mc{F}_2'})_*$
will be equivalent to $\ell_2(\Gamma).$

 The norm in the space $\mathfrak X_{\mc{F}_2'}$ is induced by a set $D_{\mc{F}_2'}$
 which in turn, is recursively defined as $\cup_{n=0}^\infty D_n.$
 The key ingredient is that the ground set $\mc{F}_2',$ which is a subset
 of $\mc{F}_2,$ is also
 defined inductively following the definition of $D_n.$ Thus in
 each step we define the set $S_n$ of the $\sigma_{\mc{F}}$-special sequences
 related to $\mc{F}_2$
 and from this set,
 the set $\mc{F}_2^n.$

 For $n=0,$ we set $S_0=\emptyset,$ $D_0=\{\pm e_n^*:\;n\in\N\}.$

 For $n=1$ we set $S_1= \cup_{j=1}^\infty F_j,$ $\mc{F}_2^1$ is defined from
 $S_1$ and $D_1$ results from $D_0\cup \mc{F}_2^1$
  after applying the operations of Definition
 \ref{def13} and taking rational convex combinations.

 Assume that $S_n,$ $\mc{F}_2^n,$ $D_n$ have been defined such that
 every $\sigma_{\mc{F}}$ special sequence $(f_1,  \dots, f_d)$ in $S_n$
 satisfies $d\le n.$ The $\sigma_{\mc{F}}$ special
 sequence $f_1, \dots, f_d$ in $S_n$ is called
 {\bf $n+1$-extendable} if for each $1\le i\le d$ there exists a
 $(18, 4j_i-3, 1, D_n, \mc{F}_2^n)$
 attracting sequence  $\chi_i=(x_k, x_k^*)_{k=1}^{4n_{j_i}-3}$,
 with $f_i=g_{\chi_i}$ (Definition
 \ref{def12}).
 Here a $(18, 4j_i-3, 1, D_n, \mc{F}_2^n)$ attracting sequence is defined as in
Definition \ref{depseq} where the norm of the underlying space is
induced by the set $D_n$ and moreover $\|x_{2k-1}\|_{D_n}\le 18$
and $\|x_{2k-1}\|_{\mc{F}_2^n}\le \frac{1}{n_{4j_i-3}^2}.$

Then we set $S_{n+1} =S_n\cup \{ (f_1, \dots, f_d):\; (f_1,
\ldots, f_{d-1})$ is a  $n+1\text{-extendable}$ $\sigma_{\mc{F}}$
special sequence$\}$.

Next we define $\mc{F}_2^{n+1}$ from $S_{n+1}$ in the usual manner
and then $D_{n+1}$ from $D_n \cup \mc{F}_2^{n+1}$ as before.

This completes the inductive definition. We set
$\mc{F}_2'=\cup_n\mc{F}_2^n$ and $D_{\mc{F}_2'} = \cup_n D_n.$

It is easy to see that for every $b=(f_n)_n$ such that $b^*\in
\overline{\mc{F}_2'}^{w^*}$ the sequence $(f_n)_n$ satisfies the
properties of Proposition \ref{aa30} and this yields that indeed
$\mathfrak X_{\mc{F}_2'}^*/ (\mathfrak X_{\mc{F}_2'})_*$ is
isomorphic to $\ell^2(\Gamma).$

\section{A nonseparable HI space with no reflexive subspace}
\label{sec4}

In this section we proceed to the construction of a nonseparable
HI space containing no reflexive subspace. The general scheme we
shall follow is similar to the one used for the definition of
$\eqs_{\mc{F}_2}$. However there are two major differences. The
first concerns  saturation methods. In the present construction we
shall use the operations $(\mc{S}_{n_j},\frac{1}{m_j})_j$ for
appropriate sequences $(m_j)_j$, $(n_j)_j$. The James Tree space
which will play the role of  the ground space is also different
from $JT_{\mc{F}_2}$. Indeed the ground set $\mc{F}_s'$ is built
on a family $(F_j)_j$ which is related to the Schreier families
$(\mc{S}_{n_{4j-3}})_j$. Furthermore in $\mc{F}_s'$ we connect the
$\sigma_{\mc{F}}$ special functionals with the use of the Schreier
operation instead of taking $\ell_2$ sums as in $\mc{F}_2$.
Finally, $\mc{F}_s'$ is defined recursively as we did in the
previous variant $\eqs_{\mc{F}_2'}$ of $\eqs_{\mc{F}_2}$. The
spaces $(\eqs_{\mc{F}_s'})_*$, $\eqs_{\mc{F}_s'}$ share the same
properties with  $(\eqs_{\mc{F}_2'})_*$, $\eqs_{\mc{F}_2'}$. The
difference occurs between $\eqs_{\mc{F}_2'}^*$ and
$(\eqs_{\mc{F}_s'})^*$. Indeed, as we have seen
$(\eqs_{\mc{F}_2'})^*/(\eqs_{\mc{F}_2'})_*$ is isomorphic to
$\ell_2(\Gamma)$, while as it will be shown
$(\eqs_{\mc{F}_s'})^*/(\eqs_{\mc{F}_s'})_*$ is isomorphic to
$c_0(\Gamma)$ with $\#\Gamma$ equal to the continuum. The later
actually yields all the desired properties for
$(\eqs_{\mc{F}_s'})^*$. Namely it is HI and it does not contain
any reflexive subspace.

We recall the definition of $(\mc{S}_n)_n$, the first infinite sequence
of the Schreier families. The first
Schreier family $\mc{S}_{1}$ is the following
\[
\mc{S}_{1}=\{F\subset \mathbb{N}:\; \#F\le \min F\}\cup\{\emptyset\}.
\]
For $n\geq 1$ the   definition goes as follows
\[
\mc{S}_{n+1}=\left\{F=\bigcup\limits_{i=1}^{d} F_{i}:\;
F_{i}\in\mc{S}_{n}\,\,F_{i}<F_{i+1},\,\,\textrm{for all $i$ and }d
\le \min F_1 \right\}.
\]
Each $\mc{S}_{n}$ is, as  can be easily verified by induction, compact, hereditary and spreading.

A finite sequence $(E_1,E_2,\ldots ,E_k)$ of successive subsets of
$\mathbb{N}$ is said to be $\mc{S}_{n}$ admissible, $n\in\N$, if
$\{\min E_{i}:i=1,\ldots, k\}\in\mc{S}_{n}$. A finite sequence
$(f_1,f_2,\ldots,f_k)$ of vectors in $c_{00}$ is said to be
$\mc{S}_{n}$ admissible if the sequence $(\supp f_1,\supp
f_2,\ldots,\supp f_k)$ is $\mc{S}_{n}$ admissible.

 We fix two sequences of integers \seq{m}{j} and \seq{n}{j}
 defined as follows:
 \begin{itemize}
 \item  $m_1=2$  and $m_{j+1}=m_j^{m_j}$.
 \item $n_1=1$,  and $n_{j+1}=2^{2m_{j+1}}n_j$.
 \end{itemize}

 \begin{definition}\label{bscc}
 {\bf (basic special convex combinations)}
 Let $\varepsilon >0$ and $j\in\mathbb{N}$, $j>1$.
 A convex combination $\sum\limits_{k\in F}a_ke_k$ of the basis
 \seq{e}{k} is said to be an $(\varepsilon,j)$
 basic special convex combination
 ($(\varepsilon,j)$ B.S.C.C.) if
 \begin{enumerate}
  \item $F\in \mc{S}_{n_j}$
  \item For every $P\in \mc{S}_{2\log_2(m_j)(n_{j-1}+1)}$ we have
  that   $\sum\limits_{k\in P}a_k<\varepsilon$.
  \item The sequence  $(a_k)_{k\in F}$ is a non increasing
  sequence of positive reals.
 \end{enumerate}
 \end{definition}

 \begin{remark}\label{rem8} The basic special convex combinations
 have been used implicitly in \cite{AD}, their  exact definition
 was given in \cite{AMT} while they have systematically studied in
 \cite{AT1}.
 \end{remark}

 \begin{definition}\label{scc}
 {\bf (special convex combinations)}
 Let $\varepsilon >0$, $j\in\mathbb{N}$ with $j>1$ and let
 \seq{x}{k} be a block sequence of the standard basis.
 A convex combination $\sum\limits_{k\in F}a_kx_k$ of the sequence
 \seq{x}{k} is said to be an $(\varepsilon,j)$
 special convex combination
 ($(\varepsilon,j)$ S.C.C.) of \seq{x}{k} if $\sum\limits_{k\in F}a_ke_{t_k}$
 (where $t_k=\min\supp x_k$ for each $k$) is an
 $(\varepsilon,j)$ basic special convex combination.

 Moreover, if   $\sum\limits_{k\in F}a_kx_k$ is a S.C.C. in a Banach space
 $(X,\|\;\|)$ such that
 $\Vert x_k\Vert \le 1$ for all $k$
 and $\Vert \sum\limits_{k\in F}a_kx_k\Vert \ge\df{1}{2}$ we say that
 $\sum\limits_{k\in F}a_kx_k$  is a seminormalized
 $(\varepsilon,j)$ special convex combination of
 $(x_k)_{k\in\mathbb{N}}$.
 \end{definition}

 \begin{definition}\label{def18}
 We set $F_{0}=\{\pm e_n^*:\;n\in\N\}$ while for $j=1,2,\ldots$ we
 set $F_j=\{\frac{1}{m_{4j-3}^2}\sum\limits_{i\in I}\pm e_i^*:\;
 I\in \mc{S}_{n_{4j-3}}\}\cup\{0\}$.
 We also set $F=\bigcup\limits_{j=0}^{\infty}F_j$.
 \end{definition}

 Let's observe that the sequence $\mc{F}=(F_j)_{j=0}^{\infty}$
 is a JTG family. The $\sigma_{\mc{F}}$ special sequences corresponding to
 this family are defined exactly as in Definition \ref{def4}.

  \begin{definition}\label{def20}
 {\bf ($\sigma$ coding,  special sequences and attractor sequences)}
  Let $\Q_s$ denote the set of all finite
 sequences $(\phi_1,\phi_2,\ldots,\phi_d)$ such that
 $\phi_i\in \co$, $\phi_i\neq 0$ with $\phi_i(n)\in \Q$ for all $i,n$ and
 $\phi_1<\phi_2<\cdots<\phi_d$.
 We fix a pair $\Omega_1,\Omega_2$ of disjoint infinite subsets of
 $\N$.
 From the fact that $\Q_s$ is
 countable we are able to define a Gowers-Maurey type injective
 coding function
 $\sigma:\Q_s\to \{2j:\;j\in\Omega_2\}$ such that
 $m_{\sigma(\phi_1,\phi_2,\ldots,\phi_d)}>\max\{\frac{1}{|\phi_i(e_l)|}:\;l\in\supp
 \phi_i,\;i=1,\ldots,d\}\cdot\max\supp \phi_d$.
 Also, let $(\Lambda_i)_{i\in\N}$ be a sequence of pairwise
 disjoint infinite subsets of $\N$ with $\min \Lambda_i>m_i$.
 \begin{enumerate}
 \item[(A)]
 A finite  sequence $(f_i)_{i=1}^{d}$ is said to be a
 {\bf $\mc{S}_{n_{4j-1}}$ special sequence}  provided that
 \begin{enumerate}
 \item[(i)]  $(f_1,f_2,\ldots,f_{d})\in\Q_s$ and
 $(f_1,f_2,\ldots,f_{d})$ is a $\mc{S}_{n_{4j-1}}$ admissible
 sequence,
 $f_i\in D_G$
 for $i=1,2,\ldots,n_{4j-1}$.
 \item[(ii)] $w(f_1)=m_{2k}$ with $k\in\Omega_1$, $m_{2k}^{1/2}>n_{4j-1}$
  and   for each $1\le i<d$, $w(f_{i+1})=m_{\sigma(f_1,\ldots,f_{i})}$.
 \end{enumerate}
 \item[(B)]
 A finite  sequence $(f_i)_{i=1}^{d}$ is said to be a
 {\bf $\mc{S}_{n_{4j-3}}$ attractor sequence}  provided that
 \begin{enumerate}
 \item[(i)]  $(f_1,f_2,\ldots,f_{d})\in\Q_s$ and
 $(f_1,f_2,\ldots,f_{d})$ is a $\mc{S}_{n_{4j-3}}$ admissible sequence.
 \item[(ii)] $w(f_1)=m_{2k}$ with $k\in\Omega_1$, $m_{2k}^{1/2}>n_{4j-3}$
  and $w(f_{2i+1})=m_{\sigma(f_1,\ldots,f_{2i})}$
 for each $1\le i<\frac{d}{2}$.
 \item[(iii)]  $f_{2i}=e^*_{l_{2i}}$ for some
  $l_{2i}\in\Lambda_{\sigma(f_1,\ldots,f_{2i-1})}$, for
  $i=1,\ldots,\frac{d}{2}$.
 \end{enumerate}
 \end{enumerate}
  \end{definition}

 \begin{definition}\label{def14}
  {\bf (The space $\eqs_{\mc{F}_s'}$)}
 In order to define the norming set $D$ of the space
 $\eqs_{\mc{F}_s'}$ we shall
 inductively define four sequences of subsets of $c_{00}(\N)$, denoted as
 \seq{K}{n}, \seq{\tau}{n}, \seq{G}{n}, \seq{D}{n}.
 \end{definition}

 We set $K_0=F$ ($K_0^0=F$, $K_0^j=\emptyset$, $j=1,2,\ldots$),
  $G_0=F$, $\tau_0=\emptyset$ and $D_0=\convq(F)$.
 Suppose that $K_{n-1}$, $\tau_{n-1}$, $G_{n-1}$ and $D_{n-1}$
 have been defined. The inductive properties of
  \seq{K}{n}, \seq{\tau}{n}, \seq{G}{n}, \seq{D}{n} are included
  in the inductive definition.
  We set
  \[K_n^{2j}=K_{n-1}^{2j}\cup\{\frac{1}{m_{2j}}\sum\limits_{i=1}^df_i:\;f_1<\cdots<f_d\mbox{
  is }\mc{S}_{n_{2j}}\mbox{ admissible},\;f_i\in D_{n-1}\}\]
  \begin{eqnarray*}
  K_n^{4j-3}=K_{n-1}^{4j-3}&\cup&\{\pm
   E(\frac{1}{m_{4j-3}}\sum\limits_{i=1}^df_i):\;(f_1,\ldots,f_d)\mbox{
  is a }\mc{S}_{n_{4j-3}} \mbox{ attractor}\\
  & &\mbox{sequence, } f_i\in K_{n-1}\mbox{ and }E \mbox{ is an interval of $\N$}\}
  \end{eqnarray*}
  \begin{eqnarray*}
  K_n^{4j-1}=K_{n-1}^{4j-1}&\cup&\{\pm E(\frac{1}{m_{4j-1}}\sum\limits_{i=1}^df_i):\;(f_1,\ldots,f_d)\mbox{
  is a }\mc{S}_{n_{4j-1}}\mbox{ special}\\
  & & \mbox{sequence, }f_i\in K_{n-1}\mbox{ and } E \mbox{ is an interval of $\N$}\}
  \end{eqnarray*}
  \[K_n^0=F\].

  We set $K_n=\bigcup\limits_{j=0}^{\infty}K_n^j$.

  In order to define $\tau_n$ we need the following definition.
 \begin{definition}
 {\bf ($(D_{n-1},j)$ exact functionals)}
 A functional $f\in F$ is said to be $(D_{n-1},j)$ exact if $f\in F_j$
 and there exists $x\in \co$ with $\|x\|_{D_{n-1}}\le 1000$,
 $\ran(x)\subset\ran(f)$, $f(x)=1$  such that for every
 $i\neq j$, we have that $\|x\|_{F_i}\le \frac{1000}{m_{4i-3}^2}$  if
 $i<j$ while $\|x\|_{F_i}\le 1000 \frac{m_{4j-3}^2}{m_{4i-3}^2}$ if
 $i>j$.
 \end{definition}
 We set
 \begin{eqnarray*}
 \tau_n &=&\{\pm E(\sum\limits_{i=1}^d\phi_i):\; d\le n, E\mbox{ is an
 interval, }(\phi_i)_{i=1}^d \mbox{ is }\sigma_{\mc{F}}\mbox{ special}\\
 & & \mbox{and each }\phi_i \mbox{ is }(D_{n-1},\ind(\phi_i))
 \mbox{ exact}\}.
 \end{eqnarray*}
 We recall that for $\Phi=\pm E(\sum\limits_{i=1}^d\phi_i)\in
 \tau_n$, $\ind(\Phi)=\{\ind(\phi_i):\; E\cap\ran \phi_i\neq
 \emptyset\}$. We set
 \begin{eqnarray*}G_n&=&\{\sum\limits_{i=1}^d\e_i\Phi_i:\;\Phi_i\in\tau_n,\;\e_i\in\{-1,1\}, \;\min\supp
 \Phi_i\ge d,\\
 & & (\ind(\Phi_i))_{i=1}^d\mbox{ are pairwise
 disjoint}\}
 \end{eqnarray*}
 We set $D_n=\convq(K_n\cup G_n\cup D_{n-1})$.

 We finally set $D=\bigcup\limits_{n=0}^\infty D_n$. We also set
 $\tau=\bigcup\limits_{n=0}^\infty \tau_n$,
 $\mc{F}_s'=\bigcup\limits_{n=0}^\infty G_n$,
 $K=\bigcup\limits_{n=0}^\infty K_n$. We set
 $K^j=\bigcup\limits_{n=1}^{\infty}K_n^j$ for $j=1,2,\ldots$. For
 $f\in K^j$ we write $w(f)=m_j$. We notice that $w(f)$ is not
 necessarily uniquely determined.

 We also need the following definition.
 \begin{definition}\label{def16}
  {\bf ($(D,j)$ exact functionals)}
 A functional $f\in F$ is said to be $(D,j)$ exact if $f\in F_j$
 and there exists $x\in \co$ with $(\|x\|_{D}=)\|x\|\le 1000$,
 $\ran(x)\subset\ran(f)$, $f(x)=1$  such that for every
 $i\neq j$, we have that $\|x\|_{F_i}\le \frac{1000}{m_{4i-3}^2}$  if
 $i<j$ while $\|x\|_{F_i}\le 1000 \frac{m_{4j-3}^2}{m_{4i-3}^2}$ if
 $i>j$.
 \end{definition}

 \begin{remarks}\label{rem10}
 \begin{enumerate}
 \item[(i)] If the functional $\phi$ is $(D_n,j)$ exact then it is
 also $(D_k,j)$ exact for all $k\le n$.
 \item[(ii)] Let \seq{\phi}{i} be a $\sigma_{\mc{F}}$ special sequence
 such that each $\phi_i$ is $(D,\ind(\phi_i))$  exact.
 Then each $\phi_i$ is $(D_n,\ind(\phi_i))$  exact for all $n$ and
 $\sum\limits_{i=1}^d\phi_i\in \tau_n\subset \tau$ for all $n\ge
 d$. It follows that
 $\sum\limits_{i=1}^{\infty}\phi_i\in \overline{\tau}^{w^*}\subset
 \overline{\mc{F}_s'}^{w^*}\subset\overline{D}^{w^*}=B_{\eqs_{\mc{F}_s'}^*}$.
 \item[(iii)] Let \seq{\phi}{i} be a $\sigma_{\mc{F}}$ special sequence
 such that $\sum\limits_{i=1}^d\phi_i\in\tau$ for all $d$. In this
 case we call the  $\sigma_{\mc{F}}$ special sequence \seq{\phi}{i}
 survivor and the functional $\Phi=\sum\limits_{i=1}^{\infty}\phi_i$ a survivor
 $\sigma_{\mc{F}}$ special functional.
 Then each $\phi_i$ is $(D_n,j_i)$  exact (where $j_i=\ind(\phi_i)$) for all $n$,
 thus for each $n$ there exists $x_{i,n}$ with
 $\|x_{i,n}\|_{D_{n}}\le 1000$,
 $\ran(x_{i,n})\subset\ran(\phi_i)$, $\phi_i(x_{i,n})=1$
 and  such that for every
 $k\neq j_i$, we have that $\|x_{i,n}\|_{F_k}\le
 \frac{1000}{m_{4k-3}^2}$  if
 $k<j$ while $\|x_{i,n}\|_{F_k}\le 1000
 \frac{m_{4j-3}^2}{m_{4k-3}^2}$ if
 $k>j$. Taking a subsequence of $(x_{i,n})_{n\in\N}$ norm
 converging to some $x_i$ it is easily checked that
 $\|x_i\|\le 1000$, $\phi_i(x_i)=1$ while
 $\|x_i\|_{F_k}\le
 \frac{1000}{m_{4k-3}^2}$  for
 $k<j$ and $\|x_i\|_{F_k}\le 1000
 \frac{m_{4j-3}^2}{m_{4k-3}^2}$ for
 $k>j$.

 A  sequence \seq{x}{i} satisfying the above property is called a
 sequence witnessing that the $\sigma_{\mc{F}}$ special sequence
 \seq{\phi}{i} (or the special functional
 $\Phi=\sum\limits_{i=1}^{\infty}\phi_i$)
  is survivor.
 \end{enumerate}
 \end{remarks}

 \begin{lemma}\label{lem24}
 The norming set $D$ of the space $\eqs_{\mc{F}_s'}$ is the minimal
 subset of \co satisfying the following conditions:
 \begin{enumerate}
 \item[(i)] $\mc{F}_s'\subset D$.
 \item[(ii)] $D$ is closed in the
 $(\mc{S}_{n_{2j}},\frac{1}{m_{2j}})$ operations.
 \item[(iii)] $D$ is closed in the
 $(\mc{S}_{n_{4j-1}},\frac{1}{m_{4j-1}})$ operations on
 $\mc{S}_{n_{4j-1}}$ special sequences.
 \item[(iv)] $D$ is closed in the
 $(\mc{S}_{n_{4j-3}},\frac{1}{m_{4j-3}})$ operations on
 $\mc{S}_{n_{4j-3}}$ special sequences.
 \item[(v)] $D$ is symmetric, closed in the restrictions of its
 elements on intervals of $\N$ and rationally convex.
 \end{enumerate}
 \end{lemma}

It is easily proved that the Schauder basis \seq{e}{n} of the
space $\eqs_{\mc{F}_s'}$ is boundedly complete and that
$\eqs_{\mc{F}_s'}$ is an asymptotic $\ell_1$ space. Since the
space $JT_{\mc{F}_s'}$ is $c_0$ saturated
 (see Remark \ref{brem7} where we use the notation $JT_{\mc{F}_{\tau,s}}$ for such a space)
  we get the following.

 \begin{proposition}\label{prop23} The identity operator $I:\eqs_{\mc{F}_s'}\to
 JT_{\mc{F}_{s}}$ is strictly singular.
 \end{proposition}

 \begin{remark}\label{rem9}
 Applying the methods of \cite{AT1} and taking into account that
 the identity operator $I:\eqs_{\mc{F}_s'}\to
 JT_{\mc{F}_{s}}$ is strictly singular we may prove the
 following. For every $\e>0$ and $j>1$ every block subspace of
 $\eqs_{\mc{F}_s'}$ contains a vector $x$ which is a seminormalized
 $(\e,j)$ S.C.C. with $\|x\|_{\mc{F}_s'}<\e$.
 \end{remark}

 \begin{definition}\label{defexact}
  {\bf (exact pairs in $\eqs_{\mc{F}_s'}$)}
A pair $(x,f)$ with $x\in c_{00}$ and $f\in K$  is said to be a
$(12,j,\theta)$ exact pair, where $j\in\N$, if the following
conditions are satisfied:
\begin{enumerate}
\item[(i)]\,\, $1\leq \|x\|\leq 12$, $f(x)=\theta$ and
$\ran(f)=\ran(x)$. \item[(ii)]\,\, For every $g\in K$ with
$w(g)=m_{i}$ and $i< j$, we have that $|g(x)|\leq 24/m_{i}$.
\item[(iii)] For every sequence  $(\phi_i)_i$ in $K$ with
$m_j<w(\phi_1)<w(\phi_2)<\cdots$ we have that
$\sum\limits_i|\phi_i(x)|\le 12/m_{j}$.
\end{enumerate}
\end{definition}
%Let us observe that for every $j>1$, if $x=\sum_{k\in
%F}a_{k}e_{k}$ is an $(1/m_{2j}^{2},2j)$ basic s.c.c. and
%$f=\frac{1}{m_{2j}}\sum_{k\in F}e_{k}^{*}$ then from Proposition
%??? it follows easily that the pair $(m_{2j}x,f)$ is a $(12,2j)$
%exact pair.

\begin{proposition}\label{prop31}
For every $j\in\N$, $\e>0$ and every block subspace $Z$ of
$\eqs_{\mc{F}_s'}$, there exists a $(12,2j,1)$ exact pair $(z,f)$
with $z\in Z$ and $\|z\|_{\mc{F}_s'}<\e$.
\end{proposition}
\begin{proof}[\bf Proof.]
Since the identity operator $I:\eqs_{\mc{F}_s'}\to
 JT_{\mc{F}_{s}}$ is strictly singular we may assume, passing to a
 block subspace of $Z$, that $\|z\|_{\mc{F}_s'}\le
 \frac{\e}{12}\|z\|$ for every $z\in Z$.

 Let \seq{x}{k} be a block sequence in
$Z$ such that \seq{x}{k} is a $(2,\frac{1}{m_{2j}})$ R.I.S.
and each $x_k$ is a seminormalized
 $(\frac{1}{m_{j_k}},j_k)$ S.C.C.
 Passing to a subsequence we may assume that $(b^*(x_k))_{k\in\N}$ converges for
every $\sigma_{\mc{F}}$
branch $b$. We set $z_k=x_{2k-1}-x_{2k}$. Then \seq{z}{k} is a $(4,\frac{1}{m_{2j}})$
R.I.S. such that $b^*(z_k)\to 0$ for every branch $b$.

We recall that each $g\in \mc{F}_s'$ has the form
$g=\sum\limits_{i=1}^d\e_i\Phi_i^*$  with $\e_i\in\{-1,1\}$,
$(\Phi_i^*)_{i=1}^d\in\tau$ with $\min\supp x_i^*\ge d$ and
$(\ind(x_i^*))_{i=1}^d$ pairwise disjoint. We may assume,
replacing \seq{z}{k} by an appropriate subsequence,  that for
every $g\in G$ we have that the set $\{\min\supp
z_k:\;|g(z_k)|>\frac{1}{m_{2j}}\}$ belongs to $\mc{S}_2$, the
second Schreier family.

It follows now from Proposition 6.2 of \cite{AT1}
that if $z=\sum_{k\in F}a_{k}z_{k}$ is a $(1/m_{2j}^{2},2j)$
special convex combination of $\seq{z}{k}$
 and $f$ is of the form $f=1/m_{2j}\sum_{k\in F}f_{k}$ where
 $f_k\in K$ with
 $f_{k}(z_{k})=1$ and $\ran(f_k)=\ran(z_k)$ then
 $(z,f)$ is the desired $(12,2j,1)$ exact pair.
\end{proof}

 \begin{definition}\label{depseqfs}
 {\bf (dependent sequences and attracting sequences in $\eqs_{\mc{F}_s'}$)}
 \begin{enumerate}
  \item[(A)]
 A double sequence $(x_k,x_k^*)_{k=1}^{d}$ is said to be a
 $(C,4j-1,\theta)$ dependent sequence  (for $C>1$,
 $j\in\mathbb{N}$,     and
 $0\le\theta\le 1$) if there exists a
 sequence $(2j_k)_{k=1}^{d}$ of even integers such that the
 following conditions are fulfilled:
 \begin{enumerate}
 \item[(i)] $(x^*_k)_{k=1}^{d}$ is a $\mc{S}_{n_{4j-1}}$ special
 sequence with
 $w(x^*_{k})=m_{2j_{k}}$ for each $k$.
 \item[(ii)] Each $(x_{k},x_{k}^*)$ is a
 $(C,2j_{k},\theta)$ exact pair.
 \item[(iii)]  Setting $t_k=\min\supp x_k$, we have that
  $t_1>m_{2j}$ and $\{t_1,\ldots,t_d\}$ is a maximal element of
  $\mc{S}_{n_{4j-1}}$. (Observe, for later use, that
 Remark 3.18  of \cite{AT1} yields that there exist
 $(a_k)_{k=1}^d$ such that $\sum\limits_{k=1}^da_ke_{t_k}$
 is a $(\frac{1}{m_{4j-1}^2},4j-1)$ basic special convex combination).
 \end{enumerate}
 \item[(B)]
 A double sequence $(x_k,x_k^*)_{k=1}^{d}$ is said to be a
 $(C,4j-3,\theta)$ attracting sequence (for $C>1$,
 $j\in\mathbb{N}$,     and
 $0\le\theta\le 1$) if there exists a
 sequence $(2j_k)_{k=1}^{d}$ of even integers such that the
 following conditions are fulfilled:
 \begin{enumerate}
 \item[(i)] $(x^*_k)_{k=1}^{d}$ is a $\mc{S}_{n_{4j-3}}$ attractor sequence
 with  $w(x^*_{2k-1})=m_{2j_{2k-1}}$ and
 $x^*_{2k}=e^*_{l_{2k}}$ where $l_{2k}\in \Lambda_{2j_{2k}}$
  for all $k\le d/{2}$.
 \item[(ii)] $x_{2k}=e_{l_{2k}}$.
 \item[(iii)]  Setting $t_k=\min\supp x_k$, we have that
  $t_1>m_{2j}$ and $\{t_1,\ldots,t_d\}$ is a maximal element of
  $\mc{S}_{n_{4j-3}}$. (Observe that
 Remark 3.18  of \cite{AT1} yields that there exist
 $(a_k)_{k=1}^d$ such that $\sum\limits_{k=1}^da_ke_{t_k}$
 is a $(\frac{1}{m_{4j-3}^2},4j-3)$ basic special convex combination).
 \item[(iv)] Each $(x_{2k-1},x_{2k-1}^*)$ is a
 $(C,2j_{2k-1},\theta)$ exact pair.
 \end{enumerate}
 \end{enumerate}
 \end{definition}

  \begin{proposition}\label{prop24}
 The space $\eqs_{\mc{F}_s'}$ is reflexively saturated and
 Hereditarily Indecomposable.
 \end{proposition}
 \begin{proof}[\bf Proof.]
The proof that $\eqs_{\mc{F}_s'}$ is reflexively saturated is
 a consequence of the fact that the identity operator $I:\eqs_{\mc{F}_s'}\to
 JT_{\mc{F}_{s}}$ is strictly singular and its proof is identical
 to that of Proposition \ref{prop15}.

In order to show that the space $\eqs_{\mc{F}_s'}$ is Hereditarily
Indecomposable we consider  a pair of block subspaces $Y$ and $Z$
and $\delta>0$. We choose $j$ such that
$m_{4j-1}>\frac{192}{\delta}$.

Using Proposition \ref{prop31} we may choose a $(12,4j-1,1)$
dependent sequence $(x_k,x_k^*)_{k=1}^d$ such that
$\|x_{2k-1}\|_{\mc{F}_{s}}<\frac{2}{m_{4j-1}^2}$ for all $k$ while
$x_{2k-1}\in Y$ if $k$ is odd and $x_{2k-1}\in Z$ if $k$ is even.
>From the observation in Definition \ref{depseqfs}(A)(iii) there
exist $(a_k)_{k=1}^d$ such that $\sum\limits_{k=1}^da_ke_{t_k}$ is
a $(\frac{1}{m_{4j-1}^2},4j-1)$ basic special convex combination
(where $t_k=\min\supp x_k$). A variant of Proposition \ref{prop11}
(i) in terms of the space $\eqs_{\mc{F}_s'}$ (using Proposition
6.2 of \cite{AT1}) yields that
$\|\sum\limits_{k=1}^d(-1)^{k+1}a_{k}x_{k}\|\le
\frac{96}{m_{4j-1}^2}$. On the other hand the functional
$f=\frac{1}{m_{4j-1}}\sum\limits_{k=1}^dx_k^*$ belongs to the
norming set $D$ of the space $\eqs_{\mc{F}_s'}$ and estimating
$f(\sum\limits_{k}a_{k}x_{k})$ we get that
$\|\sum\limits_{k=1}^da_{k}x_{k}\|\ge \frac{1}{m_{2j-1}}$.

   Setting $y=\sum\limits_{k\;\text{odd}}a_{k}x_{k}$ and
   $z=\sum\limits_{k\;\text{even}}a_{k}x_{k}$ we have that $y\in Y$ and $z\in
   Z$ while from the above inequalities we get that $\|y-z\|\le
   \delta \|y+z\|$. Therefore
 $\eqs_{\mc{F}_s'}$ is a Hereditarily Indecomposable space.
 \end{proof}

 \begin{proposition}\label{prop25}
 The dual space $\eqs_{\mc{F}_s'}^*$  is the norm closed
 linear span of the $w^*$ closure of $\mc{F}_s'$ i.e.
 \[   \eqs_{\mc{F}_s'}^*=\overline{\spann}(\overline{\mc{F}_s'}^{w^*}).\]
 \end{proposition}
 \begin{proof}[\bf Proof.] Assume the contrary. Then using arguments similar
 to those of the proof of Proposition \ref{prop14} we
 may choose a $x^*\in \eqs_{\mc{F}_s'}^*$ with $\|x^*\|=1$,  and a
 block sequence \seq{x}{k} in $\eqs_{\mc{F}_s'}$ with $x^*(x_k)>1$
 and $\|x_k\|\le 2$  such that
 $x_k \stackrel{w}{\longrightarrow}0$ in $JT_{\mc{F}_{s}}$.
 Observe that the action of $x^*$ ensures that every convex
 combination of \seq{x}{k} has norm greater than 1.

 We may choose a convex block sequence \seq{y}{k} of \seq{x}{k}
 with $\|y_k\|_{\mc{F}_{s}}<\frac{\e}{2}$ where
 $\e=\frac{1}{m_4}$.
 We select a block sequence \seq{z}{k} of \seq{y}{k} such that
 each $z_k$ is a convex combination  of \seq{y}{k}
 and such that \seq{z}{k} is
 $(4,\e)$ R.I.S.
 This is possible if we consider each $z_k$ to be a
 $(\frac{1}{m_{i_k}},i_k)$ S.C.C. and $m_{i_{k+1}}\e>\max\supp
 z_k$ for an appropriate increasing sequence of integers
 \seq{i}{k}.
 We then consider $x=\sum a_k z_k$, an $(\e,4)$ S.C.C.
 of \seq{z}{k}.  A variant of Proposition 6.2(1a) of \cite{AT1}
 yields that $\|x\|\le \frac{20}{m_4}<1$. On the other hand, since
 $x$ is a convex combination of \seq{x}{k} we get that $\|x\|>1$, a
 contradiction.
 \end{proof}

  \begin{definition} \label{sepfs} Let $(x_n)_{n\in\N}$ be a bounded block
sequence in $JT_{\mc{F}_s}$ and $\e >0.$ We say that
$(x_n)_{n\in\N}$ is {\em $\e$-separated} if for every $\phi\in
\cup_{j\in \mathbb{N}}F_j$
\[ \# \{n: |\phi(x_n)|\ge\e \}\le 1. \]

In addition, we say that $(x_n)_{n\in\N}$ is {\em separated} if
for every $L\in [\N]$ and $\e >0$ there exists an $M\in [L]$ such
that $(x_n)_{n\in M}$ is $\e$-separated.
\end{definition}

 \begin{lemma} \label{sepfslem} Let $(x_n)_{n\in\N}$ be a weakly null
separated sequence in $JT_{\mc{F}_s}.$ Then for every $\e
>0,$ there exists an $L\in [\N]$ such that for all $y^* \in
\overline{\tau}^{w^*},$
\[\# \{ n\in L: |y^*(x_n)| \ge \e \}\le 2. \]
\end{lemma}
 The proof of the above lemma is similar to that of Lemma
 \ref{aa2}.

 \begin{lemma}\label{lem34}
 Let \seq{z}{k} be a block sequence in $\eqs_{\mc{F}_s'}$ such that each
 $z_k$ is a $(\frac{1}{m_{2j_k}},2j_k)$  special
 convex combination of a normalized block sequence,
 where \seq{j}{k} is strictly increasing. Then
 the sequence \seq{z}{k} is separated.
 \end{lemma}
 \begin{proof}[\bf Proof.]
 Given $\e>0$ and $L\in [\N]$ we have to find an $M\in[L]$ such
 that for every $\phi\in\bigcup\limits_{j\in N}F_j$ we have that
 $|\phi(z_k)|>\e$ for at most one $k\in M$. For simplicity in our
 notation we may assume, passing to a subsequence, that
 $\frac{1}{m_{2j_1}}<\e$ and $\max\supp z_{k-1}<\e m_{2j_k}$ for
 each $k$.

 Now let $\phi\in\bigcup\limits_{j\in N}F_j$. Then  $\phi$ takes the form
 $\phi=\frac{1}{m_{4j-3}^2}\sum\limits_{i\in I}\pm e_i^*$ with $I\in
 \mc{S}_{n_{4j-3}}$ for some $j$. Let $k_0$  such that
 $2j_{k_0}<4j-3<2j_{k_0+1}$.

  We have that
 $|\phi(z_k)|\le \frac{1}{m_{4j-3}^2}\sum\limits_{i\in
 I}|e_i^*(z_k)|<
 \frac{1}{m_{4j-3}}\frac{1}{m_{2j_{k_0}}}\#\supp(z_k)<\e$ for every
 $k<k_0$.
 Also for $k>k_0$ we get that
 $|\phi(z_k)|\le \frac{2}{m_{2j_k}}<\e$.
 Thus the subsequence we have selected is $\e$-separated and this
 finishes the proof of the lemma.
 \end{proof}

 \begin{remark}\label{remark}
  Let's  observe, for later use, that easy
 modifications of the previous proof yield that for a sequence
 \seq{z}{k} as above the sequence $(z_{2k-1}-z_{2k})_{k\in\N}$ is also
 separated.
 \end{remark}

 \begin{lemma}\label{lem35}
 Let \seq{z}{k} be a weakly null separated sequence in
 $\eqs_{\mc{F}_s'}$. Then for every $\e>0$ there exists $L\in[\N]$
 such that for every $\phi\in \mc{F}_s'$,
 $\{\min\supp z_k:\; k\in L,\;|\phi(z_k)|>\e\}\in\mc{S}_2$.
 \end{lemma}
 \begin{proof}[\bf Proof.]
 The proof is almost identical to the proof of Lemma 10.9 of
 \cite{AT1}. For the sake of completeness we include the proof
 here.

   Using Lemma \ref{sepfslem} we construct a
 sequence \seq{L}{k} of infinite subsets of the natural numbers
 such that
 the following conditions hold
 \begin{enumerate}
 \renewcommand{\labelenumi}{(\roman{enumi})}
 \item  $\min\supp x_{l_1} \ge 3$.
 \item  $l_k=\min L_k\not\in L_{k+1}$  and $L_{k+1}\subset L_k$
 for each $k\in\mathbb{N}$.
 \item For each $k\in\mathbb{N}$ if $p_k=\max\supp x_{l_k}$ then for every
 segment $s\in \overline{\tau}^{w^*}$ we have that
 \[ \#\{n\in L_{k+1}:\; |s^*(x_n)|>\df{\varepsilon}{p_k}\}\le 2.\]
 \end{enumerate}
 We set $L=\{l_1,l_2,l_3,\ldots \}$ and we claim that  the set
 $L$ satisfies the required condition.

 Indeed, let  $\phi=\sum\limits_{i=1}^d\varepsilon_is_i^*\in \mc{F}_{s}$
 where
 segments $s_1,s_2,\ldots,s_d$ are in ${\overline{\tau}}^{w^*}$, have  pairwise
 disjoint sets of indices, $d\le\min s_i$ and $\varepsilon_i\in \{-1,1\}$
  for each
 $i=1,2,\ldots ,d$.
 We set
 \[   l_{i_0}=\min\{n\in L:\;\supp x_n\cap\supp \phi\neq\emptyset\}.\]
 Observe that $d\le p_{i_0}$.
 We set \[  F=\{n\in L:\;|\phi(x_n)|>\varepsilon\}.\]
 We have that $F\subset \{l_{i_0},l_{i_0+1},l_{i_0+2},\ldots \}$
 thus $F\setminus \{l_{i_0}\}\subset L_{i_0+1}$.
 Also (iii) yields that  for each $i=1,2,\ldots,d$ the set
 $F_i=\{n\in L_{i_0+1}:\; |s_i^*(x_n)|>\df{\varepsilon}{p_{i_0}}\}$
 has at most two elements.

 We observe that
 \[F\setminus \{l_{i_0}\}\subset \bigcup\limits_{i=1}^dF_i.\]
 Indeed if $n\in L_{i_0+1}$ and $n\not\in\bigcup\limits_{i=1}^dF_i$
 then by our inductive construction
 $|s_i^*(x_n)|\le\df{\varepsilon}{p_{i_0}}$
 for each $i=1,2,\ldots ,d$  thus
 $|\sum\limits_{i=1}^d\varepsilon_is_i^*(x_n)|\le d\df{\varepsilon}{p_{i_0}}$
 and it follows that
 $|\phi(x_n)|\le\varepsilon$ therefore $n\not\in F$.

 We conclude that
 $\#(F\setminus \{l_{i_0}\})\le 2d$. Also $\min\supp x_n>p_{i_0}\ge d$
 for each $n\in F\setminus \{l_{i_0}\}$ hence the set
 $\left\{\min\supp x_n:\; n\in F\setminus \{l_{i_0}\}\right\}$
 is the union of two  sets belonging to the first Schreier family
 $\mc{S}_1$. Since $\min\supp x_{l_{i_0}}\ge 3$ the set
 $\{\min\supp x_n:\; n\in F\}$ is the union of three sets of $\mc{S}_1$
 and its minimum is greater or equal to 3.
 It follows that
 \[ \{\min\supp x_n:\; n\in F\}\in \mc{S}_2\]
 which completes the proof of the Lemma.
 \end{proof}

 \begin{lemma}\label{lem33}
 For every block subspace $Z$
  of $(\eqs_{\mc{F}_s'})_*$ and  and
 $j>1$, $\e>0$ there exists a $(60,2j,1)$ exact pair $(z,z^*)$ such
 that
 $\dist(z^*,Z)<\e$ and $\|z\|_{\mc{F}_s'}<\frac{2}{m_{2j}}$.
 \end{lemma}
 \begin{proof}[\bf Proof.]
 As in the proof of Theorem 8.3 of \cite{AT1} we may select
 a block sequence \seq{z}{k} and a sequence \seq{f}{k} in $D$
 such that
 \begin{enumerate}
 \item[(i)] Each $z_k$ is a $(\frac{1}{m_{2j_k}},2j_k)$ S.C.C. of
 a normalized block sequence and the sequence $(j_k)_k$ is
 strictly increasing.
 \item[(ii)] $f_k(z_k)>\frac{1}{3}$ and $\ran f_k=\ran z_k$.
 \item[(iii)] $\dist(f_k,Z)<\frac{1}{2^k}$.
 \end{enumerate}
  We assume that the sequence \seq{z}{k}
  is weakly Cauchy. Then the sequence $(z_{2k-1}-z_{2k})_{k\in
  \N}$ is weakly null while
  from Remark \ref{remark} this sequence is
  separated.
 From Lemma \ref{lem35}, for
 every $\e>0$ there exists $L\in[\N]$
 such that for every $\phi\in \mc{F}_s'$,
 $\{\min\supp (z_{2k-1}-z_{2k})
 :\; k\in L,\;|\phi(z_{2k-1}-z_{2k})|>\e\}\in\mc{S}_2$.

  The rest of the proof follows the argument of Theorem 8.3 of
  the Memoirs monograph \cite{AT1}.
 \end{proof}

 \begin{proposition}\label{prop26}
 Every infinite dimensional
 subspace of $(\eqs_{\mc{F}_s'})_*$ has nonseparable second
 dual. In particular the space $(\eqs_{\mc{F}_s'})_*$ contains no
 reflexive subspace.
 \end{proposition}
 \begin{proof}[\bf Proof.]
  Using Lemma \ref{lem33} for every block subspace and every $j$
 we may select,   similarly to Lemma \ref{lem20},
 a $(60,4j-3,1)$ attracting sequence $\chi=(x_k,x_k^*)_{k=1}^d$ with
 $\sum\limits_k\|x_{2k-1}\|_{\mc{F}_s'}<\frac{1}{m_{4j-3}^2}$ and
 $\sum\limits_k\dist(x_{2k-1}^*,Z)<\frac{1}{m_{4j-3}}$.
 We recall at this point (see Definition \ref{depseqfs}(B)(iii))
 that there exists a sequence
 $(a_k)_{k=1}^d$ such that $\sum\limits_{k=1}^da_ke_{t_k}$  is a
$(\frac{1}{m_{4j-3}^2},4j-3)$ basic special convex combination
(where $t_k=\min\supp x_k$).

 We set
$F_{\chi}=-\frac{1}{m_{4j-3}^2}\sum\limits_kx_{2k-1}^*$ and
$g_{\chi}=\frac{1}{m_{4j-3}^2}\sum\limits_kx_{2k}^*$.
 From the fact that $(x_k^*)_{k=1}^d$ is a $\mc{S}_{n_{4j-3}}$
  special sequence we have that
 $\|\frac{1}{m_{2j-1}}(x_1^*+x_2^*+\cdots+x_d^*)\|\le 1$
 and since
  $g_{\chi}-F_{\chi}=\frac{1}{m_{2j-1}}
(\frac{1}{m_{2j-1}}(x_1^*+x_2^*+\cdots+x_d^*))$ we get that
$\|g_{\chi}-F_{\chi}\|\le\frac{1}{m_{2j-1}}$. We also have that
$\dist(F_{\chi},Z)<\frac{1}{m_{4j-3}}$.

Similarly to Proposition 7.5 of \cite{AT1} and to Proposition
\ref{prop11} of the present paper, we may prove that
$\|\sum\limits_{k=1}^d(-1)^ka_kx_k\|\le\frac{300}{m_{4j-3}^2}$.
Observe also that $g_{\chi}(\sum\limits_{k=1}^d(-1)^ka_kx_k)
=\frac{1}{m_{4j-3}^2}\sum\limits_ka_{2k}\ge\frac{1}{3m_{4j-3}^2}$.
From these inequalities it follows that there exists
$1\le\theta_{\chi}\le 900$ such that $g_{\chi}(d_{\chi})=1$ and
$\|d_{\chi}\|\le 900$ where
$d_{\chi}=3\theta_{\chi}m_{4j-3}^2\sum\limits_{k=1}^d(-1)^ka_kx_k$.
It is also easily checked that
$\|d_{\chi}\|_{F_i}\le\frac{1000}{m_{4i-3}^2}$ for $i<j$ while
$\|d_{\chi}\|_{F_i}\le 1000\frac{m_{4j-3}^2}{m_{4i-3}^2}$ if
$i>j$. Thus the vector $d_{\chi}$ witnesses that the functional
$g_{\chi}$ is $(D,j)$ exact.

Using arguments similar to those of Theorem \ref{th3}, for a given
block subspace $Z$ of $(\eqs_{\mc{F}_s'})_*$ we construct a family
$(\chi_a)_{a\in\mc{D}}$ ($\mc{D}$ is the dyadic tree) of dependent
sequences with properties analogous to (i),(ii), (iii) of Theorem
\ref{th3} and such that for every $a\in\mc{D}$ the functional
$g_{\chi_a}$ is $(D,j_a)$ exact. It follows that for every branch
$b$ of the dyadic tree the sum $\sum\limits_{a\in b}g_{\chi_a}$
converges in the $w^*$ topology to a survivor $\sigma_{\mc{F}}$
special functional $g_b\in B_{\eqs_{\mc{F}_s'}^*}$ and there
exists $z_b\in Z^{**}$ with $\|z_b-g_b\|<\frac{1}{1152}$. Then, as
in the proof of Theorem \ref{th3} we obtain that $Z^{**}$ is
nonseparable.
\end{proof}

 \begin{proposition}\label{prop32}
 The space $(\eqs_{\mc{F}_s'})_*$ is Hereditarily Indecomposable.
 \end{proposition}
 \begin{proof}[\bf Proof.]
 Let $Y,Z$ be a pair of block subspaces of $(\eqs_{\mc{F}_s'})_*$.
 For every $j>1$,
 using Lemma \ref{lem33}, we are able to construct a
 $(60,4j-1,1)$ dependent sequence $(x_k,x_k^*)_{k=1}^{n_{4j-1}}$ such that
 $\|x_k\|_{\mc{F}_s'}<\frac{1}{m_{4j-1}^2}$
 while $\sum\dist(x_{2k-1}^*,Y)<$
 and   $\sum\dist(x_{2k}^*,Z)<  $ for each $k$.
 From the observation in Definition \ref{depseqfs}(A)(iii)
 there exist
 $(a_k)_{k=1}^d$ such that $\sum\limits_{k=1}^da_ke_{t_k}$
 is a $(\frac{1}{m_{4j-1}^2},4j-1)$ basic special convex combination.
  As in Proposition \ref{prop24} we get
  that $\|\sum\limits_{k=1}^d(-1)^{k+1}a_{k}x_{k}\|\le
  \frac{480}{m_{4j-1}^2}$.

 We set $h_Y=\frac{1}{m_{4j-1}}\sum\limits_{k\;\text{odd}}x^*_{k}$
 and
 $h_Z=\frac{1}{m_{4j-1}}\sum\limits_{k\;\text{even}}x^*_{k}$.
 The functional
$h_Y+h_Z=\frac{1}{m_{4j-1}}\sum\limits_{k=1}^dx_k^*$ belongs to
the norming set $D$ hence $\|h_Y+h_Z\|\le 1$.
 On the other hand the action of $h_Y-h_Z$ to the vector
 $\sum\limits_{k=1}^d(-1)^{k+1}a_{k}x_{k}$ yields that
 $\|h_Y-h_Z\|\ge\frac{m_{4j-1}}{480}$.

  From the above estimates and since $\dist(h_Y,Y)<1$
  and $\dist(h_Z,Z)<1$ we may choose $f_Y\in Y$ and $f_Z\in Z$
  with
 $\|f_Y-f_Z\|\ge (\frac{m_{4j-1}}{1440}-\frac{2}{3})\|f_Y+f_Z\|$.
 Since this can be done for arbitrary large $j$ we obtain that
 $(\eqs_{\mc{F}_s'})_*$ is Hereditarily Indecomposable.
 \end{proof}

 \begin{proposition}\label{prop27}
 The quotient space $\eqs_{\mc{F}_s'}^*/(\eqs_{\mc{F}_s'})_*$ is
 isomorphic to $c_{0}(\Gamma)$ where the set $\Gamma$ coincides
 with the set of all survivor $\sigma_{\mc{F}}$ special sequences.
 \end{proposition}
 \begin{proof}[\bf Proof.] As follows from \ref{prop25} the
 quotient space $\eqs_{\mc{F}_s'}^*/(\eqs_{\mc{F}_s'})_*$ is
 generated in norm by the classes of the elements of the set
 $\overline{\mc{F}_s'}^{w^*}$.
 Since clearly
 \begin{eqnarray*} \overline{\mc{F}_s'}^{w^*}&=&F\cup\{
 \sum\limits_{i=1}^d\e_i\Phi_i:\;\Phi_i\in\tau,\;\e_i\in\{-1,1\}, \;\min\supp
 \Phi_i\ge d,\\
 && (\ind(\Phi_i))_{i=1}^d\mbox{ are pairwise
 disjoint}\}
 \end{eqnarray*}
 we get that
 \[
 \eqs_{\mc{F}_s'}^*
 =\overline{\spann}(\{e_n^*:\;n\in\N\}\cup\{\Phi:\;
   \Phi
 \mbox{ is a survivor }\sigma_{\mc{F}}\mbox{ special
 functional}\})
 \]
Thus
 $\eqs_{\mc{F}_s'}^*/(\eqs_{\mc{F}_s'})_*
 =\overline{\spann}\{\Phi+(\eqs_{\mc{F}_s'})_*:\; \Phi
 \mbox{ is a survivor }\sigma_{\mc{F}}\mbox{ special functional}\}$
 To prove that this space is isomorphic to $c_0(\Gamma)$
 we shall  show that for
 every choice $(\Phi)_{i=1}^d$ of pairwise different survivor
 $\sigma_{\mc{F}}$ special functionals and every choice of signs
 $(\e_i)_{i=1}^d$ we have that
 \begin{equation} \label{eq26}   \frac{1}{2000}\le
 \|\sum\limits_{i=1}^d\e_i(\Phi_i+(\eqs_{\mc{F}_s'})_*)\| \le 1.
 \end{equation}

 We have that
 $\|\sum\limits_{i=1}^d\e_i(\Phi_i+(\eqs_{\mc{F}_s'})_*)\|
 =\lim\limits_k\|\sum\limits_{i=1}^d\e_i(E_k\Phi_i)\|$
 where for each $k$,
 $E_k=\{k,k+1,\ldots\}$.
 The right part of  inequality \eqref{eq26} follows directly, since for
 all but finite $k$ the functional
 $E_k(\sum\limits_{i=1}^d\e_i\Phi_i)$  belongs to
 $\overline{\mc{F}_s'}^{w^*}\subset B_{\eqs_{\mc{F}_s'}^*}$.

 For each $i=1,\ldots,d$ let
 $\Phi_i=\sum\limits_{l=1}^{\infty}\phi^i_l$ with
 $(\phi^i_l)_{l\in\N}$ a survivor $\sigma_{\mc{F}}$
 special sequence and
 and let $(x^i_l)$ be
 a sequence witnessing this fact (see Remark \ref{rem10} (iii)).
 We choose $k_0$ such that $(\ind(E_{k_0}\Phi_i))_{i=1}^d$ are
 pairwise disjoint and $\min(\bigcup\limits_{i=1}^d
 \ind(E_{k_0}\Phi_i))=r_0$ with $m_{2r_0-1}>10^{10}$.

 Let $k\ge k_0$. We choose $t$ such that $\ran(\phi^1_t)\subset
 E_k$ and let $\ind(\phi^1_t)=l_0$. We get that
 $\sum\limits_{i=2}^d|\Phi_i(x^1_t)|\le
 \sum\limits_{r=r_0}^{l_0-1}\frac{1000}{m_{4r-3}^2}+
 \sum\limits_{r=l_0+1}^{\infty}\frac{1000m_{2l_0-1}^2}{m_{4r-3}^2}<\frac{1}{2}$.
 We thus get that
 \[ \|E_k(\sum\limits_{i=1}^d\e_i\Phi_i)\| \ge
 \frac{1}{1000}(\Phi_1(x^1_t)-\sum\limits_{i=2}^d|\Phi_i(x^1_t)|)>
 \frac{1}{1000}(1-\frac{1}{2})=\frac{1}{2000}.\]

 The proof of the proposition is complete.
 \end{proof}

 \begin{theorem}\label{th5}
 There exists a Banach space $\eqs_{\mc{F}_s'}$ satisfying
  the following
 properties:
 \begin{enumerate}
 \item[(i)] The space $\eqs_{\mc{F}_s'}$ is an asymptotic $\ell_1$
 space with a boundedly complete Schauder basis \seq{e}{n} and is
 Hereditarily Indecomposable and reflexively saturated.
 \item[(ii)] The predual space
 $(\eqs_{\mc{F}_s'})_*=\overline{\spann}\{e_n^*:\;n\in\N\}$ is
 Hereditarily Indecomposable and each infinite dimensional
 subspace  of $(\eqs_{\mc{F}_s'})_*$ has nonseparable second
 dual. In particular the space $(\eqs_{\mc{F}_s'})_*$ contains no
 reflexive subspace.
 \item[(iii)] The dual space $\eqs_{\mc{F}_s'}^*$ is nonseparable,
 Hereditarily Indecomposable and contains no reflexive subspace.
 \item[(iv)] Every bounded linear operator $T:X\to X$
 where $X=(\eqs_{\mc{F}_s'})_*$ or $X=\eqs_{\mc{F}_s'}$
 or $X=\eqs_{\mc{F}_s'}^*$ takes the form $T=\lambda I+W$ with $W$
 a weakly compact operator. In particular each
 $T:\eqs_{\mc{F}_s'}^* \to \eqs_{\mc{F}_s'}^*$ is of the form
 $T=Q^*+K$ with $Q:\eqs_{\mc{F}_s'}\to\eqs_{\mc{F}_s'}$ and $K$ a
compact operator, hence $T=\lambda I+R$ with  $R$ an operator with separable range.
 \end{enumerate}
 \end{theorem}
 \begin{proof}[\bf Proof.]
 As we have observed the Schauder basis \seq{e}{n} of
 $\eqs_{\mc{F}_s'}$ is  boundedly complete and
 $\eqs_{\mc{F}_s'}$ is asymptotic $\ell_1$. In Proposition \ref{prop24}
 we have shown that  $\eqs_{\mc{F}_s'}$ is reflexively saturated and
 Hereditarily Indecomposable. The facts that $(\eqs_{\mc{F}_s'})_*$
 is Hereditarily Indecomposable and that every subspace of it has
 nonseparable second dual have been shown in Proposition
 \ref{prop26} and Proposition \ref{prop32}.

 From the facts that the quotient space $\eqs_{\mc{F}_s'}^*/(\eqs_{\mc{F}_s'})_*$ is
 isomorphic to $c_{0}(\Gamma)$ and $(\eqs_{\mc{F}_s'})_*$
 is Hereditarily Indecomposable and taking into account that
 $\eqs_{\mc{F}_s'}$, being a Hereditarily Indecomposable space,
 contains no isomorphic copy of $\ell_1$ we get that the dual space
 $\eqs_{\mc{F}_s'}^*$ is also Hereditarily Indecomposable
 (Corollary 1.5 of \cite{AT1}).
 Since $\eqs_{\mc{F}_s'}^*$ is  Hereditarily Indecomposable and
 contains a subspace (which is $(\eqs_{\mc{F}_s'})_*$ with no
 reflexive subspace we conclude that
 $\eqs_{\mc{F}_s'}^*$ also does not have any reflexive subspace.

 Using similar arguments to those of the proof of
 Theorems \ref{th7}, \ref{aa12} and Corollary \ref{cor6} we may prove
 that every bounded linear operator
  $T:(\eqs_{\mc{F}_s'})_*\to (\eqs_{\mc{F}_s'})_*$ and every bounded
  linear operator
 $T:\eqs_{\mc{F}_s'}\to \eqs_{\mc{F}_s'}$  takes the form $T=\lambda I+W$ with $W$
 a strictly singular and weakly compact operator.
 Since $\eqs_{\mc{F}_s'}$ contains no isomorphic copy of $\ell_1$
 and $\eqs_{\mc{F}_s'}^{**}$ is isomorphic to $\eqs_{\mc{F}_s'}\oplus\ell_1(\Gamma)$
 Proposition 1.7 of \cite{AT1} yields that every bounded linear
 operator $T:\eqs_{\mc{F}_s'}^*\to\eqs_{\mc{F}_s'}^*$ is of the form
 $T=Q^*+K$ with $Q:\eqs_{\mc{F}_s'}\to \eqs_{\mc{F}_s'}$
  and $K$ a compact operator. From the
 form of the operators of $\eqs_{\mc{F}_s'}$ we have mentioned
 before we conclude that $T$ takes the form $T=\lambda I+R$ with
 $R$ a weakly compact operator and hence of separable range.
 \end{proof}

\begin{remark}
It is worth mentioning that the key ingredient to obtain
$\eqs_{\mc{F}_s'}^*/(\eqs_{\mc{F}_s'})_*$ isomorphic to
$c_0(\Gamma)$ which actually yields the HI property of
$\eqs_{\mc{F}_s'}^*$ is that in the ground set $\mc{F}_s'$ we
connect the $\sigma_{\mc{F}}$ special functionals using the
Schreier operation. This forces us to work with the saturation
families $(\mc{S}_{n_j},\frac{1}{m_j})_j$ instead of
$(\mc{A}_{n_j},\frac{1}{m_j})_j$. The reason for this is that
working with $\mc{F}_s'$ built on $(F_j)_j$ with
$F_j=\big\{\frac{1}{m_{4j-3}^2}\sum\limits_{i\in I}\pm
e_i^*:\;\#(I)\le\frac{n_{4j-3}}{2}\big\}$ the extension with
attractors of this ground set $\mc{F}_s'$ based on
$(\mc{A}_{n_j},\frac{1}{m_j})_j$ is not strongly strictly
singular.

 However there exists an alternative way of connecting the
 $\sigma_{\mc{F}}$ special functionals lying between the Schreier
 operation and the $\ell_2$ sums. This yields the James Tree
 space $JT_{\mc{F}_{2,s}}$ defined and studied in Appendix B. It
 is easy to check that the corresponding HI extension with
 attractors $\eqs_{\mc{F}_{2,s}'}$ of $JT_{\mc{F}_{2,s}'}$ is a
 strictly singular one either we work on
 in the frame of $(\mc{A}_{n_j},\frac{1}{m_j})_j$ or of
 $(\mc{S}_{n_j},\frac{1}{m_j})_j$.
 It is  open   whether the corresponding space
 $\eqs_{\mc{F}_{2,s}'}^*$  contains $\ell_2$ or not. If
 it does not contain $\ell_2$ then  $\eqs_{\mc{F}_{2,s}'}^*$
 will be also  a nonseparable
 HI space not containing any reflexive subspace
 with the additional property that
 $\eqs_{\mc{F}_{2,s}'}^*/(\eqs_{\mc{F}_{2,s}'})_*$ is isomorphic
 to $\ell_2(\Gamma)$.
\end{remark}

\section{A HJT space with  unconditionally and reflexively saturated dual}
 This section concerns the definition of the space
 $\eqs_{\mc{F}_2}^{us}$ namely a separable space with a boundedly
 complete basis which is reflexive and unconditionally saturated
 and its predual $(\eqs_{\mc{F}_2}^{us})_*$ is  HJT space hence it
 does not contain any reflexive subspace. This construction starts
 with the ground set $\mc{F}_2$ used in Section 4. In the
 extensions we use only attractors for which we eliminate a sufficient
 part of their conditional structure. The proof of the property
 that $\eqs_{\mc{F}_2}^{us}$ is unconditionally saturated follows
 the arguments of \cite{AM},\cite{AT2} while the HJT property of
 $(\eqs_{\mc{F}_2}^{us})_*$ results from the remaining part of the
 conditional structure of the attractors. We additionally  show
 that $(\eqs_{\mc{F}_2}^{us})_*$  is HI.

 Let $\mathbf{Q}$ be the set of all finitely supported
 scalar sequences with
 rational coordinates, of maximum modulus 1  and  nonempty support.
 We  set
 \begin{eqnarray*}
 \mathbf{Q_s}& = & \{(x_1,f_1,\ldots,x_n,f_n):\; x_i,f_i\in
 \Q,\;i=1,\ldots,n\\
 & & \ran(x_i)\cup\ran(f_i)<\ran(x_{i+1})\cup\ran(f_{i+1}),\;
 i=1,\ldots,n-1\}.
 \end{eqnarray*}
 For $\phi=(x_1,f_1,\ldots,x_n,f_n)\in \mathbf{Q_s}$ and $l\le n$
 we denote by $\phi_l$ the sequence $(x_1,f_1,\ldots,x_l,f_l)$.
 We consider an injective coding function
 $\sigma:\mathbf{Q_s}\to\{2j:\; j\in\mathbb{N}\}$
  such that for every
 $\phi=(x_1,f_1,\ldots,x_n,f_n)\in\mathbf{Q_s}$
 \begin{eqnarray*}
 \sigma(x_1,f_1,\ldots,x_{n-1},f_{n-1})<
 \sigma(x_1,f_1,\ldots,x_{n},f_{n}) \\ \mbox{ and }
 \max\{\ran(x_n)\cup\ran(f_n)\}\le
 m_{\sigma(\phi)}^{\frac{1}{2}}.
 \end{eqnarray*}

 The norming set
 $D^{us}$ of the space $\eqs_{\mc{F}_{2}}^{us}$ will be defined
 as $D^{us}=\bigcup\limits_{n=0}^{\infty}D_n$ after defining
 inductively two sequences $(K_n)_{n=0}^{\infty}$, $(D_n)_{n=0}^{\infty}$
 of subsets of \co with $D_n=\conv_{\mathbb{Q}}(K_n)$.

 Let $\mc{F}_{2}$ be the set defined in the beginning of the third
 section.
 We set \[K_0=\mc{F}_{2}\;\;\text{ and }
 \; D_0=\convq(K_0).\]
 Assume that $K_{n-1}$ and $D_{n-1}$ have been defined.
 Then for each $j\in \N$ we set
 \begin{equation*}
 K_n^{2j}=K_{n-1}^{2j}\cup\{\df{1}{m_{2j}}\sum\limits_{i=1}^df_i:\;
 f_1<\cdots<f_d,\; f_i\in D_{n-1},\; d\le n_{2j}\}.
 \end{equation*}
 For fixed $j\in\N$ we consider the collection of all sequences
 $\phi=(x_i,f_i)_{i=1}^{n_{4j-3}}$ satisfying the following
 conditions:
 \begin{enumerate}
 \item[(i)] $x_1=e_{l_1}$ and $f_1=e_{l_1}^*$ for some $l_1\in
 \Lambda_{2j_1}$ where $j_1$ is an integer with
 $m_{2j_1}^{1/2}>n_{4j-3}$.
 \item[(ii)] For $1\le i\le n_{4j-3}/2$, $f_{2i}\in
 K_{n-1}^{\sigma(\phi_{2i-1})}$ and
 $\|x_{2i}\|_{K_{n-1}}\le\frac{18}{m_{\sigma(\phi_{2i-1})}}$.
 \item[(iii)] For $1\le i< n_{4j-3}/2$, $x_{2i+1}=e_{l_{2i+1}}$ and
 $f_{2i+1}=e_{l_{2i+1}}^*$ for some $l_{2i+1}\in
 \Lambda_{\sigma(\phi_{2i})}$.
 \end{enumerate}
 For every  $\phi$ satisfying (i),(ii) and (iii) we define the set
  \begin{eqnarray*}
   K_{n,\phi}^{4j-3} & = & \big\{ \df{\pm 1}{m_{4j-3}}
   E\big(\sum\limits_{i=1}^{n_{4j-3}/2}(\lambda_{f'_{2i}}f_{2i-1}+f'_{2i})\big):\;
   E\text{ is an interval of }\N,\\
  & &
    f'_{2i}\in
   K_{n-1}^{\sigma(\phi_{2i-1})},\;\;
   \lambda_{f'_{2i}}=f'_{2i}(m_{\sigma(\phi_{2i-1})}x_{2i}),\\
   &&
  (x_{2i-1},f_{2i-1},x_{2i},f'_{2i})_{i=1}^{n_{4j-3}/2}\in \mathbf{Q_s}
    \big\}.
 \end{eqnarray*}
 We define
 \[ K_n^{4j-3}=\cup\{K^{4j-3}_{n,\phi}:\;\phi \text{ satisfies
 conditions (i), (ii), (iii)} \}\cup K_{n-1}^{4j-3}   \]

 and we set
 \[  K_n=\bigcup\limits_j (K_n^{2j}\cup K_{n}^{4j-3})\; \text{ and
 }\;D_n=\convq(K_n). \]

 We finally set \[K^{us}=\bigcup\limits_{n=0}^{\infty}K_n
 \qquad\text{and}\qquad
         D^{us}=\bigcup\limits_{n=0}^{\infty}D_n \]

 The space $\eqs^{us}_{\mc{F}_{2}}$ is the completion of the space
 $(c_{00},\|\;\|_{D^{us}})$ where
 \[ \norm[x]_{D^{us}}=\sup\{f(x):\; f\in D^{us}\}.\]

 Using the same arguments as those in Proposition
 \ref{propa1} we get the following.

 \begin{lemma}\label{lem37}
 The identity operator $I:\eqs^{us}_{\mc{F}_{2}}\to
 JT_{\mc{F}_{2}}$ is strongly strictly singular (Definition \ref{aa11}).
 \end{lemma}

 \begin{definition}
 The sequence
 \[\phi=(x_1,f_1,x_2,f_2,x_3,f_3,x_4,f_4,\ldots,x_{n_{4j-3}},f_{n_{4j-3}})
 \in\mathbf{Q_s}\]
 is said to be a $n_{4j-3}$ attracting sequence provided that
 \begin{enumerate}
  \item[(i)] $x_1=e_{l_1}$ and $f_1=e_{l_1}^*$ for some $l_1\in
 \Lambda_{2j_1}$ where $j_1$ is an integer with
 $m_{2j_1}^{1/2}>n_{4j-3}$.
 \item[(ii)] For $1\le i\le n_{4j-3}/2$,
 $(m_{\sigma(\phi_{2i-1})}x_{2i},f_{2i})$
 is a $(18,\sigma(\phi_{2i-1}),1)$ exact pair (Definition
 \ref{def19}) while
 $\sum\limits_{i=1}^{n_{4j-3}/2}\|x_{2i}\|_{\mc{F}_2}<\frac{1}{n_{4j-3}}$.
  \item[(iii)] For $1\le i< n_{4j-3}/2$, $x_{2i+1}=e_{l_{2i+1}}$ and
 $f_{2i+1}=e_{l_{2i+1}}^*$ for some $l_{2i+1}\in
 \Lambda_{\sigma(\phi_{2i})}$.
 \end{enumerate}
 We consider the vectors $d_{\phi}$ in $\eqs^{us}_{\mc{F}_{2}}$,
 and $g_{\phi},F_{\phi}$ in $(\eqs^{us}_{\mc{F}_{2}})_*$ as they are
 defined in Definition \ref{def12}. Let also notice, for later use,
  that the
 analogue of Lemma \ref{lem21} remains valid.
 \end{definition}

  \begin{lemma}\label{lem36}
  For every block subspace $Z$ of the predual space
 $(\eqs^{us}_{\mc{F}_{2}})_*$ and every $j\in N$ there exists
  a
 $n_{4j-3}$ attracting sequence\\
 $\phi=(x_1,f_1,x_2,f_2,\ldots,x_{n_{4j-3}},f_{n_{4j-3}})$ with
 $\sum\limits_{i=1}^{n_{4j-3}/2}\dist(f_{2i},Z)<\frac{1}{m_{4j-3}^2}$.
 \end{lemma}
 \begin{proof}[\bf Proof.]
 Since the identity $I:\eqs^{us}_{\mc{F}_{2}}\to
 JT_{\mc{F}_{2}}$ is strongly strictly singular (\ref{lem37}), we
 may construct, using the analogue of Lemma \ref{lem25} in terms
 of $\eqs^{us}_{\mc{F}_{2}}$, the desired attracting sequence.
 \end{proof}

 \begin{lemma}\label{lem38}
 Let $\chi=(x_{2k},x_{2k}^*)_{k=1}^{n_{4j-3}/2}$ be a $(18,4j-3,1)$
 attracting sequence such that
 $\sum\limits_{k=1}^{n_{4j-3}/2}\|x_{2k-1}\|_{\mc{F}_{2}}<\frac{1}{m_{4j-3}^2}$.
 Then for every branch $b$ such that $j\not\in\ind(b)$ we have
 that
 $|b^*(d_{\chi})|<\frac{3}{m_{4j-3}}$.
 \end{lemma}
 \begin{proof}[\bf Proof.]
 Let $b=(f_1,f_2,f_3,\ldots)$ be a branch (i.e. \seq{f}{i} is a
 $\sigma_{\mc{F}}$ special sequence).
 We recall that
 $d_{\chi}=\frac{m_{4j-3}^2}{n_{4j-3}}\sum\limits_{k=1}^{n_{4j-3}}(-1)^kx_k$
 and we set\\
 $d_1=-\frac{m_{4j-3}^2}{n_{4j-3}}\sum\limits_{k=1}^{n_{4j-3}/2}x_{2k-1}$
 and
 $d_2=\frac{m_{4j-3}^2}{n_{4j-3}}\sum\limits_{k=1}^{n_{4j-3}/2}x_{2k}$.

 Our assumption
 $\sum\limits_{k=1}^{n_{4j-3}/2}\|x_{2k-1}\|_{\mc{F}_{2}}<\frac{1}{m_{4j-3}^2}$
 yields that
 \begin{equation}\label{eq22}
 |b^*(d_1)|\le
 \frac{m_{4j-3}^2}{n_{4j-3}}\sum\limits_{k=1}^{n_{4j-3}/2}|b^*(x_{2k-1})|
 \le
 \frac{m_{4j-3}^2}{n_{4j-3}}\sum\limits_{k=1}^{n_{4j-3}/2}\|x_{2k-1}\|_{\mc{F}_{2}}
 <\frac{1}{n_{4j-3}}.
 \end{equation}

 We decompose $b^*$ as $b^*=x^*+y^*$ with
 $\ind(x^*)\subset\{1,\ldots,j-1\}$ and
 $\ind(y^*)\subset\{j+1,j+2,\ldots\}$.
 We recall  that an $f\in F$ with $\inf(f)=l$ is of the form
 $f=\frac{1}{m_{4l-3}^2}\sum\limits_{i\in \supp(f)}\pm e_i^*$ with
 $\supp(f)\le n_{4l-3}/2$.
  Thus
 $\supp(x^*)\le
 \frac{n_1}{2}+\frac{n_5}{2}+\cdots+\frac{n_{4j-7}}{2}<n_{4j-4}$.
 Hence
 \begin{equation}\label{eq23}
 |x^*(d_2)|\le
 \frac{m_{4j-3}^2}{n_{4j-3}}n_{4j-4}<\frac{1}{m_{4j-3}}.
 \end{equation}
 On the other hand we have that $\|y^*\|_{\infty}\le
 \frac{1}{m_{4j+1}}$, therefore
 \begin{equation}\label{eq24}
 |y^*(d_2)|\le
  \frac{1}{m_{4j+1}}\cdot\frac{m_{4j-3}^2}{n_{4j-3}}\cdot\frac{n_{4j-3}}{2}
  <\frac{1}{m_{4j-3}}.
 \end{equation}

 From \eqref{eq22},\eqref{eq23} and  \eqref{eq24} we obtain that
 $|b^*(d_{\chi})|<\frac{3}{m_{4j-3}}$.
 \end{proof}

 \begin{proposition}\label{prop33}
 The space $(\eqs^{us}_{\mc{F}_{2}})_*$ is Hereditarily
 Indecomposable.
 \end{proposition}
 \begin{proof}[\bf Proof.]
 Let $Z_1$, $Z_2$ be a pair of block subspaces in $(\eqs^{us}_{\mc{F}_{2}})_*$
 and let $0<\delta<1$.

 We may inductively construct, using Lemma \ref{lem36},
 a sequence $(\chi_r)_{r\in\N}$ such
 that the following conditions are satisfied.
 \begin{enumerate}
 \item[(i)] Each $\chi_r=(x^r_k,(x^r_k)^*)_{k=1}^{n_{4j_r-3}}$ is
 a $(18,4j_r-3,1)$ attracting sequence with
 $\sum\limits_{k=1}^{n_{2j_r-1}}\|x^r_{2k-1}\|_{\mc{F}_{2}}<\frac{1}{m_{4j_r-3}^2}$
 and additionally \\$\dist(F_{\chi_r},Z_1)<\frac{1}{m_{2j_r-1}}$ if
 $r$ is odd, while $\dist(F_{\chi_r},Z_2)<\frac{1}{m_{2j_r-1}}$ if
 $r$ is even.
 \item[(ii)] $(d_{\chi_r})_{r\in\N}$ is a block sequence.
 \item[(iii)] For $r>1$,
 $j_r=\sigma_{\mc{F}}(g_{\chi_1},\ldots,g_{\chi_{r-1}})$.
 \end{enumerate}

 \begin{claim}
 The sequence $(d_{\chi_{2r-1}}-d_{\chi_{2r}})_{r\in\N}$ is
 a weakly null sequence in $\eqs_{\mc{F}_{2}}$.
 \end{claim}
 \begin{proof}[\bf Proof of the claim.]
 From the analogue of Proposition \ref{prop14}
 the space $\eqs_{\mc{F}_{2}}^*$ is the
 closed linear span of the pointwise closure
 $\overline{\mc{F}_{2}}^{w^*}$ of the set $\mc{F}_{2}$.
 From the observation after Theorem \ref{th8} we have that
 \begin{eqnarray*}
 \overline{\mc{F}_{2}}^{w^*}&=&F_0\cup
 \big\{\sum\limits_{i=1}^{\infty}a_ix_i^*:\;\sum\limits_{i=1}^{\infty}a_i^2\le1,\;
 (x_i^*)_{i=1}^{d}\mbox{ are $\sigma_{\mc{F}}$ special functionals}\\
  && \mbox{with }(\ind(x_i^*))_{i=1}^d
 \mbox{ pairwise disjoint }\min\supp x_i^*\ge d\big\}.
 \end{eqnarray*}
 Thus it is enough to show that
 $b^*(d_{\chi_{2r-1}}-d_{\chi_{2r}})\stackrel{r\to
 \infty}{\longrightarrow}0$ for every branch $b$.

 Let $b$ be an arbitrary branch. If
 $b=(g_{\chi_1},g_{\chi_2},g_{\chi_3},g_{\chi_4},\ldots)$
 from  we obtain that
 $g(d_{\chi_{2r-1}}-d_{\chi_{2r}})
 =g_{\chi_{2r-1}}(d_{\chi_{2r-1}})-g_{\chi_{2r}}(d_{\chi_{2r}})
 =\frac{1}{2}-\frac{1}{2}=0$ for every $r$.
 If $b\neq (g_{\chi_1},g_{\chi_2},g_{\chi_3},g_{\chi_4},\ldots)$
 the injectivity of the coding function $\sigma_{\mc{F}}$ yields that
 there exist $r_0\in\N$ such that $j_r\not\in \ind(b^*)$ for all
 $r> 2r_0$.  Hence, for $r>r_0$, Lemma \ref{lem38} yields that
 $|b(d_{\chi_{2r-1}}-d_{\chi_{2r}})|\le|b(d_{\chi_{2r-1}})|+|b(d_{\chi_{2r}})|
 <\frac{3}{m_{4j_{2r-1}-3}}+\frac{3}{m_{4j_{2r}-3}}<\frac{4}{m_{4j_{2r-1}-3}}$
 and therefore
 $b^*(d_{\chi_{2r-1}}-d_{\chi_{2r}})\stackrel{r\to
 \infty}{\longrightarrow}0$.

 The proof of the claim is complete.
 \end{proof}

 It follows from the claim that there exists a convex combination
 of the sequence $(d_{\chi_{2r-1}}-d_{\chi_{2r}})_{r\in\N}$ with
 norm less than $\frac{\delta}{3}$; let $(a_r)_{r=1}^d$
 be nonnegative reals with $\sum\limits_{r=1}^da_r=1$ such that
 $\|\sum\limits_{r=1}^da_r(d_{\chi_{2r-1}}-d_{\chi_{2r}})\|
 <\frac{\delta}{3}$.

 We set $g=\sum\limits_{r=1}^{2d}g_{\chi_r}$ and
 $g'=\sum\limits_{r=1}^{d}(g_{\chi_{2r-1}}-g_{\chi_{2r}})$.
 Since $g\in \mc{F}_{2}$ we have that $\|g\|\le1$.
 On the other hand
 \begin{eqnarray*}
 \|g'\|&\ge& \frac{g'\big(\sum\limits_{r=1}^da_r(d_{\chi_{2r-1}}
 -d_{\chi_{2r}})\big)}{\|\sum\limits_{r=1}^da_r(d_{\chi_{2r-1}}-d_{\chi_{2r}})\|}\\
  &>&\frac{\sum\limits_{r=1}^da_r(g_{\chi_{2r-1}}-g_{\chi_{2r}})
 (d_{\chi_{2r-1}}-d_{\chi_{2r}})}{\frac{\delta}{3}}
  =\frac{\sum\limits_{r=1}^d
 a_r(\frac{1}{2}+\frac{1}{2})}{\frac{\delta}{3}}=\frac{3}{\delta}.
 \end{eqnarray*}

 For each $r\le 2d$ with $r$ odd we select $z_r^*\in Z_1$ such that
 $\|z_r^*-F_{\chi_r}\|<\frac{1}{m_{4j_r-3}}$,
 while  for $r$ even we
 select $z_r^*\in Z_2$ such that
 $\|z_r^*-F_{\chi_r}\|<\frac{1}{m_{4j_r-3}}$.
 We set
 \[F_1=\sum\limits_{r=1}^dz_{2r-1}^*(\in Z_1)
 \quad\mbox{and}\quad F_2=\sum\limits_{r=1}^dz_{2r}^*(\in Z_2).\]

 From our choice of $z_r^*$ and the analogue of
 Lemma \ref{lem21} we get that
 \begin{equation}\label{eq25}
 \sum\limits_{r=1}^{2d}\|z_r-g_{\chi_r}\|
 \le\sum\limits_{r=1}^{2d}(\|z_r-F_{\chi_r}\|+\|F_{\chi_r}-g_{\chi_r}\|)
 \le\sum\limits_{r=1}^{2d}(\frac{1}{m_{4j_r-3}}+\frac{1}{m_{4j_r-3}})<1.
 \end{equation}

 From \eqref{eq25} we obtain that $\|(F_1+F_2)-g\|<1$ and
 $\|(F_1-F_2)-g'\|<1$.
 Thus, the facts that $\|g\|\le 1$ and $\|g'\|> \frac{3}{\delta}$
 yield that
 $\|F_1+F_2\|<2$ and $\|F_1-F_2\|>
 \frac{3}{\delta}-1>\frac{2}{\delta}$,
 therefore
 $\|F_1+F_2\|<\delta\|F_1+F_2\|$. The proof of the proposition is
 complete.
 \end{proof}

 \begin{proposition}\label{prop34}
 The space $(\eqs^{us}_{\mc{F}_{2}})_*$ is HJT. In particular
 the space $(\eqs^{us}_{\mc{F}_{2}})_*$
  contains no reflexive subspace
  and every infinite dimensional subspace $Z$ of
  $(\eqs^{us}_{\mc{F}_{2}})_*$
 has nonseparable second dual
  $Z^{**}$.
 \end{proposition}
 \begin{proof}[\bf Proof.]
 The proof is identical to that of Theorem \ref{th3}.
 \end{proof}

 \begin{theorem}
 The space $\eqs^{us}_{\mc{F}_{2}}$ has the following properties
 \begin{enumerate}
 \item[(i)] Every subspace $Y$ of $\eqs^{us}_{\mc{F}_{2}}$
 contains a further subspace $Z$ which is reflexive and has an
 unconditional basis.
 \item[(ii)] The predual $(\eqs^{us}_{\mc{F}_{2}})_*$ of the space
 $\eqs^{us}_{\mc{F}_{2}}$
 is Hereditarily
 Indecomposable and has no reflexive subspace.
 \end{enumerate}
 \end{theorem}
 \begin{proof}[\bf Proof.]
  First,  the identity operator $I:\eqs^{us}_{\mc{F}_{2}}\to
 JT_{\mc{F}_{2}}$, being
 strongly strictly singular, is strictly singular
 (Proposition \ref{propsss}). Therefore the space
 $\eqs^{us}_{\mc{F}_{2}}$ is reflexively saturated.
  Let $Z$ be an arbitrary block subspace of the space
 $\eqs^{us}_{\mc{F}_{2}}$. From the fact that
 $I:\eqs^{us}_{\mc{F}_{2}}\to JT_{\mc{F}_{2}}$
 is strictly singular we may choose a block sequence
 \seq{z}{k} in $Z$ with $\|z_k\|=1$ and
 $\sum\limits_{k=1}^{\infty}\|z_k\|_{\mc{F}_{2}}<\frac{1}{16}$.
 We may prove that \seq{z}{k} is an unconditional basic sequence
 following the procedure used in the proof of Proposition 3.6 of
 \cite{AT2}.

  The facts that the space $(\eqs^{us}_{\mc{F}_{2}})_*$ is
 Hereditarily
 Indecomposable and has no reflexive subspace have been proved in
 Propositions
 \ref{prop33} and \ref{prop34}.
 \end{proof}

 Defining the norming set of the present section
 using $\mc{F}_s$ (instead of $\mc{F}_2$)
 in the first inductive step (namely in the definition of $K_0$)
 and using the saturation methods $(\mc{S}_{n_j},\frac{1}{m_j})_j$
 (instead of $(\mc{A}_{n_j},\frac{1}{m_j})_j$) we produce a Banach
 space $\mathfrak X_{\mc{F}_s}^{us}$ which is unconditionally
 saturated while its predual and its dual share similar
 properties with the space $\eqs_{\mc{F}_s}$ of Section \ref{sec4}.
 Namely we have the following.

 \begin{theorem} There exists a Banach space $\eqs_{\mc{F}_s}^{us}$
 with the properties:
 \begin{enumerate}
 \item[(i)] The predual $(\eqs_{\mc{F}_s}^{us})_*$ of $\eqs_{\mc{F}_s}^{us}$
  is HI and every
 infinite dimensional subspace of $(\eqs_{\mc{F}_s}^{us})_*$ has
 nonseparable second dual. In particular
 $(\eqs_{\mc{F}_s}^{us})_*$ contains no reflexive subspace.
 \item[(ii)] The space $\eqs_{\mc{F}_s}^{us}$ is unconditionally
 and reflexively saturated.
 \item[(iii)] The dual space $(\eqs_{\mc{F}_s}^{us})^*$ is
 nonseparable HI and contains no reflexive subspace.
 \end{enumerate}
 \end{theorem}

%%%%%%%%%%%%%%%%%%%%%%%%%%%%%%%%%

  \begin{appendix}
 \section{The auxiliary space and the basic inequality}

 The basic inequality is the main tool in providing upper bounds
 for the action of functionals on certain vectors of $X_G$. It has
 appeared in several variants in previous works like \cite{AT1},
 \cite{ALT}, \cite{ArTo}. In this section we present another
 variant which mainly concerns the case of strongly strictly
 singular extensions and in particular we provide the proof of
 Proposition \ref{prop13} stated in Section 1.
  The proof of the present variant follows the
 same lines as the previous ones.

 \begin{definition}\label{tree}
   {\bf The tree $T_f$
 of a functional $f\in W$.}
 Let $f\in D$. By a  tree of $f$ (or tree corresponding to
 the analysis of $f$) we mean a finite family
 $T_f=(f_a)_{a\in\mc{A}}$ indexed by a finite tree $\mc{A}$ with a unique
 root $0\in\mc{A}$ such that the following conditions are satisfied:
 \begin{enumerate}
 \item[1.] $f_0=f$ and $f_a\in D$ for all $a\in \mc{A}$.
 \item[2.] An $a\in\mc{A}$ is maximal if and only if $f_a\in G$.
 \item[3.] For every $a\in \mc{A}$ which is not maximal,
 denoting by $S_a$ the set of the
 immediate successors of $a$, exactly one of the following  holds:
 \begin{enumerate}
  \item[(a)] $S_a=\{\beta_1,\ldots,\beta_d\}$ with $f_{\beta_1}<\cdots<f_{\beta_d}$
 and there exists $j\in\mathbb{N}$ such that
 $d\le n_j$  and
 $f_a=\frac{1}{m_j}\sum\limits_{i=1}^df_{\beta_i}$ (recall that
 in this case we say that $f_a$ is of type I).
 \item[(b)] $S_a=\{\beta_1,\ldots,\beta_d\}$
  and there exists a family of positive rationals
 $\{r_{\beta_i}:\;i=1,\ldots,d\}$
 with $\sum\limits_{i=1}^dr_{\beta_i}=1$
 such that $f_a=\sum\limits_{i=1}^dr_{\beta_i}f_{\beta_i}$.
 Moreover for all $i=1,\ldots,d,$\, $\ran f_{\beta_i}\subset \ran f_a$.
 (recall that in this case we say that $f_a$ is of type II).
 \end{enumerate}
 \end{enumerate}
 \end{definition}
 It is obvious that every $f\in D$
 has a tree which is not necessarily unique.

 \begin{definition}{\bf (The auxiliary space $T_{j_0}$)} \label{auxi}
 Let $j_0>1$ be fixed. We set $C_{j_0}=\{\sum\limits_{i\in F}\pm
 e_i^*:\; \#(F)\le n_{j_0-1}\}$.

  The auxiliary space $T_{j_0}$
 is the completion of $(\co,\|\;\|_{D_{j_0}})$ where the  norming
 set $D_{j_0}$ is defined  to be the minimal subset of \co
 which (i) Contains $C_{j_0}$. (ii) It is closed under
 $(\mc{A}_{5n_j},\frac{1}{m_j})$ operations for all $j\in\N$.
 (iii) It is rationally convex.

 Observe that the Schauder basis \seq{e}{n} of $T_{j_0}$ is
 1-unconditional.
 \end{definition}

 \begin{remark}\label{rem5}
 Let $D_{j_0}'$ be the
 minimal subset of \co which (i) Contains $C_{j_0}$. (ii) Is closed
 under
 $(\mc{A}_{5n_j},\frac{1}{m_j})$ operations for all $j\in\N$.
 We notice that each $f\in D_{j_0}'$ has a tree $(f_a)_{a\in\mc{A}}$
 in which for $a\in \mc{A}$ which is not maximal, $f$ is the result
 of an $(\mc{A}_{5n_j},\frac{1}{m_j})$ operation (for some $j$) of
 the functionals $(f_{\beta})_{\beta\in S_a}$.

  It can be shown
 that $D_{j_0}'$ is also a norming set for
 the space $T_{j_0}$ and that for every $j\in\N$ we have that
 $\conv_{\mathbb{Q}}\{f\in D_{j_0}: w(f)=m_j\}
 =\conv_{\mathbb{Q}}\{f\in D_{j_0}': w(f)=m_j\}$.
 For proofs of similar results in a different context we refer to
  \cite{AT1} (Lemma 3.5).
  \end{remark}

  \begin{lemma} \label{lem14}
 Let $j_0\in \N$ and $f\in D_{j_0}'$. Then for every family
 $k_1<k_2<\ldots<k_{n_{j_0}}$ we have that
  \begin{equation}\label{eq13}
 |f(\frac{1}{n_{j_0}}\sum\limits_{l=1}^{n_{j_0}}e_{k_{l}})| \leq
 \begin{cases}
 \frac{2}{m_i\cdot m_{j_0}},\quad &\text{if}\,\,\,w(f)=m_{i},\;i<j_0\\
 \frac{1}{m_i},\quad &\text{if}\,\,\,w(f)= m_{i},\;i\ge j_0
 \end{cases}
 \end{equation}
  In particular
 $\|\frac{1}{n_{j_0}}\sum\limits_{l=1}^{n_{j_0}}e_{k_{l}}\|_{D_{j_0}}
 \le\frac{1}{m_{j_0}}$.

 If we additionally assume that the functional $f$ admits a tree
 $(f_{\alpha})_{\alpha\in\mc{A}}$ such that
  $w(f_{\alpha})\not= m_{j_0}$
 for every $\alpha\in\mc{A}$, then  we have that
  \begin{equation}\label{eq14}
 |f(\frac{1}{n_{j_0}}\sum\limits_{l=1}^{n_{j_0}}e_{k_{l}})| \leq
 \begin{cases}
 \frac{2}{m_i\cdot m_{j_0}^2},\quad &\text{if}\,\,\,w(f)=m_{i},\;i<j_0\\
 \frac{1}{m_i},\quad &\text{if}\,\,\,w(f)= m_{i},\;i> j_0
 \end{cases}
 \le \frac{1}{m_{j_0}^2}.
 \end{equation}
 \end{lemma}
 \begin{proof}[\bf Proof.]
 We  first prove the following claim.
 \begin{claim} Let $h\in D_{j_0}'$. Then
 \begin{enumerate}
 \item[(i)]  $\#\{k:\;
 |h(e_k)|>\frac{1}{m_{j_0}}\}<(5n_{j_0-1})^{\log_2(m_{j_0})}$.
 \item[(ii)] If the functional $h$ has a tree $(h_a)_{a\in\mc{A}}$
 with $w(h_a)\neq m_{j_0}$ for each $a\in\mc{A}$ then\\
 $\#\{k:\;
 |h(e_k)|>\frac{1}{m_{j_0}^2}\}<(5n_{j_0-1})^{2\log_2(m_{j_0})}$.
  \end{enumerate}
 \end{claim}
 \begin{proof}[\bf Proof of the claim.]
 We shall prove only part (i) of the claim, as the proof of (ii)
 is similar. Let $(h_a)_{a\in\mc{A}}$ be a tree of $h$ and let $n$
 be its height (i.e. the length of its maximal branch). We may
 assume that $|h(e_k)|>\frac{1}{m_{j_0}}$ for all $k\in\supp h$.
 Let $h=h_0,h_1,\ldots,h_n$ be a maximal branch (then $h_n\in
 C_{j_0}$) and let $k\in \supp h_n$. Then
 $\frac{1}{m_{j_0}}<|h(e_k)|=\prod\limits_{l=0}^{n-1}\frac{1}{w(h_l)}
 \le\frac{1}{2^n}$, hence $n\le \log_2(m_{j_0})-1$.

 On the other hand, since $|h(e_k)|>\frac{1}{m_{j_0}}$ for all
 $k\in\supp h$, each $h_a$ with $a$ non maximal is a result of an
 $(\mc{A}_{5n_j},\frac{1}{m_j})$ operation for
 $j\le j_0-1$.
 An inductive argument yields that for $i\le n$ the cardinality of
 the set $\{h_a:\; |a|=i\}$ is less or equal to $(5n_{j_0-1})^i$.
 The facts that $n\le\log_2(m_{j_0})-1$ and that each element of
 $g\in C_{j_0}$ has $\#(\supp(g))\le n_{j_0-1}$ yield that
 $\#(\supp(h)) \le n_{j_0-1}(5n_{j_0-1})^{\log_2(m_{j_0})-1}<
 (5n_{j_0-1})^{\log_2(m_{j_0})}$.

 The proof of the claim is complete.
 \end{proof}

 We  pass to the proof of the lemma.
 The case $w(f)=m_i$, $i\ge j_0$ is straightforward. Let $f\in
 D_{j_0}'$ with $w(f)=m_i$, $i<j$. Then
 $f=\frac{1}{m_i}\sum\limits_{t=1}^{d}f_t$ where $f_1<\cdots<f_d$
 belong to $D_{j_0}'$ and $d\le n_i$.

 For $t=1,\ldots,d$ we set $H_t=\{k:\;
 |f_t(e_k)|>\frac{1}{m_{j_0}}\}$. Part (i) of the claim yields that
 $\#(H_t) <(5n_{j_0-1})^{\log_2(m_{j_0})}$. Thus, setting
 $H=\bigcup\limits_{t=1}^dH_t$, we get that
 $\#(H)<d(5n_{j_0-1})^{\log_2(m_{j_0})}\le(5n_{j_0-1})^{\log_2(m_{j_0})+1}$.
 Therefore
 \begin{eqnarray*}
 |f(\frac{1}{n_{j_0}}\sum\limits_{l=1}^{n_{j_0}}e_{k_{l}})|& \le&
 \frac{1}{m_i}\Big(\big|(\sum\limits_{t=1}^df_t)_{|H}
 (\frac{1}{n_{j_0}}\sum\limits_{l=1}^{n_{j_0}}e_{k_{l}})\big|\Big)\\
 & &+   \frac{1}{m_i}\Big(\big|(\sum\limits_{t=1}^df_t)_{|(\N\setminus H)}
  (\frac{1}{n_{j_0}}\sum\limits_{l=1}^{n_{j_0}}e_{k_{l}})\big|\Big)\\
   & \le &
   \frac{1}{m_i}\#(H)\frac{1}{n_{j_0}}+\frac{1}{m_i}\frac{1}{m_{j_0}}<\frac{2}{m_im_{j_0}}.
 \end{eqnarray*}

  The second part is proved similarly  by using  part (ii) of the
  claim.
 \end{proof}

 \begin{proposition}{\bf (The basic inequality)} \label{prop17}
 Let $(x_k)_{k\in \mathbb{N}}$ be a $(C,\e)$ R.I.S. in $X_G$ and
 $j_0>1$ such that for every $g\in G$
 the set $\{k:\; |g(x_k)|>\e\}$ has cardinality at most
 $n_{j_0-1}$. Let
 $({\lambda}_k)_{k\in\mathbb{N}}\in c_{00}$ be a sequence of
 scalars.
 Then for every $f\in D$ of type I
 we can find $g_1$, such that either $g_1=h_1$ or
 $g_1=e_t^*+h_1$ with $t\not\in \supp h_1$ where
 $h_1\in {\conv}_{\mathbb{Q}}\{h\in D_{j_0}':\;w(h)=w(f)\}$
 and $g_2\in c_{00}(\N)$ with $\|g_2\|_{\infty}\le \varepsilon$
 with $g_1,g_2$ having nonnegative coordinates and
 such that
 \begin{equation}
 |f(\sum \lambda_kx_k)|\le C(g_1+g_2)(\sum |\lambda_k|e_k).
 \end{equation}

 If we additionally assume that for  every
 $h\in D$ with $w(h)=m_{j_0}$
 and every interval  $E$ of the natural numbers we have that
 \begin{equation} \label{eq12}
   |h(\sum\limits_{k\in E} \lambda_kx_k)|
   \le C(\max\limits_{k\in E}|\lambda_k|
     +\varepsilon \sum\limits_{k\in E}|\lambda_k|)
 \end{equation}
 then, if $w(f)\neq m_{j_0}$,
   $h_1$ may be selected satisfying additionally
 the following property:
 $h_1=\sum\limits r_l\tilde{h_l}$ with $r_l\in\mathbb{Q}^+$, $\sum r_l=1$ and
 for each $l$ the functional $\tilde{h_l}$ belongs to $D_{j_0}'$
 with $w(\tilde{h_l})=w(f)$
  and admits  a tree
 $T_{\tilde{h_l}}=(f^l_a)_{a\in \mc{C}_l}$ with $w(f^l_a)\neq m_{j_0}$ for
 all $a\in \mc{C_l}$.
 \end{proposition}
  \begin{proof}[\bf Proof.]
  The proof in the general case (where \eqref{eq12} is not
  assumed) and in the special case (where we assume \eqref{eq12})
  is actually the same. We shall give the proof only in the
  special case.
 The proof in the
 general case arises by omitting any reference to distinguishing cases
 whether a functional has weight $m_{j_0}$ or not and treating
  the functionals with $w(f)=m_{j_0}$ as for any other $j$.

  We fix a tree $T_f=(f_a)_{a\in \mc{A}}$ of $f$.
  Before passing to the proof we adopt some useful notation and
 state two lemmas. Their proofs can be found in \cite{AT1} (Lemmas
 4.4 and 4.5).

 \begin{definition}\label{def8}
  For each  $k\in\mathbb{N}$
 we define the set $\mc{A}_k$ as follows:
 \begin{eqnarray*}
 \mc{A}_k & = & \Big\{a\in \mc{A}\,\mbox{ such that }
 f_a \mbox{ is not of type } II \mbox{ and } \\
  & & (i)\;\;  \ran f_{a}\cap\ran x_k\neq
 \emptyset\\
 & & (ii)\;\;\forall\; \gamma<a\mbox{ if }f_{\gamma } \mbox{ is of type }
 I\mbox{ then }w(f_{\gamma})\neq m_{j_0}\\
  & & (iii)\;\; \forall\; \beta\le a\;\;
 \mbox{ if }\beta\in S_\gamma \mbox{ and } f_\gamma  \mbox{ is of type }I \\
 & & \mbox{ then } \ran f_{\beta}\cap\ran x_k=\ran f_\gamma \cap\ran x_k \\
 & & (iv)\;\; \mbox{if } w(f_a)\neq m_{j_0} \mbox{ then for all }\beta\in S_a \\
 & &\;\;
  \ran f_{\beta}\cap\ran x_k \subsetneqq \ran f_a\cap\ran x_k \Big\}
 \end{eqnarray*}
 \end{definition}

 The next lemma describes the properties of the set $\mc{A}_k$.
 \begin{lemma}\label{lem12}
  For every $k\in\mathbb{N}$ we have the following:
 \begin{enumerate}
 \renewcommand{\labelenumi}{(\roman{enumi})}
 \item[(i)] If $a\in\mc{A}$ and $f_a$ is of type $II$ then
 $a\not\in\mc{A}_k$.\\
 (Hence $\mc{A}_k\subset \{a\in\mc{A}:\;f_a\mbox{ is of type }I\mbox{ or }
 f_a\in G\}$.)
 \item[(ii)] If $a\in\mc{A}_k$, then for every  $\beta<a$
 if $f_\beta $ is of type  $I$ then $w(f_\beta )\neq m_{j_0}$.
 \item[(iii)] If $\mc{A}_k$ is not a singleton then its members are incomparable
 members of the tree $\mc{A}$. Moreover
 if $a_1,a_2$  are two different elements of $\mc{A}_k$
 and $\beta$ is the (necessarily uniquely determined) maximal element of
 $\mc{A}$  satisfying $\beta<a_1$ and $\beta<a_2$
 then $f_\beta $ is of type $II$.
 \item[(iv)] If $a\in\mc{A}$ is such that
 $\supp f_a \cap\ran x_k \neq \emptyset$
 and $\gamma\not\in \mc{A}_k$ for all $\gamma<a$ then there
 exists $\beta\in \mc{A}_k$  with $a\le\beta$. In particular if
 $\supp f \cap\ran x_k \neq \emptyset$ then
 $\mc{A}_k\neq\emptyset$.
 \end{enumerate}
 \end{lemma}

 \begin{definition}\label{def9}
 For every $a\in \mc{A}$ we define
 $D_a=\bigcup_{\beta\ge a}\{k:\;\beta\in \mc{A}_k\}$.
 \end{definition}

 \begin{lemma}\label{lem13}
 According to the notation  above we have the following:
 \begin{enumerate}
 \item[(i)] If $\supp f\cap\ran x_k\neq \emptyset$ then
 $k\in D_0$   (recall that 0
 denotes the unique root of $\mc{A}$ and $f=f_0$).
 Hence $f(\sum\lambda_k x_k)=f(\sum_{k\in D_0}\lambda_k
 x_k)$.
 \item[(ii)] If $f_a$ is of type $I$ with $w(f_a)=m_{j_0}$ then
 $D_a$ is an interval of $\mathbb{N}$.
 \item[(iii)] If $f_a$ is of type $I$ with $w(f_a)\neq m_{j_0}$ then
 \[\left\{\{k\}:\;k\in D_a\setminus \bigcup\limits_{\beta\in S_a}D_\beta \right\}
 \cup \left\{D_\beta :\; \beta\in S_a\right\}  \] is a family of successive
 subsets of $\mathbb{N}$. Moreover
 for every  $k\in D_a\setminus \bigcup_{\beta\in
 S_a}D_\beta$ (i.e. for $k$ such that $a\in\mc{A}_k$) such that
 $\supp f_a\cap\ran x_k\neq \emptyset$
 there exists a $\beta\in S_a$ such that either
 $\min\supp x_k\le \max\supp f_\beta <\max\supp x_k$ or
 $\min\supp x_k< \min\supp f_\beta \le\max\supp x_k$.
 \item[(iv)] If $f_a$ is of type $II$, $\beta\in S_a$ and
 $k\in D_a \setminus D_{\beta}$ then $\supp f_{\beta}\cap \ran x_k=\emptyset$
 and hence $f_{\beta}(x_k)=0.$
 \end{enumerate}
 \end{lemma}

 Recall that
 we have fixed a tree $(f_a)_{a\in\mc{A}}$ for the given $f$.
   We construct two families
 $(g_a^1)_{a\in\mc{A}}$ and $(g_a^2)_{a\in\mc{A}}$ such that the
 following conditions are fulfilled.

 \begin{enumerate}
 \renewcommand{\labelenumi}{(\roman{enumi})}
 \item[(i)] For every $a\in\mc{A}$ such that $f_a$ is not of type  $II$,
 $g^1_a=h_a$ or $g^1_a=e_{k_a}^*+h_a$ with
 $t_a\not\in \supp h_a$, where $h_a\in \conv_{\mathbb{Q}}(D_{j_0}')$ and
 $g_a^2\in c_{00}(\N)$ with ${\|g_a^2\|}_{\infty}\le \varepsilon$.
 \item[(ii)] For every $a\in\mc{A}$, $\supp g^1_a\subset D_a$ and
 $\supp g^2_a\subset D_a$ and the functionals
 $g_a^1$, $g_a^2$ have nonnegative  coordinates.
 \item[(iii)] For $a\in\mc{A}$ with $f_a\in G$ and $D_a\neq \emptyset$ we have
 that $g^1_a\in C_{j_0}$.
 \item[(iv)] For $f_a$ of type  $II$ with
  $f=\sum_{\beta\in S_a}r_\beta f_\beta $
 (where $r_\beta \in\mathbb{Q}^+$ for every $\beta\in S_a$ and
 $\sum_{\beta\in S_a}r_\beta =1$) we have
  $g^1_a=\sum_{\beta\in S_a}r_\beta g^1_\beta $
 and $g^2_a=\sum_{\beta\in S_a}r_\beta g^2_\beta $.
 \item[(v)] For $f_a$ of type $I$ with $w(f)=m_{j_0}$ we have
 $g^1_a=e_{k_a}^*$ where $k_a\in D_a$ is such that
 $|\lambda_{k_a}|=\max_{k\in D_a}|\lambda_k|$ and
 $g^2_a=\sum_{k\in D_a}\varepsilon e^*_{k}$.
 \item[(vi)] For $f_a$ of type $I$ with $w(f)=m_j$ for  $j\neq j_0$
 we have $g^1_a=h_a$ or $g^1_a=e^*_{k_a}+h_a$
 with $h_a\in \conv_{\mathbb{Q}}\{h\in D_{j_0}':\;w(h)=m_j\}$
 and $k_a\not\in\supp h_a$.
 \item[(vii)] For every $a\in \mc{A}$ the following inequality holds:
 \[ |f_a(\sum\limits_{k\in D_a}\lambda_kx_k)|\le C(g^1_a+g^2_a)
 (\sum\limits_{k\in D_a}|\lambda_k|e_{k}).   \]
 \end{enumerate}
 When the construction of $(g_a^1)_{a\in\mc{A}}$ and $(g_a^2)_{a\in\mc{A}}$
 has been accomplished, we set
 $g_1=g^1_0$ and $g_2=g^2_0$ (where $0$ is the root of $\mc{A}$ and
 $f=f_0$) and
 we observe that these are the desired functionals. To show that
 such  $(g_a^1)_{a\in\mc{A}}$ and $(g_a^2)_{a\in\mc{A}}$ exist we use finite
 induction starting with $a\in\mc{A}$ which are maximal
 and in the  general inductive step we assume that $g_\beta^1$,
 $g_\beta^2$ have been defined for all $\beta>a$ satisfying the
 inductive assumptions and
 we define $g_a^1$ and $g_a^2$.

 {\bf   $\boldsymbol 1\stackrel{st}{=}$ inductive step}\\
 Let $a\in\mc{A}$ which is maximal. Then $f_a\in G$.
 If $D_a=\emptyset$ we define $g^1_a=0$ and $g^2_a=0$.
 If $D_a\neq\emptyset$ we set
 \[ E_a=\{k\in D_a:\;|f_a(x_k)|>\varepsilon\} \mbox{\; and \;}
 F_a=D_a\setminus E_a.\]
  From our assumption we have that $\#(E_a)\le n_{j_0-1}$  and we define
 \[ g_a^1=\sum\limits_{k\in E_a}e_k^*  \mbox{\; and \;}
 g_a^2=\sum\limits_{k\in F_a}\varepsilon e_k^*.   \]
 We observe that $g_a^1\in C_{j_0}$ and $\|g^2_a\|_{\infty}\le\e$.
 Inequality (vii) is easily checked (see   Proposition 4.3 of \cite{AT1}).

 {\bf   General inductive step}\\
 Let $a\in\mc{A}$ and suppose that
 $g_{\gamma}^1$ and $g_{\gamma}^2$  have been defined for every
 $\gamma>a$ satisfying the inductive assumptions.
 If $D_a=\emptyset$ we set $g^1_a=0$
 and $g^2_a=0$.
 In the remainder of the proof we assume that
 $D_a\neq \emptyset$. We consider the following three cases:

 {\bf   $\boldsymbol 1\stackrel{st}{=}$ case }
 The functional $f_a$ is of type $II$.\\
 Let $f_a=\sum_{\beta\in S_a}r_\beta f_\beta $
 where $r_\beta \in \mathbb{Q}^+$
 are such that $\sum_{\beta\in S_a}r_\beta =1$.
 In this case,  we have that
 $D_a=\bigcup_{\beta\in S_a}D_\beta$.  We define
 \[ g^1_a=\sum\limits_{\beta\in S_a}r_\beta  g^1_\beta  \mbox{\; and \;}
 g^2_a=\sum\limits_{\beta\in S_a}r_\beta  g^2_\beta .\]
 For the proof of  inequality (vii) see  Proposition 4.3 of \cite{AT1}.

 {\bf   $\boldsymbol 2\stackrel{nd}{=}$ case }
  The functional $f_a$ is of type $I$ with $w(f)=m_{j_0}$. \\
 In this case  $D_a$ is an interval of the natural numbers (Lemma \ref{lem13}(ii)).
 Let  $k_a\in D_a$ be such that $|\lambda_{k_a}|=\max_{k\in
 D_a}|\lambda_k|$.
 We define  \[g^1_a=e_{k_a}^* \mbox{\;  and \;}
 g^2_a=\sum\limits_{k\in D_a}\varepsilon e^*_{k}. \]
 Inequality (vii) is easily established.

 {\bf   $\boldsymbol 3\stackrel{rd}{=}$ case }
 The functional $f_a$ is of type $I$
 with $w(f)=m_j$ for $j\neq j_0$.\\
 Then $f_a=\frac{1}{m_j}\sum_{\beta\in S_a}f_\beta $ and the family
 $\{f_\beta :\; \beta\in S_a\}$ is a family of successive
 functionals with $\#(S_a)\le n_j$.
 We set
 \begin{eqnarray*}
  E_a & = & \{k:\;a\in \mc{A}_k\mbox{ and }  \supp f_a\cap\ran
 x_k\neq\emptyset\}\\
     & &(=\{k\in D_a\setminus \bigcup_{\beta\in
     S_a}D_\beta:\;\supp f_a\cap\ran  x_k\neq\emptyset\}).
 \end{eqnarray*}
  We  consider the following partition of $E_a$.
 \[E_a^2=\{k\in E_a:\; m_{j_{k+1}}\le m_j\}\mbox{\; and \;}
 E_a^1=E_a\setminus E_a^2.\]

 We define  \[ g^2_a=\sum\limits_{k\in E_a^2}\varepsilon e_{k}^*
                             +\sum\limits_{\beta\in S_a}g^2_\beta. \]
 Observe that ${\|g^2_a\|}_{\infty}\le \varepsilon$.
 Let $E_a^1=\{k_1<k_2<\cdots <k_l\}$.
 From the  definition of $E^1_a$ we get that
 $m_j<m_{j_{k_2}}<\cdots <m_{j_{k_l}}$.
  We set
 \[k_a=k_1 \mbox{\; and \;} g^1_a=e_{k_a}^*+h_a \mbox{\; where \;}
 h_a=\frac{1}{m_j} ( \sum\limits_{i=2}^le_{k_i}^*
                 +  \sum\limits_{\beta\in S_a}g^1_\beta   )  \]
 (The term $e_{k_a}^*$ does not appear if $E_a^1=\emptyset$.)\\
  For the verification of inequality (vii) see  Proposition 4.3 of \cite{AT1}.

 It remains to show that $h_a\in\conv_{\Q}\{h\in
 D_{j_0}':\;w(h)=m_j\}$
 By the second part of Lemma \ref{lem13}(iii), for
 every $k\in E_a$ there exists an element of the set
 $N=\{\min\supp f_\beta,\max\supp f_\beta:\; \beta\in S_a\}$
 belonging to $\ran x_k$.
 Hence $\#(E^1_a)\le\#(E_a)\le 2n_j$.

 We next show that $h_a\in \conv_{\mathbb{Q}}\{g\in D_{j_0}':\;
 w(g)=m_j\}$. We first examine the case that
 for every $\beta\in S_a$ the functional $f_\beta $ is not of type $II$.
  Then for every $\beta\in S_a$ one of the following holds:
 \begin{enumerate}
 \renewcommand{\labelenumi}{(\roman{enumi})}
 \item $f_\beta \in G$. In this case $g^1_\beta \in C_{j_0}$ (by the first
 inductive step).
 \item $f_\beta $ is of type $I$ with $w(f_\beta )=m_{j_0}$. In this case
 $g^1_\beta =e_{k_\beta }^*\in D'_{j_0}$.
 \item $f_\beta $ is of type $I$ with $w(f_\beta )=m_j$ for $j\neq j_0$. In
 this case $g^1_\beta =e_{k_\beta }^*+h_\beta$
  (or $g^1_\beta =h_\beta$)
 where $h_\beta \in  \conv_{\mathbb{Q}}(D'_{j_0})$
 and $k_\beta \not\in \supp h_\beta$ .
 We set $E_\beta^1=\{n\in\mathbb{N}:\;n<k_\beta \}$,
 $E_\beta^2=\{n\in\mathbb{N}:\;n>k_\beta \}$
  and $h_\beta^1=E_\beta^1h_\beta$,
 $h_\beta^2=E_\beta^2h_\beta $.
 The functionals $h_\beta^1$, $e_{t_\beta }^*$,
  $h_\beta^2$ are successive
 and belong to  $D_{j_0}=\conv_{\mathbb{Q}}(D'_{j_0})$.
 \end{enumerate}
 We set
 \begin{eqnarray*}
 T_a^1 & = & \{\beta\in S_a:\; f_\beta \in G\}\\
 T_a^2 & = & \{\beta\in S_a:\; f_\beta \mbox{ of type }
 I \mbox{ and } w(f_\beta )= m_{j_0}\} \\
 T_a^3 & = & \{\beta\in S_a:\; f_\beta \mbox{ of type } I
 \mbox{ and } w(f_\beta )\neq m_{j_0}\}.
 \end{eqnarray*}

 The family of successive (see Lemma \ref{lem13}(iii))
  functionals of $D_{j_0}$,
 \begin{eqnarray*}
 \{e_{k_i}^*:\;i=2,\ldots ,l\}\cup
 \{g^1_\beta :\;\beta\in T_a^1\}\cup
 \{g^1_\beta :\;\beta\in T_a^2\}\cup    \\
 \cup \{h^1_\beta :\;\beta\in T_a^3\}\cup
 \{e_{k_\beta }^*:\;\beta\in T_a^3\}\cup
 \{h^2_\beta :\;\beta\in T_a^3\}
 \end{eqnarray*}
 has cardinality $\le 5n_j$, thus we get that $h_a\in D_{j_0}$
 with $w(h_a)=m_j$. Therefore from Remark \ref{rem5} we get that
 \[h_a\in\conv_{\mathbb{Q}}\{h\in D_{j_0}':\;w(h)=m_j\}.\]

 For the case that for some $\beta\in S_a$ the functional $f_\beta $
 is of type $II$ see \cite{AT1} Proposition 4.3.
 \end{proof}

 \begin{proof}[\bf Proof of Proposition \ref{prop13}.]
 The proof is an application of the basic inequality
 (Proposition \ref{prop17}) and Lemma \ref{lem14}. Indeed, let
 $f\in D$ with $w(f)=m_i$. Proposition \ref{prop17} yields the
 existence of a functional $h_1$ with
 $h_1\in {\conv}_{\mathbb{Q}}\{h\in D_{j_0}':\;w(h)=m_i\}$,
  a $t\in\N$
 and a $h_2\in c_{00}(\N)$ with $\|h_2\|_{\infty}\le\e$,
   such that
 \[
 |f(\frac{1}{n_{j_0}}\sum\limits_{k=1}^{n_{j_0}}x_{k})| \leq
 C(e_t^*+h_1+h_2)\big(\frac{1}{n_{j_0}}
 \sum\limits_{k=1}^{n_{j_0}}e_k\big).\]

  If $i\ge j_0$ we get that
 $|f(\frac{1}{n_{j_0}}\sum_{k=1}^{n_{j_0}}x_{k})| \leq
 C(\frac{1}{n_{j_0}}+\frac{1}{m_i}+\e)<\frac{C}{n_{j_0}}+\frac{C}{m_i}
 +C\e$.
 If $i< j_0$, using Lemma \ref{lem14} we get that
 $|f(\frac{1}{n_{j_0}}\sum_{k=1}^{n_{j_0}}x_{k})| \leq
 C(\frac{1}{n_{j_0}}+\frac{2}{m_i\cdot
 m_{j_0}}+\e)<\frac{3C}{m_i\cdot m_{j_0}}$.

 In order to prove 2) let $(b_k)_{k=1}^{n_{j_0}}$ be scalars
 with $|b_k|\le 1$  such that \eqref{eq16} is satisfied.
 Then condition \eqref{eq12} of the basic inequality is satisfied
 for the linear combination
 $\frac{1}{n_{j_0}}\sum_{k=1}^{n_{j_0}}b_kx_k$
 and thus  for every $f\in D$ with $w(f)=m_i$, $i\neq j_0$,
 there exist a $t\in\N$ and $h_1,h_2\in c_{00}(\N)$
 with $h_1,h_2$ having nonnegative coordinates and
 $\|h_2\|_{\infty}\le\e$
 such that
  \begin{eqnarray*}\label{eq17}
 |f(\frac{1}{n_{j_0}}\sum\limits_{k=1}^{n_{j_0}}b_kx_{k})| &\leq&
 C(e_t^*+h_1+h_2)\big(\frac{1}{n_{j_0}}\sum\limits_{k=1}^{n_{j_0}}|b_k|e_k\big)\\
 & \le&
 C(e_t^*+h_1+h_2)\big(\frac{1}{n_{j_0}}\sum\limits_{k=1}^{n_{j_0}}e_k\big)
 \end{eqnarray*}
  with  $h_1$ being  a rational convex combination
 $h_1=\sum r_l\tilde{h_l}$ and
 for each $l$ the functional $\tilde{h_l}$ belongs to $D_{j_0}'$
 with $w(\tilde{h_l})=m_i$ and has a tree
 $T_{\tilde{h_l}}=(f^l_a)_{a\in \mc{C}_l}$ with $w(f^l_a)\neq m_{j_0}$ for
 all $a\in \mc{C_l}$.
  Using the second part of Lemma \ref{lem14} we deduce that
 \[|f(\frac{1}{n_{j_0}}\sum\limits_{k=1}^{n_{j_0}}b_kx_{k})| \leq
 C(\frac{1}{n_{j_0}}+\frac{1}{m_{j_0}^2}+\e)<\frac{4C}{m_{j_0}^2}.\]
 For $f\in D$ with $w(f)=m_{j_0}$ it follows from condition \eqref{eq16}
 that\\
 $|f(\frac{1}{n_{j_0}}\sum_{k=1}^{n_{j_0}}b_kx_{k})| \le
 \frac{C}{n_{j_0}}(1+\frac{2}{m_{j_0}^2}
 n_{j_0})<\frac{4C}{m_{j_0}^2}$.
   \end{proof}

 \section{The James tree spaces $JT_{\mc{F}_{2,s}}$, $JT_{\mc{F}_{2}}$
 and $JT_{\mc{F}_s}$}

 In this part we continue the study of the James Tree spaces
 initialized in Section 3. We give a slightly different definition
 of JTG sequences and then we define
 the space   $JT_{\mc{F}_2}$ exactly as in section 3. We also
 define the spaces
 $JT_{\mc{F}_{s}}$, $JT_{\mc{F}_{2,s}}$.
 We prove that $JT_{\mc{F}_2}$ is $\ell_2$ saturated while
 $JT_{\mc{F}_{s}}$, $JT_{\mc{F}_{2,s}}$ are $c_0$ saturated. We
 also give an example of $JT_{\mc{F}_2}$, defined for a precise
 family $(F_j)_j$ such that the basis \seq{e}{n} of the space is
 normalized weakly null and for every subsequence $(e_n)_{n\in
 M}$, $M\in [\N]$ the space $X_M=\overline{\spann}\{e_n:\;n\in
 M\}$ has nonseparable dual. As we have mentioned before the study
 of the James Tree spaces does not require techniques related to HI
 constructions.

 \begin{definition}\label{bdef21}
  {\bf (JTG families)}
 A sequence  $(F_j)_{j=0}^{\infty}$  of subsets of
 $c_{00}(\N)$ is said to be a \emph{James Tree Generating family}
  (JTG
 family)  provided that it satisfies the following conditions:
 \begin{enumerate}
 \item[(A)] $F_0=\{\pm e_n^*:\;n\in\N\}$
 and each $F_j$ is nonempty, countable, symmetric, compact in the
 topology of pointwise convergence and closed under restrictions to
 intervals of $\N$.
 \item[(B)] Setting ${\tau}_j=\sup\limits\{\|f\|_{\infty}:\; f\in
 F_j\}$, the sequence \seq{\tau}{j}
 is strictly
 decreasing and $\sum\limits_{j=1}^{\infty}{\tau}_j\le 1$.
 \item[(C)] For every block sequence \seq{x}{k} of $c_{00}(\N)$,
 every $j=0,1,2,\ldots$ and every $\delta>0$ there exists a vector
 $x\in\spann\{x_k:\; k\in\N\}$ such that
 \[  \delta\cdot \sup\{f(x):\; f\in
 \bigcup\limits_{i=0}^{\infty}F_i\}>
 \sup\{f(x):\; f\in F_j\}.  \]
 \end{enumerate}

 We set
 $F=\bigcup\limits_{j=0}^{\infty}F_j$. The set $F$ defines a norm
 $\|\;\|_F$
 on $c_{00}(\N)$ by the rule
 \[  \|x\|_F=\sup\{f(x):\;f\in F\}.\]
 The space $Y_F$ is the completion of the space
 $(c_{00}(\N),\|\;\|_F)$.
  \end{definition}

 \begin{examples}\label{bexa1}
 We provide some examples of JTG families.
 \begin{enumerate}
 \item[(i)] The first example is what we call the
  Maurey-Rosenthal JTG
 family, related to the
 first construction
  of a normalized
 weakly null sequence with no unconditional subsequence
  (\cite{MR}). In particular the norming set for their example
  is the set $F=\bigcup\limits_{j=0}^{\infty}F_j$
  together with the $\sigma_{\mc{F}}$
  special functionals
  resulting from the family $F$. We proceed defining the sets
   $(F_j)_{j=0}^{\infty}$.

  Let \seq{k}{j} be a strictly
  increasing sequence of integers such
 that
  \[\sum\limits_{j=1}^{\infty}\sum\limits_{n\neq j}\min\big\{\frac{\sqrt{k_n}}{\sqrt{k_j}},
 \frac{\sqrt{k_j}}{\sqrt{k_n}}\big\} \le1.\]

 We set $F_0=\{\pm e_n^*:\;n\in\N\}$ while for $j=1,2,\ldots$ we set
 \[ F_j=\big\{ \frac{1}{\sqrt{k_j}}\big(\sum\limits_{i\in
 F}\pm e_i^*\big):\; \emptyset\neq F\subset \N,\;\#(F)\le
 k_j\big\}\cup\big\{0\big\}.\]
 The above conditions (1) and (2) for the sequence
 \seq{k}{j} easily yield that $(F_j)_{j=0}^{\infty}$ is a JTG
 family.
 \item[(ii)] The second example is the family introduced in
 Section \ref{secr}. For completeness we recall its definition.
 Let $\seq{m}{j}$ and $\seq{n}{j}$ defined as follows:
 \begin{itemize}
 \item  $m_1=2$  and $m_{j+1}=m_j^5$.
 \item $n_1=4$,  and $n_{j+1}=(5n_j)^{s_j}$  where  $s_j=\log_2 m_{j+1}^3$.
 \end{itemize}
 We set $F_0=\{\pm e_n^*:\;n\in\N\}$ and for
  $j=1,2,\ldots$ we set
 \[ F_j=\big\{\frac{1}{m_{2j-1}^2}\sum\limits_{i\in I}\pm e^*_i:\;
  \#(I)\le \frac{n_{2j-1}}{2} \big\}\cup\big\{0\big\}. \]

  We shall show that
    the sequence $(F_j)_{j=0}^{\infty}$ is a JTG family.
 Conditions (A), (B) of Definition \ref{bdef21} are obviously satisfied.
  Suppose that
 condition (C)  fails. Then for some $j\in\N$,
 there exists a block sequence \seq{x}{k} in \co with
 $\|x_k\|_F=1$ and a $\delta>0$ such that
 $\delta \|\sum a_kx_k\|_F \le
 \|\sum a_kx_k\|_{F_j}$ for every
 sequence of scalars $\seq{a}{k}\in \co$. We observe that
 $\|x_k\|_{\infty}\ge \frac{2\delta m_{2j-1}^2}{n_{2j-1}}$
 for all $k$. Indeed,
  if $\|x_k\|_{\infty}< \frac{2\delta
 m_{2j-1}^2}{n_{2j-1}}$ then for every $f\in F_j$,
 $f=\frac{1}{m_{2j-1}^2}\sum\limits_{i\in I}\pm e^*_i$, with
 $\#(I)\le \frac{n_{2j-1}}{2}$ we would have that
 $|f(x_k)|\le \frac{1}{m_{2j-1}^2}
 \sum\limits_{i\in I}|e_i^*(x_k)|<\frac{1}{m_{2j-1}^2}\frac{n_{2j-1}}{2}
 \frac{2\delta  m_{2j-1}^2}{n_{2j-1}}=\delta$ which yields that
 $\|x_k\|_{F_j}<\delta$, a contradiction.
 Hence for each $k$ we may select a $t_k\in\supp x_k$ such that
 $|e_{t_k}^*(x_k)|\ge \frac{2\delta m_{2j-1}^2}{n_{2j-1}}$.
 Since the sequence  $(\frac{n_{2i-1}}{m_{2i-1}^2})_{i\in\N}$ increases to
 infinity we may choose a $j'\in\N$ such that
 $\delta^2\frac{n_{2j'-1}}{m_{2j'-1}^2}>(\frac{n_{2j-1}}{m_{2j-1}^2})^2$.
 We consider the vector $y=\sum\limits_{k=1}^{n_{2j'-1}/2}x_k$.
 We have that $\delta\|y\|_F\ge \delta \frac{1}{m_{2j'-1}^2}
 \sum\limits_{k=1}^{n_{2j'-1}/2}|e_{t_k}^*(x
 _k)|\ge
 \delta
  \frac{1}{m_{2j'-1}^2}\frac{n_{2j'-1}}{2}\frac{2\delta m_{2j-1}^2}{n_{2j-1}}
 >\frac{1}{m_{2j-1}^2}{n_{2j-1}}$.
 On the other hand $\|y\|_{F_j}\le
 \frac{1}{m_{2j-1}^2}\frac{n_{2j-1}}{2}$,
 a contradiction.

% We shall use this JTG family in our constructions in section
% \ref{bsec3} in order to construct the ground set for the spaces
% $\eqs_{\mc{F}_{2,s}}$ and $\eqs_{\mc{F}_{2,s}}^{us}$.

 \item[(iii)] Another example of a JTG family has been given in
 Definition \ref{def18} and we have used it to define the
 ground set for the space $\eqs_{\mc{F}_s}$.
 \end{enumerate}
 \end{examples}

 \begin{remarks}\label{brem12}
 \begin{enumerate}
 \item[(i)] The standard basis \seq{e}{n} of $c_{00}(\N)$ is a
 normalized bimonotone Schauder basis of the space $Y_F$.
 \item[(ii)] The set $F$ is compact in the topology of pointwise
 convergence. Indeed, let $\seq{f}{n}$ be a sequence in $F$. There
 are two cases. Either there exists $j_0\in\N$ such that
  the set $F_{j_0}$
 contains a subsequence of $\seq{f}{n}$ in which case
  the compactness of
 $F_{j_0}$ yields the existence of a further subsequence
 converging pointwise to some $f\in F_{j_0}$, otherwise
 if no such $j_0$
 exists, then we may find a subsequence $(f_{k_n})_{n\in\N}$
 and a strictly increasing sequence $\seq{i}{n}$
 of integers with $f_{k_n}\in F_{i_n}$ and thus
 $\|f_{k_n}\|_{\infty}\le {\tau}_{i_n}$ for all $n$. Since condition
 (B) of Definition \ref{bdef21} yields that ${\tau}_n\to 0$ we get that
 $f_{k_n}\stackrel{p}{\to}0\in F$.
 \item[(iii)] The fact that $F$ is countable and compact yields
 that the space $(C(F),\|\;\|_{\infty})$ is $c_0$ saturated \cite{BP}.
 It
 follows that the space $Y_F$ is also $c_0$ saturated, since $Y_F$ is
 isometric to a subspace of $(C(F),\|\;\|_{\infty})$.
 \item[(iv)] For each $j$ we consider the
 seminorm $\|\;\|_{F_j}:\co\to\R$
 defined by $\|x\|_{F_j}=\sup\{|f(x)|:\; f\in F_j\}$. In general
 $\|\;\|_{F_j}$ is not  a norm.
 Defining $Y_{F_j}$ to be the completion of the space
 $(\co,\|\;\|_{F_j})$, condition (C) of Definition \ref{bdef21}
 is equivalent to saying that the
 identity operator $I:Y_F\to Y_{F_j}$ is strictly singular.

 Furthermore, observe
 that setting $H_j=\cup_{i=0}^jF_i$ the identity
 operator $I:Y_F\to Y_{H_j}$ ($Y_{H_j}$ is similarly defined)
 is also strictly singular. Indeed,
 let \seq{x}{k} be a block sequence of $\co$ and let $\delta>0$.
  We choose
 a block sequence $(x^0_k)_{k\in\N}$ of \seq{x}{k}
 with $\|x^0_k\|_F=1$ and
 $\sum\limits_{k=1}^{\infty}\|x^0_k\|_{F_1}<\delta.$
 Then for every $x\in\spann\{x^0_k:\; k\in\N\}$ we have that
 $\delta \|x\|_F\ge\|x\|_{F_0}$. We then select a block sequence
 $(x^1_k)_{k\in\N}$ of $(x^0_k)_{k\in\N}$ such that
 $\delta \|x\|_F\ge\|x\|_{F_1}$ for every $x\in\spann\{x^1_k:\;
 k\in\N\}$. Following this procedure, after $j+1$ steps we may
 select a block sequence $(x_k^j)_{k\in\N}$ of \seq{x}{k}
 such that $\delta \|x\|_F\ge\|x\|_{F_i}$ for $i=1,\ldots,j$ and
 thus $\delta \|x\|_F\ge\|x\|_{H_j}$ for every $x\in\spann\{x^j_k:\;
 k\in\N\}$.
 \end{enumerate}
 \end{remarks}

 Next using the $\sigma_{\mc{F}}$ coding
  defined in Definition
 \ref{def4} we introduce the $\sigma_{\mc{F}}$ special sequences
 and functionals in the same manner as in Definition \ref{def222}.
 For a $\sigma_{\mc{F}}$ special functional $x^*$ the index
 $\ind(x^*)$ has the analogous meaning. Finally we denote by
 $\mathscr{S}$ the set of all finitely supported
  $\sigma_{\mc{F}}$ special
 functionals.

 The next proposition is an immediate consequence of the above
 definition and describes the tree-like interference of two
 $\sigma_{\mc{F}}$ special sequences.
 \begin{proposition} \label{btrpr}
 Let $(f_i)_i$, $(h_i)_i$ be two distinct $\sigma_{\mc{F}}$ special
 sequences.
 Then $\ind(f_i)\neq\ind(h_j)$ for $i\neq j$
 while there exists $i_0$ such that $f_i=h_i$ for all $i<i_0$ and
  $\ind(f_i)\neq \ind(h_i)$
 for  $i>i_0$.
 \end{proposition}

% A tree structure $(\mc{T}_{\sigma_{\mc{F}}},\preceq)$
% is naturally induced by the family of the
% $\sigma_{\mc{F}}$ special sequences. We  set
% \[\mc{T}_{\sigma_{\mc{F}}}
% =\{(f_1,\ldots,f_d):\mbox{ finite }\sigma_{\mc{F}}\mbox{ special
% sequence}\}.\]
% The  relation $\preceq$ is defined by the rule
% \[(f_1,\ldots,f_d)\preceq (f'_1,\ldots,f'_{d'})
% \iff d\le d' \mbox{ and }f_i=f'_i \mbox{ for }i=1,\ldots,d.\]
%
% An $f\in F\setminus F_0$ appears finitely many times as a final
% member of an element of $\mc{T}_{\sigma_{\mc{F}}}$.
%  First the sequence $(f)$ of
% length $1$ is a member of $\mc{T}$. The set $I=\{n:\;f\in F_n\}$
% is finite due to condition (B) of Definition \ref{bdef21}.
%  Since $\sigma_{\mc{F}}$ is injective the
% set
% $P_f=\{(f_1,\ldots,f_d)\in\mc{W}:\;\sigma_{\mc{F}}(f_1,\ldots,f_d)\in I
% \mbox{ and } f_d<f\}$ is finite. The set of all members of
% $\mc{T}$ having $f$ as a final member is the set
% \[ \{(f_1,\ldots,f_d,f):\;(f_1,\ldots,f_d)\in P_f\}\cup\{(f)\}.\]

 \begin{definition}\label{bdef17}
  {\bf (The norming sets  $\mc{F}_{2,s}$, $\mc{F}_2$, $\mc{F}_{s}$)}
  Let $(F_j)_{j=0}^{\infty}$ be a JTG family.
 We set
 \begin{eqnarray*}
 \mc{F}_2 & = & F_0\cup \big\{\sum\limits_{k=1}^da_kx^*_k:\; a_k\in\Q,\;
 \sum\limits_{k=1}^da_k^2\le 1,\;x_k^*\in
 \mathscr{S}\cup \bigcup\limits_{i=1}^{\infty}F_i,\;k=1,\ldots,d
  \\
   & & \mbox{ with }
   (\ind(x_k^*))_{k=1}^d\mbox{
   pairwise disjoint} \big\}
 \end{eqnarray*}
 \begin{eqnarray*}
 \mc{F}_{2,s} &  & F_0\cup \big\{\sum\limits_{k=1}^da_kx^*_k:\; a_k\in\Q,\;
 \sum\limits_{k=1}^da_k^2\le 1,\;x_k^*\in
 \mathscr{S}\cup \bigcup\limits_{i=1}^{\infty}F_i,\;k=1,\ldots,d
 \\    & & \mbox{ with } (\ind(x_k^*))_{k=1}^d\mbox{
   pairwise disjoint and }\min\supp x_k^*\ge d \big\},
 \end{eqnarray*}
 and
 \begin{eqnarray*}
 \mc{F}_{s}&=&F_0\cup\{
 \sum\limits_{k=1}^d\e_kx_k^*:\;\e_1,\ldots,\e_d\in\{-1,1\},
 \;x_k^*\in
  \mathscr{S}\cup \bigcup\limits_{i=1}^{\infty}F_i,\;k=1,\ldots,d  \\
 & & \mbox{  with }(\ind(x_i^*))_{i=1}^d\mbox{ pairwise disjoint  and
 }\min\supp x_i^*\ge d,\;d\in\N\}.
 \end{eqnarray*}

 The space $JT_{\mc{F}_{2,s}}$ is defined as the
 completion of the space $(\co,\|\;\|_{\mc{F}_{2,s}})$,
  the space $JT_{\mc{F}_2}$ is defined to be the
 completion of the space $(\co,\|\;\|_{\mc{F}_{2}})$
 while $JT_{\mc{F}_s}$  the completion of
 $(\co,\|\;\|_{\mc{F}_s})$
 (where $\|x\|_{\mc{F}_*}
 =\sup\{f(x):\; f\in \mc{F}\}$ for $x\in\co$,
 for either $\mc{F}_*=\mc{F}_2$ or
 $\mc{F}_*=\mc{F}_{2,s}$ or
 $\mc{F}_*=\mc{F}_{s}$).
 For a functional $f\in \mc{F}_*\setminus F_0$ of the form
 $f=\sum\limits_{k=1}^la_kx^*_k$ the set of its indices $\ind(f)$ is
 defined to be the set $\ind(f)=\bigcup\limits_{k=1}^l\ind(x_k^*)$.
 \end{definition}

 \begin{remark}\label{brem11}
 The standard basis \seq{e}{n} of \co is a normalized bimonotone
 Schauder basis for the space $JT_{\mc{F}_*}$.

 Let's observe that the only difference between the definition of
 $\mc{F}_{2,s}$ and that of $\mc{F}_2$ is the way we connect the
 $\sigma_{\mc{F}}$ special functionals.
 In the case of $\mc{F}_2$ the $\sigma_{\mc{F}}$ special functionals are
 connected more freely than in $\mc{F}_{2,s}$
  and obviously $\mc{F}_{2,s}\subset \mc{F}_2$.
 This difference leads the spaces $JT_{\mc{F}_{2,s}}$ and $JT_{\mc{F}_2}$
 to have extremely different structures. We study the structure of
 these two spaces as well as the structure of $JT_{\mc{F}_s}$.
   Namely we have the following theorem.
  \end{remark}

 \begin{theorem}\label{bth1}
 \begin{enumerate}
 \item[(i)] The space $JT_{\mc{F}_{2,s}}$ is $c_0$ saturated.
 \item[(ii)] The space $JT_{\mc{F}_2}$  is $\ell_2$ saturated.
 \item[(iii)] The space $JT_{\mc{F}_{s}}$ is $c_0$ saturated.
 \end{enumerate}
 \end{theorem}

 Proposition \ref{btrpr} yields
 that the set of all finite $\sigma_{\mc{F}}$ special sequences
 is naturally endowed with a tree structure.
  The set of  infinite branches of this tree structure
  is identified with the set of all infinite
 $\sigma_{\mc{F}}$ special sequences. For such a branch
 $b=(f_1,f_2,\ldots)$ the functional
 $b^*=\lim\limits_d\sum\limits_{i=1}^df_i$ (where the
 limit is taken in the pointwise topology) is a cluster point of
 the sets  $\mc{F}_2$, $\mc{F}_{2,s}$ and $\mc{F}_s$
  and hence belongs to the unit balls of the dual spaces $JT_{\mc{F}_2}^*$
  $JT_{\mc{F}_{2,s}}^*$ and $JT_{\mc{F}_s}^*$.
 Let also point out that a $\sigma_{\mc{F}}$
 special functional $x^*$ is either finite or takes the
 form $Eb^*$ for some branch $b$ and some interval $E$.
 Furthermore, it is easy to check that the set $\{Ex^*:\; E\mbox{
 interval, }x^*\;\sigma_{\mc{F}}\mbox{ special functional}\}$ is closed
 in the pointwise topology.

 Our main goal in this section is to prove Theorem \ref{bth1}.
 Many of the Lemmas used
 in proving this theorem
 are common for  $JT_{\mc{F}_{2,s}}$, $JT_{\mc{F}_2}$ and
 $JT_{\mc{F}_s}$. For this reason it is convenient
 to use the symbol $\mc{F}_*$ when stating or proving a property which is
 valid for
 $\mc{F}_*=\mc{F}_{2,s}$,  $\mc{F}_*=\mc{F}_2$ and
 $\mc{F}_*=\mc{F}_s$.

 \begin{lemma} \label{blem1}
 The identity operator $I:JT_{\mc{F}_*}\to Y_F$ is strictly
 singular.
 \end{lemma}
 \begin{proof}[\bf Proof.]
 Assume  the contrary. Then there exists a block subspace $Y$ of
 $JT_{\mc{F}_*}$ such that the identity operator
 $I:(Y,\|\;\|_{\mc{F}_*})\to  (Y,\|\;\|_{F})$ is an isomorphism.
 Since $Y_F$ is $c_0$ saturated (Remark \ref{brem12} (iii)) we may
 assume that $(Y,\|\;\|_{\mc{F}_*})$ is
 is spanned by a block basis which is equivalent to the standard
 basis of $c_0$. Using
 property (C) of Definition \ref{bdef21} and Remark \ref{brem12} (iv)
 we inductively choose a normalized  block sequence \seq{x}{n} in
 $(Y,\|\;\|_{\mc{F}_*})$ and a
 strictly increasing sequence \seq{j}{n} of integers such that for
 some $\delta$ determined by the isomorphism, the following hold:
 \begin{enumerate}
 \item[(i)] $\|x_n\|_{F_{j_n}}>\delta$.
 \item[(ii)]
 $\|x_{n+1}\|_{\bigcup\limits_{k=1}^{j_n}F_k}<\delta$.
 \end{enumerate}
  From (i) and (ii) and the definition of each $\mc{F}_*$ we
  easily get that $\|x_1+\cdots+x_n\|\stackrel{n}{\to}\infty$.
 This is a contradiction since \seq{x}{n}, being a normalized
 block basis of a sequence equivalent to the standard basis of
 $c_0$, is also equivalent to the standard basis of $c_0$.
 \end{proof}

 The following lemma, although it refers exclusively to the functional
 $b^*$, its proof is crucially depended on the fact that in
 $\mc{F}_*$ we connect the special functionals under certain
 norms.
 A similar result is also obtained in \cite{AT1} (Lemma 10.6).

 \begin{lemma} \label{blem4}
 Let \seq{x}{n} be a bounded block sequence in $JT_{\mc{F}_*}$.
  Then there
 exists an $L\in [\N]$ such that for every branch $b$ the limit
 $\lim\limits_{n\in L}b^*(x_n)$ exists. In particular, if
 the sequence \seq{x}{n} is
  seminormalized (i.e. $\inf\|x_n\|_{\mc{F}_*}>0$) and
 $L=\{l_1<l_2<l_3<\cdots\}$ then the sequence
 $y_n=\frac{x_{l_{2n-1}}-x_{l_{2n}}}{\|x_{l_{2n-1}}-x_{l_{2n}}\|}$
 satisfies $\|y_n\|_{\mc{F}_*}=1$ and
 $\lim\limits_nb^*(y_n)=0$ for every branch $b$.
 \end{lemma}
 \begin{proof}[\bf Proof.] We first prove the following claim.
 \begin{claim}
 For every $\e>0$ and $M\in[\N]$ there exists $L\in [M]$ and a
 finite collection of branches $\{b_1,\ldots,b_l\}$ such that for
 every branch $b$ with $b\not\in\{b_1,\ldots,b_l\}$ we have that
 $\limsup\limits_{n\in L} |b^*(x_n)|\le \e$.
 \end{claim}
 \begin{proof}[\bf Proof of the claim.]
 Assume the contrary. Then we may inductively construct a sequence
 $M_1\supset M_2\supset M_3\cdots $ of infinite subsets of $\N$
 and a sequence $b_1,b_2,b_3,\ldots$ of pairwise different
 branches satisfying $|b_i^*(x_n)|>\e$ for all $n\in M_i$.

 We set $C=\sup\limits_n\|x_n\|_{\mc{F}_*}$ and we consider
 $k>\frac{C}{\e}$. Since the branches $b_1,b_2,\ldots,b_{k^2}$ are
 pairwise different we may choose an infinite interval
  $E$ with $\min E\ge k^2$
 such
 that the functionals $(Eb_i^*)_{i=1}^{k^2}$ have disjoint
 indices. We also consider any $n\in M_{k^2}$ with $\supp
 x_n\subset E$ and we set $\e_i=\sgn b_i^*(x_n)$ for
 $i=1,2,\ldots,k^2$. Then the functional
 $f=\sum\limits_{i=1}^{k^2}\frac{\e_i}{k}b_i^*$ belongs to
 $B_{JT_{\mc{F}_*}^*}$.
 Therefore
 \[ \|x_n\|_{\mc{F}_*} \ge f(x_n)
 =\sum\limits_{i=1}^{k^2}\frac{1}{k}|b_i^*(x_n)|\ge
 \sum\limits_{i=1}^{k^2}\frac{1}{k}\cdot\e=k\cdot\e>C,\]
 a contradiction completing the proof of the claim.
 \end{proof}

 Using the claim we inductively select a sequence $L_1\supset L_2\supset
 L_3\supset\cdots$ of infinite subsets of $\N$ and a sequence
 $B_1,B_2,B_3,\ldots$ of finite collections of branches such that for
 every branch $b\not\in B_i$ we have that $|b^*(x_n)|<\frac{1}{i}$
 for all $n\in L_i$.  We then choose a diagonal set $L_0$ of the
 nested sequence \seq{L}{i}. Then for every branch $b$ not
 belonging to $B=\bigcup\limits_{i=1}^{\infty}B_i$ we have that
 $\lim\limits_{n\in L_0}b^*(x_n)=0$. Since the set $B$ is countable,
 we may choose, using a diagonalization argument, an $L\in [L_0]$
 such that the sequence $(b^*(x_n))_{n\in L}$ converges for every
 $b\in F$. The set $L$ clearly satisfies the conclusion of the
 lemma.
 \end{proof}

  Combining Lemma \ref{blem1} and Lemma \ref{blem4} we get the
  following.

 \begin{corollary} \label{bcor5}
 Every block subspace of $JT_{\mc{F}_*}$ contains a block sequence
 \seq{y}{n} such that $\|y_n\|_{\mc{F}_*}=1$, $\|y_n\|_F\stackrel{n\to
 \infty}{\longrightarrow}0$ and $b^*(y_n)\stackrel{n\to
 \infty}{\longrightarrow}0$ for every branch $b$.
 \end{corollary}

 \begin{lemma}\label{blem8}
 Let $Y$ be a block subspace of $JT_{\mc{F}_*}$ and let $\e>0$.
 Then there exists
  a finitely supported
 vector $y\in Y$ such that $\|y\|_{\mc{F}_*}=1$
 and $|x^*(y)|<\e$ for every
 $\sigma_{\mc{F}}$ special functional $x^*$.
 \end{lemma}
 \begin{proof}[\bf Proof.]
 Assume the contrary. Then there exists a block subspace $Y$ of
 $JT_{\mc{F}_*}$ and an $\e>0$ such that
 \begin{equation} \label{beq9}
   \e\cdot\|y\|_{\mc{F}_*}
   \le \sup\{|x^*(y)|:\;x^*\mbox{ is a }\sigma_{\mc{F}} \mbox{ special
 functional}\}
 \end{equation}
 for every $y\in Y$. Let $q>\frac{8}{\e^2}$.
 From Corollary \ref{bcor5} we may
 select a block sequence \seq{y}{n} in $Y$ such that
 $\|y_n\|_{\mc{F}_*}=1$,
 $\|y_n\|_F\stackrel{n\to\infty}{\longrightarrow}0$ and
 $b^*(y_n)\stackrel{n\to\infty}{\longrightarrow}0$ for every
 branch $b$.
 Observe also that \seq{y}{n} is a separated sequence (Definition
 \ref{aa1}) hence from Lemma \ref{aa2} we may assume passing to a
 subsequence  that for every $\sigma_{\mc{F}}$ special functional
 $x^*$ we have that
 $|x^*(y_n)|\ge\frac{1}{q^2}$ for at most two $y_n$.
 (Although Lemma \ref{aa2} is stated for $JT_{\mc{F}_2}$ with
 small modifications in the proof remains valid for either
 $\mc{F}_s$ or  $\mc{F}_{2,s}$.)

 We set $t_1=1$. From \eqref{beq9} there exists a $\sigma_{\mc{F}}$ special
  functional $y_1^*$ with $\ran
 y_1^*\subset\ran y_{t_1}$ such that $|y_1^*(y_{t_1})|>\frac{\e}{2}$.
 Setting $d_1=\max\ind y_1^*$ we select $t_2$ such that
 $\|y_{t_2}\|_F<\frac{\e}{4d_1}$. Let $z_2^*$ be a $\sigma_{\mc{F}}$
 special
 functional with $\ran z_2^*\subset\ran y_{t_2}$ such that
 $|z_2^*(y_{t_2})|>\frac{3\e}{4}$. We write $z_2^*=x_2^*+y_2^*$
 with $\ind x_2^*\subset\{1,\ldots,d_1\}$ and
 $\ind y_2^*\subset\{d_1+1,\ldots\}$. We have that
 $|x_2^*(y_{t_2})|\le d_1\|y_{t_2}\|_F<\frac{\e}{4}$  and thus
 $|y_2^*(y_{t_2})|>\frac{\e}{2}$. Following this procedure we
 select a finite collection $(t_n)_{n=1}^{q^2}$ of integers and a finite
 sequence of $\sigma_{\mc{F}}$
 special functionals $(y_n^*)_{n=1}^{q^2}$ with
 $\ran y_n^*\subset \ran y_{t_n}$ and $|y_n^*(y_{t_n})|>\frac{\e}{2}$
 such that the sets of indices $(\ind(y_n^*))_{n=1}^{q^2}$ are
 pairwise disjoint. We consider the vector
 $y=y_{t_1}+y_{t_2}+\ldots+y_{t_{q^2}}$.

 We set $\e_i=\sgn y_n^*(y_{t_n})$ for $i=1,\ldots,q^2$.
 The functional $f=\sum\limits_{n=1}^{q^2}\frac{\e_n}{q}y_n^*$
 belongs to $\mc{F}_{2,s}$ ($\subset \mc{F}_2$), while
 $qf\in\mc{F}_{s}$. Therefore
 \begin{equation}\label{beq10}
 \|y\|_{\mc{F}_*}\ge f(y_1+y_2+\ldots+y_{q^2})
 \ge
 \frac{1}{q}\sum\limits_{n=1}^{q^2}|y_n^*(y_{t_n})|
 >q\frac{\e}{2}>\frac{4}{\e}.
 \end{equation}

 It is enough to show that $\sup\{|x^*(y)|:\;x^*\mbox{ is a $\sigma_{\mc{F}}$ special
 functional}\}\le 3$ so as to derive a contradiction with
 \eqref{beq9} and \eqref{beq10}.
 Let $x^*$ be a $\sigma_{\mc{F}}$ special  functional.
 Then from our assumptions  that $|x^*(y_n)|>\frac{1}{q^2}$ for at most
 two $y_n$,

 \[ |x^*(y)|\le 2 + (q^2-2)\frac{1}{q^2}<3 .\]
 The proof of the lemma is complete.
 \end{proof}

 \begin{theorem} \label{bth2}
 Let $Y$ be a subspace  of either $JT_{\mc{F}_{2,s}}$
 or of $JT_{\mc{F}_s}$.
 Then for every $\e>0$, there
 exists a subspace of $Y$ which $1+\e$ isomorphic to $c_0$.
 \end{theorem}
 \begin{proof}[\bf Proof.]
 Let $Y$ be a block subspace  of $JT_{\mc{F}_{2,s}}$
 or of $JT_{\mc{F}_s}$
  and let $\e>0$.
 Using Lemma \ref{blem8} we may inductively select a
 normalized block sequence
 \seq{y}{n} in $Y$ such that, setting $d_n=\max\supp y_n$ for each
 $n$ and $d_0=1$, $|x^*(y_n)|<\frac{\e}{2^nd_{n-1}}$ for every $\sigma_{\mc{F}}$
 special functional $x^*$.

 We claim that \seq{y}{n} is $1+\e$ isomorphic to the standard
 basis of $c_0$. Indeed, let $(\beta_n)_{n=1}^N$ be a sequence of
 scalars. We shall show that $\max\limits_{1\le n\le N}|\beta_n|
 \le \|\sum\limits_{n=1}^N\beta_ny_n\|_{\mc{F}_*} \le (1+\e)
 \max\limits_{1\le n\le  N}|\beta_n|$ for either
 $\mc{F}_*=\mc{F}_{2,s}$ or $\mc{F}_*=\mc{F}_s$.
 We may assume that $\max\limits_{1\le n\le N}|\beta_n|=1$. The
 left inequality follows directly from the bimonotonicity of the
 Schauder basis \seq{e}{n} of $JT_{\mc{F}_*}$.

 To see the right inequality we consider an arbitrary $g\in
 \mc{F}_*$.
 Then there exist $d\in \N$, $(x_i^*)_{i=1}^d$
 in $\mathscr{S}\cup(\bigcup\limits_{i=1}^{\infty}F_i)$
 with
 $(\ind(x_i^*))_{i=1}^d$ pairwise disjoint
 and $\min\supp x_i^*\ge d$,
 such that
 $g=\sum\limits_{i=1}^da_ix_i^*$
 with $\sum\limits_{i=1}^d a_i^2\le 1$ in the case
 $\mc{F}_*=\mc{F}_{2,s}$, while $a_i\in\{-1,1\}$ in the case
 $\mc{F}_*=\mc{F}_s$.
  Let $n_0$ be the minimum integer $n$ such that $d\le d_n$.
 Since $\min\supp g\ge d>d_{n_0-1}$ we get that $g(y_n)=0$
 for $n<n_0$. In either case we get that
 \begin{eqnarray*}
  g(\sum\limits_{n=1}^N\beta_ny_n)&\le&
  |g(y_{n_0})|+\sum\limits_{n=n_0+1}^N|g(y_n)|
 \le 1+\sum\limits_{n=n_0+1}^{N}\sum\limits_{i=1}^d|x_i^*(y_n)|\\
 && <1+\sum\limits_{n=n_0+1}^{N} d\frac{\e}{2^nd_{n-1}}<
 1+ \sum\limits_{n=n_0+1}^{N} \frac{\e}{2^n}<1+\e.
 \end{eqnarray*}
 The proof of the theorem is complete.
 \end{proof}

  \begin{lemma} \label{blem2}
 For every $x\in\co$ and every $\e>0$ there exists $d\in\N$
 (denoted  $d=d(x,\e)$)
 such
 that for every $g\in \mc{F}_2\setminus F_0$ with
 $\ind(g)\cap\{1,\ldots,d\}=\emptyset$ we have that $|g(x)|<\e$.
 \end{lemma}
 \begin{proof}[\bf Proof.]
 Let $C=\|x\|_{\ell_1}$ be the $\ell_1$ norm of the vector $x$.
 We choose $d\in\N$ such that
 $\sum\limits_{l=d+1}^{\infty}{\tau}_l^2<(\frac{\e}{C})^2$.
 Now let $g=\sum\limits_{i=1}^ka_ix_i^*\in \mc{F}_2$,
  such that
 $\ind(g)\cap\{1,\ldots,d\}=\emptyset$. Each $x_i^*$
 takes the form
 $x_i^*=\sum\limits_{j=1}^{r_i}x_{i,j}^*$,
 where for each $i$ either $r_i=1$ and
 $x_{i,1}\in\bigcup\limits_{i=1}^{\infty}F_i$ or
 $(x_{i,j})_{j=1}^{r_i}$ is a $\sigma_{\mc{F}}$  special sequence,
 and the indices $(\ind(x_{i,j}^*))_{i,j}$ are pairwise different
 elements of
 $\{d+1,d+2,\ldots\}$. We get that
 \begin{eqnarray*}
 |g(x)|& \le & \sum\limits_{i=1}^k|a_i|\cdot|x_i^*(x)|\le
   \big(\sum\limits_{i=1}^k|a_i|^2)\big)^{1/2}
   \big(\sum\limits_{i=1}^k|x_i^*(x)|^2\big)^{1/2} \\
  & \le &
  1\cdot \big(\sum\limits_{i=1}^k\|x\|^2_{\ell_1}\|x_i^*\|^2_{\infty} \big)^{1/2}
  \le C
  \big(\sum\limits_{l=d+1}^{\infty}{\tau}_l^2\big)^{1/2}
  <\e.
 \end{eqnarray*}
 \end{proof}

 \begin{theorem} \label{bthm9}
 For every subspace $Y$ of $JT_{\mc{F}_2}$ and every $\e>0$ there exists a
 subspace of $Y$ which is $1+\e$ isomorphic to $\ell_2$.
 \end{theorem}
 \begin{proof}[\bf Proof.]
 Let $Y$ be a block subspace of $JT_{\mc{F}_2}$ and let $\e>0$. We choose a
 sequence  \seq{\e}{n} of positive reals satisfying
 $\sum\limits_{n=1}^{\infty}\e_n<\frac{\e}{2}$.
 We shall produce a block sequence \seq{x}{n} in $Y$ and a
 strictly increasing sequence of integers \seq{d}{n} such that
 \begin{enumerate}
 \item[(i)] $\|x_n\|_{\mc{F}_2}=1.$
 \item[(ii)] For every $\sigma_{\mc{F}}$ special functional $x^*$  we have that
 $|x^*(x_n)|<\frac{\e_n}{3d_{n-1}}$.
 \item[(iii)] If $g=\sum\limits_{i=1}^ka_iy_i^*\in \mc{F}_2$
 is such that
 $\ind(g)\cap\{1,2,\ldots,d_n\}=\emptyset$, then $|g(x_n)|<\e_n$.
 \end{enumerate}
  The construction is inductive. We choose an arbitrary finitely
  supported vector $x_1\in Y$ with $\|x_1\|_{\mc{F}_2}=1$ and we set
  $d_1=d(x_1,\e_1)$ (see  the notation in the statement of Lemma
  \ref{blem2}). Then, using Lemma \ref{blem8}
   we select a vector
  $x_2\in Y\cap\co$ with $x_1<x_2$ such that $\|x_2\|_{\mc{F}_2}=1$ and
  $|x^*(x_2)|<\frac{\e_2}{3d_1}$
  for every $\sigma_{\mc{F}}$ special functional $x^*$.
   We set $d_2=d(x_2,\e_2)$.
  It is clear how the inductive construction proceeds.
  We shall show that for every sequence of scalars $(\beta_n)_{n=1}^N$
  we have that
  \begin{equation}\label{beq11}
  (1-\e)\big(\sum\limits_{n=1}^{N}\beta_n^2\big)^{1/2}
  \le \|\sum\limits_{n=1}^N\beta_nx_n\|_{\mc{F}_2}
  \le (1+\e) \big(\sum\limits_{n=1}^{N}\beta_n^2\big)^{1/2}.
  \end{equation}
  We may assume that $\sum\limits_{n=1}^{N}\beta_n^2=1$.

 We first show the left hand inequality of \eqref{beq11}.
  For each $n$ we choose
 $g_n\in \mc{F}_2$, $g_n=\sum\limits_{i=1}^{l_n}a_{n,i}g_{n,i}$,
 with $\ran g_n\subset \ran  x_n$ and
 \begin{equation} \label{beq2}
 g_n(x_n)>1-\frac{\e}{3}.
 \end{equation}
  For each $(n,i)$ we write the  functional
 $g_{n,i}$ as the sum of  three successive functionals,
 $g_{n,i}=x_{n,i}^*+y_{n,i}^*+z_{n,i}^*$ such that
 $\ind(x_{n,i}^*)\subset \{1,\ldots,d_{n-1}\}$,
 $\ind(y_{n,i}^*)\subset \{d_{n-1}+1,\ldots,d_{n}\}$ and
 $\ind(z_{n,i}^*)\subset \{d_n+1,\ldots\}$.
 From the choice of the vector $x_n$ and the definition
 $x_{n,i}^*$ we get that
 \begin{equation} \label{beq3}
 |\sum\limits_{i=1}^{l_n}a_{n,i}x^*_{n,i}(x_n)|\le
 \sum\limits_{i=1}^{l_n}|x^*_{n,i}(x_n)|\le d_{n-1}\cdot
 \|x_n\|_F<d_{n-1}\frac{\e_n}{3d_{n-1}}<\frac{\e}{3}.
 \end{equation}
 The definition of the number $d_n$ yields also that
 \begin{equation} \label{beq4}
 |\sum\limits_{i=1}^{l_n}a_{n,i}z^*_{n,i}(x_n)|<\frac{\e}{3}.
 \end{equation}
 From \eqref{beq2}, \eqref{beq3} and \eqref{beq4} we get that
 $g'_n(x_n)>1-\e$ where the functional
 $g_n'=\sum\limits_{i=1}^{l_n}a_{n,i}y^*_{n,i}$, belongs to $\mc{F}_2$,
 satisfies $\ran(g_n')\subset\ran(g_n)\subset
 \ran(x_n)$ and $\ind(g_n')\subset
 \{d_{n-1}+1,\ldots,d_n\}$.
 We consider the functional
 $g=\sum\limits_{n=1}^N\beta_ng_n'
 =\sum\limits_{n=1}^N\sum\limits_{i=1}^{l_n}\beta_na_{n,i}y_{n,i}^*$.
 Since $\sum\limits_n\sum\limits_i(\beta_na_{n,i})^2\le 1$ and the
 sets  $(\ind(y^*_{n,i}))_{n,i}$ are pairwise disjoint we
 get that $g\in \overline{\mc{F}_2}^p\subset B_{JT_{\mc{F}_2}^*}$.
 Therefore
 \[\|\sum\limits_{n=1}^N\beta_nx_n\|_{\mc{F}_2}\ge
 g(\sum\limits_{n=1}^N\beta_nx_n) \sum\limits_{n=1}^N\beta_n^2g_n'(x_n)>1-\e.\]

 We next show the right hand inequality of \eqref{beq11}.
 Let $(\beta_n)_{n=1}^N$ be any sequence of scalars such that
 $\sum\limits_{n=1}^N\beta_n^2\le 1$.
 We consider an arbitrary
 $f\in \mc{F}_2$,
 and we shall show that $f(\sum\limits_{n=1}^N\beta_nx_n)\le 1+\e$.
 Let $f=\sum\limits_{i=1}^ka_ix_i^*$, where $(x_i^*)_{i=1}^k$
 belong to $\mathscr{S}\cup(\bigcup\limits_{i=1}^{\infty}F_i)$
 with pairwise disjoint sets of
 indices and $\sum\limits_{i=1}^ka_i^2\le 1$.

 We partition the set $\{1,2,\ldots,k\}$  in the following manner.
 We set
 \[  A_1=\big\{i\in\{1,2,\ldots,k\}:\;
 \ind(x_i^*)\cap\{1,2,\ldots,d_1\}\neq\emptyset\big\}.  \]
 If $A_1,\ldots,A_{n-1}$ have been defined we set
 \[  A_n=\big\{i\in\{1,2,\ldots,k\}:\;
 \ind(x_i^*)\cap\{1,2,\ldots,d_n\}\neq\emptyset\big\}\setminus
 \bigcup\limits_{i=1}^{n-1}A_i.  \]
 Finally we set
 $A_{N+1}=\{1,2,\ldots,k\}\setminus
 \bigcup\limits_{i=1}^{N}A_i.$

 The sets $(A_n)_{n=1}^{N+1}$ are pairwise disjoint and
 $\#\big(\bigcup\limits_{i=1}^{n}A_i\big)\le d_n$ for
 $n=1,2,\ldots,N$.
 We set
 \[  f_{A_n}=\sum\limits_{i\in A_n}a_ix_i^*,\qquad
 i=1,2,\ldots,N+1.\]
 It is clear that $\|f_{A_n}\|_{JT_{\mc{F}_2}^*}\le\big(\sum\limits_{i\in
 A_n}a_i^2\big)^{1/2}$ for each $n$.

 Let $n\in\{1,2,\ldots,N\}$ be fixed. We have that
 \begin{equation}\label{beq6}
  |f(\beta_nx_n)|\le\sum\limits_{l=1}^{n-1}|f_{A_l}(x_n)|+
                    |f_{A_n}(\beta_nx_n)|+
                    |\sum\limits_{l=n+1}^{N+1}f_{A_l}(x_n)|.
 \end{equation}
 From condition (ii) we get that
 \begin{equation}\label{beq7}
  \sum\limits_{l=1}^{n-1}|f_{A_l}(x_n)|
  \le\sum\limits_{i\in\cup_{l=1}^{n-1}A_l}|x_i^*(x_n)|
  <d_{n-1}\cdot\frac{\e_n}{d_{n-1}}=\e_n.
  \end{equation}
 On the other hand\\
 $\ind(\sum\limits_{l=n+1}^{N+1}f_{A_l})\cap\{1,\ldots,d_n\}=\ind(\sum\limits_{l=n+1}^{N+1}\sum\limits_{i\in
  A_l}a_ix_i^*)\cap\{1,\ldots,d_n\}=\emptyset$ and thus condition
  (iii) yields that
 \begin{equation}\label{beq8}
 |\sum\limits_{l=n+1}^{N+1}f_{A_l}(x_n)|<\e_n.
 \end{equation}
 Inequalities \eqref{beq6},\eqref{beq7} and \eqref{beq8} yield that
 \[|f(\beta_nx_n)|<|f_{A_n}(\beta_nx_n)|+2\e_n.\]
 Therefore
 \begin{eqnarray*}
 |f\big(\sum\limits_{n=1}^N\beta_nx_n\big)| & \le &
 \sum\limits_{n=1}^N|f(\beta_nx_n)|\le
 \sum\limits_{n=1}^N\big(|f_{A_n}(\beta_nx_n)|+2\e_n\big)\\
  & \le &  \sum\limits_{n=1}^N|\beta_n|
 |f_{A_n}(x_n)|+2\sum\limits_{n=1}^N \e_n\\
  &<&
 \big(\sum\limits_{n=1}^N|\beta_n|^2\big)^{1/2}
 \big(\sum\limits_{n=1}^N|f_{A_n}(x_n)|^2\big)^{1/2}+\e \\
  & \le &  1\cdot \big(\sum\limits_{n=1}^N\sum\limits_{i\in
 A_n}a_i^2\big)^{1/2}+\e \le 1+\e.
 \end{eqnarray*}

  The proof  of the theorem is complete.
 \end{proof}

 \begin{proposition}\label{bprop10}
 The dual space $JT_{\mc{F}_*}^*$ is equal to the closed linear span of the
 set containing $(e_n^*)_{n\in\N}$ and $b^*$ for every branch $b$,
 \[
 JT_{\mc{F}_*}^*=\overline{\spann}(\{e_n^*:\;n\in\N\}\cup\{b^*:\;b\mbox{
 is a $\sigma_{\mc{F}}$ branch}\}).  \]
 Moreover the Schauder basis
 \seq{e}{n} of the space $JT_{\mc{F}_*}$ is weakly null.
 \end{proposition}
 \begin{proof}[\bf Proof.]
 Since the space $JT_{\mc{F}_*}$ is $c_0$ saturated
 for $\mc{F}_*=\mc{F}_{2,s}$ or $\mc{F}_*=\mc{F}_s$
 (Theorem \ref{bth2}) or
 $\ell_2$ saturated (for $\mc{F}_*=\mc{F}_2$)
 it contains no isomorphic copy of $\ell_1$. Haydon's
 theorem yields that the unit ball of $JT_{\mc{F}_*}^*$ is the norm closed
 convex hull of its extreme points. Since the set $\mc{F}_*$ is the
 norming set of the space $JT_{\mc{F}_*}$  we have that
 $B_{JT_{\mc{F}_*}^*}=\overline{\conv(\mc{F}_*)}^{w^*}$ hence
 $\Ext(B_{JT_{\mc{F}_*}^*})\subset \overline{\mc{F}_*}^{w^*}$.
 We thus get that
 $JT_{\mc{F}_*}^*=\overline{\spann}(\overline{\mc{F}_*}^{w^*})$.

 We observe that
 \begin{eqnarray*}
 \overline{\mc{F}_{2,s}}^{w^*}&=&F_0\cup
 \big\{\sum\limits_{i=1}^{d}a_ix_i^*:\;\sum\limits_{i=1}^{d}a_i^2\le1,\;
 (x_i^*)_{i=1}^{d}\mbox{ are $\sigma_{\mc{F}}$ special functionals}\\
  && \mbox{with }(\ind(x_i^*))_{i=1}^d
 \mbox{ pairwise disjoint and }\min\supp x_i^*\ge d\big\},
 \end{eqnarray*}
 \begin{eqnarray*}
 \overline{\mc{F}_2}^{w^*}&=& F_0\cup
 \big\{\sum\limits_{i=1}^{\infty}a_ix_i^*:
 \;\sum\limits_{i=1}^{\infty}a_i^2\le1,\;
 (x_i^*)_{i=1}^{\infty}\mbox{ are $\sigma_{\mc{F}}$ special functionals }\\
 && \mbox{ with }(\ind(x_i^*))_{i=1}^{\infty}
 \mbox{ pairwise disjoint}\big\}.
 \end{eqnarray*}
 and
 \begin{eqnarray*}
 \overline{\mc{F}_s}^{w^*}&=&F_0\cup
 \big\{\sum\limits_{i=1}^{d}\e_ix_i^*:\;\e_i\in\{-1,1\},\;
 (x_i^*)_{i=1}^{d}\mbox{ are $\sigma_{\mc{F}}$ special functionals}\\
  && \mbox{with }(\ind(x_i^*))_{i=1}^d
 \mbox{ pairwise disjoint and }\min\supp x_i^*\ge d\big\},
 \end{eqnarray*}
 The first and third equality follow easily. For the second the
 arguments are similar to Lemma 8.4.5 of \cite{Fa}.

 The first part of the proposition
 for the cases $\mc{F}_*=\mc{F}_{2,s}$ or $\mc{F}_*=\mc{F}_{s}$
 follows directly
 while for
 the case $\mc{F}_*=\mc{F}_2$ it is enough to observe that
  $\|\sum\limits_{i=1}^{\infty}a_ix_i^*\|_{JT_{\mc{F}_2}^*}
 \le\big(\sum\limits_{i=1}^{\infty}a_i^2\big)^{1/2}$
 for every $g=\sum\limits_{i=1}^{\infty}a_ix_i^*\in \mc{F}_2$.
 Therefore
 \[
 JT_{\mc{F}_*}^*=\overline{\spann}(\{e_n^*:\;n\in\N\}\cup\{b^*:\;b\mbox{
 branch}\}).  \]

 From the first part of the proposition, to show that
 the basis \seq{e}{n} is weakly null, it is enough to show that
 $b^*(e_n)\stackrel{n\to\infty}{\longrightarrow}0$ for every
 branch $b$. But if $b=(f_1,f_2,f_3,\ldots)$ is an arbitrary
 branch then the sequence $k_n=\ind(f_n)$ is strictly increasing
 and hence,
 since $\|f_n\|_{\infty}\le {\tau}_{k_n}$, the conclusion follows.
 \end{proof}

 \begin{remark}\label{brem7}
  Let $(F_j)_{j=0}^{\infty}$ be a JTG family  (Definition \ref{bdef21}).
 If $\tau$ is a subfamily of the family of finite
 $\sigma_{\mc{F}}$ special functionals such that $F\subset \tau$ and
 $Ex^*\in\tau$ for every $x^*\in \tau$ and  interval $E$ of $\N$,
 then,
 setting
 \begin{eqnarray*}
 \mc{F}_{\tau,s}&=&\{\sum\limits_{i=1}^d\e_ix_i^*:\;\e_1,\ldots,\e_d\in\{-1,1\},
 \;x_1^*,\ldots,x_d^*\in \tau\\
 & & \mbox{  with }\ind(x_i^*)_{i=1}^d\mbox{ pairwise disjoint  and
 }\min\supp x_i^*\ge d,\;d\in\N\}.
 \end{eqnarray*}
 the space $JT_{\mc{F}_{\tau,s}}$, which is defined to be the completion of
  $(\co,\|\;\|_{\mc{F}_{\tau,s}})$, is also $c_0$
 saturated.
 \end{remark}

% In the sequel we shall denote by $(\mc{D},<)$ the dyadic tree and
% we also denote by $<_{lex}$ the lexicographical order of it.
%  The dyadic tree $\mc{D}$ with
% the lexicographical ordering is isomorphic to $\N$
%  in its natural
% ordering. Observe also that
% for $a\in\mc{D}$,
% $\{\beta\in\mc{D}:\;\beta<a\}\subset\{\beta\in\mc{D}:\;\beta<_{lex}a\}.$

 \begin{theorem}\label{bth6} There exists
 a Banach space $X$ with a weakly null
 Schauder basis \seq{e}{n} such that $X$ is $\ell_2$ saturated
 ($c_0$ saturated)  and
 for every $M\in [\N]$ the space
 $X_M=\overline{\spann}\{e_n:\;n\in M\}$ has nonseparable dual.
 \end{theorem}

 A similar result has been also obtained by E. Odell in \cite{O}
 using a different approach.

 \begin{proof}[\bf Proof.]
 Let $(F_j)_{j=0}^{\infty}$ be the Maurey-Rosenthal JTG family
 (Example \ref{bexa1} (i)).
 As we have seen  the space $X=JT_{\mc{F}_2}$ (Definition \ref{bdef17})
 has a normalized weakly null Schauder basis \seq{e}{n}
 (Proposition \ref{bprop10}) and it is $\ell_2$ saturated
 (Theorem \ref{bthm9}).

 Let now $M\in [\N]$.
 We inductively construct
 $(x_a,f_a,j_a)_{a\in\mc{D}}$, where $\mc{D}$ is the dyadic tree
 and the induction runs on the lexicographical order of
 $\mc{D}$, such that the following conditions are satisfied:
 \begin{enumerate}
 \item[(i)] For every $a\in \mc{D}$ there exists $F_a\subset M$
 with $\#(F_a)=k_{j_a}$ such that
 $x_a=\frac{1}{\sqrt{k_{j_a}}}\sum\limits_{i\in F_a}e_i$ and
 $f_a=\frac{1}{\sqrt{k_{j_a}}}\sum\limits_{i\in F_a}e_i^*$.
 \item[(ii)] $j_{\emptyset}\in \Xi_1$ with  $j_{\emptyset}\ge 2$
 while for $a\in \mc{D}$, $a\neq \emptyset$,
 $j_a=\sigma_{\mc{F}}((f_{\beta})_{\beta<a})$.
 \item[(iii)] If $a<_{lex}\beta$ then $F_a<F_\beta$.
 \end{enumerate}
 Our construction yields that for every branch $b$ of the
 dyadic tree the sequence $(f_a)_{a\in b}$ is a $\sigma_{\mc{F}}$ special
 sequence. Hence the $w^*$ sum $g_b=\sum\limits_{a\in b}f_a$ is
 a member of $\overline{\mathscr{S}}^{w^*}$ and thus it belongs to
 the unit ball of $JT_{\mc{F}_2}$.
  We shall show that
 $\|g_b|_{X_M}-g_{b'}|_{X_M}\|_{X_M^*}\ge\frac{1}{2}$ for  infinite branches
$b\neq  b'$   of the dyadic tree.

 We first
 observe that for every $a\in\mc{D}$ and $f\in F_j$
 we have that $|f(x_a)|\le \min\{\frac{\sqrt{k_j}}{\sqrt{k_{j_a}}},
  \frac{\sqrt{k_{j_a}}}{\sqrt{k_j}}\}$.
 Thus, if $g=\sum\limits_{i=1}^da_ix_i^*\in\mc{F}_2$ (Definition
 \ref{bdef17}), then
 \[ |g(x_a)|\le \sum\limits_{i=1}^d|x_i^*(x_a)|\le
 \sum\limits_{j<j_a}\frac{\sqrt{k_j}}{\sqrt{k_{j_a}}}+1+
 \sum\limits_{j>j_a}\frac{\sqrt{k_{j_a}}}{\sqrt{k_j}}\le 1+1=2.\]
 We conclude that
  $x_a\in X_M$ with $\|x_a\|_{\mc{F}_2}\le 2$.

 Therefore, if $b\neq
 b'$ are infinite branches of $\mc{D}$ and $a\in b\setminus b'$
 then
  \[\|g_b|_{X_M}-g_{b'}|_{X_M}\|_{X_M^*}\ge
  \frac{(g_b-g_{b'})(x_a)}{\|x_a\|_{\mc{F}_2}}
  \ge\frac{f_a(x_a)}{2}=\frac{1}{2}.\]

  The $c_0$ saturated space of the statement is
  the space  $X=JT_{\mc{F}_{2,s}}$ and the proof is the same.
 \end{proof}

  \end{appendix}

 \end{document}